\documentclass[12pt]{amsart}
\usepackage{amsfonts,amssymb,amsmath,amscd,amsthm, amstext,bm,color}

\usepackage[colorlinks, linkcolor=blue, citecolor=blue, urlcolor=blue,  pagebackref,hypertexnames=false]{hyperref}
\usepackage{enumerate}
\usepackage{array}
\usepackage{fancyhdr}
\usepackage{geometry}
\usepackage{graphicx}
\usepackage{epstopdf}
\usepackage{mathrsfs}
\usepackage{txfonts}
\usepackage{tikz}
\usetikzlibrary{arrows.meta, positioning, decorations.pathreplacing}
\usepackage{url}

\usepackage{comment}
\allowdisplaybreaks
\numberwithin{equation}{section}
\theoremstyle{plain}
\newtheorem{theorem} {Theorem} [section]
\newtheorem{thmA} {Theorem}

\newtheorem{proposition} [theorem]{Proposition}
\newtheorem{lemma} [theorem]{Lemma}
\newtheorem{corollary} [theorem]{Corollary}
\newtheorem{remark} [theorem]{Remark}
\newtheorem{example} [theorem]{Example}

\newcommand{\R} {\mathbb R}

\title[Reverse square-function estimates]
{Aperture and kernel dependence of reverse estimates for one-parameter and Zygmund square functions}

\author{Ji Li}
\address{Ji Li, School of Mathematical and Physical Sciences,
Macquarie University, NSW 2109, Australia}
\email{ji.li@mq.edu.au}

\author{Liangchuan Wu}
\address{Liangchuan Wu, School of Mathematical Sciences, Anhui University, Hefei 230601, P.R.~China}
\email{wuliangchuan@ahu.edu.cn}

\date{}

\begin{document}

\begin{abstract}

In this paper, we discover a previously overlooked dependence on the kernel and aperture
in $L^p$ estimates for the square functions of the type considered by Chang--Wilson--Wolff and subsequent 
work.

In the one-parameter setting on $\R^d$, these constants lie between orders $\sqrt p$ and $p$
for nonzero real-valued smooth compactly supported radial kernels with integral zero.
A support-distance criterion characterizes the possible support radii of
normalized Calder\'on dual kernels and ensures the order $\sqrt p$ when
the support distance is sufficiently small relative to the aperture.
At any fixed aperture, different kernels have best reverse constants of orders $\sqrt p$ and $p$.
We also construct a fixed kernel for which the order changes from $p$ at small apertures
to $\sqrt p$ at large apertures.

For the area function of R. Fefferman and Pipher associated with Zygmund dilations,
we prove sharp order-$p$ reverse estimates under support conditions on the kernels.
Removing one or both conditions leads to upper bounds of orders $p^{3/2}$ and $p^2$, respectively.
At any fixed aperture, these bounds are sharp over the corresponding kernel classes.
For every fixed kernel pair, the best constant has order $p$ at all sufficiently large apertures;
we also construct a pair for which it has order $p^2$ at all sufficiently small apertures.

The proofs use the one-parameter and product Haar estimates from our joint work with Jill Pipher on
Hoeffding inequalities \cite{LPW}.
For the Zygmund estimates, we adapt the intermediate square function method to the coupled scales
and construct dual kernels with supports suited to the aperture.
As an application, we prove the sharp order-$p^2$ bound anticipated by
R. Fefferman and Pipher for their Fourier multiplier class.
\end{abstract}



\maketitle
\tableofcontents

\section{Introduction}\label{sec:intro}

The exponential-square theorem of Chang, Wilson and Wolff \cite[Theorem~3.2]{CWW}
sharpens the John--Nirenberg inequality \cite{JN} for functions with
bounded Poisson area function. To recall the estimate, let $u(y,t)=P_t*f(y)$
be the Poisson extension of a function $f$ on $\R^d$, where $d\ge1$.
For $\alpha>0$ and $x\in\R^d$, consider the cone
$$
\Gamma_\alpha(x)=\left\{(y,t)\in\R^d\times(0,\infty):|y-x|<\alpha t\right\}
$$
and the associated area function
$$
A_\alpha^\nabla f(x)
=\bigg(\int_{\Gamma_\alpha(x)}|t\nabla_yu(y,t)|^2\,\frac{dy\,dt}{t^{d+1}}\bigg)^{1/2}.
$$
Here $\nabla_y$ is the spatial gradient. If
$0<\|A_\alpha^\nabla f\|_{L^\infty(\R^d)}<\infty$, then
$$
\frac1{|Q|}\int_Q
\exp\left( c_{d,\alpha} \frac{|f(x)-f_Q|^2}  {\|A_\alpha^\nabla f\|_{L^\infty(\R^d)}^2} \right)\,dx
\le C_{d,\alpha}
$$
for every cube $Q\subset\R^d$, where $f_Q=|Q|^{-1}\int_Qf$.
In terms of the local $L^p$ means, this estimate is equivalent to
$$
\bigg(\frac1{|Q|}\int_Q|f(x)-f_Q|^p\,dx\bigg)^{1/p}
\le C_{d,\alpha}\sqrt p\,\|A_\alpha^\nabla f\|_{L^\infty(\R^d)},\qquad 2\le p<\infty,
$$
with a constant independent of $p$ and $Q$. At $p=2$ this already implies
$f\in\mathrm{BMO}(\R^d)$. For a general BMO function, the John--Nirenberg estimate is
$$
\bigg(\frac1{|Q|}\int_Q|f(x)-f_Q|^p\,dx\bigg)^{1/p}
\le C_dp\,\|f\|_{\mathrm{BMO}(\R^d)},\qquad 2\le p<\infty,
$$
where $\|f\|_{\mathrm{BMO}(\R^d)}=\sup_Q|Q|^{-1}\int_Q|f-f_Q|$.

Chang, Wilson and Wolff \cite{CWW} established the exponential-square estimate
in their study of weighted inequalities and Schr\"odinger operators.
R.~Fefferman and Pipher \cite{FP} applied Chang--Wilson--Wolff estimates
to sharp weighted inequalities and multiparameter operators.
For the one-parameter theory of these estimates, see Wilson \cite[Chapters~3--4]{Wilson08}.
Chen, Duong, Wu and Yan \cite{CDWY} proved exponential-square estimates
for square functions associated with operators and applied them to spectral
multipliers and Schr\"odinger operators.

Grafakos, Honz\'ik and Seeger \cite{GHS} used the Gaussian decay in the
dyadic Chang--Wilson--Wolff inequality \cite[Corollary~3.1]{CWW} to prove
the optimal $\sqrt{\log(N+1)}$ bound for maximal functions associated
with $N$ uniformly controlled Mikhlin--H\"ormander multipliers.
Chen, Lin, Wu and Yan \cite{CLWY} obtained the same dependence on $N$
for families of spectral multipliers associated with a class of
nonnegative self-adjoint operators on doubling metric measure spaces.
Guo, Roos, Seeger and Yung \cite{GRSY} applied the dyadic inequality to
maximal families of Hilbert transforms along homogeneous curves.

The local moment estimate suggests a global question. For fixed
$\alpha>0$, consider the reverse inequality
$$
    \|f\|_{L^p(\R^d)}\le C_p\|A_\alpha^\nabla f\|_{L^p(\R^d)},\qquad 2\le p<\infty.
$$
Does its best constant $C_p$ also have order $\sqrt p$ as $p\to\infty$?

For dyadic atomic square functions under the almost orthogonality assumptions in \cite{LPW},
we proved the equivalence of the Chang--Wilson--Wolff type exponential-square estimate,
the reverse $L^p$ estimate with factor $\sqrt p$, and a Chernoff--Hoeffding
inequality for the atomic sums, analogous to the classical concentration inequality
for sums of independent bounded random variables \cite{Hoeffding,Vershynin}.
A basic example is the classical dyadic martingale square function on $\R$,
whose atomic summands are multiples of orthogonal Haar functions.
Wu \cite{WuJGA} established a corresponding equivalence
for Poisson-semigroup area functions of nonnegative self-adjoint operators
on doubling spaces of infinite measure, under Gaussian bounds and
H\"older regularity of the heat kernels.
For $L=-\Delta$ on $\R^d$, the area function is formed from $t^2\partial_t^2u$,
where $u=e^{-t\sqrt{-\Delta}}f$ is the Poisson extension.

For $A_\alpha^\nabla$, the best reverse constant has order $\sqrt p$ at every fixed
aperture $\alpha>0$ by Proposition~\ref{sp:p:cww}. The proof uses a reproducing
formula for the Poisson gradient in which the dual kernel can be supported in
an arbitrarily small ball, as in \cite[Lemma~3.1]{CWW}.

Even for this classical area function, the aperture affects boundedness at the $L^\infty$ endpoint.
Chang, Wilson and Wolff \cite[p.~239]{CWW} report an example of P.~Jones:
$A_{\gamma_1}^\nabla f\in L^\infty$ but $A_{\gamma_2}^\nabla f\notin L^\infty$
for some $0<\gamma_1<\gamma_2$.
To investigate the dependence on the aperture beyond the Poisson model,
we consider area functions formed from smooth compactly supported radial kernels with integral zero.
For each fixed $1<p<\infty$, changing the aperture gives equivalent $L^p$ norms of these area functions,
but the comparison constants may depend on $p$.
{\bf Does this equivalence at each fixed $p$ preserve the optimal order of the reverse constant as $p\to\infty$?}

Surprisingly, {\bf the optimal growth can depend on both the aperture and the kernel.}
This subtlety has been overlooked in the literature until now. For example, in Theorem 4.2 in \cite{FP}, the square function bounded $p^{3\over2}$ was obtained using an iteration scheme. 
But this did not take into account the difference in the bound for different one parameter square function associated with various kernels.

In this article, we study this dependence both for one-parameter area functions on $\R^d$
and for the area function associated with Zygmund dilations on $\R^3$.
Here one-parameter refers to the common dilation scale $t$ in all spatial directions, regardless of the dimension $d$.
For these one-parameter area functions, we first prove in every dimension
that the best reverse constant lies between orders $\sqrt p$ and $p$.
A support-distance condition relating the kernel and the aperture
ensures the optimal order $\sqrt p$.
In particular, this order holds at all sufficiently large apertures
for every fixed kernel; see Theorem~\ref{sp:t:general}.
Both the kernel and the aperture can change this order:
\begin{itemize}
\item \emph{A fixed aperture and different kernels.}
At any fixed aperture, including aperture $\alpha=1$, there are two fixed kernels
for which the optimal reverse constants have orders $\sqrt p$ and $p$, respectively;
see Theorem~\ref{sp:t:kernels}.

\smallskip

\item \emph{A fixed kernel and different apertures.}
There is a fixed kernel for which the optimal reverse constant has order $p$
at all sufficiently small apertures and order $\sqrt p$ at all sufficiently large apertures;
see Theorem~\ref{sp:t:aperture}.
\end{itemize}

For product area functions with independent dilation parameters, the best
reverse constant is the product of the one-parameter constants;
see Proposition~\ref{sp:p:product}.
The same questions become more difficult when the scales in different
coordinates depend on one another.
We then turn to the area function of R. Fefferman and Pipher \cite{FP} associated with the Zygmund dilations
$$
(x_1,x_2,x_3)\longmapsto(sx_1,tx_2,stx_3),\qquad s,t>0,
$$
where the scale in the third coordinate is the product of the other two.
Here the area function is formed from a kernel on $\R$ and another kernel on $\R^2$.
We prove a sharp order-$p$ reverse estimate at every positive aperture
under two support conditions at the origin; see Theorem~\ref{t:primal}.
For every fixed kernel pair, the optimal reverse constant has order $p$
at all sufficiently large apertures.
There is a fixed pair for which it has order $p^2$ at all sufficiently small apertures;
see Theorem~\ref{t:large-ap} and Corollary~\ref{c:aperture-change}.
As an application, we prove the sharp order-$p^2$ bound anticipated by R. Fefferman and Pipher
for their Fourier multiplier class in \cite[p.~360]{FP}; see Theorem~\ref{t:multipliers}.

The one-parameter results provide the first-variable reverse estimate
and the test functions for the Zygmund lower bounds.

For the remaining variables, we keep the first-variable square integral in an intermediate square function.
Its comparison with the Zygmund area function uses the coupled scales $t,st$
and dual kernels with support suited to the aperture.

We first prove the one-parameter results in Sections~\ref{sp:sec:dual}--\ref{sp:sec:dependence},
then develop the Zygmund estimates and their multiplier application.

\subsection{One-parameter reverse estimates}\label{sp:sub:results}

Let $\varphi\in C_c^\infty(\R^d)$ be nonzero, real-valued and radial, with
$
\int_{\R^d}\varphi(x)\,dx=0.
$
These are the basic assumptions on the kernel.
In dimension one, radiality means evenness.
For $t>0$, write $\varphi_t(x)=t^{-d}\varphi(x/t)$.
Using the cones $\Gamma_\alpha(x)$ defined above, we consider the area function
\begin{equation}\label{sp:e:area}
\mathcal A_{\varphi,\alpha}f(x)^2
=\int_{\Gamma_\alpha(x)}
|f*\varphi_t(y)|^2\,\frac{dy\,dt}{t^{d+1}}.
\end{equation}
For $2\le p<\infty$, the best reverse constant is
\begin{equation}\label{sp:e:best}
\mathfrak c_{p,d}^{\varphi}(\alpha)
=\sup_{0\ne f\in L^p(\R^d)}
\frac{\|f\|_{L^p(\R^d)}}{\|\mathcal A_{\varphi,\alpha}f\|_{L^p(\R^d)}}.
\end{equation}
Our aim is to determine the growth of $\mathfrak c_{p,d}^{\varphi}(\alpha)$
as $p\to\infty$, with $d$, $\varphi$ and $\alpha$ fixed.
The constants in the upper and lower bounds for
$\mathfrak c_{p,d}^{\varphi}(\alpha)$ may depend on $d$, $\varphi$ and $\alpha$,
but are independent of $p$.

We use the Fourier transform
$$
\widehat f(\xi)=\int_{\R^d}f(x)e^{-2\pi ix\cdot\xi}\,dx.
$$
For $f\in L^2(\R^d)$, our starting point is the Calder\'on reproducing formula
\begin{equation}\label{e:intro-reproduce}
f=\int_0^\infty(f*\varphi_t)*\Phi_t\,\frac{dt}{t}
\qquad\text{in }L^2(\R^d),
\end{equation}
where $\Phi_t(x)=t^{-d}\Phi(x/t)$ and $\Phi\in C_c^\infty(\R^d)$ is real and radial,
with integral zero.
Here $\varphi$ is the given kernel in the area function, whereas $\Phi$ is chosen
to reconstruct $f$ from the convolutions $f*\varphi_t$.
For \eqref{e:intro-reproduce}, the required normalization is
\begin{equation}\label{e:intro-dual}
\int_0^\infty\widehat\varphi(re_1)\widehat\Phi(re_1)\,\frac{dr}{r}=1.
\end{equation}
A suitable positive multiple of $\varphi$ always satisfies \eqref{e:intro-dual}.
For the order-$\sqrt p$ estimate, we seek a normalized dual whose support is small enough relative to the aperture.
The question is how small this support can be, since it determines how far the values
$f*\varphi_t$ enter the reconstruction at each scale.

Take $\Phi=-\Delta\chi$ with $\chi\in C_c^\infty(\R^d)$ real and radial.
Then $\int_{\R^d}\Phi=0$ without enlarging the support.
In \eqref{e:intro-dual}, the Laplacian contributes $4\pi^2r^2$.
Comparison with the radial measure in polar coordinates gives
\begin{equation}\label{e:first-cond-why}
r^2\,\frac{dr}{r}=r\,dr=r^{2-d}r^{d-1}\,dr.
\end{equation}
This leads to the function
$$
(|D|^{2-d}\varphi)(x)
=\int_{\R^d}|\xi|^{2-d}\widehat\varphi(\xi)e^{2\pi ix\cdot\xi}\,d\xi
$$
and the identity
\begin{equation}\label{e:intro-pair}
\int_0^\infty\widehat\varphi(re_1)\widehat{-\Delta\chi}(re_1)\,\frac{dr}{r}
=\frac{4\pi^2}{|\mathbb S^{d-1}|}
\int_{\R^d}(|D|^{2-d}\varphi)(x)\chi(x)\,dx.
\end{equation}
By \eqref{e:first-cond-why}, boundedness of $\widehat\varphi$ near zero
suffices for convergence at zero in \eqref{e:intro-pair}.
For $d=1$, the function in this pairing is $|D|\varphi$; for $d=2$, it is $\varphi$.

Proposition~\ref{sp:p:local-dual} proves that there is a real radial
$\Phi\in C_c^\infty(B(0,\rho))$ with integral zero satisfying
\eqref{e:intro-dual} if and only if
\begin{equation}\label{e:intro-dual-radius}
\operatorname{dist}\bigl(0,\operatorname{supp}(|D|^{2-d}\varphi)\bigr)<\rho.
\end{equation}
Thus the support distance is exactly the infimum of the possible support radii
of these normalized dual kernels.
The necessity applies to every dual kernel in this class:
each can be written as $-\Delta\chi$ with $\chi$ supported in the same ball,
so the pairing in \eqref{e:intro-pair} cannot be normalized if
$|D|^{2-d}\varphi$ vanishes throughout the permitted ball.

At scale $t$, a dual kernel supported in $B(0,\rho)$ contributes only within
distance $\rho t$ of the integration point, while the area function uses
$|x-y|<\alpha t$.
Thus the fixed kernel imposes a lower bound on $\rho$, while the cone comparison imposes
an upper bound determined by the aperture.
We obtain the order $\sqrt p$ when these two requirements are compatible.
This leads to the following result.

\begin{thmA}\label{sp:t:general}
Let $\varphi$ satisfy the basic assumptions.
For every $\alpha>0$ and $2\le p<\infty$,
$$
c_{d,\varphi,\alpha}\sqrt p
\le \mathfrak c_{p,d}^{\varphi}(\alpha)
\le C_{d,\varphi,\alpha}p.
$$
Moreover, there is a constant $c_d>0$, depending only on $d$, such that if
$
\operatorname{dist}\bigl(0,\operatorname{supp}(|D|^{2-d}\varphi)\bigr)<c_d\alpha,
$
then
$$
\mathfrak c_{p,d}^{\varphi}(\alpha)
\simeq_{d,\varphi,\alpha}\sqrt p,
\qquad 2\le p<\infty.
$$
Consequently, $\mathfrak c_{p,d}^{\varphi}(\alpha)$ has order $\sqrt p$
at every positive aperture if $0\in\operatorname{supp}(|D|^{2-d}\varphi)$,
and at all sufficiently large apertures for every fixed kernel.
\end{thmA}

Under the distance condition in Theorem~\ref{sp:t:general}, we may choose
$$
\operatorname{dist}\bigl(0,\operatorname{supp}(|D|^{2-d}\varphi)\bigr)
<\rho<c_d\alpha.
$$
The first inequality permits normalization; the second ensures that the
square sum of the atomic coefficients is controlled pointwise by
$\bigl(\mathcal A_{\varphi,\alpha}f\bigr)^2$, as in \eqref{sp:e:discrete-area}.
The Haar estimate then supplies the factor $\sqrt p$.
This is a sufficient condition for the reverse estimate; the exact
characterization in \eqref{e:intro-dual-radius} concerns the support of the dual kernel.
The choice of a local reconstruction scale suited to a given aperture also
appears in the discrete square-function comparison of \cite[Section~3]{CDWY}.
Here Proposition~\ref{sp:p:local-dual} identifies the restriction imposed
by the fixed kernel $\varphi$.

At a fixed aperture, both orders occur for different kernels.
\begin{thmA}\label{sp:t:kernels}
For every $d\ge1$ and every fixed $\alpha_0>0$, there are kernels
$\varphi^{(0)},\varphi^{(1)}$ satisfying the basic assumptions such that
$$
\mathfrak c_{p,d}^{\varphi^{(0)}}(\alpha_0)
\simeq_{d,\varphi^{(0)},\alpha_0}\sqrt p,
\qquad
\mathfrak c_{p,d}^{\varphi^{(1)}}(\alpha_0)
\simeq_{d,\varphi^{(1)},\alpha_0}p,
\qquad 2\le p<\infty.
$$
Both kernels are independent of $p$.
In particular, both orders occur at aperture $\alpha_0=1$.
\end{thmA}

The kernels in Theorem~\ref{sp:t:kernels} can be normalized so that
$$
\int_0^\infty|\widehat{\varphi^{(j)}}(re_1)|^2\,\frac{dr}{r}=1,
\qquad j=0,1.
$$
Their best reverse constants then agree at $p=2$ by \eqref{sp:e:energy},
although their growth as $p\to\infty$ is different.
The construction in Subsection~\ref{sp:sec:kernels} also ensures that
these kernels are not dilates of one another.

There also exists a single kernel for which changing the aperture changes the optimal order.

\begin{thmA}\label{sp:t:aperture}
For every $d\ge1$, there exist a kernel $\varphi$ satisfying the basic assumptions
and $0<\alpha_-<\alpha_+$ such that, for every $2\le p<\infty$,
$$
\mathfrak c_{p,d}^{\varphi}(\alpha)\simeq_{d,\varphi,\alpha}p
\quad (0<\alpha\le\alpha_-),
$$
whereas
$$
\mathfrak c_{p,d}^{\varphi}(\alpha)\simeq_{d,\varphi,\alpha}\sqrt p
\quad (\alpha\ge\alpha_+).
$$
The kernel $\varphi$ and the apertures $\alpha_-$ and $\alpha_+$
are independent of $p$.
\end{thmA}

The general lower bound in Theorem~\ref{sp:t:general} follows from lacunary sums.
We consider
$$
\mathcal P_N(x)=\sum_{k=1}^N\cos(2\pi2^kx_1),
\qquad N=\lceil p\rceil.
$$
This choice is motivated by the wave-packet example in Guth's notes
\cite[Example~1.6]{GuthDecoupling}:
terms at different frequencies are all close to $1$ on a small common region.
Here $\mathcal P_N\ge N/2$ when $0\le x_1\le2^{-N-4}$, so
$$
\|\mathcal P_N\|_{L^p([0,1]^d)}
\ge\frac N2\,2^{-(N+4)/p}\ge cp.
$$
For the area function, the term at frequency $2^k$ is largest at scales comparable to $2^{-k}$.
The decay away from these scales allows us to bound the squared area function by $CN$:
$$
\|\mathcal A_{\varphi,\alpha}\mathcal P_N\|_{L^\infty(\R^d)}
\le C_{d,\varphi,\alpha}\sqrt N
\le C_{d,\varphi,\alpha}\sqrt p.
$$
After localization, the test functions in Lemma~\ref{sp:l:lower}
have norm ratio at least $c\sqrt p$;
see also our work \cite[Lemmas~6.2 and~6.3]{LPWCWW}
for the trigonometric sums and localization.
This proves the optimality of the order $\sqrt p$ under the distance condition.

For an arbitrary aperture $\alpha$, choose a sufficiently large $\beta\ge\alpha$
so that the distance condition holds. Lemma~\ref{sp:l:aperture} then yields
$$
\|f\|_p
\le C_{d,\varphi,\beta}\sqrt p\,\|\mathcal A_{\varphi,\beta}f\|_p
\le C_{d,\varphi,\alpha,\beta}p\,\|\mathcal A_{\varphi,\alpha}f\|_p.
$$
The first factor $\sqrt p$ comes from the Haar estimate; the second comes
from the comparison of apertures.
The examples below show that the resulting order $p$ can be optimal.

For the order-$p$ examples in Theorems~\ref{sp:t:kernels} and~\ref{sp:t:aperture},
the natural test function is a truncated logarithm.
The function $\log|u|$ belongs to $\mathrm{BMO}(\R)$, and
$$
\bigg(\int_0^1|\log u|^p\,du\bigg)^{1/p}\simeq p,\qquad p\ge2.
$$
This is a classical example of the sharp dependence on $p$ in the John--Nirenberg inequality;
see also Slavin and Vasyunin \cite[Section~8]{SVBMO} for logarithmic extremizers in sharp BMO estimates.
The additional task here is to keep its area-function norm bounded.
For this purpose, we construct a kernel whose one-dimensional integral $\kappa$ satisfies
$$
U(u):=(\kappa*\log|\cdot|)(u)=0\qquad(|u|<a).
$$
Since $\int\kappa=0$, the untruncated logarithm satisfies
$$
(\log|\cdot|*\kappa_t)(y)=U(y/t).
$$
Thus these convolution values vanish in the cone $|y|<\alpha t$ when
$\alpha<a$. The proof controls the errors from truncating the logarithm
and, in higher dimensions, from localization in the other variables.
The construction produces a fixed kernel $\varphi$ and a compactly supported
function $h$, both independent of $p$, such that, for every $0<\alpha<a$,
$$
\|h\|_{L^p(\R^d)}\ge cp,
\qquad
\|\mathcal A_{\varphi,\alpha}h\|_{L^p(\R^d)}\le C_{d,\varphi,\alpha},
\qquad p\ge2;
$$
see Lemma~\ref{sp:l:log-test}.
Thus the reverse constant is at least $c_{d,\varphi,\alpha}p$.
Combined with the general upper bound, this proves that its order is exactly $p$.
At sufficiently large apertures, the same kernel satisfies the distance condition,
and its best reverse constant has order $\sqrt p$.

The restriction on local reconstruction and the logarithmic example concern the same kernel structure.
Indeed, $U''=2\pi^2|D|\kappa$, so the cancellation $U=0$ on $(-a,a)$ also implies
$|D|\kappa=0$ there. By Proposition~\ref{sp:p:marginal},
$$
\operatorname{dist}\bigl(0,\operatorname{supp}(|D|^{2-d}\varphi)\bigr)\ge a.
$$
Consequently, no normalized dual kernel from Proposition~\ref{sp:p:local-dual}
can be supported in $B(0,a)$.
For these kernels, the restriction encountered in local reconstruction is accompanied
by an order-$p$ lower bound at small apertures.
Thus the dependence on the kernel and the aperture affects the optimal estimate itself.

At a fixed aperture $\alpha_0$, choose the logarithmic kernel so that its small-aperture range
contains $\alpha_0$, and choose a second kernel $\varphi^{(0)}$ with
$0\in\operatorname{supp}(|D|^{2-d}\varphi^{(0)})$.
These are the two kernels in Theorem~\ref{sp:t:kernels}.

For the fixed logarithmic example, the area function is bounded at small apertures
and unbounded at sufficiently large apertures.
Despite the small-aperture bound, $h-h_Q$ has no exponential-square integrability
on a cube meeting its logarithmic singularity; see Remark~\ref{sp:r:log-endpoint}.
Moreover, at two suitable fixed apertures, the ratio of its larger-aperture and smaller-aperture
$L^p$ norms has order $\sqrt p$; see Corollary~\ref{sp:c:endpoint}.
The additional $\sqrt p$ in the general reverse bound is therefore forced by this fixed function.
Its reverse norm ratios at the smaller and larger apertures have orders $p$ and $\sqrt p$, respectively;
see Remark~\ref{sp:r:reconstruction}.

We also construct, in every dimension, a kernel satisfying
$(|D|^{2-d}\varphi)(0)\ne0$ for which the class of functions
with bounded conical area function depends on the aperture;
see Proposition~\ref{sp:p:bounded-area}.
For this kernel, the best reverse constant nevertheless has order
$\sqrt p$ at every fixed positive aperture.

For the one-parameter Poisson and heat area functions, and the spatial-gradient Poisson area function
in \cite{CWW}, the best reverse constants have order $\sqrt p$ at every fixed positive aperture;
see Propositions~\ref{sp:p:semigroups} and~\ref{sp:p:cww}.
The local reconstruction is explained in Remark~\ref{sp:r:cww-dual}
and Subsections~\ref{sp:sub:semigroups}--\ref{sp:sub:cww}.
For these classical kernels, the relevant local pairing can be normalized
with arbitrarily small support. Hence every positive aperture admits
the local reconstruction needed for the $\sqrt p$ estimate.

For independent dilation parameters, the one-parameter constants combine
without any further loss. Let $\varphi_i$ satisfy the basic assumptions on
$\R^{d_i}$, fix $\alpha_i>0$, $i=1,2$, and let $\mathfrak c_p^{\mathrm{prod}}$
denote the best reverse constant for the product area function
in \eqref{sp:e:product-area}. Proposition~\ref{sp:p:product} proves
$$
\mathfrak c_p^{\mathrm{prod}}
=\mathfrak c_{p,d_1}^{\varphi_1}(\alpha_1)
 \mathfrak c_{p,d_2}^{\varphi_2}(\alpha_2),\qquad 2\le p<\infty.
$$
The upper bound follows by successive Hilbert-valued one-parameter estimates,
and tensor products of test functions prove the lower bound.
Thus the orders $p$, $p^{3/2}$ and $p^2$ already occur for product area functions.
For Zygmund dilations, the third-coordinate scale depends on both dilation
parameters, and this independent iteration must be replaced by an argument
that uses their relation.

\subsection{Reverse estimates for Zygmund dilations}\label{sub:results}

The Zygmund dilations on $\R^3$ are
$$
\delta_{s,t}(x_1,x_2,x_3)=(sx_1,tx_2,stx_3),\qquad s,t>0.
$$
These dilations arise in the study of differentiation over rectangles
whose side lengths are related.
The strong maximal theorem of Jessen, Marcinkiewicz and Zygmund \cite{JMZ}
concerns independent side lengths; C\'ordoba \cite{Cor} treated families
that include the side lengths $s,t,st$, while Soria \cite{Sor} showed
limitations for more general families.
The associated singular integrals and weighted estimates were studied by
R.~Fefferman and Pipher \cite{FP}.
Han, Li, Lin and Tan \cite{HLLT} use Littlewood--Paley estimates and a Calder\'on identity
to study Zygmund singular integrals, and
Hyt\"onen, Li, Martikainen and Vuorinen \cite{HLMV} develop dyadic multiresolution analysis
and weighted estimates for these dilations.
The scales $s,t,st$ form a basic nonproduct model: there are two independent
parameters, but each one changes the scale in the third coordinate.

For the associated Lusin area function, we prove sharp order-$p$ reverse estimates
under two support conditions at the origin, then remove one or both conditions.
The dependence on the aperture and the kernels takes the following forms.
\begin{itemize}
\item \emph{A fixed aperture and different kernel pairs.}
At any fixed aperture, including aperture one, suitable kernel pairs have optimal reverse constants
of orders $p$, $p^{3/2}$ and $p^2$;
see Theorems~\ref{t:primal} and \ref{t:sharp-cases}.

\smallskip

\item \emph{A fixed kernel pair and different apertures.}
For every fixed kernel pair, the optimal reverse constant has order $p$
at all sufficiently large apertures.
There is a fixed pair for which it has order $p^2$ at all sufficiently small apertures;
see Theorem~\ref{t:large-ap} and Corollary~\ref{c:aperture-change}.
\end{itemize}

\smallskip

\noindent\textbf{Notation.}\  \   
We write $x=(x_1,x_2,x_3)$, $y=(y_1,y_2,y_3)$ and $\xi=(\xi_1,\xi_2,\xi_3)$ for points in $\R^3$.
For the last two coordinates, we use
$
\pmb x=(x_2,x_3),\ \pmb y=(y_2,y_3),\ \pmb\xi=(\xi_2,\xi_3)\in\R^2.
$
We also write $\pmb v=(v_2,v_3)\in\R^2$, with $d\pmb y=dy_2\,dy_3$ and $d\pmb v=dv_2\,dv_3$.

Following R. Fefferman and Pipher \cite[Definition~2.8]{FP}, we consider an even kernel
$\psi_1\in C_c^\infty(\R)$ and a radial kernel $\psi_2\in C_c^\infty(\R^2)$.
Both kernels are real-valued and nonzero, with
$$
\int_\R\psi_1(x_1)\,dx_1=0,\qquad
\int_{\R^2}\psi_2(x_2,x_3)\,dx_2\,dx_3=0.
$$
These are the basic assumptions on the Zygmund kernel pair throughout the compactly supported case.
Write $\psi=\psi_1\otimes\psi_2$ and
$$
\psi_{s,t}(x)=\frac1{s^2t^2}
\psi_1\left(\frac{x_1}{s}\right)
\psi_2\left(\frac{x_2}{t},\frac{x_3}{st}\right),\qquad s,t>0.
$$
For $\alpha>0$, define
$$
S_{z,\alpha}^{\psi}f(x)^2
=\int_0^\infty\int_0^\infty
\int_{\substack{|u-x_1|<\alpha s\\
|v_2-x_2|<\alpha t\\ |v_3-x_3|<\alpha st}}
|(f*\psi_{s,t})(u,v_2,v_3)|^2
\,\frac{du\,dv_2\,dv_3\,dt\,ds}{s^3t^3}.
$$
For $2\le p<\infty$, the best reverse constant is
\begin{equation}\label{e:reverse-constant}
\mathfrak C_p^\psi(\alpha)
=\sup_{0\ne f\in L^p(\R^3)}
\frac{\|f\|_{L^p(\R^3)}}
{\|S_{z,\alpha}^{\psi}f\|_{L^p(\R^3)}}.
\end{equation}
We study the growth of $\mathfrak C_p^\psi(\alpha)$ as $p\to\infty$, with $\psi$ and $\alpha$ fixed.
The constants in its upper and lower bounds may depend on $\psi$ and $\alpha$,
but are independent of $p$.

The following theorem proves the sharp order-$p$ estimate under two support conditions at the origin.

\begin{thmA}\label{t:primal}
Assume
\begin{equation}\label{e:supp-origin}
0\in\operatorname{supp}(|D|\psi_1),
 \qquad 
0\in\operatorname{supp}\psi_2.
\end{equation}
Then, for every $\alpha>0$, every $2\le p<\infty$,
and every $f\in L^p(\R^3)$,
$$
 \left\|   f \right\|_{L^p(\R^3)}
\le C_{\psi,\alpha}p\,
 \left\|   S_{z,\alpha}^{\psi}f \right\|_{L^p(\R^3)}.
$$
The constant is independent of $p$ and $f$, and the factor $p$ is optimal.
\end{thmA}

{
Following the intermediate square function method of \cite[Sections~4--5]{LPW},
we first apply Theorem~\ref{sp:t:general} in $x_1$ with the area function
$$
S_{1,\alpha}^{\psi_1}f(x)
=\mathcal A_{\psi_1,\alpha}\bigl(f(\cdot,x_2,x_3)\bigr)(x_1).
$$
The first condition in \eqref{e:supp-origin} gives
$\|f\|_p\le C_{\psi_1,\alpha}\sqrt p\,\|S_{1,\alpha}^{\psi_1}f\|_p$.
The ratio of the two remaining scales $t$ and $st$ depends on $s$.
We therefore keep the square integral in $x_1$ when comparing $S_{1,\alpha}^{\psi_1}f$
with $S_{z,\alpha}^{\psi}f$. Under the second condition in \eqref{e:supp-origin},
this comparison costs a further factor $\sqrt p$.
The two estimates are proved in Subsections~\ref{sub:first-estimate} and~\ref{sub:remaining-estimate}.

At a fixed aperture, weaker support conditions suffice. Write
$$
d_1(\psi_1)=\operatorname{dist}\bigl(0,\operatorname{supp}(|D|\psi_1)\bigr),
\qquad
d_2(\psi_2)=\operatorname{dist}\bigl(0,\operatorname{supp}\psi_2\bigr).
$$
By Proposition~\ref{p:support-distance}, the first estimate holds if $d_1(\psi_1)<c_*\alpha$,
and the second if $d_2(\psi_2)<c_*\alpha$, with an absolute constant $c_*>0$.
For a fixed kernel pair, both inequalities hold at all sufficiently large apertures.
}

\begin{thmA}\label{t:large-ap}
Let $\psi_1,\psi_2$ be any fixed kernel pair satisfying the basic assumptions.
There is $\alpha_\psi>0$ such that, for  
$2\le p<\infty$, $\alpha\ge\alpha_\psi$ and $f\in L^p(\R^3)$,
$$
\|f\|_{L^p(\R^3)}
\le C_{\psi,\alpha}p\,\|S_{z,\alpha}^{\psi}f\|_{L^p(\R^3)}.
$$
The constant is independent of $p$ and $f$, and the factor $p$ is optimal.
\end{thmA}

For every kernel pair satisfying the basic assumptions,
Proposition~\ref{p:lower-bound} proves
$$
c_{\psi,\alpha}p\le\mathfrak C_p^\psi(\alpha),\qquad \forall\, \alpha>0,\  \   2\le p<\infty.
$$
This proves the sharpness in Theorems~\ref{t:primal} and \ref{t:large-ap}.

We now fix any $\alpha>0$ and remove the first distance condition
$d_1(\psi_1)<c_*\alpha$, while still assuming $d_2(\psi_2)<c_*\alpha$. 

\begin{thmA}\label{t:first-free}
Assume $d_2(\psi_2)<c_*\alpha$.
Then, for every $2\le p<\infty$ and $f\in L^p(\R^3)$,
\begin{equation}\label{e:first-free}
\left\|f\right\|_{L^p(\R^3)}
\le C_{\psi,\alpha}p^{3/2}\,
\left\|S_{z,\alpha}^{\psi}f\right\|_{L^p(\R^3)}.
\end{equation}
The constant is independent of $p$ and $f$.
No support condition on $|D|\psi_1$ is required.
In particular, \eqref{e:first-free} holds for every $\alpha>0$
when $0\in\operatorname{supp}\psi_2$.
\end{thmA}

We instead remove the second distance condition $d_2(\psi_2)<c_*\alpha$.
In (i) of the following theorem, we assume $d_1(\psi_1)<c_*\alpha$;
in (ii), we remove this condition as well.

\begin{thmA}\label{t:second-free}
For every $\alpha>0$, $2\le p<\infty$ and $f\in L^p(\R^3)$,
\begin{itemize}
\item[(i)] if $d_1(\psi_1)<c_*\alpha$, then
$
\|f\|_{L^p(\R^3)}
\le C_{\psi,\alpha}p^{3/2}\,
\|S_{z,\alpha}^{\psi}f\|_{L^p(\R^3)};
$
\item[(ii)] without either additional support condition, we have
$
\|f\|_{L^p(\R^3)}
\le C_{\psi,\alpha}p^2\,
\|S_{z,\alpha}^{\psi}f\|_{L^p(\R^3)},
$
\end{itemize}
where $C_{\psi,\alpha}$ is independent of $p$ and $f$.
In particular, (i) holds for every $\alpha>0$ when
$0\in\operatorname{supp}(|D|\psi_1)$, with no support condition on $\psi_2$.
\end{thmA}

{
To prove Theorem~\ref{t:second-free}, we use two dual kernels supported in strips.
Each kernel satisfies the cone restriction in the narrow direction of its strip.
Only the long direction requires a change of aperture (see Section~\ref{sec:second-free}).
}

\begin{samepage}
The bounds in Theorems~\ref{t:first-free} and \ref{t:second-free} are optimal on their respective kernel classes.
The next theorem constructs fixed kernel pairs with the stated orders.

\begin{thmA}\label{t:sharp-cases}
For every $a>0$, there are three kernel pairs satisfying the basic
assumptions, one for each case below.
Each pair is fixed independently of $\alpha$ and $p$:
\begin{itemize}
\item[(i)] $d_1(\psi_1)\ge a$ and $d_2(\psi_2)=0$, with
$$
c_{\psi,\alpha}p^{3/2}\le\mathfrak C_p^\psi(\alpha)
\le C_{\psi,\alpha}p^{3/2};
$$
\item[(ii)] $d_1(\psi_1)=0$ and $d_2(\psi_2)\ge a$, with
$$
c_{\psi,\alpha}p^{3/2}\le\mathfrak C_p^\psi(\alpha)
\le C_{\psi,\alpha}p^{3/2};
$$
\item[(iii)] $d_1(\psi_1),d_2(\psi_2)\ge a$, with
$$
c_{\psi,\alpha}p^2\le\mathfrak C_p^\psi(\alpha)
\le C_{\psi,\alpha}p^2.
$$
\end{itemize}
All estimates hold for $0<\alpha<a$ and $2\le p<\infty$.
The positive constants are independent of $p$.
\end{thmA}
\end{samepage}

Combining (iii) of Theorem~\ref{t:sharp-cases} with
Theorem~\ref{t:large-ap}, we obtain a fixed kernel pair
and  $0<\alpha_-<\alpha_+$ such that
$$
\mathfrak C_p^\psi(\alpha)\simeq_{\psi,\alpha}
\begin{cases}
p^2,&0<\alpha\le\alpha_-,\\
p,&\alpha\ge\alpha_+,
\end{cases}
\qquad 2\le p<\infty;
$$
see Corollary~\ref{c:aperture-change}.

Table~\ref{tab:main-results} summarizes these reverse estimates
and shows how the bounds for the two intermediate comparisons combine.
All norms are in $L^p(\R^3)$, and the displayed bounds are up to
constants independent of $p$.
For $\mathfrak C_p^\psi(\alpha)$, the first row gives the optimal order
for every kernel pair satisfying its conditions.
For each fixed aperture, the kernel pairs in Theorem~\ref{t:sharp-cases}
show that the orders in the last three rows are sharp
on the corresponding kernel classes.

\begin{table}[htbp]
\centering
\small
\setlength{\tabcolsep}{4pt}
\renewcommand{\arraystretch}{1.2}
\begin{tabular}{@{}>{\raggedright\arraybackslash}p{.28\textwidth}
>{\centering\arraybackslash}p{.145\textwidth}
>{\centering\arraybackslash}p{.145\textwidth}
>{\centering\arraybackslash}p{.10\textwidth}
>{\raggedright\arraybackslash}p{.25\textwidth}@{}}
\hline
\textbf{Distance conditions}
& $\displaystyle\frac{\|f\|_p}{\|S_{1,\alpha}^{\psi_1}f\|_p}$
& $\displaystyle\frac{\|S_{1,\alpha}^{\psi_1}f\|_p}{\|S_{z,\alpha}^{\psi}f\|_p}$
& $\mathfrak C_p^\psi(\alpha)$
& \textbf{Upper bound; sharpness}\\
\hline
\noalign{\vskip6pt}
$d_1(\psi_1)<c_*\alpha$,\par $d_2(\psi_2)<c_*\alpha$
& $\sqrt p$ & $\sqrt p$ & $p$
& Proposition~\ref{p:support-distance};\par Proposition~\ref{p:lower-bound}\\
\noalign{\vskip5pt}
No condition on $d_1(\psi_1)$;\par $d_2(\psi_2)<c_*\alpha$
& $p$ & $\sqrt p$ & $p^{3/2}$
& Theorem~\ref{t:first-free};\par (i) of Theorem~\ref{t:sharp-cases}\\
\noalign{\vskip5pt}
$d_1(\psi_1)<c_*\alpha$;\par no condition on $d_2(\psi_2)$
& $\sqrt p$ & $p$ & $p^{3/2}$
& (i) of Theorem~\ref{t:second-free};\par (ii) of Theorem~\ref{t:sharp-cases}\\
\noalign{\vskip5pt}
Neither distance condition is required
& $p$ & $p$ & $p^2$
& (ii) of Theorem~\ref{t:second-free};\par (iii) of Theorem~\ref{t:sharp-cases}\\
\noalign{\vskip 6pt}
\hline
\end{tabular}
\caption{Reverse estimates and sharpness for any fixed $\alpha>0$ and $2\le p<\infty$.
All kernels satisfy the basic assumptions.}
\label{tab:main-results}
\end{table}

{
For the lower bounds, we dilate a test function in the third coordinate.
In the limit the planar kernel $\psi_2$ reduces to its marginal
$\kappa(u)=\int_\R\psi_2(u,v)\,dv$, and the two one-dimensional norm ratios multiply
(Lemma~\ref{l:lower-transfer}). Thus for every fixed kernel pair and every $\alpha>0$,
$$
\mathfrak C_p^\psi(\alpha)
\ge \frac{1}{\sqrt{2\alpha}}\,
\mathfrak c_{p,1}^{\psi_1}(\alpha)\mathfrak c_{p,1}^{\kappa}(\alpha),
\qquad 2\le p<\infty.
$$
The logarithmic examples in Subsection~\ref{sub:sharp-kernels} give the lower bounds
$p^{3/2}$ and $p^2$ by this formula. Under either of the two conditions on $\psi_2$ and
$\alpha$ in Proposition~\ref{p:constants}, the product also gives an upper bound
up to a constant independent of $p$.
}

We also prove sharp reverse estimates for the Poisson and heat
area functions under Zygmund dilations.
These classical semigroup models have kernels that are not compactly supported.
Let
$$
b_{\mathrm P}(r)=re^{-r},
\qquad
b_{\mathrm H}(r)=r^2e^{-r^2},
\qquad r\ge0.
$$
For $s,t>0$ and $f\in\mathcal S(\R^3)$, define
$$
\widehat{Q_{s,t}^{\nu}f}(\xi)
=b_{\nu}(s|\xi_1|)
b_{\nu}\left(t\sqrt{\xi_2^2+s^2\xi_3^2}\right)\widehat f(\xi),\qquad \nu\in\{\mathrm P,\mathrm H\}.
$$
The spatial scales are $s,t,st$.
 The integrable kernels of $Q_{s,t}^{\nu}$ also define these operators on $L^p(\R^3)$,
$1\le p\le\infty$; see Section~\ref{sec:semigroups}.
For $\alpha>0$, define
\begin{equation}\label{e:sg-area}
S_{z,\alpha}^{\nu}f(x)^2
=\int_0^\infty\int_0^\infty
\int_{\substack{|u-x_1|<\alpha s\\
|v_2-x_2|<\alpha t\\|v_3-x_3|<\alpha st}}
|Q_{s,t}^{\nu}f(u,v_2,v_3)|^2
\,\frac{du\,dv_2\,dv_3\,dt\,ds}{s^3t^3},\qquad \nu\in\{\mathrm P,\mathrm H\}.
\end{equation}

\begin{thmA}\label{t:semigroup}
For each $\nu\in\{\mathrm P,\mathrm H\}$ and every $\alpha>0$,
there are positive constants $c_{\nu,\alpha},C_{\nu,\alpha}$ such that
for $2\le p<\infty$ and $f\in L^p(\R^3)$,
$$
\|f\|_{L^p(\R^3)}
\le C_{\nu,\alpha}p\,\|S_{z,\alpha}^{\nu}f\|_{L^p(\R^3)}.
$$
The order $p$ is sharp: for every $p\ge2$, there is
$f_p\in\mathcal S(\R^3)$ with
$$
\operatorname{supp}\widehat f_p\Subset\{\xi_1\xi_2\xi_3\ne0\},\qquad
0<\|S_{z,\alpha}^{\nu}f_p\|_{L^p(\R^3)}<\infty,\qquad
\|f_p\|_{L^p(\R^3)}\ge c_{\nu,\alpha}p\,\|S_{z,\alpha}^{\nu}f_p\|_{L^p(\R^3)}.
$$
Both constants are independent of $p$.
\end{thmA}

The proof in Section~\ref{sec:semigroups} uses compactly supported dual kernels
for the reverse estimate and the one-dimensional tests from Proposition~\ref{sp:p:semigroups}
for sharpness, via third-variable dilation.

\subsection{Sharp Zygmund multiplier estimates}\label{sub:mult}

R. Fefferman and Pipher introduced a class of Fourier multipliers
associated with Zygmund dilations in \cite{FP}.
They improved the $O(p^3)$ bound from the general three-parameter
Marcinkiewicz theorem to $O(p^{5/2})$ as $p\to\infty$;
see \cite[Theorem~4.1]{FP}.
Since the Zygmund dilations have two independent parameters,
they expected an $O(p^2)$ bound, as for two-parameter product
singular integrals; see \cite[p.~360]{FP}.

 We are able to improve the $O(p^{5/2})$ bound to $O(p^2)$ bound 
and show that $p^2$ is optimal for a fixed multiplier.
Following \cite[p.~349]{FP}, in our coordinates, let
$$
A^{x_1}
=\left\{\xi\in\R^3:
1/2<|\xi_1|\le1,\quad 1/2<|(\xi_2,\xi_3)|\le1\right\}.
$$
For a positive integer $N$, define
\begin{equation}\label{e:mult-snorm}
\|m\|_{\mathcal M_{Z,x_1}^N}
=\max_{0\le a,b,c\le N}\sup_{r,q>0}\sup_{\xi\in A^{x_1}}
\left|\partial_1^a\partial_2^b\partial_3^c
[m(r\xi_1,q\xi_2,rq\xi_3)]\right|.
\end{equation}
The derivatives are taken with respect to $\xi$, with $r,q$ fixed.

\begin{thmA}\label{t:multipliers}
There is a positive integer $N_0$ such that, if
$
m\in C^{3N_0}(\{\xi:\xi_1\ne0,\ (\xi_2,\xi_3)\ne(0,0)\}),
$ with $\|m\|_{\mathcal M_{Z,x_1}^{N_0}}<\infty,
$
then the Fourier multiplier $\widehat{T_mf}=m\widehat f$ extends to $L^p(\R^3)$ and
$$
\|T_m\|_{L^p\to L^p}\le Cp^2\|m\|_{\mathcal M_{Z,x_1}^{N_0}},
\qquad 2\le p<\infty.
$$
The constant is independent of $p$ and $m$, and the factor $p^2$ is optimal.
One may take $N_0=12$.
\end{thmA}

The improvement from $p^{5/2}$ to $p^2$ uses the freedom to choose
the auxiliary square-function kernels in the multiplier proof.
We choose convolution squares, fixed independently of $m$ and $p$,
that satisfy both support conditions in \eqref{e:supp-origin}.

{Theorem~\ref{t:primal} bounds $\|T_mf\|_{L^p}$ by its area-function norm
with a factor $p$.
Following \cite[Theorem~4.4 and Lemmas~4.6--4.7]{FP},
we estimate this norm through an auxiliary weighted square function,
whose upper $L^p$ bound contributes another factor $p$. Thus
$$
\underbrace{p}_{\text{reverse area estimate}}
\;\times\;
\underbrace{p}_{\text{weighted square-function estimate}}
=p^2.
$$
See Section~\ref{sec:mult} for the proof and the dependence on \eqref{e:mult-snorm}.}

{For the lower bound we use the dilation-invariant symbol}
$$
m_\star(\xi)=
\begin{cases}
\displaystyle\frac{\xi_1\xi_2}{((\xi_1\xi_2)^2+\xi_3^2)^{1/2}},
 &(\xi_1\xi_2,\xi_3)\ne(0,0),\\
0, &(\xi_1\xi_2,\xi_3)=(0,0).
\end{cases}
$$
{It satisfies $\|m_\star\|_{\mathcal M_{Z,x_1}^N}<\infty$ for every positive integer $N$, and}
$$
cp^2\le\|T_{m_\star}\|_{L^p\to L^p}\le C_\star p^2,
\qquad 2\le p<\infty,
$$
with $c,C_\star$ independent of $p$.

To prove $\|T_{m_\star}\|_{L^p\to L^p}\ge cp^2$, we use the limit
$$
m_\star(\xi_1,\xi_2,\xi_3/L)
\to\operatorname{sgn}\xi_1\operatorname{sgn}\xi_2
\qquad(L\to\infty),\qquad \xi_1\xi_2\ne0.
$$
The limiting symbol is the multiplier of $-H_1H_2$,
where $H_j$ is the Hilbert transform in the $j$th variable.
Each Hilbert transform has operator norm of order $p$,
so their product has norm of order $p^2$.
Dilating test functions in the third variable and using Fatou's lemma,
we obtain the same lower bound for $T_{m_\star}$;
see Section~\ref{sec:mult}.

{
This paper is organized as follows.
We prove Theorem~\ref{sp:t:general} in Section~\ref{sp:sec:reverse} using the dual-kernel
criterion of Section~\ref{sp:sec:dual}. The proofs of Theorems~\ref{sp:t:kernels}
and~\ref{sp:t:aperture} are in Subsections~\ref{sp:sec:kernels} and~\ref{sp:sec:apertures}, respectively.
Sections~\ref{sec:dyadic} and~\ref{sec:atoms} develop the estimates used to prove
Theorem~\ref{t:primal} in Section~\ref{sec:remaining}.
We prove Theorems~\ref{t:large-ap} and~\ref{t:first-free} in Subsections~\ref{sub:large-ap}
and~\ref{sub:first-free}, Theorem~\ref{t:second-free} in Section~\ref{sec:second-free},
and Theorem~\ref{t:sharp-cases} in Subsection~\ref{sub:sharp-kernels}.
The proofs of Theorems~\ref{t:semigroup} and~\ref{t:multipliers} are in
Sections~\ref{sec:semigroups} and~\ref{sec:mult}, respectively.
}

\section{Calder\'on dual kernels and the support distance}\label{sp:sec:dual}

{Throughout this section, let $\varphi\in C_c^\infty(\R^d)$ be nonzero, real and radial,
with $\int_{\R^d}\varphi=0$. For this fixed kernel, we determine how small the support
of a dual kernel $\Phi$ in the Calder\'on reproducing formula \eqref{sp:e:reproduce} can be.}

\subsection{\texorpdfstring{$L^p$}{Lp} continuity and approximation}

Write $v_d=|B(0,1)|$ and $\sigma_{d-1}=|\mathbb S^{d-1}|$,
with $\sigma_0=2$.
{Let}
\begin{equation}\label{sp:e:energy-constant}
c_\varphi=\int_0^\infty|\widehat\varphi(re_1)|^2\,\frac{dr}{r}.
\end{equation}
{Since $\varphi$ has integral zero and is even, $\widehat\varphi(0)=0$ and
$\nabla\widehat\varphi(0)=0$. Hence $\widehat\varphi(\xi)=O(|\xi|^2)$ near zero.
Together with rapid decay at infinity, this proves $c_\varphi<\infty$.
Radiality and $\varphi\ne0$ imply $c_\varphi>0$.}

\begin{lemma}\label{sp:l:continuity}
For $f\in L^2(\R^d)$,
\begin{equation}\label{sp:e:energy}
\int_0^\infty\left\|f*\varphi_t\right\|_{L^2(\R^d)}^2\,\frac{dt}{t}=c_\varphi\left\|f\right\|_{L^2(\R^d)}^2,
\qquad
\left\|\mathcal A_{\varphi,\alpha}f\right\|_{L^2(\R^d)}^2=v_d\alpha^dc_\varphi\left\|f\right\|_{L^2(\R^d)}^2.
\end{equation}
For every $1<p<\infty$, there is $B_{p,d,\varphi,\alpha}<\infty$ such that
$$
\left\|\mathcal A_{\varphi,\alpha}f\right\|_{L^p(\R^d)}\le B_{p,d,\varphi,\alpha}\left\|f\right\|_{L^p(\R^d)}.
$$
In particular, $f_j\to f$ in $L^p$ implies
$\mathcal A_{\varphi,\alpha}f_j\to\mathcal A_{\varphi,\alpha}f$ in $L^p$.
\end{lemma}

\begin{proof}
For the first identity, we use Plancherel's theorem and the change of variables $r=t|\xi|$:
\begin{align*}
\int_0^\infty\left\|f*\varphi_t\right\|_{L^2(\R^d)}^2\,\frac{dt}{t}
= &\int_{\R^d}|\widehat f(\xi)|^2
\bigg(\int_0^\infty|\widehat\varphi(t\xi)|^2\,\frac{dt}{t}\bigg)\,d\xi
= c_\varphi\left\|f\right\|_{L^2(\R^d)}^2.
\end{align*}
For the area-function identity, we integrate first in $x$:
\begin{align*}
\left\|\mathcal A_{\varphi,\alpha}f\right\|_{L^2(\R^d)}^2
= &\int_0^\infty\int_{\R^d}|f*\varphi_t(y)|^2
\bigg(\int_{\R^d}\mathbf1_{\{|x-y|<\alpha t\}}\,dx\bigg)
\,\frac{dy\,dt}{t^{d+1}}\\
= &v_d\alpha^d\int_0^\infty\left\|f*\varphi_t\right\|_{L^2(\R^d)}^2\,\frac{dt}{t}
=v_d\alpha^dc_\varphi\left\|f\right\|_{L^2(\R^d)}^2.
\end{align*}

For the $L^p$ estimate, {write} $y=x+tu$, so that $|y-x|<\alpha t$ becomes $|u|<\alpha$
and $dy\,dt/t^{d+1}=du\,dt/t$.
Define the convolution kernel
$$
K(x;u,t)=t^{-d}\varphi\left(\frac{x}{t}+u\right),
\qquad |u|<\alpha,\quad t>0.
$$
For $f\in L^p(\R^d)$, we have
$$
(f*\varphi_t)(x+tu)
=\int_{\R^d}K(x-z;u,t)f(z)\,dz.
$$
With the substitution $r=|x|/t$, for $x\ne0$ and $|\gamma|=0,1$,  we obtain
\begin{align*}
\int_{|u|<\alpha}\int_0^\infty
|\partial_x^\gamma K(x;u,t)|^2\,\frac{dt}{t}\,du
&= |x|^{-2d-2|\gamma|}\int_{|u|<\alpha}\int_0^\infty
r^{2d+2|\gamma|-1}
\left|(\partial^\gamma\varphi)\left(r\frac{x}{|x|}+u\right)\right|^2\,dr\,du\\
&\le C_{d,\varphi,\alpha}|x|^{-2d-2|\gamma|}.
\end{align*}
Indeed, if $\operatorname{supp}\varphi\subset B(0,A)$, the integrand vanishes for $r>A+\alpha$.
For $0<2|z|\le|x|$, integrate the gradient along the segment from $x$ to $x-z$.
By Minkowski's inequality,
\begin{align*}
&\bigg(\int_{|u|<\alpha}\int_0^\infty
|K(x-z;u,t)-K(x;u,t)|^2\,\frac{dt}{t}\,du\bigg)^{1/2}\\
\le &|z|\int_0^1
\bigg(\int_{|u|<\alpha}\int_0^\infty
|\nabla_xK(x-\theta z;u,t)|^2\,\frac{dt}{t}\,du\bigg)^{1/2}\,d\theta\\
\le &C_{d,\varphi,\alpha}|z|\int_0^1\frac{d\theta}{|x-\theta z|^{d+1}}
\le C_{d,\varphi,\alpha}\frac{|z|}{|x|^{d+1}}.
\end{align*}
By \eqref{sp:e:energy}, the kernel estimates above and the
$L^2(du\,dt/t)$-valued Calder\'on--Zygmund theorem, we obtain, for $1<p<\infty$ and $f\in L^p(\R^d)$,
$$
\left\|\mathcal A_{\varphi,\alpha}f\right\|_{L^p(\R^d)}
\le B_{p,d,\varphi,\alpha}\left\|f\right\|_{L^p(\R^d)};
$$
see \cite[Chapter~I, Section~6.4]{Stein}.
For continuity, apply the triangle inequality in the same $L^2$ space pointwise:
$$
|\mathcal A_{\varphi,\alpha}f-\mathcal A_{\varphi,\alpha}g|
\le\mathcal A_{\varphi,\alpha}(f-g).
$$
Consequently,
$$
\left\|\mathcal A_{\varphi,\alpha}f-\mathcal A_{\varphi,\alpha}g\right\|_{L^p(\R^d)}
\le B_{p,d,\varphi,\alpha}\left\|f-g\right\|_{L^p(\R^d)}.
$$
The proof is complete.
\end{proof}

For the lower bounds, we may restrict the supremum in \eqref{sp:e:best} to
$0\ne f\in\mathcal S_0(\R^d)$, where
$$
\mathcal S_0(\R^d)
=\left\{g\in\mathcal S(\R^d):
\operatorname{supp}\widehat g\Subset\R^d\setminus\{0\}\right\}.
$$

For $1<p<\infty$, first approximate an $L^p(\R^d)$ function by $g\in\mathcal S(\R^d)$.
Choose $\vartheta\in C_c^\infty(\R^d)$ equal to one near zero, and define, for $0<\varepsilon<1$,
$$
\widehat g_\varepsilon(\xi)
=\vartheta(\varepsilon\xi)
\bigl(1-\vartheta(\xi/\varepsilon)\bigr)\widehat g(\xi).
$$
Then $g_\varepsilon\in\mathcal S_0(\R^d)$.
Convolution with $\varepsilon^{-d}\check\vartheta(\cdot/\varepsilon)$ approximates $g$ in $L^p$.
For the low-frequency term, we have by Young's inequality
$$
\left\|g*\bigl(\varepsilon^d\check\vartheta(\varepsilon\,\cdot)\bigr)\right\|_{L^p(\R^d)}
\le\varepsilon^{d(1-1/p)}
\left\|g\right\|_{L^1(\R^d)}\left\|\check\vartheta\right\|_{L^p(\R^d)}\to0.
$$
Combining the two cutoffs, we obtain
\begin{align*}
\left\|g_\varepsilon-g\right\|_{L^p(\R^d)}
\le &\left\|g*\bigl(\varepsilon^{-d}\check\vartheta(\cdot/\varepsilon)\bigr)-g
\right\|_{L^p(\R^d)}+\left\|\check\vartheta\right\|_{L^1(\R^d)}
\left\|g*\bigl(\varepsilon^d\check\vartheta(\varepsilon\,\cdot)\bigr)\right\|_{L^p(\R^d)}
\to0.
\end{align*}
Thus $\mathcal S_0(\R^d)$ is dense in $L^p(\R^d)$.

To pass to the norm ratios, we first show that
$\|\mathcal A_{\varphi,\alpha}f\|_{L^p(\R^d)}>0$ for $f\ne0$.
Otherwise, integrating the square of the area function first in $x$, as above, gives
$$
0=v_d\alpha^d\int_0^\infty\int_{\R^d}|f*\varphi_t(y)|^2\,dy\,\frac{dt}{t}.
$$
Since $(y,t)\mapsto f*\varphi_t(y)$ is continuous, $f*\varphi_t=0$ for every $t>0$.
Taking Fourier transforms in the sense of tempered distributions, we obtain
$$
\operatorname{supp}\widehat f
\subset\bigcap_{t>0}\left\{\xi:\widehat\varphi(t\xi)=0\right\}=\{0\},
$$
where the last equality uses radiality, $\varphi\ne0$ and $\widehat\varphi(0)=0$.
Thus $f$ is a polynomial in $L^p(\R^d)$ and must vanish.

For $2\le p<\infty$ and $0\ne f\in L^p(\R^d)$,
choose $0\ne f_j\in\mathcal S_0(\R^d)$ with $f_j\to f$ in $L^p$.
By Lemma~\ref{sp:l:continuity},
$$
\frac{\|f_j\|_{L^p(\R^d)}}{\|\mathcal A_{\varphi,\alpha}f_j\|_{L^p(\R^d)}}
\to
\frac{\|f\|_{L^p(\R^d)}}{\|\mathcal A_{\varphi,\alpha}f\|_{L^p(\R^d)}}.
$$
This proves the assertion about the supremum in \eqref{sp:e:best}.

\subsection{Localized dual kernels}

{For every integer $j\geq0$, polar coordinates give
$$
\int_{|\xi|<1}|\xi|^{j+2-d}|\widehat\varphi(\xi)|\,d\xi
\leq \sigma_{d-1}\|\widehat\varphi\|_\infty\int_0^1 r^{j+1}\,dr<\infty.
$$
Since $\widehat\varphi$ decays rapidly at infinity, the Fourier integral for
$|D|^{2-d}\varphi$ and all its derivatives converge absolutely.
Thus $|D|^{2-d}\varphi$ is smooth, real and radial.}

For a real radial $\chi\in C_c^\infty(B(0,\rho))$, the function $\Phi=-\Delta\chi$
has support in $B(0,\rho)$ and integral zero.
We now determine when $\Phi$ can be normalized.

\begin{proposition}\label{sp:p:local-dual}
Fix $\rho>0$.
There is a real radial $\Phi\in C_c^\infty(B(0,\rho))$ satisfying
\begin{equation}\label{sp:e:dual}
\int_{\R^d}\Phi=0,
\qquad
\int_0^\infty\widehat\varphi(re_1)\widehat\Phi(re_1)\,\frac{dr}{r}=1
\end{equation}
if and only if $\operatorname{dist}\bigl(0,\operatorname{supp}(|D|^{2-d}\varphi)\bigr)<\rho$.
\end{proposition}

\begin{proof}
For a real radial $\chi\in C_c^\infty(\R^d)$, write the normalizing integral
in polar coordinates and then use Fourier inversion:
\begin{align}\label{sp:e:pair}
\int_0^\infty\widehat\varphi(re_1)\widehat{-\Delta\chi}(re_1)\,\frac{dr}{r}
= &4\pi^2\int_0^\infty r\widehat\varphi(re_1)\widehat\chi(re_1)\,dr\notag\\
= &\frac{4\pi^2}{\sigma_{d-1}}\int_{\R^d}
|\xi|^{2-d}\widehat\varphi(\xi)\widehat\chi(\xi)\,d\xi\notag\\
= &\frac{4\pi^2}{\sigma_{d-1}}
\int_{\R^d}(|D|^{2-d}\varphi)(x)\chi(x)\,dx.
\end{align}
{The second equality uses radiality and $r=r^{2-d}r^{d-1}$.
The last follows by Fubini's theorem and $\widehat\chi(-\xi)=\widehat\chi(\xi)$.}

Suppose first that such a $\Phi$ exists.
We shall write it as $\Phi=-\Delta\chi$, with $\chi$ real, radial and supported
in the same ball, so that \eqref{sp:e:pair} also proves necessity.
Define
$$
\chi(x)=\int_{|x|}^\infty r^{1-d}
\bigg(\int_0^r t^{d-1}\Phi(te_1)\,dt\bigg)\,dr.
$$
Choose $r_0<\rho$ with $\operatorname{supp}\Phi\subset B(0,r_0)$.
By radiality and the zero integral,
$$
\int_0^r t^{d-1}\Phi(te_1)\,dt
=\frac1{\sigma_{d-1}}\int_{\R^d}\Phi=0
\qquad\text{for }r\ge r_0,
$$
hence $\chi(x)=0$ for $|x|\ge r_0$.
Since
$$
\left|\int_0^r t^{d-1}\Phi(te_1)\,dt\right|
\le\frac{r^d}{d}\|\Phi\|_{L^\infty(\R^d)},
$$
the integral defining $\chi(0)$ is finite.
By the radiality of $\Phi$,
\begin{align*}
\chi(x)-\chi(0)
= &-\int_0^{|x|}r^{1-d}
\bigg(\int_0^r\tau^{d-1}\Phi(\tau e_1)\,d\tau\bigg)\,dr
= -\int_0^{|x|}r\bigg(\int_0^1v^{d-1}\Phi(rv e_1)\,dv\bigg)\,dr\\
= &-|x|^2\int_0^1\int_0^1tv^{d-1}\Phi(tv|x|e_1)\,dv\,dt
=-|x|^2\int_0^1\int_0^1tv^{d-1}\Phi(tv x)\,dv\,dt.
\end{align*}
Thus $\chi\in C_c^\infty(B(0,\rho))$ is real and radial as well.
For $x\ne0$, taking $r=t|x|$, we have
\begin{align*}
\nabla\chi(x)
= &-\frac{x}{|x|^d}\int_0^{|x|}r^{d-1}\Phi(re_1)\,dr
= -x\int_0^1t^{d-1}\Phi(t|x|e_1)\,dt
=-x\int_0^1t^{d-1}\Phi(tx)\,dt.
\end{align*}
Consequently,
\begin{align*}
-\Delta\chi(x)
= &\sum_{j=1}^d\partial_{x_j}
\bigg(x_j\int_0^1t^{d-1}\Phi(tx)\,dt\bigg)
= \int_0^1\bigl[dt^{d-1}\Phi(tx)+t^dx\cdot\nabla\Phi(tx)\bigr]\,dt\\
= &\int_0^1\frac{d}{dt}\bigl[t^d\Phi(tx)\bigr]\,dt=\Phi(x).
\end{align*}
Both identities also hold at $x=0$ by continuity.

The identity in \eqref{sp:e:pair} now becomes
$$
1=\int_0^\infty\widehat\varphi(re_1)\widehat{-\Delta\chi}(re_1)\,\frac{dr}{r}
=\frac{4\pi^2}{\sigma_{d-1}}\int_{\R^d}(|D|^{2-d}\varphi)\chi.
$$
Hence $|D|^{2-d}\varphi$ cannot vanish on $B(0,\rho)$,
and so $\operatorname{dist}\bigl(0,\operatorname{supp}(|D|^{2-d}\varphi)\bigr)<\rho$.

Conversely, if $\operatorname{dist}\bigl(0,\operatorname{supp}(|D|^{2-d}\varphi)\bigr)<\rho$, choose a nonnegative radial
$\zeta\in C_c^\infty(B(0,\rho))$ with
$\int_{\R^d}\zeta(|D|^{2-d}\varphi)^2>0$.
The function
$$
\Phi
=-\frac{\sigma_{d-1}\Delta\bigl(\zeta\,|D|^{2-d}\varphi\bigr)}
{4\pi^2\int_{\R^d}\zeta(|D|^{2-d}\varphi)^2}
$$
satisfies all the requirements, with normalization verified by \eqref{sp:e:pair}.
\end{proof}

{
\begin{remark}\label{sp:r:vanishing-moments}
One may ask whether \eqref{sp:e:pair} requires $\widehat\varphi$
to vanish to higher order at the origin.
The radial factor is $r\,dr$, so boundedness of $\widehat\varphi$
near zero suffices for convergence there.
In fact, evenness and mean zero already give
$\widehat\varphi(\xi)=O(|\xi|^2)$ as $\xi\to0$.

Higher-order vanishing does occur later in the paper.
In Lemma~\ref{sp:l:log-kernel}, the factor $|\xi|^{2m}$ comes from
the choice $\Omega=\Delta^m\zeta$.
In Section~\ref{sec:mult}, we choose the auxiliary kernels $\sigma_1,\sigma_2$
with the vanishing moments in \eqref{e:mult-moments}
to prove the pointwise estimate \eqref{e:mult-ptw}.
\end{remark}
}

For a dual kernel satisfying \eqref{sp:e:dual}, write $\Phi_t(x)=t^{-d}\Phi(x/t)$.
Radiality implies
$$
\int_0^\infty\widehat\varphi(t\xi)\widehat\Phi(t\xi)\,\frac{dt}{t}=1
\qquad(\xi\ne0).
$$
It follows that
\begin{equation}\label{sp:e:reproduce}
f=\int_0^\infty(f*\varphi_t)*\Phi_t\,\frac{dt}{t}
\quad\text{in }L^2(\R^d).
\end{equation}
For $\xi\ne0$ and $0<\varepsilon<R<\infty$,  it follows from Cauchy--Schwarz that
$$
\bigg|\int_\varepsilon^R\widehat\varphi(t\xi)\widehat\Phi(t\xi)\,\frac{dt}{t}\bigg|
\le\bigg(\int_0^\infty|\widehat\varphi(t\xi)|^2\,\frac{dt}{t}\bigg)^{1/2}
\bigg(\int_0^\infty|\widehat\Phi(t\xi)|^2\,\frac{dt}{t}\bigg)^{1/2}
=\sqrt{c_\varphi c_\Phi}.
$$
As $\varepsilon\to0$ and $R\to\infty$, this multiplier converges to one for $\xi\ne0$.
By Plancherel's theorem and dominated convergence,
$$
\bigg\|f-\int_\varepsilon^R(f*\varphi_t)*\Phi_t\,\frac{dt}{t}\bigg\|_{L^2(\R^d)}^2
=\int_{\R^d}|\widehat f(\xi)|^2
\bigg|1-\int_\varepsilon^R\widehat\varphi(t\xi)\widehat\Phi(t\xi)\,\frac{dt}{t}\bigg|^2\,d\xi
\to0.
$$

\begin{remark}
A normalized Calder\'on dual kernel always exists: one may take $\Phi=\varphi/c_\varphi$.
The condition $0\in\operatorname{supp}(|D|^{2-d}\varphi)$ allows us to choose such kernels
with arbitrarily small support.

For a fixed radius $\rho$, Proposition~\ref{sp:p:local-dual} requires only
$$
\operatorname{dist}\bigl(0,\operatorname{supp}(|D|^{2-d}\varphi)\bigr)<\rho.
$$
At scale $t$, the support of $\Phi_t$ is contained in $B(0,\rho t)$, whereas
the cone has cross-section $B(x,\alpha t)$.
The dyadic covering in {\it Step 1} of Proposition~\ref{sp:p:reverse} introduces
the factor $16\sqrt d$.
Taking $16\sqrt d\,\rho<\alpha$ then permits the pointwise comparison
\eqref{sp:e:discrete-area} with the area function at aperture $\alpha$.
Proposition~\ref{sp:p:reverse} therefore establishes the $\sqrt p$ estimate under the condition
$$
\operatorname{dist}\bigl(0,\operatorname{supp}(|D|^{2-d}\varphi)\bigr)<c_d\alpha.
$$
This holds for sufficiently large $\alpha$, and for every $\alpha>0$ if
$0\in\operatorname{supp}(|D|^{2-d}\varphi)$.
\end{remark}

\begin{remark}\label{sp:r:cww-dual}
The same operator occurs in the local reconstruction of Chang, Wilson and Wolff
\cite[Lemma~3.1]{CWW}.
Their area function is formed from the vector kernel $\nabla k$, where $k$ is the Poisson kernel.
They choose a real radial $K_0$ with small support such that
$$
\int_{\R^d}K_0(x)(-\Delta)^{1-d/2}k(x)\,dx\ne0.
$$
With our Fourier convention, $(-\Delta)^{1-d/2}=(2\pi)^{2-d}|D|^{2-d}$.
Indeed, for $\xi\ne0$, radiality gives
\begin{align*}
\int_0^\infty\sum_{j=1}^d
\widehat{\partial_j k}(t\xi)\widehat{\partial_jK_0}(t\xi)\,\frac{dt}{t}
= &-4\pi^2\int_0^\infty r\widehat k(re_1)\widehat K_0(re_1)\,dr
= -\frac{4\pi^2}{\sigma_{d-1}}\int_{\R^d}(|D|^{2-d}k)(x)K_0(x)\,dx.
\end{align*}
Thus their gradient construction and our scalar construction in \eqref{sp:e:pair}
use the same local pairing, applied to $k$ and $\varphi$, respectively.

For the Poisson kernel, $\widehat k(\xi)=e^{-2\pi|\xi|}$, and
$$
(|D|^{2-d}k)(0)
=\sigma_{d-1}\int_0^\infty re^{-2\pi r}\,dr>0.
$$
Hence the pairing can be normalized with a real radial
$K_0\in C_c^\infty(B(0,\rho))$ for every $\rho>0$.
Such a choice is not automatic for the more general nonnegative radial
Schwartz kernels mentioned in \cite[p.~239]{CWW}.
In dimension two, for example, $|D|^{2-d}k=k$.
If $k$ vanishes on $B(0,\rho)$, then
$$
\int_{\R^2}k(x)K_0(x)\,dx=0
\qquad\text{for every }K_0\in C_c^\infty(B(0,\rho)),
$$
which prevents normalization within this ball.
We return to the Poisson gradient area function in Subsection~\ref{sp:sub:cww}.
\end{remark}

\section{One-parameter reverse estimates}\label{sp:sec:reverse}

Our starting point is the comparison between the discrete square function
associated with an atomic decomposition and the Lusin area function,
as in our earlier work \cite{CDWY,WuJGA,LPW}.
Here the condition
$$
\operatorname{dist}\bigl(0,\operatorname{supp}(|D|^{2-d}\varphi)\bigr)<c_d\alpha
$$
allows a Calder\'on dual kernel with support small enough to obtain
the pointwise comparison \eqref{sp:e:discrete-area}.
{We expand the smooth atoms in Haar functions on dyadic cubes containing
their supports, as in \cite[Proposition~3.6]{WuJGA}. At each relative generation
we use the Haar estimate from \cite[Proposition~1.5]{LPW}, then sum the decay
from the smoothness of the atoms.}

\subsection{Haar functions and smooth atoms}

Let $\mathcal D$ be a dyadic grid of cubes in $\R^d$.
For a dyadic interval $J$, write
$$
h_J^0=|J|^{-1/2}\mathbf1_J,
\qquad
h_J^1=|J|^{-1/2}(\mathbf1_{J_-}-\mathbf1_{J_+}),
$$
where $J_-$ and $J_+$ are its two children.
For $Q=J_1\times\cdots\times J_d$ and $\eta\in\{0,1\}^d\setminus\{0\}$, define
$$
h_Q^\eta(x)=\prod_{j=1}^dh_{J_j}^{\eta_j}(x_j).
$$
These functions are constant on the children of $Q$, with
$$
\int_{\R^d}h_Q^\eta=0,
\qquad
|h_Q^\eta|=|Q|^{-1/2}\mathbf1_Q,
\qquad
\left\|h_Q^\eta\right\|_{L^2(\R^d)}=1.
$$

\begin{lemma}\label{sp:l:haar}
For every finitely supported family of scalar coefficients $\{b_{Q,\eta}\}$ and every $2\le p<\infty$,
\begin{equation}\label{sp:e:haar}
\bigg\|\sum_{Q\in\mathcal D}\sum_{\eta\ne0}b_{Q,\eta}h_Q^\eta\bigg\|_{L^p(\R^d)}
\le C_d\sqrt p\,
\bigg\|\bigg(\sum_{Q\in\mathcal D}\sum_{\eta\ne0}|b_{Q,\eta}|^2
\frac{\mathbf1_Q}{|Q|}\bigg)^{1/2}\bigg\|_{L^p(\R^d)}.
\end{equation}
\end{lemma}

\begin{proof}
First suppose that the coefficients are real, and denote
$$
A_Q=\sum_{\eta\ne0}b_{Q,\eta}h_Q^\eta.
$$
Each $A_Q$ is supported in $Q$ and constant on its children.
For every dyadic cube $R$ and every subfamily $\mathscr E\subset\{Q\in\mathcal D:Q\subseteq R\}$,
orthogonality of the Haar functions implies
\begin{align*}
\frac1{|R|}\int_R\bigg|\sum_{Q\in\mathscr E}A_Q(x)\bigg|^2\,dx
= &\frac1{|R|}\sum_{Q\in\mathscr E}\sum_{\eta\ne0}|b_{Q,\eta}|^2
= \frac1{|R|}\int_R\sum_{Q\in\mathscr E}|A_Q(x)|^2\,dx
\le \bigg\|\sum_{Q\in\mathscr E}|A_Q|^2\bigg\|_{L^\infty(R)}.
\end{align*}
Thus $\{A_Q\}$ satisfies (P1), (P2) and (P3) of \cite[Theorem~A]{LPW}, with the constant in (P3) equal to one.
For each maximal cube $R$ in the finite family, we apply
\cite[Proposition~1.5]{LPW} to $\{A_Q:Q\subseteq R\}$ and obtain
$$
\bigg\|\sum_{Q\subseteq R}A_Q\bigg\|_{L^p(R)}
\le C_d\sqrt p\,
\bigg\|\bigg(\sum_{Q\subseteq R}|A_Q|^2\bigg)^{1/2}\bigg\|_{L^p(R)}.
$$
Summing the $p$th powers over the disjoint maximal cubes, we obtain
$$
\bigg\|\sum_QA_Q\bigg\|_{L^p(\R^d)}
\le C_d\sqrt p\,\bigg\|\bigg(\sum_Q|A_Q|^2\bigg)^{1/2}\bigg\|_{L^p(\R^d)}.
$$
Finally, Cauchy--Schwarz in $\eta$ shows that
$$
|A_Q(x)|^2
\le (2^d-1)\sum_{\eta\ne0}|b_{Q,\eta}|^2|h_Q^\eta(x)|^2
=(2^d-1)\sum_{\eta\ne0}|b_{Q,\eta}|^2\frac{\mathbf1_Q(x)}{|Q|}.
$$
This proves \eqref{sp:e:haar} for real coefficients.
For complex coefficients, apply this estimate to their real and imaginary parts and use the triangle inequality.
\end{proof}

The cubes $I$ index the pieces of the Calder\'on decomposition below,
whose smooth atoms are supported in $5I$.
Here $cI$ denotes the closed concentric cube with side length $c\ell(I)$.
For the Haar expansion, we need a dyadic cube containing $5I$ and comparable in size to $I$.
Such a cube may not exist in the original grid, so we use adjacent dyadic systems.

Consider
$$
\mathcal D^{\pmb\varepsilon}
=\{2^{-k}([0,1)^d+m+(-1)^k\pmb\varepsilon/3):
k\in\mathbb Z,\ m\in\mathbb Z^d\},
\qquad \pmb\varepsilon\in\{0,1,2\}^d.
$$
The choice $\pmb\varepsilon=0$ gives the standard system $\mathcal D^0$.
Each system is nested, since the shift between consecutive scales,
measured in units of the smaller side length, is an integer vector:
$$
2(-1)^k\pmb\varepsilon/3-(-1)^{k+1}\pmb\varepsilon/3
=(-1)^k\pmb\varepsilon\in\mathbb Z^d.
$$

For each $I\in\mathcal D^0$, we can choose and fix
$Q(I)$ and $\pmb\varepsilon(I)$ such that
\begin{equation}\label{sp:e:cover}
5I\subset Q(I)\in\mathcal D^{\pmb\varepsilon(I)},
\qquad \ell(Q(I))=16\ell(I),
\qquad
\sup_Q\#\{I\in\mathcal D^0:Q(I)=Q\}\le16^d.
\end{equation}
Indeed, at side length $16\ell(I)$, the boundary points of distinct
one-dimensional grids are separated by at least $16\ell(I)/3>5\ell(I)$.
Each coordinate projection of $5I$ therefore avoids the boundaries of at least two grids.
Choose a containing interval from one of these grids in each coordinate;
their product is the required $Q(I)$.

For a fixed $Q$, the cubes $I$ with $Q(I)=Q$ are disjoint, lie in $Q$,
and have volume $16^{-d}|Q|$. Hence
$$
16^{-d}|Q|\,\#\{I\in\mathcal D^0:Q(I)=Q\}
=\sum_{\substack{I\in\mathcal D^0\\Q(I)=Q}}|I|
\le |Q|,
$$
which proves the multiplicity bound in \eqref{sp:e:cover}.

Notably, in the Zygmund arguments below, it suffices to use only two adjacent dyadic grids on $\R$;
see Lemma~\ref{l:adj-int-cover}.

We apply the Haar estimate separately in each selected grid.
For a fixed $\pmb\varepsilon$, we consider finite families
$$
\mathscr I\subset
\{I\in\mathcal D^0:\pmb\varepsilon(I)=\pmb\varepsilon\},
\qquad
Q(I)\in\mathcal D^{\pmb\varepsilon}\quad(I\in\mathscr I).
$$
The next lemma estimates the corresponding atomic sum.
Its parameter $M$ bounds the number of atoms with the same chosen supporting cube;
in \eqref{sp:e:cover}, we may take $M=16^d$.

\begin{lemma}\label{sp:l:smooth}
Let $\mathcal D$ be a dyadic system on $\R^d$, and let
$\mathscr I\subset\mathcal D^0$ be a finite family of dyadic cubes.
For each $I\in\mathscr I$, let $\mathsf a_I\in C_c^1(\R^d)$
and $Q(I)\in\mathcal D$ satisfy
$$
I\subset Q(I),
\qquad
\operatorname{supp}\mathsf a_I\subset Q(I),
\qquad
\int_{\R^d}\mathsf a_I=0.
$$
Suppose that, for some $M\ge1$ and $K>0$,
$$
\sup_{Q\in\mathcal D}\#\{I\in\mathscr I:Q(I)=Q\}\le M
$$
and
$$
\left\|\mathsf a_I\right\|_{L^\infty(\R^d)}
+\ell(Q(I))\left\|\nabla\mathsf a_I\right\|_{L^\infty(\R^d)}
\le K|Q(I)|^{-1/2},
\qquad I\in\mathscr I.
$$
Then, for all scalar coefficients $\{\lambda_I\}_{I\in\mathscr I}$
and $2\le p<\infty$,
\begin{equation}\label{sp:e:smooth}
\bigg\|\sum_{I\in\mathscr I}\lambda_I\mathsf a_I\bigg\|_{L^p(\R^d)}
\le C_dK\sqrt{Mp}\,
\bigg\|\bigg(\sum_{I\in\mathscr I}|\lambda_I|^2
\frac{\mathbf1_{Q(I)}}{|Q(I)|}\bigg)^{1/2}\bigg\|_{L^p(\R^d)},
\end{equation}
where $C_d$ depends only on $d$.
\end{lemma}

\begin{proof}
For $Q\in\mathcal D$ and $n\in\mathbb N_0$, denote
$$
Ch_n(Q)=\{P\in\mathcal D:P\subseteq Q,\ \ell(P)=2^{-n}\ell(Q)\}.
$$
Here $n$ is the generation relative to $Q$.
The derivative bound will supply a factor $2^{-n}$, so that summing over the
generations will not introduce another power of $p$.
Since $\mathsf a_I$ is supported in $Q(I)$ and has integral zero,
its Haar expansion is
$$
\mathsf a_I
=\sum_{n=0}^\infty\sum_{P\in Ch_n(Q(I))}\sum_{\eta\ne0}
\langle\mathsf a_I,h_P^\eta\rangle h_P^\eta
\quad\text{in }L^2(\R^d).
$$
For $P\in Ch_n(Q(I))$, let $c_P$ be its centre. Note that
\begin{align}\label{sp:e:haar-coeff}
|\langle\mathsf a_I,h_P^\eta\rangle|
= &\bigg|\int_P\bigl(\mathsf a_I(x)-\mathsf a_I(c_P)\bigr)h_P^\eta(x)\,dx\bigg|\notag\\
\le &\left\|\nabla\mathsf a_I\right\|_{L^\infty(\R^d)}
\int_P|x-c_P||h_P^\eta(x)|\,dx\notag\\
\le &C_dK\,\frac{\ell(P)}{\ell(Q(I))}
\bigg(\frac{|P|}{|Q(I)|}\bigg)^{1/2}
=C_dK\,2^{-n}\bigg(\frac{|P|}{|Q(I)|}\bigg)^{1/2}.
\end{align}
Consequently,
\begin{align*}
\sum_{P\in Ch_n(Q(I))}\sum_{\eta\ne0}
|\langle\mathsf a_I,h_P^\eta\rangle h_P^\eta(x)|
\le &C_dK\,2^{-n}|Q(I)|^{-1/2}
\sum_{P\in Ch_n(Q(I))}\mathbf1_P(x)
= C_dK\,2^{-n}\frac{\mathbf1_{Q(I)}(x)}{|Q(I)|^{1/2}}.
\end{align*}
Thus the Haar series converges absolutely and uniformly to its $L^2$ limit $\mathsf a_I$.
For every $2\le p<\infty$,
$$
\bigg\|\mathsf a_I-
\sum_{n=0}^{N-1}\sum_{P\in Ch_n(Q(I))}\sum_{\eta\ne0}
\langle\mathsf a_I,h_P^\eta\rangle h_P^\eta\bigg\|_{L^p(\R^d)}
\le C_dK\,2^{-N}|Q(I)|^{1/p-1/2}\to0.
$$

For fixed $n\in\mathbb N_0$, $P\in\mathcal D$ and $\eta\ne0$, define
$$
b_{P,\eta,n}
=\sum_{I\in\mathscr I:\,P\in Ch_n(Q(I))}
\lambda_I\langle\mathsf a_I,h_P^\eta\rangle.
$$
Because all the cubes $Q(I)$ belong to the same system $\mathcal D$,
the condition $P\in Ch_n(Q(I))$ determines $Q(I)$ uniquely:
it is the $n$th dyadic ancestor of $P$ in this system.
At most $M$ atoms can therefore contribute to this Haar coefficient, and
$$
\#\left\{I\in\mathscr I:P\in Ch_n(Q(I))\right\}\le M.
$$
By Cauchy--Schwarz and \eqref{sp:e:haar-coeff},
\begin{align}\label{sp:e:depth-sq}
\sum_{P\in\mathcal D}\sum_{\eta\ne0}|b_{P,\eta,n}|^2\frac{\mathbf1_P}{|P|}
\le &M\sum_{I\in\mathscr I}|\lambda_I|^2
\sum_{P\in Ch_n(Q(I))}\sum_{\eta\ne0}
|\langle\mathsf a_I,h_P^\eta\rangle|^2\frac{\mathbf1_P}{|P|}\notag\\
\le &C_dMK^2\,2^{-2n}\sum_{I\in\mathscr I}
\frac{|\lambda_I|^2}{|Q(I)|}\sum_{P\in Ch_n(Q(I))}\mathbf1_P\notag\\
= &C_dMK^2\,2^{-2n}\sum_{I\in\mathscr I}|\lambda_I|^2
\frac{\mathbf1_{Q(I)}}{|Q(I)|}.
\end{align}
Since $\mathscr I$ is finite, we may sum these expansions in $L^p$.
Applying Lemma~\ref{sp:l:haar} at each generation and then \eqref{sp:e:depth-sq}, we obtain
\begin{align*}
\bigg\|\sum_{I\in\mathscr I}\lambda_I\mathsf a_I\bigg\|_{L^p(\R^d)}
\le &\sum_{n=0}^\infty
\bigg\|\sum_{P\in\mathcal D}\sum_{\eta\ne0}b_{P,\eta,n}h_P^\eta\bigg\|_{L^p(\R^d)}\\
\le &C_d\sqrt p\sum_{n=0}^\infty
\bigg\|\bigg(\sum_{P\in\mathcal D}\sum_{\eta\ne0}
|b_{P,\eta,n}|^2\frac{\mathbf1_P}{|P|}\bigg)^{1/2}\bigg\|_{L^p(\R^d)}\\
\le &C_dK\sqrt{Mp}\sum_{n=0}^\infty2^{-n}
\bigg\|\bigg(\sum_{I\in\mathscr I}|\lambda_I|^2
\frac{\mathbf1_{Q(I)}}{|Q(I)|}\bigg)^{1/2}\bigg\|_{L^p(\R^d)}.
\end{align*}
Summing the geometric series proves \eqref{sp:e:smooth}.
\end{proof}

\subsection{The reverse estimate under the distance condition}

\begin{proposition}\label{sp:p:reverse}
There is $c_d>0$ depending only on the dimension such that, whenever
$$
\operatorname{dist}\bigl(0,\operatorname{supp}(|D|^{2-d}\varphi)\bigr)<c_d\alpha,
$$
we have
$$
\left\|f\right\|_{L^p(\R^d)}
\le C_{d,\varphi,\alpha}\sqrt p\,
\left\|\mathcal A_{\varphi,\alpha}f\right\|_{L^p(\R^d)}
$$
for every $2\le p<\infty$ and $f\in L^p(\R^d)$.
\end{proposition}

\begin{proof}
We may take $c_d=(32\sqrt d)^{-1}$; the factor $\sqrt p$ will come from Lemma~\ref{sp:l:haar}.
Fix $\alpha>0$ with $\operatorname{dist}\bigl(0,\operatorname{supp}(|D|^{2-d}\varphi)\bigr)<c_d\alpha$, let $\rho=c_d\alpha$,
and choose the dual kernel $\Phi$ from Proposition~\ref{sp:p:local-dual}.
The constants below may depend on $d,\varphi,\alpha$, but not on $p,f$ or the finite families of cubes.
First suppose $f\in L^2(\R^d)\cap L^p(\R^d)$.

\smallskip
\noindent{\it Step 1.} We decompose the Calder\'on reproducing formula for $\varphi$ into smooth atoms.

At scale $t$, the dual kernel $\Phi_t$ has support radius at most $\rho t$.
We therefore choose the spatial cubes so that $\ell(I)$ is comparable to $\rho t$.
For the standard dyadic system $\mathcal D^0$, let
$$
\mathcal W(I)=I\times[\ell(I)/\rho,2\ell(I)/\rho),
\qquad I\in\mathcal D^0.
$$
At each fixed scale the cubes $I$ partition $\R^d$, and the intervals
$[\ell(I)/\rho,2\ell(I)/\rho)$ for the different scales partition $(0,\infty)$.
Hence
$$
\R^d\times(0,\infty)=\bigcup_{I\in\mathcal D^0}\mathcal W(I)
\quad\text{(disjoint)}.
$$
Define
$$
\lambda_I(f)=\bigg(\int_{\mathcal W(I)}
|f*\varphi_t(y)|^2\,dy\,\frac{dt}{t}\bigg)^{1/2}
$$
and
$$
\mathsf a_I(f)(x)=\frac1{\lambda_I(f)}
\int_{\mathcal W(I)}(f*\varphi_t)(y)\Phi_t(x-y)\,dy\,\frac{dt}{t},
$$
with $\mathsf a_I(f)=0$ when $\lambda_I(f)=0$.
The Calder\'on reproducing formula has the atomic form
\begin{align}\label{sp:e:atomic-repro}
f(x)
= &\int_0^\infty\int_{\R^d}(f*\varphi_t)(y)\Phi_t(x-y)\,dy\,\frac{dt}{t}
= \sum_{I\in\mathcal D^0}\lambda_I(f)\mathsf a_I(f)(x),
\end{align}
with convergence in $L^2(\R^d)$ proved in {\it Step 3}.

For $(y,t)\in\mathcal W(I)$, we have $y\in I$ and $\rho t<2\ell(I)$.
Since $\operatorname{supp}\Phi\subset B(0,\rho)$,
the factor $\Phi_t(x-y)$ vanishes for $x\notin5I$.
Integrating first in $x$ and using $\int_{\R^d}\Phi=0$, we obtain
$$
\operatorname{supp}\mathsf a_I(f)\subset5I,
\qquad
\int_{\R^d}\mathsf a_I(f)=0.
$$
For a multi-index $\gamma$ with $|\gamma|\le1$, by Cauchy--Schwarz, we have 
\begin{align}\label{sp:e:atom-derivative}
|\partial^\gamma\mathsf a_I(f)(x)|
\le &\bigg(\int_{\ell(I)/\rho}^{2\ell(I)/\rho}
\left\|\partial^\gamma\Phi_t\right\|_{L^2(\R^d)}^2\,\frac{dt}{t}\bigg)^{1/2}
= \left\|\partial^\gamma\Phi\right\|_{L^2(\R^d)}
\bigg(\int_{\ell(I)/\rho}^{2\ell(I)/\rho}t^{-d-2|\gamma|}\,\frac{dt}{t}\bigg)^{1/2}\\
\le & C\ell(I)^{-d/2-|\gamma|}.\notag
\end{align}

For each $I$, use the cube $Q(I)\in\mathcal D^{\pmb\varepsilon(I)}$
fixed in \eqref{sp:e:cover} for the Haar expansion of $\mathsf a_I(f)$.
Since $\ell(Q(I))=16\ell(I)$, \eqref{sp:e:atom-derivative} gives
$$
\left\|\mathsf a_I(f)\right\|_{L^\infty(\R^d)}
+\ell(Q(I))\left\|\nabla\mathsf a_I(f)\right\|_{L^\infty(\R^d)}
\le C|Q(I)|^{-1/2}.
$$

To apply the Haar estimate separately in each system, define, for
$\pmb\varepsilon\in\{0,1,2\}^d$,
$$
\mathscr A_{\pmb\varepsilon}
=\{I\in\mathcal D^0:\pmb\varepsilon(I)=\pmb\varepsilon\}.
$$
Thus
$$
\mathcal D^0=\bigcup_{\pmb\varepsilon\in\{0,1,2\}^d}\mathscr A_{\pmb\varepsilon}
\quad\text{(disjoint)},
\qquad
Q(I)\in\mathcal D^{\pmb\varepsilon}
\quad(I\in\mathscr A_{\pmb\varepsilon}).
$$

Since the discrete square function uses $\mathbf1_{Q(I)}$, we need the cone
inclusion for every $x\in Q(I)$; see Figure~\ref{fig:dual-cone}.
For $x\in Q(I)$ and $(y,t)\in\mathcal W(I)$, \eqref{sp:e:cover} shows that
$$
|x-y|\le16\sqrt d\,\ell(I)\le16\sqrt d\,\rho t<\alpha t,
\qquad
\frac1{|Q(I)|}\le C_d\rho^{-d}t^{-d}.
$$
Since the regions $\mathcal W(I)$ are disjoint, we obtain the pointwise comparison
\begin{align}\label{sp:e:discrete-area}
\sum_{I\in\mathcal D^0}\lambda_I(f)^2\frac{\mathbf1_{Q(I)}(x)}{|Q(I)|}
\le &C_d\rho^{-d}\sum_{I\in\mathcal D^0}\mathbf1_{Q(I)}(x)
\int_{\mathcal W(I)}|f*\varphi_t(y)|^2\,\frac{dy\,dt}{t^{d+1}}\notag\\
\le &C_d\rho^{-d}\int_0^\infty\int_{|y-x|<\alpha t}
|f*\varphi_t(y)|^2\,\frac{dy\,dt}{t^{d+1}}
=C_d\rho^{-d}\mathcal A_{\varphi,\alpha}f(x)^2.
\end{align}

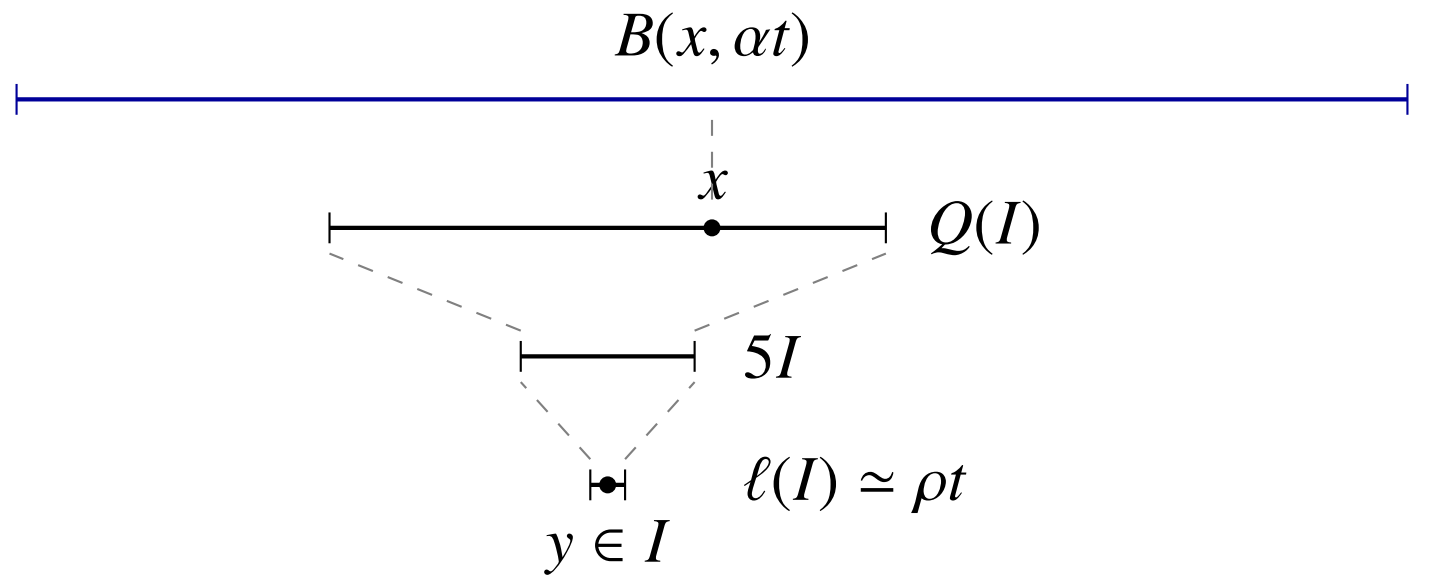
\begin{figure}[tbp]
\centering
\begin{tikzpicture}[x=1.15cm,y=0.85cm,>=Stealth]
\draw[thick,blue!60!black] (-3.4,3)--(4.6,3);
\draw[blue!60!black] (-3.4,2.88)--(-3.4,3.12) (4.6,2.88)--(4.6,3.12);
\node[above] at (0.6,3.08) {$B(x,\alpha t)$};
\draw[thick] (-1.6,2)--(1.6,2);
\draw (-1.6,1.88)--(-1.6,2.12) (1.6,1.88)--(1.6,2.12);
\node[right] at (1.7,2) {$Q(I)$};
\fill (0.6,2) circle (1.6pt);
\node[above] at (0.6,2.05) {$x$};
\draw[dashed,gray] (0.6,2.22)--(0.6,2.9);
\draw[thick] (-0.5,1)--(0.5,1);
\draw (-0.5,0.88)--(-0.5,1.12) (0.5,0.88)--(0.5,1.12);
\node[right] at (0.65,1) {$5I$};
\draw[thick] (-0.1,0)--(0.1,0);
\draw (-0.1,-0.12)--(-0.1,0.12) (0.1,-0.12)--(0.1,0.12);
\fill (0,0) circle (1.6pt);
\node[below] at (0,-0.1) {$y\in I$};
\node[right] at (0.65,0) {$\ell(I)\simeq\rho t$};
\draw[dashed,gray] (-0.1,0.2)--(-0.5,0.8);
\draw[dashed,gray] (0.1,0.2)--(0.5,0.8);
\draw[dashed,gray] (-0.5,1.2)--(-1.6,1.8);
\draw[dashed,gray] (0.5,1.2)--(1.6,1.8);
\end{tikzpicture}
\caption{A one-dimensional schematic of the spatial sets at a fixed scale $t$.
The atom is supported in $5I$, and its Haar expansion uses $Q(I)$.
Since $x$ may lie anywhere in $Q(I)$, the comparison with the area function
requires $|x-y|<\alpha t$ for every $y\in I$.
The choice $16\sqrt d\,\rho<\alpha$ ensures this inclusion in dimension $d$.}
\label{fig:dual-cone}
\end{figure}

\smallskip
\noindent{\it Step 2.} We estimate each class separately and then sum over
the $3^d$ adjacent dyadic systems.

Fix $\pmb\varepsilon\in\{0,1,2\}^d$ and a finite family $\mathscr I\subset\mathscr A_{\pmb\varepsilon}$.
All supporting cubes $Q(I)$ lie in $\mathcal D^{\pmb\varepsilon}$.
Applying Lemma~\ref{sp:l:smooth} with $M\le16^d$ and $K\le C$
as verified in \eqref{sp:e:cover} and {\it Step 1},
and then using \eqref{sp:e:discrete-area}, we obtain
\begin{align}\label{sp:e:class-reverse}
\bigg\|\sum_{I\in\mathscr I}\lambda_I(f)\mathsf a_I(f)\bigg\|_{L^p(\R^d)}
\le &C\sqrt p\,
\bigg\|\bigg(\sum_{I\in\mathscr I}\lambda_I(f)^2
\frac{\mathbf1_{Q(I)}}{|Q(I)|}\bigg)^{1/2}\bigg\|_{L^p(\R^d)}\notag\\
\le &C\sqrt p\,
\bigg\|\bigg(\sum_{I\in\mathcal D^0}\lambda_I(f)^2
\frac{\mathbf1_{Q(I)}}{|Q(I)|}\bigg)^{1/2}\bigg\|_{L^p(\R^d)}
\le C\sqrt p\,\left\|\mathcal A_{\varphi,\alpha}f\right\|_{L^p(\R^d)}.
\end{align}
The constant is independent of $\pmb\varepsilon$ and $\mathscr I$.

Now let $\widetilde{\mathscr I}\subset\mathcal D^0$ be any finite family.
Its decomposition into these classes is
$$
\widetilde{\mathscr I}
=\bigcup_{\pmb\varepsilon\in\{0,1,2\}^d}
(\widetilde{\mathscr I}\cap\mathscr A_{\pmb\varepsilon})
\quad\text{(disjoint)}.
$$
Consequently, the triangle inequality and \eqref{sp:e:class-reverse} prove
\begin{align}\label{sp:e:finite-reverse}
\bigg\|\sum_{I\in\widetilde{\mathscr I}}\lambda_I(f)\mathsf a_I(f)\bigg\|_{L^p(\R^d)}
\le &\sum_{\pmb\varepsilon\in\{0,1,2\}^d}
\bigg\|\sum_{I\in\widetilde{\mathscr I}\cap\mathscr A_{\pmb\varepsilon}}
\lambda_I(f)\mathsf a_I(f)\bigg\|_{L^p(\R^d)}\notag\\
\le &3^dC\sqrt p\,\left\|\mathcal A_{\varphi,\alpha}f\right\|_{L^p(\R^d)}
\le C_{d,\varphi,\alpha}\sqrt p\,\left\|\mathcal A_{\varphi,\alpha}f\right\|_{L^p(\R^d)}.
\end{align}

\smallskip
\noindent{\it Step 3.}
We verify \eqref{sp:e:atomic-repro} in $L^2$ and pass from \eqref{sp:e:finite-reverse} to $f$.

By \eqref{sp:e:energy},
$$
\sum_{I\in\mathcal D^0}\lambda_I(f)^2
=c_\varphi\left\|f\right\|_{L^2(\R^d)}^2<\infty,
\qquad
c_\Phi=\int_0^\infty|\widehat\Phi(re_1)|^2\,\frac{dr}{r}<\infty.
$$
Choose finite families $\mathscr I_j\subseteq\mathscr I_{j+1}$ with
$\bigcup_j\mathscr I_j=\mathcal D^0$.
For $g\in L^2(\R^d)$, \eqref{sp:e:reproduce} and evenness of $\Phi$ imply
\begin{align*}
\bigg|\bigg\langle f-\sum_{I\in\mathscr I_j}\lambda_I(f)\mathsf a_I(f),g\bigg\rangle\bigg|
= &\bigg|\sum_{I\notin\mathscr I_j}\int_{\mathcal W(I)}
(f*\varphi_t)(y)\overline{(g*\Phi_t)(y)}\,dy\,\frac{dt}{t}\bigg|\\
\le &\bigg(\sum_{I\notin\mathscr I_j}\lambda_I(f)^2\bigg)^{1/2}
\bigg(\int_0^\infty\int_{\R^d}|g*\Phi_t(y)|^2\,dy\,\frac{dt}{t}\bigg)^{1/2}\\
= &c_\Phi^{1/2}\bigg(\sum_{I\notin\mathscr I_j}\lambda_I(f)^2\bigg)^{1/2}
\left\|g\right\|_{L^2(\R^d)}.
\end{align*}
The same Cauchy--Schwarz bound proves absolute integrability of the product in this calculation.
Taking the supremum over $\|g\|_{L^2(\R^d)}=1$, we obtain
$$
\bigg\|f-\sum_{I\in\mathscr I_j}\lambda_I(f)\mathsf a_I(f)\bigg\|_{L^2(\R^d)}
\le c_\Phi^{1/2}
\bigg(\sum_{I\notin\mathscr I_j}\lambda_I(f)^2\bigg)^{1/2}\to0.
$$
This proves \eqref{sp:e:atomic-repro}, and a subsequence converges to $f$ almost everywhere.
Fatou's lemma applied to \eqref{sp:e:finite-reverse} along this subsequence gives
$$
\left\|f\right\|_{L^p(\R^d)}
\le\liminf_{k\to\infty}
\bigg\|\sum_{I\in\mathscr I_{j_k}}\lambda_I(f)\mathsf a_I(f)\bigg\|_{L^p(\R^d)}
\le C_{d,\varphi,\alpha}\sqrt p\,
\left\|\mathcal A_{\varphi,\alpha}f\right\|_{L^p(\R^d)}.
$$
For a general $f\in L^p(\R^d)$, take $f_j\in L^2(\R^d)\cap L^p(\R^d)$ with
$\|f_j-f\|_{L^p(\R^d)}\to0$.
By Lemma~\ref{sp:l:continuity},
$$
\left\|\mathcal A_{\varphi,\alpha}f_j-\mathcal A_{\varphi,\alpha}f\right\|_{L^p(\R^d)}
\le B_{p,d,\varphi,\alpha}\left\|f_j-f\right\|_{L^p(\R^d)}\to0.
$$
Letting $j\to\infty$ in the estimate for $f_j$, we obtain the reverse estimate for $f$.
\end{proof}

\subsection{Change of aperture and the general upper bound}

If the condition in Proposition~\ref{sp:p:reverse},
$$
\operatorname{dist}\bigl(0,\operatorname{supp}(|D|^{2-d}\varphi)\bigr)<c_d\alpha,
$$
fails at $\alpha$, we apply the proposition at an auxiliary aperture $\beta\ge\alpha$ with
$$
\operatorname{dist}\bigl(0,\operatorname{supp}(|D|^{2-d}\varphi)\bigr)<c_d\beta.
$$
The following lemma lets us return to aperture $\alpha$.

\begin{lemma}\label{sp:l:aperture}
For $0<\alpha\le\beta$, $2\le p<\infty$ and $f\in L^p(\R^d)$,
\begin{equation}\label{sp:e:aperture}
\left\|\mathcal A_{\varphi,\beta}f\right\|_{L^p(\R^d)}
\le C_d\left(\frac\beta\alpha\right)^{d/2}\sqrt p\,
\left\|\mathcal A_{\varphi,\alpha}f\right\|_{L^p(\R^d)},
\end{equation}
where the constant $C_d$ depends only on $d$.
\end{lemma}

The sharpness of the factor $\sqrt p$ in \eqref{sp:e:aperture}
will be proved in Corollary~\ref{sp:c:endpoint} below.

\begin{proof}[Proof of Lemma~\ref{sp:l:aperture}]
For $p=2$, the energy identity \eqref{sp:e:energy} gives
$$
\left\|\mathcal A_{\varphi,\beta}f\right\|_{L^2(\R^d)}^2
=\left(\frac\beta\alpha\right)^d
\left\|\mathcal A_{\varphi,\alpha}f\right\|_{L^2(\R^d)}^2.
$$
For $p>2$, we use $L^{p/2}$ duality, and compare the apertures
using the Hardy--Littlewood maximal operator:
$$
\mathsf Mw(x):=\sup_{B\ni x}\frac1{|B|}\int_B|w(z)|\,dz,
$$
where the supremum is over Euclidean balls $B\subset\R^d$.

Let $r=(p/2)'$.
Note that $\mathcal A_{\varphi,\beta}f\in L^p$ by Lemma~\ref{sp:l:continuity}, so
$$
\left\|\mathcal A_{\varphi,\beta}f\right\|_{L^p(\R^d)}^2
=\sup_{\substack{w\ge0\\\|w\|_{L^r(\R^d)}=1}}
\int_{\R^d}\left(\mathcal A_{\varphi,\beta}f\right)^2w.
$$
For any $w\ge0$ with $\|w\|_{L^r}=1$ and $|z-y|<\alpha t$,
$$
\int_{|x-y|<\beta t}w(x)\,dx
\le v_d(\beta t)^d\mathsf Mw(z).
$$
Averaging over $|z-y|<\alpha t$, we obtain
$$
\int_{|x-y|<\beta t}w(x)\,dx
\le\left(\frac\beta\alpha\right)^d
\int_{|z-y|<\alpha t}\mathsf Mw(z)\,dz.
$$
The maximal estimate for these weights is (see \cite{GM} for $d=1$)
$$
\left\|\mathsf Mw\right\|_{L^r(\R^d)}
\le C_d\frac r{r-1}\|w\|_{L^r(\R^d)}
\le C_dp\|w\|_{L^r(\R^d)}.
$$
Here $r/(r-1)=p/2$; taking square roots at the end will account for the
factor $\sqrt p$ in \eqref{sp:e:aperture}.
Changing the order of integration and applying H\"older's inequality, we obtain
\begin{align*}
\int_{\R^d}\left(\mathcal A_{\varphi,\beta}f\right)^2w
= &\int_0^\infty\int_{\R^d}|(f*\varphi_t)(y)|^2
\bigg(\int_{|x-y|<\beta t}w(x)\,dx\bigg)\,\frac{dy\,dt}{t^{d+1}}\\
\le &\left(\frac\beta\alpha\right)^d
\int_0^\infty\int_{\R^d}|(f*\varphi_t)(y)|^2
\bigg(\int_{|z-y|<\alpha t}\mathsf Mw(z)\,dz\bigg)\,\frac{dy\,dt}{t^{d+1}}\\
= &\left(\frac\beta\alpha\right)^d
\int_{\R^d}\left(\mathcal A_{\varphi,\alpha}f\right)^2\mathsf Mw\\
\le &\left(\frac\beta\alpha\right)^d
\left\|\mathcal A_{\varphi,\alpha}f\right\|_{L^p(\R^d)}^2
\left\|\mathsf Mw\right\|_{L^r(\R^d)}
\le C_d\left(\frac\beta\alpha\right)^dp\,
\left\|\mathcal A_{\varphi,\alpha}f\right\|_{L^p(\R^d)}^2.
\end{align*}
Taking the supremum over $w\ge0$ with $\|w\|_{L^r(\R^d)}=1$, we obtain \eqref{sp:e:aperture}.
\end{proof}

The Haar estimate contributes a factor $\sqrt p$ at a sufficiently large aperture.
Applying Lemma~\ref{sp:l:aperture} to return to the given aperture,
we obtain the following general upper bound with factor $p$.

\begin{proposition}\label{sp:p:general-upper}
For every $\alpha>0$, $2\le p<\infty$ and $f\in L^p(\R^d)$,
\begin{equation}\label{sp:e:general-upper}
\left\|f\right\|_{L^p(\R^d)}
\le C_{d,\varphi,\alpha}p\,
\left\|\mathcal A_{\varphi,\alpha}f\right\|_{L^p(\R^d)}.
\end{equation}
The constant is independent of $p$ and $f$.
No support condition on $|D|^{2-d}\varphi$ is required.
\end{proposition}

 The factor $p$ in \eqref{sp:e:general-upper} is optimal over the full kernel class;
see Lemma~\ref{sp:l:log-test} below.

\begin{proof}[Proof of Proposition~\ref{sp:p:general-upper}]
Since $|D|^{2-d}\varphi\ne0$, we have $\operatorname{dist}\bigl(0,\operatorname{supp}(|D|^{2-d}\varphi)\bigr)<\infty$.
Choose $\beta\ge\alpha$ such that
$\operatorname{dist}\bigl(0,\operatorname{supp}(|D|^{2-d}\varphi)\bigr)<c_d\beta$, with $c_d$ as in Proposition~\ref{sp:p:reverse}.
By Proposition~\ref{sp:p:reverse} and Lemma~\ref{sp:l:aperture},
\begin{align*}
\left\|f\right\|_{L^p(\R^d)}
\le &C_{d,\varphi,\beta}\sqrt p\,
\left\|\mathcal A_{\varphi,\beta}f\right\|_{L^p(\R^d)}
\le C_{d,\varphi,\beta}\left(\frac\beta\alpha\right)^{d/2}p\,
\left\|\mathcal A_{\varphi,\alpha}f\right\|_{L^p(\R^d)}
\le C_{d,\varphi,\alpha}p\,
\left\|\mathcal A_{\varphi,\alpha}f\right\|_{L^p(\R^d)}.
\end{align*}
Here $\beta$ depends only on $d$, $\varphi$ and $\alpha$, not on $p$.
\end{proof}

\subsection{A general lower bound}

To complete the proof of Theorem~\ref{sp:t:general}, we prove the lower bound $\sqrt p$
for every kernel and every aperture using localized lacunary sums.
Together with Proposition~\ref{sp:p:reverse}, this will prove the optimality of the order $\sqrt p$
whenever the support-distance condition holds.

\begin{lemma}\label{sp:l:lower}
Let $\varphi$ satisfy the basic assumptions and fix $\alpha>0$.
For every $p\ge2$, there is a nonzero $g_p\in\mathcal S_0(\R^d)$ such that
$$
\left\|g_p\right\|_{L^p(\R^d)}\ge c_{d,\varphi,\alpha}\sqrt p\,
\left\|\mathcal A_{\varphi,\alpha}g_p\right\|_{L^p(\R^d)}.
$$
\end{lemma}

\begin{proof}
The construction is motivated by Guth's wave-packet example in
\cite[Example~1.6]{GuthDecoupling}, where terms at different frequencies
are normalized to equal $1$ at the origin and stay close to $1$ nearby.
Here we use the lacunary sum
$$
\mathcal P_N(x)=\sum_{k=1}^N\cos(2\pi2^kx_1),
\qquad\text{with}\qquad N=\lceil p\rceil.
$$
Then $\mathcal P_N(x)\ge N/2$ for $0\le x_1\le2^{-N-4}$, and so
$$
\|\mathcal P_N\|_{L^p([0,1]^d)}
\ge\frac N2\,2^{-(N+4)/p}\gtrsim p,
$$
while in {\it Step 2} we will use the separated frequencies to prove
$$
\|\mathcal A_{\varphi,\alpha}\mathcal P_N\|_{L^\infty(\R^d)}
\le C_{d,\varphi,\alpha}\sqrt N
\le C_{d,\varphi,\alpha}\sqrt p.
$$
The restriction to a region of measure comparable to $2^{-N}$ costs the factor $2^{-N/p}$.
Taking $N=\lceil p\rceil$ keeps this factor bounded away from zero and leaves the ratio $N/\sqrt N\simeq\sqrt p$.
Since $\mathcal P_N\notin L^p(\R^d)$, we consider
$$
g_L(x)=\eta(x/L)\mathcal P_N(x),\qquad 1\ll L\in\mathbb N,
$$
where $\eta\in\mathcal S(\R^d)$ is a real-valued and even function, satisfying
$\operatorname{supp}\widehat\eta\Subset B(0,1/4)$ and $\eta\ge1$ on $[-1,1]^d$.
The cutoff makes $g_L$ integrable without introducing frequencies near zero.
See also our work \cite[Lemmas~6.2 and~6.3]{LPWCWW}
for the trigonometric sums and localization.

\smallskip
\noindent{\it Step 1.}
We first estimate $\|g_L\|_{L^p(\R^d)}$.
Note that the $2L-1$ intervals
$$
I_j=[j-2^{-N-4},j+2^{-N-4}],
\qquad j=-L+1,\ldots,L-1,
$$
are disjoint and lie in $[-L,L]$.
For each $j=-L+1,\ldots,L-1$,
$$
\sup_{x_1\in I_j}|2\pi2^k(x_1-j)|
\le\frac{\pi}{8}\qquad(1\le k\le N).
$$
By periodicity,
$$
\inf_{x_1\in I_j}\sum_{k=1}^N\cos(2\pi2^kx_1)
=\inf_{x_1\in I_j}\sum_{k=1}^N\cos\bigl(2\pi2^k(x_1-j)\bigr)
\ge N\cos\frac{\pi}{8}\ge\frac N2.
$$
Since $\eta(x/L)\ge1$ on $[-L,L]^d$,
$|g_L(x)|\ge N/2$ for every
$x\in\bigl(\bigcup_{j=-L+1}^{L-1}I_j\bigr)\times[-L,L]^{d-1}$.
Consequently,
\begin{equation}\label{sp:e:lacunary-size}
\|g_L\|_{L^p(\R^d)}
\ge\frac N2\bigl((2L-1)2^{-N-3}(2L)^{d-1}\bigr)^{1/p}
\ge cpL^{d/p}.
\end{equation}
Also, $g_L\in\mathcal S(\R^d)$ and
$$
\widehat g_L(\xi)
=\frac{L^d}{2}\sum_{k=1}^N
\left[\widehat\eta\bigl(L(\xi-2^ke_1)\bigr)
+\widehat\eta\bigl(L(\xi+2^ke_1)\bigr)\right].
$$
Since $L\ge1$,
$$
\operatorname{supp}\widehat g_L
\subset\bigcup_{k=1}^N
\left[B\left(2^ke_1,\frac1{4L}\right)
\cup B\left(-2^ke_1,\frac1{4L}\right)\right]
\Subset\R^d\setminus\{0\}.
$$

\smallskip
\noindent{\it Step 2.}
We next bound the area function of $\mathcal P_N$.
By $\int_{\R^d}\varphi=0$ and integration by parts for
$\varphi\in C_c^\infty(\R^d)$, respectively, we have for $r\ne0$,
$$
\widehat\varphi(re_1)
=\int_{\R^d}\varphi(z)(e^{-2\pi irz_1}-1)\,dz,
\qquad
\widehat\varphi(re_1)
=\frac1{2\pi ir}\int_{\R^d}\partial_{z_1}\varphi(z)e^{-2\pi irz_1}\,dz.
$$
Thus
\begin{align*}
|\widehat\varphi(re_1)|
\le &\min\left\{2\pi|r|\int_{\R^d}|z_1\varphi(z)|\,dz,\
\frac1{2\pi|r|}\int_{\R^d}|\partial_{z_1}\varphi(z)|\,dz\right\}
\le C_\varphi\min\{|r|,|r|^{-1}\}.
\end{align*}
The factor $\widehat\varphi(t2^ke_1)$ decays when $t$ is far from $2^{-k}$.
At each fixed $t$, the sum over all frequencies is therefore bounded independently of $N$:
$$
\sum_{k=1}^N|\widehat\varphi(t2^ke_1)|
\le C_\varphi\,
\bigg(\sum_{1\le k\le N:\,t2^k\le1}t2^k+
\sum_{1\le k\le N:\,t2^k>1}(t2^k)^{-1}\bigg)
\le C_\varphi.
$$
Each term contributes the same finite amount to the scale integral:
$$
\int_0^\infty\sum_{k=1}^N|\widehat\varphi(t2^ke_1)|\,\frac{dt}{t}
=N\int_0^\infty|\widehat\varphi(re_1)|\,\frac{dr}{r}
\le C_\varphi N.
$$
Since $\varphi$ is real-valued and even,
$$
(\mathcal P_N*\varphi_t)(y)
=\sum_{k=1}^N\widehat\varphi(t2^ke_1)\cos(2\pi2^ky_1).
$$
Consequently,
\begin{align}\label{sp:e:periodic-area}
\mathcal A_{\varphi,\alpha}\mathcal P_N(x)^2
= &\int_0^\infty\int_{|y-x|<\alpha t}|(\mathcal P_N*\varphi_t)(y)|^2
\,\frac{dy\,dt}{t^{d+1}}\notag\\
\le &v_d\alpha^d\int_0^\infty
\bigg(\sum_{k=1}^N|\widehat\varphi(t2^ke_1)|\bigg)^2\,\frac{dt}{t}\notag\\
\le &v_d\alpha^d
\bigg(\sup_{t>0}\sum_{k=1}^N|\widehat\varphi(t2^ke_1)|\bigg)
\int_0^\infty\sum_{k=1}^N|\widehat\varphi(t2^ke_1)|\,\frac{dt}{t}
\le C_{d,\varphi,\alpha}N.
\end{align}

\smallskip
\noindent{\it Step 3.}
We now use {\it Step 2} to estimate the area function of $g_L$.
Compare
$$
(g_L*\varphi_t)(y)\quad\text{with}\quad
\eta(x/L)(\mathcal P_N*\varphi_t)(y),
\qquad\text{whenever }|y-x|<\alpha t.
$$
Let $y=x+tu$ with $|u|<\alpha$. Then $dy\,dt/t^{d+1}=du\,dt/t$.
Since $g_L=\eta(\cdot/L)\mathcal P_N$, by Fourier inversion,
\begin{align}\label{sp:e:cutoff-difference}
&(g_L*\varphi_t)(x+tu)-\eta(x/L)(\mathcal P_N*\varphi_t)(x+tu)\notag\\
= &\frac12\sum_{\nu\in\{\pm2^ke_1:1\le k\le N\}}e^{2\pi i\nu\cdot(x+tu)}
\int_{\R^d}\widehat\eta(\zeta)e^{2\pi ix\cdot\zeta/L}
\left[e^{2\pi itu\cdot\zeta/L}\widehat\varphi(t(\nu+\zeta/L))
-\widehat\varphi(t\nu)\right]\,d\zeta.
\end{align}
To control the area-function error in $L^p$, we need decay in $x/L$ as well as integrability in $t$.
Choose an integer $M>d$.
On $|\xi-\nu|\le1/4$, the Fourier bounds for $\varphi$ imply
$$
\max_{|\gamma|\le2M+1}
|\partial_\xi^\gamma(\widehat\varphi(t\xi))|
\le C_{N,\varphi,M}\min\{t,t^{-2M-3}\}.
$$
For $t\le1$, use $\widehat\varphi(0)=0$ when $|\gamma|=0$ and the factor $t^{|\gamma|}$ otherwise.
For $t\ge1$, use $|\xi|\ge1$ and rapid decrease of all derivatives.
For $|\gamma|=0$, write the expression in brackets in \eqref{sp:e:cutoff-difference} as
\begin{align*}
e^{2\pi itu\cdot\zeta/L}\widehat\varphi(t(\nu+\zeta/L))-\widehat\varphi(t\nu)
= &e^{2\pi itu\cdot\zeta/L}
\bigl(\widehat\varphi(t(\nu+\zeta/L))-\widehat\varphi(t\nu)\bigr)
+\bigl(e^{2\pi itu\cdot\zeta/L}-1\bigr)\widehat\varphi(t\nu).
\end{align*}
The mean value theorem and $|e^{2\pi itu\cdot\zeta/L}-1|\le C_\alpha t|\zeta|/L$ then imply, on $|\zeta|\le1/4$,
$$
\bigl|e^{2\pi itu\cdot\zeta/L}\widehat\varphi(t(\nu+\zeta/L))-\widehat\varphi(t\nu)\bigr|
\le \frac{C_{N,d,\varphi,\alpha,M}}{L}(1+t)\min\{t,t^{-2M-3}\}.
$$
For $1\le|\gamma|\le2M$, the derivative of the constant term vanishes, and the product rule implies
\begin{align*}
&\partial_\zeta^\gamma
\bigl[e^{2\pi itu\cdot\zeta/L}\widehat\varphi(t(\nu+\zeta/L))\bigr]
= L^{-|\gamma|}e^{2\pi itu\cdot\zeta/L}
\sum_{\beta\le\gamma}\binom{\gamma}{\beta}
(2\pi it u)^{\gamma-\beta}
\left[\partial_\xi^\beta\bigl(\widehat\varphi(t\xi)\bigr)\right]_{\xi=\nu+\zeta/L}.
\end{align*}
Since $L\ge1$ and $|\gamma|\le2M$, the Fourier bounds above yield
$$
\bigg|\partial_\zeta^\gamma
\bigl[e^{2\pi itu\cdot\zeta/L}\widehat\varphi(t(\nu+\zeta/L))\bigr]\bigg|
\le\frac{C_{N,d,\varphi,\alpha,M}}{L}
(1+t)^{2M}\min\{t,t^{-2M-3}\}.
$$
Multiplying by $\widehat\eta$, applying the product rule once more, and using its fixed compact support, we obtain
\begin{align}\label{sp:e:cutoff-derivatives}
&\sum_{|\gamma|\le2M}
\bigg\|\partial_\zeta^\gamma
\left\{\widehat\eta(\zeta)
\left[e^{2\pi itu\cdot\zeta/L}\widehat\varphi(t(\nu+\zeta/L))
-\widehat\varphi(t\nu)\right]\right\}\bigg\|_{L^1(d\zeta)}
\le \frac{C_{N,d,\varphi,\alpha,M}}{L}\min\{t,t^{-2}\}.
\end{align}

 In \eqref{sp:e:cutoff-difference}, we integrate by parts using
$(1-\Delta_\zeta)^Me^{2\pi ix\cdot\zeta/L}
=(1+4\pi^2|x|^2/L^2)^Me^{2\pi ix\cdot\zeta/L}$.
After summing over the $2N$ frequencies and applying \eqref{sp:e:cutoff-derivatives},
we obtain
$$
|(g_L*\varphi_t)(x+tu)-\eta(x/L)(\mathcal P_N*\varphi_t)(x+tu)|
\le\frac{C_{N,d,\varphi,\alpha}}{L}
(1+|x|/L)^{-2M}\min\{t,t^{-2}\}.
$$
By the triangle inequality in $L^2(du\,dt/t)$ and \eqref{sp:e:periodic-area},
\begin{align*}
\mathcal A_{\varphi,\alpha}g_L(x)
\le &|\eta(x/L)|\,\mathcal A_{\varphi,\alpha}\mathcal P_N(x)\\
&+\bigg(\int_{|u|<\alpha}\int_0^\infty
|(g_L*\varphi_t)(x+tu)-\eta(x/L)(\mathcal P_N*\varphi_t)(x+tu)|^2\,\frac{dt}{t}\,du\bigg)^{1/2}\\
\le &C_{d,\varphi,\alpha}\sqrt N\,|\eta(x/L)|
+\frac{C_{N,d,\varphi,\alpha}}{L}(1+|x|/L)^{-2M}
\bigg(\int_0^1t^2\,\frac{dt}{t}+\int_1^\infty t^{-4}\,\frac{dt}{t}\bigg)^{1/2}\\
\le &C_{d,\varphi,\alpha}\sqrt N\,|\eta(x/L)|
+\frac{C_{N,d,\varphi,\alpha}}{L}(1+|x|/L)^{-2M}.
\end{align*}
The $L^p$ norms of $\eta$ and $(1+|x|)^{-2M}$ are bounded uniformly for $p\ge2$.
Consequently,
$$
\left\|\mathcal A_{\varphi,\alpha}g_L\right\|_{L^p(\R^d)}
\le\left(C_{d,\varphi,\alpha}\sqrt N+
\frac{C_{N,d,\varphi,\alpha}}{L}\right)L^{d/p}.
$$
 Choose an integer $L$ large enough that the second term in parentheses is at most one,
and let $g_p=g_L$.
Together with \eqref{sp:e:lacunary-size} and $N\le p+1$, this proves the lemma.
\end{proof}

\begin{proof}[Proof of Theorem~\ref{sp:t:general}]
By Lemma~\ref{sp:l:lower} and Proposition~\ref{sp:p:general-upper}, we have
$$
c_{d,\varphi,\alpha}\sqrt p
\le\mathfrak c_{p,d}^{\varphi}(\alpha)
\le C_{d,\varphi,\alpha}p.
$$
By Proposition~\ref{sp:p:reverse}, we also have
$$
\mathfrak c_{p,d}^{\varphi}(\alpha)
\simeq_{d,\varphi,\alpha}\sqrt p,
\qquad \operatorname{dist}\bigl(0,\operatorname{supp}(|D|^{2-d}\varphi)\bigr)<c_d\alpha,
$$
where $c_d>0$ is the constant in that proposition.
If $0\in\operatorname{supp}(|D|^{2-d}\varphi)$, this condition holds for every $\alpha>0$.
For fixed $\varphi$, it holds for sufficiently large $\alpha$.
\end{proof}

\subsection{Poisson and heat area functions}\label{sp:sub:semigroups}

We now consider the classical Poisson and heat area functions on $\R^d$.
As in our earlier work \cite{CDWY,WuJGA}, we use local reconstruction and estimates for sums of atomic functions;
see also Han~\cite{Han94} for Calder\'on reproducing formulas associated with para-accretive functions.
Although the semigroup kernels are not compactly supported, they admit dual kernels
with arbitrarily small compact support. The resulting atoms then satisfy the support requirements
in Proposition~\ref{sp:p:reverse}.

For $\nu\in\{\mathrm P,\mathrm H\}$, let
$
b_{\mathrm P}(r)=re^{-r},\ b_{\mathrm H}(r)=r^2e^{-r^2},\ r\ge0,
$
and define
\begin{equation}\label{sp:e:sg-kernels}
\widehat{\varphi_\nu^{(d)}}(\xi)=b_\nu(|\xi|),
\qquad
\varphi_{\nu,t}^{(d)}(x)=t^{-d}\varphi_\nu^{(d)}(x/t).
\end{equation}
In the Poisson case, if $\widehat{P_a}(\xi)=e^{-2\pi a|\xi|}$, then
$\varphi_{\mathrm P}^{(d)}=-a\partial_aP_a|_{a=(2\pi)^{-1}}$.
In the heat case,
$$
\varphi_{\mathrm H}^{(d)}(x)
=-\frac1{4\pi^2}\Delta\bigl(\pi^{d/2}e^{-\pi^2|x|^2}\bigr)
=\pi^{d/2}\bigl(d/2-\pi^2|x|^2\bigr)e^{-\pi^2|x|^2}.
$$
These kernels are real and radial. For every multi-index $\gamma$,
\begin{equation}\label{sp:e:sg-decay}
|\partial^\gamma\varphi_\nu^{(d)}(x)|
\le C_{d,\gamma}(1+|x|)^{-d-1-|\gamma|},
\qquad
\int_{\R^d}\varphi_\nu^{(d)}=0.
\end{equation}
In the Poisson case, the derivative bound follows by differentiating
$$
P_a(x)=\frac{\Gamma((d+1)/2)}{\pi^{(d+1)/2}}\,
\frac{a}{(a^2+|x|^2)^{(d+1)/2}}.
$$
In the heat case, all derivatives decrease rapidly.

We use the area function in \eqref{sp:e:area} and the best reverse constant in \eqref{sp:e:best}
with $\varphi=\varphi_\nu^{(d)}$.
The test functions in the following proposition will also be used for Zygmund dilations.

\begin{proposition}\label{sp:p:semigroups}
For every $d\ge1$, $\alpha>0$, $\nu\in\{\mathrm P,\mathrm H\}$ and $2\le p<\infty$,
\begin{equation}\label{sp:e:sg-reverse}
\|f\|_{L^p(\R^d)}
\le C_{d,\alpha}\sqrt p\,\|\mathcal A_{\varphi_\nu^{(d)},\alpha}f\|_{L^p(\R^d)},
\qquad f\in L^p(\R^d).
\end{equation}
For each $p$, there is a nonzero $g_p\in\mathcal S_0(\R^d)$ such that
\begin{equation}\label{sp:e:sg-lower}
\|g_p\|_{L^p(\R^d)}
\ge c_{d,\alpha}\sqrt p\,\|\mathcal A_{\varphi_\nu^{(d)},\alpha}g_p\|_{L^p(\R^d)}.
\end{equation}
Consequently,
$$
\mathfrak c_{p,d}^{\varphi_\nu^{(d)}}(\alpha)\simeq_{d,\alpha}\sqrt p.
$$
\end{proposition}

The absence of an aperture restriction has a simple explanation:
$$
(|D|^{2-d}\varphi_\nu^{(d)})(0)
=\sigma_{d-1}\int_0^\infty r b_\nu(r)\,dr
=
\begin{cases}
2\sigma_{d-1},&\nu=\mathrm P,\\[2pt]
\sigma_{d-1}/2,&\nu=\mathrm H.
\end{cases}
$$
Thus the local pairing does not vanish near the origin for either kernel.

\begin{proof}[Proof of Proposition \ref{sp:p:semigroups}]
\noindent{\it Step 1.} We first prove the energy identities and $L^p$ continuity.
Write
\begin{equation}\label{sp:e:sg-energy-constants}
c_{\mathrm P}=\int_0^\infty b_{\mathrm P}(r)^2\,\frac{dr}{r}=\frac14,
\qquad
c_{\mathrm H}=\int_0^\infty b_{\mathrm H}(r)^2\,\frac{dr}{r}=\frac18.
\end{equation}
For $f\in L^2(\R^d)$, we compute as in \eqref{sp:e:energy},
with radial Fourier profile $b_\nu$:
\begin{equation}\label{sp:e:sg-energy}
\int_0^\infty\|f*\varphi_{\nu,t}^{(d)}\|_{L^2(\R^d)}^2\,\frac{dt}{t}
=c_\nu\|f\|_{L^2(\R^d)}^2,
\qquad
\|\mathcal A_{\varphi_\nu^{(d)},\alpha}f\|_{L^2(\R^d)}^2
=v_d\alpha^dc_\nu\|f\|_{L^2(\R^d)}^2.
\end{equation}

For the upper $L^p$ estimate, use the convolution kernel
$$
K_\nu(x;u,t)=t^{-d}\varphi_\nu^{(d)}(x/t+u),
\qquad |u|<\alpha,\quad t>0.
$$
For $x\ne0$ and $|\gamma|=0,1$, substitute $r=|x|/t$ to obtain
\begin{align*}
\int_{|u|<\alpha}\int_0^\infty
|\partial_x^\gamma K_\nu(x;u,t)|^2\,\frac{dt}{t}\,du
= &|x|^{-2d-2|\gamma|}
\int_{|u|<\alpha}\int_0^\infty r^{2d+2|\gamma|-1}
\left|(\partial^\gamma\varphi_\nu^{(d)})(rx/|x|+u)\right|^2\,dr\,du\\
\le &C_{d,\alpha}|x|^{-2d-2|\gamma|}.
\end{align*}
For $r\le2(1+\alpha)$, the derivatives are bounded; for larger $r$,
the radial integrand is at most $C_{d,\alpha}r^{-3}$ by \eqref{sp:e:sg-decay}.
As in Lemma~\ref{sp:l:continuity}, we obtain, initially for $f\in L^2\cap L^p$,
\begin{equation}\label{sp:e:sg-forward}
\|\mathcal A_{\varphi_\nu^{(d)},\alpha}f\|_{L^p(\R^d)}
\le B_{p,d,\alpha}\|f\|_{L^p(\R^d)},
\qquad 1<p<\infty.
\end{equation}
For $f\in L^p$, choose $f_j\in L^2\cap L^p$ with $f_j\to f$ in $L^p$.
By Young's inequality, for each $N\ge2$,
$$
\int_{N^{-1}}^N
\|(f_j-f)*\varphi_{\nu,t}^{(d)}\|_{L^p(\R^d)}^p\,\frac{dt}{t}
\le 2\log N\,\|\varphi_\nu^{(d)}\|_1^p\|f_j-f\|_p^p\to0.
$$
Along a common subsequence, the convolutions converge almost everywhere in $(y,t)$.
By Fatou's lemma, \eqref{sp:e:sg-forward} holds for the convolution definition.
The $L^p$ continuity follows from the pointwise triangle inequality.

\smallskip
\noindent{\it Step 2.} We construct the localized dual kernels.
Fix $\rho>0$, choose a nonzero real radial
$\eta_d\in C_c^\infty(B(0,\rho/2))$, and define
\begin{equation}\label{sp:e:sg-dual}
c_{\nu,d}=\int_0^\infty r^2|\widehat\eta_d(re_1)|^2b_\nu(r)\,\frac{dr}{r}>0,
\qquad
\Phi_\nu^{(d)}=-\frac{\Delta(\eta_d*\eta_d)}{4\pi^2c_{\nu,d}}.
\end{equation}
Here $c_{\nu,d}>0$ because $b_\nu(r)>0$ for $r>0$.
The kernel $\Phi_\nu^{(d)}$ is real, radial and smooth, with
$$
\operatorname{supp}\Phi_\nu^{(d)}\subset B(0,\rho),
\qquad
\int_{\R^d}\Phi_\nu^{(d)}=0,
\qquad
\widehat{\Phi_\nu^{(d)}}(\xi)
=c_{\nu,d}^{-1}|\xi|^2|\widehat\eta_d(\xi)|^2\ge0.
$$
The definition of $c_{\nu,d}$ gives
\begin{equation}\label{sp:e:sg-norm}
\int_0^\infty b_\nu(r)\widehat{\Phi_\nu^{(d)}}(re_1)\,\frac{dr}{r}
=c_{\nu,d}^{-1}\int_0^\infty r^2|\widehat\eta_d(re_1)|^2b_\nu(r)\,\frac{dr}{r}=1.
\end{equation}
By \eqref{sp:e:sg-norm} and nonnegativity, the truncated multipliers lie in $[0,1]$
and converge to one for $\xi\ne0$.
By Plancherel's theorem and dominated convergence, as in \eqref{sp:e:reproduce}, we have
\begin{equation}\label{sp:e:sg-reconstruction}
f=\int_0^\infty(f*\varphi_{\nu,t}^{(d)})*\Phi_{\nu,t}^{(d)}\,\frac{dt}{t}
\quad\text{in }L^2(\R^d),
\qquad
\Phi_{\nu,t}^{(d)}(x)=t^{-d}\Phi_\nu^{(d)}(x/t).
\end{equation}
Also,
$$
c_{\Phi_\nu^{(d)}}=\int_0^\infty
|\widehat{\Phi_\nu^{(d)}}(re_1)|^2\,\frac{dr}{r}<\infty.
$$

\smallskip
\noindent{\it Step 3.} We prove the reverse estimate.
Choose $\rho>0$ with $16\sqrt d\,\rho<\alpha$, and first take
$f\in L^2(\R^d)\cap L^p(\R^d)$.
With $\mathcal W(I)$ defined using this $\rho$ as in {\it Step 1}
of Proposition~\ref{sp:p:reverse}, define the coefficients and atoms by
$$
\lambda_I(f)^2
=\int_{\mathcal W(I)}|(f*\varphi_{\nu,t}^{(d)})(y)|^2\,dy\,\frac{dt}{t},
$$
and
$$
\mathsf a_I(f)(x)=\frac1{\lambda_I(f)}
\int_{\mathcal W(I)}(f*\varphi_{\nu,t}^{(d)})(y)
\Phi_{\nu,t}^{(d)}(x-y)\,dy\,\frac{dt}{t},
$$
with $\mathsf a_I(f)=0$ when $\lambda_I(f)=0$.
By the support and cancellation of $\Phi_\nu^{(d)}$ and \eqref{sp:e:cover},
$$
\operatorname{supp}\mathsf a_I(f)\subset5I\subset Q(I),
\qquad
\int_{\R^d}\mathsf a_I(f)=0.
$$
For $|\gamma|\le1$, as in \eqref{sp:e:atom-derivative},
\begin{align*}
|\partial^\gamma\mathsf a_I(f)(x)|
\le &\|\partial^\gamma\Phi_\nu^{(d)}\|_2
\bigg(\int_{\ell(I)/\rho}^{2\ell(I)/\rho}t^{-d-2|\gamma|}\,\frac{dt}{t}\bigg)^{1/2}
\le C_{d,\alpha}|Q(I)|^{-1/2}\ell(Q(I))^{-|\gamma|}.
\end{align*}
The cone inclusion and disjointness of $\mathcal W(I)$ give,
as in \eqref{sp:e:discrete-area},
$$
\sum_I\lambda_I(f)^2\frac{\mathbf1_{Q(I)}}{|Q(I)|}
\le C_d\rho^{-d}\bigl(\mathcal A_{\varphi_\nu^{(d)},\alpha}f\bigr)^2.
$$

For a finite family $\mathscr I\subset\mathcal D^0$, decompose $\mathscr I$ as
$$
\mathscr I
=\bigcup_{\pmb\varepsilon\in\{0,1,2\}^d}
(\mathscr I\cap\mathscr A_{\pmb\varepsilon})
\quad\text{(disjoint)}.
$$
Within each class, all supporting cubes belong to
$\mathcal D^{\pmb\varepsilon}$, with multiplicity at most $16^d$.
Applying Lemma~\ref{sp:l:smooth} and the pointwise bound above,
then summing over $\pmb\varepsilon$ as in
\eqref{sp:e:class-reverse}--\eqref{sp:e:finite-reverse}, we obtain
$$
\bigg\|\sum_{I\in\mathscr I}\lambda_I(f)\mathsf a_I(f)\bigg\|_p
\le C_{d,\alpha}\sqrt p\,\|\mathcal A_{\varphi_\nu^{(d)},\alpha}f\|_p.
$$
By \eqref{sp:e:sg-energy},
$$
\sum_{I\in\mathcal D^0}\lambda_I(f)^2=c_\nu\|f\|_2^2<\infty.
$$
By \eqref{sp:e:sg-reconstruction} and $c_{\Phi_\nu^{(d)}}<\infty$,
the argument in {\it Step 3} of Proposition~\ref{sp:p:reverse}
gives $L^2$ convergence of the atomic sums to $f$.
Fatou's lemma and the $L^p$ continuity from {\it Step 1} then prove
\eqref{sp:e:sg-reverse} for every $f\in L^p$.

\smallskip
\noindent{\it Step 4.} We verify the lower bound.
The Fourier estimates needed in the proof of Lemma~\ref{sp:l:lower} hold for both kernels:
\begin{equation}\label{sp:e:sg-size}
|b_\nu(r)|\le C\min\{r,r^{-1}\},\qquad r>0,
\end{equation}
and, for every multi-index $\gamma$ and $M>0$,
\begin{equation}\label{sp:e:sg-symbol}
|\partial^\gamma[b_\nu(|\xi|)]|
\le C_{d,\gamma,M}
\begin{cases}
|\xi|^{1-|\gamma|},&0<|\xi|\le1,\\[2pt]
|\xi|^{-M},&|\xi|>1.
\end{cases}
\end{equation}
In the heat case, the first bound can be improved to $C_{d,\gamma}|\xi|^{2-|\gamma|}$.
These estimates follow by differentiating $re^{-r}$ and $r^2e^{-r^2}$.
In particular, for $0<c<C<\infty$,
\begin{equation}\label{sp:e:sg-annular}
\sup_{c\le|\xi|\le C}
|\partial_\xi^\gamma[b_\nu(t|\xi|)]|
\le C_{d,c,C,\gamma,M}\min\{t,t^{-M}\},\qquad t>0.
\end{equation}
For small $t$, use the factor $t^{|\gamma|}$ and the first line of \eqref{sp:e:sg-symbol};
for large $t$, use the second line with exponent $M+|\gamma|$.

By \eqref{sp:e:sg-size}, the periodic estimate \eqref{sp:e:periodic-area}
holds for $\varphi_\nu^{(d)}$.
Choose an integer $M_0>d$. Applying \eqref{sp:e:sg-annular} with
$|\gamma|\le2M_0+1$ and decay exponent $2M_0+3$, we obtain
\eqref{sp:e:cutoff-derivatives} for these kernels.
The localization argument in {\it Step 3} of Lemma~\ref{sp:l:lower}
then gives a nonzero $g_p\in\mathcal S_0(\R^d)$ satisfying
\eqref{sp:e:sg-lower}.
Moreover, $0<\|\mathcal A_{\varphi_\nu^{(d)},\alpha}g_p\|_p<\infty$
by \eqref{sp:e:sg-forward} and \eqref{sp:e:sg-reverse}.
\end{proof}

\subsection{The Poisson area function of Chang, Wilson and Wolff}\label{sp:sub:cww}

The kernel $\varphi_{\mathrm P}^{(d)}$ in \eqref{sp:e:sg-kernels}
comes from a derivative in the semigroup parameter.
For comparison with \cite{CWW}, we now use the spatial gradient of the Poisson extension.
Let
$$
P(x)=\frac{\Gamma((d+1)/2)}{\pi^{(d+1)/2}(1+|x|^2)^{(d+1)/2}},
\qquad P_t(x)=t^{-d}P(x/t),
\qquad \psi_j=\partial_jP,\quad 1\le j\le d.
$$
Then $(\psi_j)_t=t\partial_jP_t$, and the area function of \cite[pp.~238--239]{CWW} is
\begin{equation}\label{sp:e:cww-area}
A_\alpha^\nabla f(x)^2
=\int_{\Gamma_\alpha(x)}|t\nabla_y(P_t*f)(y)|^2\,\frac{dy\,dt}{t^{d+1}}
=\int_{\Gamma_\alpha(x)}\sum_{j=1}^d|f*(\psi_j)_t(y)|^2\,\frac{dy\,dt}{t^{d+1}}.
\end{equation}
The superscript distinguishes this vector-kernel area function from the scalar one above.

\begin{proposition}\label{sp:p:cww}
For every $d\ge1$, $\alpha>0$ and $2\le p<\infty$,
$$
c_{d,\alpha}\sqrt p
\le \sup_{0\ne f\in L^p(\R^d)}
\frac{\|f\|_{L^p(\R^d)}}{\|A_\alpha^\nabla f\|_{L^p(\R^d)}}
\le C_{d,\alpha}\sqrt p.
$$
\end{proposition}

\begin{proof}
We use all components to reconstruct $f$, obtaining scalar atoms to which Lemma~\ref{sp:l:smooth} applies.
Since
$$
\widehat{\psi_j}(\xi)=2\pi i\xi_j e^{-2\pi|\xi|},
\qquad
\int_0^\infty\sum_{j=1}^d|\widehat{\psi_j}(t\xi)|^2\,\frac{dt}{t}
=\int_0^\infty r^2e^{-2r}\,\frac{dr}{r}=\frac14
\quad(\xi\ne0),
$$
the calculation of \eqref{sp:e:energy}, summed over the components, gives
$$
\int_0^\infty\sum_{j=1}^d\|f*(\psi_j)_t\|_2^2\,\frac{dt}{t}=\frac14\|f\|_2^2,
\qquad
\|A_\alpha^\nabla f\|_2^2=\frac{v_d\alpha^d}{4}\|f\|_2^2.
$$
The spatial kernels satisfy
$$
|\partial^\gamma\psi_j(x)|
\le C_{d,\gamma}(1+|x|)^{-d-2-|\gamma|}.
$$
Repeating the kernel calculation in {\it Step 1} of Proposition~\ref{sp:p:semigroups},
now with values in
$
L^2\bigl(B(0,\alpha)\times(0,\infty),du\,dt/t;\mathbb C^d\bigr),
$
we obtain, for $1<p<\infty$,
$$
\|A_\alpha^\nabla f\|_{L^p(\R^d)}\le B_{p,d,\alpha}\|f\|_{L^p(\R^d)}.
$$
The $L^p$ continuity follows as in that proposition.

Choose $16\sqrt d\,\rho<\alpha$ and a nonzero real radial $\eta\in C_c^\infty(B(0,\rho/2))$.
To construct the dual vector kernel, let
$$
a_d=\int_0^\infty r|\widehat\eta(re_1)|^2e^{-2\pi r}\,dr>0,
\qquad
q_j=-\frac{\partial_j(\eta*\eta)}{4\pi^2a_d},\quad 1\le j\le d.
$$
Each $q_j$ has support in $B(0,\rho)$ and integral zero.
Moreover,
$$
\sum_{j=1}^d\widehat{\psi_j}(\xi)\widehat q_j(\xi)
=a_d^{-1}|\xi|^2e^{-2\pi|\xi|}|\widehat\eta(\xi)|^2\ge0,
$$
and hence
$$
\int_0^\infty\sum_{j=1}^d
\widehat{\psi_j}(t\xi)\widehat q_j(t\xi)\,\frac{dt}{t}
=a_d^{-1}\int_0^\infty re^{-2\pi r}|\widehat\eta(re_1)|^2\,dr=1
\quad(\xi\ne0).
$$
The truncated multipliers lie in $[0,1]$ and converge to one.
By Plancherel's theorem and dominated convergence, as in \eqref{sp:e:sg-reconstruction}, we have
$$
f=\int_0^\infty\sum_{j=1}^d(f*(\psi_j)_t)*(q_j)_t\,\frac{dt}{t}
\quad\text{in }L^2(\R^d).
$$
By radiality of $\eta$, the finite constant
$$
c_q=\int_0^\infty\sum_{j=1}^d|\widehat q_j(t\xi)|^2\,\frac{dt}{t}
=\frac1{4\pi^2a_d^2}\int_0^\infty r|\widehat\eta(re_1)|^4\,dr
$$
is independent of $\xi\ne0$.

For $f\in L^2\cap L^p$, define the regions $\mathcal W(I)$ using this $\rho$
as in {\it Step 1} of Proposition~\ref{sp:p:reverse}.
With the covering cubes $Q(I)$ fixed in \eqref{sp:e:cover}, let
$$
\lambda_I(f)^2
=\int_{\mathcal W(I)}\sum_{j=1}^d|f*(\psi_j)_t(y)|^2\,dy\,\frac{dt}{t},
$$
and
$$
\mathsf a_I(f)(x)
=\frac1{\lambda_I(f)}
\int_{\mathcal W(I)}\sum_{j=1}^d(f*(\psi_j)_t)(y)(q_j)_t(x-y)\,dy\,\frac{dt}{t},
$$
with $\mathsf a_I(f)=0$ when $\lambda_I(f)=0$.
Cauchy--Schwarz in $j,y,t$ gives, for $|\gamma|\le1$,
$$
|\partial^\gamma\mathsf a_I(f)(x)|
\le\bigg(\sum_{j=1}^d\|\partial^\gamma q_j\|_2^2
\int_{\ell(I)/\rho}^{2\ell(I)/\rho}t^{-d-2|\gamma|}\,\frac{dt}{t}\bigg)^{1/2}
\le C_{d,\alpha}|Q(I)|^{-1/2}\ell(Q(I))^{-|\gamma|}.
$$
These scalar atoms have integral zero and support in $5I\subset Q(I)$.
Also,
$$
\sum_I\lambda_I(f)^2\frac{\mathbf1_{Q(I)}}{|Q(I)|}
\le C_d\rho^{-d}(A_\alpha^\nabla f)^2,
\qquad
\sum_I\lambda_I(f)^2=\frac14\|f\|_2^2<\infty.
$$
For a finite family $\mathscr I\subset\mathcal D^0$, apply
Lemma~\ref{sp:l:smooth} to each class $\mathscr I\cap\mathscr A_{\pmb\varepsilon}$,
whose supporting cubes lie in $\mathcal D^{\pmb\varepsilon}$
with multiplicity at most $16^d$.
Summing over $\pmb\varepsilon$, as in
\eqref{sp:e:class-reverse}--\eqref{sp:e:finite-reverse}, we obtain
$$
\bigg\|\sum_{I\in\mathscr I}\lambda_I(f)\mathsf a_I(f)\bigg\|_p
\le C_{d,\alpha}\sqrt p\,\|A_\alpha^\nabla f\|_p.
$$
Using the square-summability above and the energy constant $c_q$
for $q_j(-\cdot)$, we obtain $L^2$ convergence of the atomic sums to $f$
as in {\it Step 3} of Proposition~\ref{sp:p:reverse}.
Fatou's lemma and $L^p$ continuity then give the upper bound for the reverse constant.

For the lower bound, take
$\mathcal P_N(x)=\sum_{k=1}^N\cos(2\pi2^kx_1)$ as in Lemma~\ref{sp:l:lower}.
Then
$$
(\mathcal P_N*(\psi_1)_t)(y)
=-\sum_{k=1}^N b_{\mathrm P}(2\pi t2^k)\sin(2\pi2^ky_1),
\qquad
\mathcal P_N*(\psi_j)_t=0\quad(j\ge2).
$$
As in \eqref{sp:e:periodic-area}, summing the geometric series gives
$\|A_\alpha^\nabla\mathcal P_N\|_\infty\le C_{d,\alpha}\sqrt N$.
Differentiating $2\pi it\xi_j e^{-2\pi t|\xi|}$ gives, for
$0<c<C<\infty$, every multi-index $\gamma$ and $M>0$,
$$
\sup_{c\le|\xi|\le C}|\partial_\xi^\gamma[\widehat{\psi_j}(t\xi)]|
\le C_{d,c,C,\gamma,M}\min\{t,t^{-M}\}.
$$
Let $\eta_0$ denote the cutoff $\eta$ in Lemma~\ref{sp:l:lower}, {and let}
$g_L=\eta_0(\cdot/L)\mathcal P_N$.
We repeat the localization argument in {\it Step 3} of that lemma
for the vector kernel $(\psi_1,\ldots,\psi_d)$ and obtain
$$
\|A_\alpha^\nabla g_L\|_p
\le\bigl(C_{d,\alpha}\sqrt N+C_{N,d,\alpha}/L\bigr)L^{d/p}.
$$
Choose $N=\lceil p\rceil$ and then $L$ sufficiently large.
Together with \eqref{sp:e:lacunary-size}, this proves the lower bound.
\end{proof}

\section{Kernel and aperture dependence}\label{sp:sec:dependence}

We prove Theorems~\ref{sp:t:kernels} and \ref{sp:t:aperture}, first varying the kernel at a fixed aperture,
then varying the aperture for a fixed kernel.
Both arguments use the logarithmic example constructed below.

\subsection{Kernel dependence at a fixed aperture}\label{sp:sec:kernels}

For the order-$\sqrt p$ example, choose a nonzero real radial
$\sigma\in C_c^\infty(\R^d)$ with $\int\sigma=0$, and let
$$
\varphi^{(0)}=\sigma*\sigma.
$$
This kernel satisfies the basic assumptions.
Since $\widehat\sigma$ is real,
$$
(|D|^{2-d}\varphi^{(0)})(0)
=\int_{\R^d}|\xi|^{2-d}|\widehat\sigma(\xi)|^2\,d\xi>0.
$$
Thus $0\in\operatorname{supp}(|D|^{2-d}\varphi^{(0)})$, and Theorem~\ref{sp:t:general} implies
\begin{equation}\label{sp:e:good-kernel}
\mathfrak c_{p,d}^{\varphi^{(0)}}(\alpha)
\simeq_{d,\varphi^{(0)},\alpha}\sqrt p
\qquad(\alpha>0,\ p\ge2).
\end{equation}

\subsubsection*{Radial kernels and their one-dimensional integrals}

The logarithmic example will have a singularity in one coordinate.
We first integrate the kernel in the other coordinates.

For $d\ge2$, define the one-dimensional kernel
$$
\kappa(u)=\int_{\R^{d-1}}\varphi(u,z)\,dz;
$$
for $d=1$, let $\kappa=\varphi$.
It is real, even, smooth and compactly supported, with integral zero.
Since $\widehat\kappa(r)=\widehat\varphi(re_1)$ and $\varphi$ is radial and nonzero,
$\kappa$ is also nonzero.

Radial Fourier inversion expresses $|D|^{2-d}\varphi$ as a spherical average of $|D|\kappa$;
we will verify this in \eqref{sp:e:spherical-fourier}.
The following lemma determines when such an average vanishes on a ball.

\begin{lemma}\label{sp:l:spherical}
Let $q\ge2$ be an integer and let $U$ be a continuous even function on $\R$.
For every $a>0$, we have
$$
\int_{\mathbb S^{q-1}}U(x\cdot\theta)\,d\sigma(\theta)=0
\qquad(|x|<a)
$$
if and only if $U$ vanishes on $(-a,a)$.
\end{lemma}

\begin{proof}
Suppose that the spherical integral vanishes on $B(0,a)$.
For $r>0$, using spherical coordinates, evenness of $U$ and the substitution $v=r^2s^2$, we have
\begin{align*}
\int_{\mathbb S^{q-1}}U(re_1\cdot\theta)\,d\sigma(\theta)
= &C_q\int_{-1}^1U(rs)(1-s^2)^{(q-3)/2}\,ds\\
= &C_qr^{2-q}\int_0^{r^2}
U(\sqrt v)v^{-1/2}(r^2-v)^{(q-3)/2}\,dv,
\end{align*}
where $C_q>0$.
Let $\nu=(q-1)/2$ and $H(v)=U(\sqrt v)v^{-1/2}$.
The function $H$ is locally integrable on $[0,a^2)$, and our hypothesis says that
$$
\int_0^sH(v)(s-v)^{\nu-1}\,dv=0
\qquad(0<s<a^2).
$$
Choose an integer $k>\nu$ and fix $0<w<a^2$.
Multiplying by $(w-s)^{k-\nu-1}$ and integrating in $s$, we obtain
\begin{align*}
0
= &\int_0^w(w-s)^{k-\nu-1}
\int_0^sH(v)(s-v)^{\nu-1}\,dv\,ds\\
= &\int_0^wH(v)\int_v^w(s-v)^{\nu-1}(w-s)^{k-\nu-1}\,ds\,dv\\
= &B(\nu,k-\nu)\int_0^wH(v)(w-v)^{k-1}\,dv.
\end{align*}
For the last equality, use $s=v+(w-v)\tau$, $0<\tau<1$.
All integrals are absolutely convergent: the same calculation with $|H|$
is finite since $H\in L^1(0,w)$.
Differentiating $k$ times in $w$, in the sense of distributions, we obtain
$$
B(\nu,k-\nu)(k-1)!\,H(w)=0\qquad(0<w<a^2).
$$
Thus $U(u)=0$ for almost every $0<u<a$, and by continuity and evenness,
$U$ vanishes throughout $(-a,a)$.
The converse follows directly from $|x\cdot\theta|\le|x|<a$.
\end{proof}

\begin{proposition}\label{sp:p:marginal}
For every $\rho>0$, the function $|D|^{2-d}\varphi$ vanishes on $B(0,\rho)$
if and only if $|D|\kappa$ vanishes on $(-\rho,\rho)$.
In particular,
\begin{equation}\label{sp:e:integrated-distance}
\operatorname{dist}\bigl(0,\operatorname{supp}(|D|^{2-d}\varphi)\bigr)=\operatorname{dist}(0,\operatorname{supp}(|D|\kappa)).
\end{equation}
\end{proposition}

\begin{proof}
For $d=1$ the assertion is the definition.
For $d\ge2$, write the absolutely convergent Fourier integral in polar coordinates and use radiality:
\begin{align}\label{sp:e:spherical-fourier}
(|D|^{2-d}\varphi)(x)
= &\int_{\mathbb S^{d-1}}\int_0^\infty
r\widehat\varphi(re_1)e^{2\pi irx\cdot\theta}\,dr\,d\sigma(\theta)
= \frac12\int_{\mathbb S^{d-1}}(|D|\kappa)(x\cdot\theta)\,d\sigma(\theta).
\end{align}
For the second equality, use $\widehat\kappa(r)=\widehat\varphi(re_1)$ and pair the directions $\theta$ and $-\theta$.
Since $|D|\kappa$ is smooth and even, both assertions follow by applying
Lemma~\ref{sp:l:spherical} to \eqref{sp:e:spherical-fourier}.
\end{proof}

\subsubsection*{A logarithmic example in every dimension}

The logarithm is the standard BMO example with exponential decay of its distribution:
$$
\bigl|\{u\in(0,1):|\log u|>\tau\}\bigr|=e^{-\tau},\qquad \tau\ge0.
$$
Its $L^p(0,1)$ norm grows like $p$, as allowed by the John--Nirenberg inequality;
see also \cite[Section~8]{SVBMO}.
We use a smooth truncation of $\log|u|$ and construct a kernel whose one-dimensional
integral $\kappa$ satisfies $\kappa*\log|\cdot|=0$ near zero.
This cancellation is what makes the area-function norm bounded at small apertures.

{If $d$ is even, we construct the kernel in $\R^d$.
If $d$ is odd, we construct it in $\R^{d+1}$ and integrate out the last variable.
In either case the construction takes place in an even dimension $2m$,
where $\Delta^m\log|X|=c_m\delta_0$.
For $\Omega=\Delta^m\zeta$ with $\zeta$ radial and supported away from zero,
$\Omega*\log|\cdot|=c_m\zeta$ vanishes near zero.
The resulting kernel $\varphi$ has the same one-dimensional marginal $\kappa$ as $\Omega$.
Lemma~\ref{sp:l:spherical} then gives $\kappa*\log|\cdot|=0$ near zero.}

\begin{lemma}\label{sp:l:log-kernel}
Fix $0<a<A$ and an integer $m$ with $2m\ge d$.
Choose a nonzero radial function
$$
0\le\zeta\in C_c^\infty(\{X\in\R^{2m}:a<|X|<A\}),
\qquad \Omega=\Delta^m\zeta.
$$
When $2m>d$, write $X=(x,z)\in\R^d\times\R^{2m-d}$. Define
\begin{equation}\label{sp:e:bad-kernel}
\varphi(x)=\int_{\R^{2m-d}}\Omega(x,z)\,dz,
\qquad
\kappa(u)=\int_{\R^{2m-1}}\Omega(u,Z)\,dZ.
\end{equation}
If $2m=d$, the first integral means $\varphi=\Omega$.
Then $\varphi$ satisfies the basic assumptions, $\operatorname{supp}\varphi\subset B(0,A)$,
and $\operatorname{dist}\bigl(0,\operatorname{supp}(|D|^{2-d}\varphi)\bigr)\ge a$.
Moreover,
\begin{equation}\label{sp:e:log-vanish}
U(u):=\int_\R\kappa(v)\log|u-v|\,dv=0
\qquad(|u|<a).
\end{equation}
\end{lemma}

\begin{proof}
By construction, $\varphi\in C_c^\infty(B(0,A))$ is real and radial, with integral zero.
For $\xi\in\R^d$,
$$
\widehat\varphi(\xi)
=(-4\pi^2)^m|\xi|^{2m}\widehat\zeta(\xi,0).
$$
Since $\widehat\zeta(0)>0$, we have $\varphi\ne0$; similarly, $\kappa\ne0$.
Also, $\kappa(u)=\int_{\R^{d-1}}\varphi(u,z)\,dz$.

We use the distributional identity
\begin{equation}\label{sp:e:log-fundamental}
\Delta^m\log|X|
=c_m\delta_0,
\qquad
c_m=(-1)^{m-1}2^{2m-1}\pi^m(m-1)!.
\end{equation}
For $m=1$, this is $\Delta\log|X|=2\pi\delta_0$.
For $m\ge2$, by radial differentiation, we have
$$
\Delta\log|X|=(2m-2)|X|^{-2},
\qquad
\Delta|X|^\beta=\beta(\beta+2m-2)|X|^{\beta-2}\quad(X\ne0),
$$
and hence
$$
\Delta^{m-1}\log|X|
=(-1)^{m-2}2^{2m-3}(m-1)!(m-2)!\,|X|^{2-2m}.
$$
These iterations hold distributionally, since the boundary terms on
$|X|=\varepsilon$ tend to zero as $\varepsilon\to0$.
For the last Laplacian, we compute
$$
\int_{|X|=\varepsilon}\partial_r|X|^{2-2m}\,d\sigma(X)
=(2-2m)\sigma_{2m-1},
\qquad
\sigma_{2m-1}=\frac{2\pi^m}{(m-1)!}.
$$
Thus $\Delta|X|^{2-2m}=(2-2m)\sigma_{2m-1}\delta_0$,
which proves \eqref{sp:e:log-fundamental}.
Since $\Omega=\Delta^m\zeta$, integration by parts then yields
\begin{equation}\label{sp:e:radial-log}
\int_{\R^{2m}}\Omega(Y)\log|X-Y|\,dY=c_m\zeta(X).
\end{equation}

By rotation,
$$
\frac1{\sigma_{2m-1}}\int_{\mathbb S^{2m-1}}\log|V\cdot\theta|\,d\sigma(\theta)
=\log|V|+b_m
\qquad(V\ne0),
$$
where $b_m$ is finite since $\log|\theta_1|$ is integrable on the sphere.
Moreover,
$$
\int_{\mathbb S^{2m-1}}|\log|V\cdot\theta||\,d\sigma(\theta)
\le \sigma_{2m-1}|\log|V||+C_m,
\qquad V\ne0.
$$
By Fubini's theorem and radiality of $\Omega$, we have,
for each $\theta\in\mathbb S^{2m-1}$,
\begin{align*}
\int_{\R^{2m}}\Omega(Y)\log|(X-Y)\cdot\theta|\,dY
= &\int_\R\int_{\R^{2m-1}}\Omega(v,Z)
\log|X\cdot\theta-v|\,dZ\,dv\\
= &\int_\R\kappa(v)\log|X\cdot\theta-v|\,dv
=U(X\cdot\theta).
\end{align*}
Averaging in $\theta$ and using $\int\Omega=0$ and \eqref{sp:e:radial-log}, we obtain
$$
c_m\zeta(X)
=\frac1{\sigma_{2m-1}}\int_{\mathbb S^{2m-1}}U(X\cdot\theta)\,d\sigma(\theta).
$$
Since $\zeta=0$ on $B(0,a)$ and $U=\kappa*\log|\cdot|$ is smooth and even,
Lemma~\ref{sp:l:spherical} proves \eqref{sp:e:log-vanish}.

Finally, using $\mathcal F(\operatorname{p.v.}(1/u))(\xi)=-i\pi\operatorname{sgn}\xi$, we have
$$
U''=2\pi^2|D|\kappa.
$$
Hence $|D|\kappa=0$ on $(-a,a)$, and Proposition~\ref{sp:p:marginal} yields
$\operatorname{dist}\bigl(0,\operatorname{supp}(|D|^{2-d}\varphi)\bigr)\ge a$.
\end{proof}

\begin{remark}
For $m=1$, the cancellation in \eqref{sp:e:log-vanish} can also be read
from the classical logarithmic potential on an interval.
By \cite[Example~1.11]{Saff}, after scaling, we have
\begin{equation}\label{e:log}
\frac1\pi\int_{-r}^r
\frac{\log|u-v|}{\sqrt{r^2-v^2}}\,dv
=\log(r/2),
\qquad r>0,\quad |u|<r.
 \end{equation}
The right-hand side is independent of $u$.

With $v=r\cos\theta$, this gives
$$
\int_0^{2\pi}\log|u-r\cos\theta|\,d\theta
=2\pi\log(r/2),\qquad |u|<r.
$$
For $\Omega=\Delta\zeta$ on $\R^2$ and $|u|<a$, we obtain by polar coordinates
and integration by parts
\begin{align*}
\int_\R\kappa(v)\log|u-v|\,dv
&=\int_{\R^2}\Delta\zeta(X)\log(|X|/2)\,dX
=\int_{\R^2}\zeta(X)\Delta\log(|X|/2)\,dX=0,
\end{align*}
because $\log|X|$ is harmonic away from zero and
$0\notin\operatorname{supp}\zeta$.
\end{remark}

We next construct a fixed function $h$ with $\|h\|_p\ge cp$ and
$\|\mathcal A_{\varphi,\alpha}h\|_p\le C$ for $0<\alpha<a$.
Here we use $U=0$ on $(-a,a)$; the condition $|D|\kappa=0$ alone
would only make the even function $U$ constant there.

Before truncation, $\int_\R\kappa=0$ implies
$$
\int_\R\kappa(v)\log|y-tv|\,dv
=\int_\R\kappa(v)\log|y/t-v|\,dv
=U(y/t),\qquad t>0.
$$
For $\alpha<a$, this convolution is zero throughout the cone $|y|<\alpha t$.
In the proof below, we truncate the logarithm and estimate the area function
at all base points.

\begin{lemma}\label{sp:l:log-test}
{Fix $d\ge1$, $0<a<A$ and an integer $m$ with $2m\ge d$.
Let $\varphi$ be the kernel constructed in Lemma~\ref{sp:l:log-kernel}.}
There is a nonzero compactly supported function
$h\in\bigcap_{2\le p<\infty}L^p(\R^d)$, independent of $p$ and $\alpha$, such that
$$
\left\|h\right\|_{L^p(\R^d)}\ge cp
\qquad(p\ge2),
$$
and, for every $0<\alpha<a$,
\begin{equation}\label{sp:e:log-area}
\left\|\mathcal A_{\varphi,\alpha}h\right\|_{L^\infty(\R^d)}<\infty,
\qquad
0<\left\|\mathcal A_{\varphi,\alpha}h\right\|_{L^p(\R^d)}\le C_{d,\varphi,\alpha}
\quad(p\ge2).
\end{equation}
{For each fixed $0<\alpha<a$ and $2\le p<\infty$, an $L^p$ approximation
of $h$ gives a function $g_p\in\mathcal S_0(\R^d)$ satisfying}
$$
\left\|g_p\right\|_{L^p(\R^d)}\ge c'p,
\qquad
0<\left\|\mathcal A_{\varphi,\alpha}g_p\right\|_{L^p(\R^d)}\le C'_{d,\varphi,\alpha}.
$$
\end{lemma}

\begin{proof}
\smallskip
\noindent{\it Step 1.} We begin with a one-dimensional logarithmic test function.

Fix $\vartheta\in C_c^\infty(\R)$ with $0\le\vartheta\le1$,
$\vartheta\equiv1$ on $[-1,1]$, and
$\operatorname{supp}\vartheta\subset[-2,2]$. Let
$$
h_1(u)=\vartheta(u)\log|u|
\qquad {\rm and}\qquad
\kappa_s(u)=s^{-1}\kappa(u/s).
$$
Then
\begin{equation}\label{sp:e:log-size}
\|h_1\|_{L^p(\R)}^p
\ge2\int_0^{e^{-p}}|\log u|^p\,du
\ge2p^pe^{-p}.
\end{equation}
For the upper bound, the part with $1\le|u|\le2$ is bounded, while
$$
\int_0^1|\log u|^p\,du
=\int_0^\infty v^pe^{-v}\,dv
\le2(2p/e)^p.
$$
Hence $\|h_1\|_{L^p(\R)}\simeq p$.

Fix $0<\alpha\le a$.
We next prove $\|\mathcal A_{\kappa,\alpha}h_1\|_{L^\infty(\R)}\le C$.
The argument uses only that $\kappa\in C_c^\infty(\R)$ is even, has integral zero,
and that $U=\kappa*\log|\cdot|$ vanishes on $(-\alpha,\alpha)$. Choose $A\ge1$ such that $\operatorname{supp}\kappa\subset[-A,A]$ and let
 $s_0=[4(A+\alpha)]^{-1}$.
With $v=u+bs$, we have
\begin{align*}
\mathcal A_{\kappa,\alpha}h_1(u)^2
&= \int_0^\infty\int_{|v-u|<\alpha s}
|(h_1*\kappa_s)(v)|^2\,\frac{dv\,ds}{s^2}\\
&= \int_{-\alpha}^{\alpha}\int_0^{s_0}
|(h_1*\kappa_s)(u+bs)|^2\,\frac{ds}{s}\,db+\int_{-\alpha}^{\alpha}\int_{s_0}^\infty
|(h_1*\kappa_s)(u+bs)|^2\,\frac{ds}{s}\,db\\
&=:Term_1+Term_2.
\end{align*}
We first consider $Term_2$.
Note that
\begin{equation}\label{sp:e:log-conv-size}
\|h_1*\kappa_s\|_{L^\infty(\R)}
\le\frac{\|h_1\|_{L^1(\R)}\|\kappa\|_{L^\infty(\R)}}{s}
\le \frac{4(\log 2)\|\kappa\|_{L^\infty(\R)}}{s}=:{C_\kappa\over s},
\qquad s>0.
\end{equation}
It follows that
\[
Term_2
\le 2\alpha C_\kappa^2\int_{s_0}^{\infty}\frac{ds}{s^3}
=\frac{\alpha C_\kappa^2}{s_0^2}.
\]
This bound is independent of $u$.

We next estimate $Term_1$. We consider two cases.

\smallskip
\noindent{\it Case 1.} Suppose that $|u|<1/2$.
For $|b|<\alpha$, $0<s<s_0$ and $w\in\operatorname{supp}\kappa$, we have
\[
|u+s(b-w)|
\le |u|+s(A+\alpha)<\frac34.
\]
Since $\vartheta\equiv1$ on $[-1,1]$, we obtain
\begin{align}\label{sp:e:log-conv}
(h_1*\kappa_s)(u+bs)
= &\int_\R\kappa(w)h_1(u+s(b-w))\,dw\notag\\
= &\int_\R\kappa(w)
\{\log s+\log|b+u/s-w|\}\,dw
= U(b+u/s),
\end{align}
where the last equality uses $\int_\R\kappa=0$.

To estimate $U(b+u/s)$, we first prove the decay of $U$.
For $|v|\ge2A$ and $w\in\operatorname{supp}\kappa$,
$$
|v-\theta w|\ge|v|/2,\qquad 0\le\theta\le1.
$$
By $\int_\R\kappa=\int_\R w\kappa(w)\,dw=0$
and Taylor's formula,
\begin{align*}
|U(v)|
= \bigg|\int_\R\kappa(w)
\bigg(\log|v-w|-\log|v|+\frac wv\bigg)\,dw\bigg|
\leq \int_\R|\kappa(w)|w^2
\int_0^1\frac{1-\theta}{|v-\theta w|^2}\,d\theta\,dw
\le\frac C{|v|^2}.
\end{align*}
On the other hand,
$$
\sup_{|v|\le2A}|U(v)|
\le\|\kappa\|_{L^\infty(\R)}
\int_{-3A}^{3A}|\log|w||\,dw<\infty.
$$
This, combined with the decay above, implies
$
|U(v)|\le C(1+|v|)^{-2},
\ v\in\R.
$

Since $U$ is smooth and $U(b)=0$ for $|b|<\alpha$, for $0<\tau\le1$ we have
$$
|U(b\pm\tau)|
=|U(b\pm\tau)-U(b)|
\le\tau\sup_{|v|\le\alpha+1}|U'(v)|
\le C_{\kappa,\alpha}\tau.
$$
For $\tau\ge1$, the inequality
$1+\tau\le(1+\alpha)(1+|b\pm\tau|)$ implies
$$
|U(b\pm\tau)|
\le\frac{C_\kappa}{(1+|b\pm\tau|)^2}
\le\frac{C_\kappa(1+\alpha)^2}{(1+\tau)^2}
\le C_{\kappa,\alpha}\tau^{-2}.
$$
If $u=0$, then \eqref{sp:e:log-conv} gives
$Term_1=0$, since $U(b)=0$ for $|b|<\alpha$.
If $0<|u|<1/2$, the change of variables $\tau=|u|/s$ gives
\begin{align*}
Term_1
&=\int_{-\alpha}^{\alpha}\int_{|u|/s_0}^{\infty}
  |U(b+\operatorname{sgn}(u)\tau)|^2\,\frac{d\tau}{\tau}\,db
\le C_{\kappa,\alpha}
  \left(\int_0^1 \tau\,d\tau
        +\int_1^\infty \frac{d\tau}{\tau^5}\right)
\le C_{\kappa,\alpha}.
\end{align*}

\smallskip
\noindent{\it Case 2.} Suppose that $|u|\ge1/2$.
For $|b|<\alpha$, $0<s<s_0$ and $w\in\operatorname{supp}\kappa$,
\[
|u+\theta s(b-w)|
\ge |u|-s(A+\alpha)\ge\frac14,
\qquad 0\le\theta\le1.
\]
Since $\int_\R\kappa(w)\,dw=\int_\R w\kappa(w)\,dw=0$,
Taylor's formula gives
\begin{align*}
|(h_1*\kappa_s)(u+bs)|
&=\left|\int_\R\kappa(w)
  \bigl[h_1(u+s(b-w))-h_1(u)-s(b-w)h_1'(u)\bigr]\,dw\right|\\
&\le \frac{s^2}{2}\sup_{|v|\ge1/4}|h_1''(v)|
  \int_\R |\kappa(w)|(b-w)^2\,dw\\
&\le C_{\kappa,\alpha}s^2.
\end{align*}
Hence
\[
Term_1
\le C_{\kappa,\alpha}\int_0^{s_0}s^3\,ds
\le C_{\kappa,\alpha}.
\]

Combining the bounds for $\mathrm{Term}_1$ and $\mathrm{Term}_2$, we obtain
\begin{equation}\label{sp:e:one-log-area}
\|\mathcal A_{\kappa,\alpha}h_1\|_{L^\infty(\R)}
\le C_{\kappa,\alpha}.
\end{equation}

\smallskip
\noindent{\it Step 2.} We extend the test function to $\R^d$.

If $d=1$, take $h=h_1$. Since $\varphi=\kappa$ in this case,
the required bounds follow from \eqref{sp:e:log-size}
and \eqref{sp:e:one-log-area}.

Consider $d\ge2$. Fix $\chi\in C_c^\infty(\R^{d-1})$ with
$0\le\chi\le1$ and $\chi=1$ on $[-1,1]^{d-1}$. Set
\[
h(x_1,x')=h_1(x_1)\chi(x').
\]
Since $\|\chi\|_{L^p(\R^{d-1})}\ge1$, \eqref{sp:e:log-size} gives
$
\|h\|_{L^p(\R^d)}
=\|h_1\|_{L^p(\R)}\|\chi\|_{L^p(\R^{d-1})}
\ge cp
$ for $2\le p<\infty$.

It remains to bound $\mathcal A_{\varphi,\alpha}h(x)$ uniformly in $x$.
We split the integral at $t=1$.

\smallskip
\noindent{\it Case 1.} $t\ge1$.
Since $h\in L^1(\R^d)$, we have
$
|(h*\varphi_t)(y)|
\le t^{-d}\|h\|_{L^1(\R^d)}\|\varphi\|_{L^\infty(\R^d)}.
$
Using $|B(x,\alpha t)|=v_d(\alpha t)^d$, we obtain
\[
\int_1^\infty\int_{|y-x|<\alpha t}
|(h*\varphi_t)(y)|^2\,\frac{dy\,dt}{t^{d+1}}
\le C_{d,\varphi,\alpha}\int_1^\infty t^{-2d}\,\frac{dt}{t}
<\infty.
\]

\smallskip
\noindent{\it Case 2.} $0<t<1$.
Write $y=(y_1,y')$ and recall that
$\kappa(z_1)=\int_{\R^{d-1}}\varphi(z_1,z')\,dz'$.
Adding and subtracting $\chi(y')$ in the convolution gives
\begin{equation}\label{sp:e:log-product}
(h*\varphi_t)(y)
=\chi(y')(h_1*\kappa_t)(y_1)+E_t(y),
\end{equation}
where
\[
E_t(y)=\int_{\R^d}\varphi(z)h_1(y_1-tz_1)
       \bigl(\chi(y'-tz')-\chi(y')\bigr)\,dz.
\]

We first estimate the term containing $h_1*\kappa_t$.
For fixed $y_1$ and $t$, the set of $y'$ such that
$|(y_1,y')-x|<\alpha t$ has measure at most
$v_{d-1}(\alpha t)^{d-1}$. Since $|\chi|\le1$, it follows that
\begin{align*}
\int_0^1\int_{|y-x|<\alpha t}
  |\chi(y')(h_1*\kappa_t)(y_1)|^2\,\frac{dy\,dt}{t^{d+1}}
&\le v_{d-1}\alpha^{d-1}
  \int_0^1\int_{|y_1-x_1|<\alpha t}
  |(h_1*\kappa_t)(y_1)|^2\,\frac{dy_1\,dt}{t^2}\\
&\le v_{d-1}\alpha^{d-1}
  \mathcal A_{\kappa,\alpha}h_1(x_1)^2\\
&\le C_{d,\varphi,\alpha},
\end{align*}
where the last inequality follows from \eqref{sp:e:one-log-area}.

We next estimate $E_t$. By the mean value theorem,
$
|\chi(y'-tz')-\chi(y')|\le t|z'|\|\nabla\chi\|_\infty.
$
Since $\operatorname{supp}\varphi\subset B(0,A)$, this gives
\[
|E_t(y)|
\le Ct\int_{-A}^A|h_1(y_1-tv)|\,dv
=C\int_{y_1-At}^{y_1+At}|h_1(u)|\,du.
\]
Let $J=[y_1-At,y_1+At]$. On $|u|<t$, we have $h_1(u)=\log|u|$.
On $|u|\ge t$, the support of $h_1$ and the bound
$0\le\vartheta\le1$ give $|h_1(u)|\le|\log t|+\log 2$. Hence
\begin{align*}
\int_J|h_1(u)|\,du
&\le 2\int_0^t(-\log u)\,du+|J|(|\log t|+\log 2)
\le C_A t(1+|\log t|).
\end{align*}
Thus $|E_t(y)|\le Ct(1+|\log t|)$ uniformly in $y$, and
\[
\int_0^1\int_{|y-x|<\alpha t}|E_t(y)|^2\,\frac{dy\,dt}{t^{d+1}}
\le C_{d,\varphi,\alpha}\int_0^1 t(1+|\log t|)^2\,dt
<\infty.
\]
Applying $|a+b|^2\le2|a|^2+2|b|^2$ in \eqref{sp:e:log-product}
gives a uniform bound for the integral over $0<t<1$.
Combining the two cases, we obtain
\[
\|\mathcal A_{\varphi,\alpha}h\|_{L^\infty(\R^d)}
\le C_{d,\varphi,\alpha}.
\]

\smallskip
\noindent{\it Step 3.} We prove \eqref{sp:e:log-area} and choose $g_p\in\mathcal S_0(\R^d)$.

Since $h\in L^2(\R^d)$, we obtain from \eqref{sp:e:energy} and the $L^\infty$ bound above
\begin{align*}
\left\|\mathcal A_{\varphi,\alpha}h\right\|_{L^p(\R^d)}
\le &\left\|\mathcal A_{\varphi,\alpha}h\right\|_{L^\infty(\R^d)}^{1-2/p}
\left\|\mathcal A_{\varphi,\alpha}h\right\|_{L^2(\R^d)}^{2/p}
\le \max\bigl\{\|\mathcal A_{\varphi,\alpha}h\|_{L^\infty(\R^d)},
\|\mathcal A_{\varphi,\alpha}h\|_{L^2(\R^d)}\bigr\}
\le C_{d,\varphi,\alpha}.
\end{align*}
The area function has positive $L^2$ norm, so its $L^p$ norm is positive.
This proves \eqref{sp:e:log-area}.

For each fixed $p$, choose $g_p\in\mathcal S_0(\R^d)$ by the approximation
following Lemma~\ref{sp:l:continuity}, so that
$$
\left\|g_p-h\right\|_{L^p(\R^d)}
\le\min\{cp/2,B_{p,d,\varphi,\alpha}^{-1}\}.
$$
It follows that
$$
\left\|g_p\right\|_{L^p(\R^d)}
\ge\left\|h\right\|_{L^p(\R^d)}-\left\|g_p-h\right\|_{L^p(\R^d)}
\ge cp/2
$$
and
\begin{align*}
\left\|\mathcal A_{\varphi,\alpha}g_p\right\|_{L^p(\R^d)}
\le &\left\|\mathcal A_{\varphi,\alpha}h\right\|_{L^p(\R^d)}
+B_{p,d,\varphi,\alpha}\left\|g_p-h\right\|_{L^p(\R^d)}
\le C_{d,\varphi,\alpha}+1.
\end{align*}
The proof is complete.
\end{proof}

\begin{remark}\label{sp:r:log-endpoint}
The function $h$ in Lemma~\ref{sp:l:log-test} also explains the distinction between
John--Nirenberg and Chang--Wilson--Wolff integrability.
On $Q=(0,1)^d$, we have $h(x)=\log x_1$ and $h_Q=-1$.
Consequently,
$$
\frac{|\{x\in Q:|h(x)-h_Q|>\tau\}|}{|Q|}=e^{-(\tau+1)},\qquad \tau\ge1,
$$
and, for every $c>0$,
$$
\int_Q e^{c|h-h_Q|^2}\,dx
=\int_0^\infty e^{c(v-1)^2-v}\,dv=\infty.
$$
Nevertheless, $\mathcal A_{\varphi,\alpha}h\in L^\infty$ for $0<\alpha<a$ by \eqref{sp:e:log-area}.
Thus boundedness of these small-aperture area functions does not force exponential-square integrability.
\end{remark}

\begin{proof}[Proof of Theorem~\ref{sp:t:kernels}]
Use $\varphi^{(0)}$ from \eqref{sp:e:good-kernel}.
For the second kernel, choose $a>\alpha_0$, $A>a$, and $m=\lceil d/2\rceil$,
and let $\varphi^{(1)}$ be the kernel $\varphi$ in \eqref{sp:e:bad-kernel}.
Lemma~\ref{sp:l:log-test} shows that
$$
\mathfrak c_{p,d}^{\varphi^{(1)}}(\alpha_0)
\ge\frac{\|h\|_{L^p(\R^d)}}{\|\mathcal A_{\varphi^{(1)},\alpha_0}h\|_{L^p(\R^d)}}
\ge c_{d,\varphi^{(1)},\alpha_0}p.
$$
Proposition~\ref{sp:p:general-upper} proves the matching upper bound.
Both kernels are independent of $p$.
\end{proof}

\begin{remark}
The kernels in Theorem~\ref{sp:t:kernels} may both be normalized by $c_\varphi=1$.
Indeed, replace $\varphi$ by $\varphi/\sqrt{c_\varphi}$.
This changes the best reverse constant by the fixed factor $\sqrt{c_\varphi}$,
and does not change its order in $p$.
For these normalized kernels, \eqref{sp:e:energy} gives
$$
\mathfrak c_{2,d}^{\varphi^{(0)}}(\alpha_0)
=\mathfrak c_{2,d}^{\varphi^{(1)}}(\alpha_0)
=(v_d\alpha_0^d)^{-1/2}.
$$
Thus the different orders in $p$ persist under the same $L^2$ normalization.
The support distances of $|D|^{2-d}\varphi^{(0)}$ and $|D|^{2-d}\varphi^{(1)}$
are respectively zero and positive.
Since dilation multiplies this distance by the dilation factor,
these two kernels are not dilates of one another; see Subsection~\ref{sp:sec:apertures}.
\end{remark}

\subsection{Aperture dependence for a fixed kernel}\label{sp:sec:apertures}

For the fixed kernel $\varphi$ and the number $a>0$ in Lemma~\ref{sp:l:log-kernel},
by \eqref{sp:e:log-area} and Proposition~\ref{sp:p:general-upper}, we have
$$
\mathfrak c_{p,d}^{\varphi}(\alpha)
\simeq_{d,\varphi,\alpha}p,
\qquad 0<\alpha<a.
$$
By Theorem~\ref{sp:t:general}, we also have
$$
\mathfrak c_{p,d}^{\varphi}(\alpha)
\simeq_{d,\varphi,\alpha}\sqrt p,
\qquad \operatorname{dist}\bigl(0,\operatorname{supp}(|D|^{2-d}\varphi)\bigr)<c_d\alpha.
$$

\begin{proof}[Proof of Theorem~\ref{sp:t:aperture}]
Fix $0<a<A$, let $m=\lceil d/2\rceil$, and choose $\varphi$ from Lemma~\ref{sp:l:log-kernel}.
Define
$$
\alpha_-=a/2,
\qquad
\alpha_+=\max\{a,2\operatorname{dist}\bigl(0,\operatorname{supp}(|D|^{2-d}\varphi)\bigr)/c_d\}.
$$
Then $0<\alpha_-<\alpha_+<\infty$, independently of $p$.
For $0<\alpha\le\alpha_-$, we have $\alpha<a$; for $\alpha\ge\alpha_+$,
we have $\operatorname{dist}\bigl(0,\operatorname{supp}(|D|^{2-d}\varphi)\bigr)<c_d\alpha$.
The two estimates above therefore apply in the respective ranges,
with constants depending on $d,\varphi,\alpha$ but not on $p$.
\end{proof}

\subsubsection*{Kernel dilation and change of aperture}

At a fixed aperture, both orders can also be obtained from two dilates of one kernel.
For $r>0$, define
$$
\varphi^{[r]}(x)=r^{-d}\varphi(x/r).
$$
Then $(\varphi^{[r]})_t=\varphi_{rt}$.
Changing variables by $s=rt$ in the cone integral, we obtain
\begin{align}\label{sp:e:kernel-dilation}
\mathcal A_{\varphi^{[r]},\alpha}f(x)^2
= &\int_0^\infty\int_{|y-x|<\alpha t}
|f*\varphi_{rt}(y)|^2\,\frac{dy\,dt}{t^{d+1}}\notag\\
= &r^d\int_0^\infty\int_{|y-x|<(\alpha/r)s}
|f*\varphi_s(y)|^2\,\frac{dy\,ds}{s^{d+1}}
=r^d\mathcal A_{\varphi,\alpha/r}f(x)^2.
\end{align}
Taking the supremum of the corresponding norm ratios, we obtain
\begin{equation}\label{sp:e:constant-dilation}
\mathfrak c_{p,d}^{\varphi^{[r]}}(\alpha)
=r^{-d/2}\mathfrak c_{p,d}^{\varphi}(\alpha/r).
\end{equation}
Also, Fourier inversion and a change of variables in \eqref{sp:e:energy-constant} imply
$$
(|D|^{2-d}\varphi^{[r]})(x)=r^{-2}(|D|^{2-d}\varphi)(x/r),
\qquad
c_{\varphi^{[r]}}=c_\varphi.
$$
Consequently,
$$
\operatorname{dist}\bigl(0,\operatorname{supp}(|D|^{2-d}\varphi^{[r]})\bigr)
=r\operatorname{dist}\bigl(0,\operatorname{supp}(|D|^{2-d}\varphi)\bigr).
$$
Thus $\operatorname{dist}\bigl(0,\operatorname{supp}(|D|^{2-d}\varphi^{[r]})\bigr)<c_d\alpha$ is equivalent to
$\operatorname{dist}\bigl(0,\operatorname{supp}(|D|^{2-d}\varphi)\bigr)<c_d\alpha/r$, in agreement with \eqref{sp:e:constant-dilation}.

For example, fix $\alpha_0>0$ and take the kernel $\varphi$ and the number $a>0$
from Lemma~\ref{sp:l:log-kernel}.
Choose fixed $r_0,r_1>0$ with
$$
r_0\operatorname{dist}\bigl(0,\operatorname{supp}(|D|^{2-d}\varphi)\bigr)<c_d\alpha_0
\qquad{\rm and}\qquad
r_1>\alpha_0/a.
$$
Apply Theorem~\ref{sp:t:general} at aperture $\alpha_0/r_0$, and
Lemma~\ref{sp:l:log-test} and Proposition~\ref{sp:p:general-upper} at aperture $\alpha_0/r_1<a$.
By \eqref{sp:e:constant-dilation}, we obtain that for $2\le p<\infty$,
$$
\mathfrak c_{p,d}^{\varphi^{[r_0]}}(\alpha_0)
\simeq_{d,\varphi,\alpha_0,r_0}\sqrt p
\qquad{\rm and}\qquad
\mathfrak c_{p,d}^{\varphi^{[r_1]}}(\alpha_0)
\simeq_{d,\varphi,\alpha_0,r_1}p.
$$
These two dilates have the same value of $c_\varphi$ and are independent of $p$.

\subsection{Bounded conical area functions and aperture comparison}\label{sp:sub:bounded-area}

For a fixed kernel $\varphi$ satisfying the basic assumptions, consider the class
of functions with bounded conical area function,
\begin{equation}\label{sp:e:bounded-class}
\mathcal X_{\varphi,\alpha}
=\{f\in L^2(\R^d):\mathcal A_{\varphi,\alpha}f\in L^\infty(\R^d)\},
\qquad \alpha>0.
\end{equation}
Since $\Gamma_\alpha(x)\subset\Gamma_\beta(x)$ for $0<\alpha<\beta$,
we have
$
\mathcal X_{\varphi,\beta}\subseteq\mathcal X_{\varphi,\alpha}.
$
The following corollary shows that this inclusion can be strict.

For the fixed function in Lemma~\ref{sp:l:log-test}, the area function belongs to $L^\infty$
at small apertures but is unbounded at sufficiently large apertures.
The same function shows that the factor $\sqrt p$ in Lemma~\ref{sp:l:aperture} is necessary in general,
even for a fixed kernel and two fixed apertures.

\begin{corollary}\label{sp:c:endpoint}
Let $\varphi,h,a$ be as in Lemmas~\ref{sp:l:log-kernel} and \ref{sp:l:log-test}.
For every $0<\alpha<a$ and every $\beta$ with $\operatorname{dist}\bigl(0,\operatorname{supp}(|D|^{2-d}\varphi)\bigr)<c_d\beta$, we have
$$
\mathcal A_{\varphi,\alpha}h\in L^\infty(\R^d),
\qquad
\mathcal A_{\varphi,\beta}h\notin L^\infty(\R^d).
$$
Moreover, for this fixed pair of apertures,
\begin{equation}\label{sp:e:sharp-aperture}
\frac{\|\mathcal A_{\varphi,\beta}h\|_{L^p(\R^d)}}
{\|\mathcal A_{\varphi,\alpha}h\|_{L^p(\R^d)}}
\simeq_{d,\varphi,\alpha,\beta}\sqrt p,
\qquad p\ge2.
\end{equation}
In particular, the factor $\sqrt p$ in \eqref{sp:e:aperture} is optimal in general.
\end{corollary}

\begin{proof}
Since $\operatorname{dist}\bigl(0,\operatorname{supp}(|D|^{2-d}\varphi)\bigr)\ge a$ and $c_d<1$, we have $\beta>\alpha$.
Applying Proposition~\ref{sp:p:reverse} and Lemma~\ref{sp:l:aperture} to the function $h$
in Lemma~\ref{sp:l:log-test}, we obtain
$$
cp\le\left\|h\right\|_{L^p(\R^d)}
\le C\sqrt p\,\left\|\mathcal A_{\varphi,\beta}h\right\|_{L^p(\R^d)},
$$
and
$$
\left\|\mathcal A_{\varphi,\beta}h\right\|_{L^p(\R^d)}
\le C\sqrt p\,\left\|\mathcal A_{\varphi,\alpha}h\right\|_{L^p(\R^d)}
\le C\sqrt p.
$$
Consequently,
$$
\frac{\|\mathcal A_{\varphi,\beta}h\|_{L^p(\R^d)}}
{\|\mathcal A_{\varphi,\alpha}h\|_{L^p(\R^d)}}
\simeq\sqrt p,\qquad 2\le p<\infty.
$$
The kernel $\varphi$, the function $h$ and the apertures $\alpha,\beta$ are fixed,
and all comparison constants are independent of $p$.

Moreover,
$
\mathcal A_{\varphi,\alpha}h\in L^\infty(\R^d),
$ and $
\mathcal A_{\varphi,\beta}h\notin L^\infty(\R^d).
$

The first assertion was proved in Lemma~\ref{sp:l:log-test}.
By \eqref{sp:e:energy}, $\mathcal A_{\varphi,\beta}h\in L^2(\R^d)$.
If $\mathcal A_{\varphi,\beta}h$ also belonged to $L^\infty$, then
$$
\left\|\mathcal A_{\varphi,\beta}h\right\|_{L^p(\R^d)}
\le
\left\|\mathcal A_{\varphi,\beta}h\right\|_{L^\infty(\R^d)}^{1-2/p}
\left\|\mathcal A_{\varphi,\beta}h\right\|_{L^2(\R^d)}^{2/p}
\le C,
$$
contrary to the lower bound
$\|\mathcal A_{\varphi,\beta}h\|_{L^p(\R^d)}\ge c\sqrt p$ above.
\end{proof}

\begin{remark}\label{sp:r:reconstruction}
For the fixed kernel $\varphi$, function $h$ and apertures $\alpha,\beta$ in
Corollary~\ref{sp:c:endpoint},
$$
\alpha<a\le\operatorname{dist}\bigl(0,\operatorname{supp}(|D|^{2-d}\varphi)\bigr)<c_d\beta.
$$
By Proposition~\ref{sp:p:local-dual}, no normalized real radial dual kernel with integral zero
can be supported in $B(0,\alpha)$.
At aperture $\beta$, we may choose a dual with support small enough for the cone comparison.

The proof above, Lemma~\ref{sp:l:log-test} and Proposition~\ref{sp:p:general-upper} show that
$$
\left\|h\right\|_{L^p(\R^d)}\simeq p,
\qquad
\left\|\mathcal A_{\varphi,\alpha}h\right\|_{L^p(\R^d)}\simeq1,
\qquad
\left\|\mathcal A_{\varphi,\beta}h\right\|_{L^p(\R^d)}\simeq\sqrt p,
\qquad p\ge2,
$$
with constants independent of $p$.
Thus the same fixed function has reverse norm ratios of orders $p$ and $\sqrt p$
at the two apertures, respectively.
In particular, the order-$\sqrt p$ reverse estimate at aperture $\alpha$ is ruled out
by these norm ratios, independently of the choice of a reproducing formula.
\end{remark}

For the compactly supported kernel $\varphi$, our example exhibits the same change of boundedness
as the example of P.~Jones reported in \cite[p.~239]{CWW}.
Moreover, the norm ratio in \eqref{sp:e:sharp-aperture} has the exact order $\sqrt p$.

For the spatial-gradient Poisson area function in the Jones example,
the best reverse constant has order $\sqrt p$ at every positive aperture;
see Proposition~\ref{sp:p:cww}.
The following proposition shows that the class in \eqref{sp:e:bounded-class}
can also depend on the aperture under the origin-support condition for our compactly supported kernels.

\begin{proposition}\label{sp:p:bounded-area}
Fix $d\ge1$ and $0<\alpha<\beta$.
There are a kernel $\varphi$ satisfying the basic assumptions with
$$
(|D|^{2-d}\varphi)(0)\ne0
$$
and a nonzero bounded compactly supported function $f$ on $\R^d$ such that
$$
\mathcal A_{\varphi,\alpha}f\in L^\infty(\R^d),
\qquad
\mathcal A_{\varphi,\beta}f\notin L^\infty(\R^d).
$$
\end{proposition}

\begin{proof}
\noindent{\it Step 1.} We construct a one-dimensional example.
Choose
$$
\alpha<a<b<\beta,
\qquad
0\ne\zeta\in C_c^\infty((a,b)),\quad \zeta\ge0,
$$
{and let}
$$
F(u)=\zeta(u)-\zeta(-u),\qquad \kappa=F',\qquad g=\mathbf1_{(0,1)}.
$$
Then $\kappa$ satisfies the basic assumptions, and
\begin{equation}\label{sp:e:interval-convolution}
(g*\kappa_t)(y)=F(y/t)-F((y-1)/t).
\end{equation}
For $|z|<\alpha$, the condition $F(x/t+z)\ne0$ requires
$$
a-\alpha<|x|/t<b+\alpha.
$$
Thus, for $x\ne0$,
$$
\int_0^\infty|F(x/t+z)|^2\,\frac{dt}{t}
\le\|F\|_\infty^2\log\frac{b+\alpha}{a-\alpha};
$$
for $x=0$, the integral is zero.
Applying this bound to $x$ and $x-1$, we obtain
\begin{align}\label{sp:e:interval-small}
\mathcal A_{\kappa,\alpha}g(x)^2
= \int_{-\alpha}^{\alpha}\int_0^\infty
|F(x/t+z)-F((x-1)/t+z)|^2\,\frac{dt}{t}\,dz
\le 8\alpha\|F\|_\infty^2\log\frac{b+\alpha}{a-\alpha}.
\end{align}

At aperture $\beta$, $(tz,t)\in\Gamma_\beta(x)$ for $|z|<b$ whenever
$$
\frac{2|x|}{\beta-b}<t<\frac1{2b}.
$$
On this range, $F(z-1/t)=0$. Hence
$$
\mathcal A_{\kappa,\beta}g(x)^2
\ge\int_{2|x|/(\beta-b)}^{1/(2b)}\int_{-b}^b|F(z)|^2\,dz\,\frac{dt}{t}
=\|F\|_2^2\log\frac{\beta-b}{4b|x|}
$$
for $0<|x|<(\beta-b)/(4b)$.
At aperture $\alpha$, each term contributes over a scale interval of uniformly bounded $dt/t$ measure.
At aperture $\beta$, this measure grows like $\log(1/|x|)$ near the jump of $g$.

Using the identity $2\pi^2|D|\kappa=(\kappa*\log|\cdot|)''$
from the proof of Lemma~\ref{sp:l:log-kernel} and integrating by parts, we have
\begin{equation}\label{sp:e:interval-origin}
(|D|\kappa)(0)
=-\frac1{2\pi^2}\int_\R\frac{\kappa(v)}{v^2}\,dv
=-\frac1{\pi^2}\int_\R\frac{F(v)}{v^3}\,dv<0.
\end{equation}
Both integrals are ordinary integrals, since $\kappa$ and $F$ vanish near zero.
This proves the proposition for $d=1$.

\smallskip
\noindent{\it Step 2.} For $d\ge2$, we construct a radial kernel with one-dimensional integral $\kappa$.
Choose an odd integer $D=2m+1\ge d$, and for $r>0$ define
$$
q_0(r)=\kappa(r),\qquad
q_{j+1}(r)=-\frac{q_j'(r)}{2\pi r},\qquad 0\le j<m.
$$
Each $q_j$ is smooth and supported in $(a,b)$; extend it by zero at $r=0$.
For $r\ge0$,
$$
\int_{\R^2}q_{j+1}\bigl(\sqrt{r^2+|z|^2}\bigr)\,dz
=2\pi\int_r^\infty q_{j+1}(s)s\,ds=q_j(r).
$$
Consequently, $\Omega(X)=q_m(|X|)$ on $\R^D$ has one-dimensional integral $\kappa$.
{Define}
$$
\varphi(x)=\int_{\R^{D-d}}\Omega(x,z)\,dz,
$$
with $\varphi=\Omega$ when $D=d$.
The kernel $\varphi$ is real, radial and smooth, is supported in $B(0,b)$, and satisfies
\begin{equation}\label{sp:e:interval-marginal}
\int_{\R^{d-1}}\varphi(u,z)\,dz=\kappa(u).
\end{equation}
Thus $\varphi\ne0$ and $\int_{\R^d}\varphi=\int_\R\kappa=0$.
By \eqref{sp:e:spherical-fourier} and \eqref{sp:e:interval-origin},
$$
(|D|^{2-d}\varphi)(0)=\frac{\sigma_{d-1}}2(|D|\kappa)(0)\ne0.
$$

\smallskip
\noindent{\it Step 3.} We localize the one-dimensional example in the other variables.
Choose $0\le\eta\in C_c^\infty(\R^{d-1})$ with $\eta=1$ on $B(0,2)$, {and let}
$$
f(x)=g(x_1)\eta(x'),\qquad x'=(x_2,\ldots,x_d).
$$
For $0<t<1$, write
$$
(f*\varphi_t)(y)=\eta(y')(g*\kappa_t)(y_1)+E_t(y).
$$
By \eqref{sp:e:interval-marginal},
\begin{align*}
|E_t(y)|
= \bigg|\int_{\R^d}g(y_1-tv_1)
\bigl(\eta(y'-tv')-\eta(y')\bigr)\varphi(v)\,dv\bigg|
\le t\|\nabla\eta\|_\infty\int_{\R^d}|v'|\,|\varphi(v)|\,dv\le Ct.
\end{align*}
For $t\ge1$, we use $|f*\varphi_t|\le t^{-d}\|f\|_1\|\varphi\|_\infty$.
Changing variables by $y=x+tz$ in the cone integral, we obtain from \eqref{sp:e:interval-small}
$$
\mathcal A_{\varphi,\alpha}f(x)^2
\le C\mathcal A_{\kappa,\alpha}g(x_1)^2
+C\int_0^1t^2\,\frac{dt}{t}
+C\int_1^\infty t^{-2d}\,\frac{dt}{t}\le C.
$$

For the large-aperture estimate, choose $\delta,\varepsilon>0$ and $t_0>0$ such that
$$
(b+\delta)^2+\varepsilon^2<\beta^2,
\qquad
t_0<\min\{(2b)^{-1},(b+\varepsilon)^{-1}\}.
$$
If $|x'|<1$, $|z_1|<b$, $|z'|<\varepsilon$ and
$$
\frac{|x_1|}{\delta}<t<t_0,
$$
then $(y,t)\in\Gamma_\beta(x)$, where $y=(tz_1,x'+tz')$.
Moreover, $|x'+t(z'-v')|<2$ for $v\in\operatorname{supp}\varphi$, so the cutoff is one
throughout the convolution. By \eqref{sp:e:interval-convolution},
$$
(f*\varphi_t)(tz_1,x'+tz')=(g*\kappa_t)(tz_1)=F(z_1).
$$
Integrating over these $z_1,z'$ and $t$, we obtain
$$
\mathcal A_{\varphi,\beta}f(x)^2
\ge v_{d-1}\varepsilon^{d-1}\|F\|_2^2
\int_{|x_1|/\delta}^{t_0}\frac{dt}{t}
=c\log\frac{\delta t_0}{|x_1|},
\qquad 0<|x_1|<\delta t_0,\quad |x'|<1.
$$
Thus $\mathcal A_{\varphi,\beta}f$ is essentially unbounded.
\end{proof}

For the kernel in Proposition~\ref{sp:p:bounded-area}, we therefore have
$$
\mathcal X_{\varphi,\beta}\subsetneqq\mathcal X_{\varphi,\alpha},
\qquad
\mathfrak c_{p,d}^{\varphi}(\gamma)\simeq_{d,\varphi,\gamma}\sqrt p
\quad(\gamma>0,\ p\ge2),
$$
where the reverse estimate follows from Theorem~\ref{sp:t:general}.
Thus the class of functions with bounded conical area function may change with the aperture
even when the optimal order of the reverse constant does not.

\subsection{Product area functions}\label{sp:sub:product}

When the dilation parameters act on independent variables, the best reverse
constants multiply. We first explain why passing to Hilbert-valued functions
introduces no loss in the reverse constant.

Let $\varphi_i$ satisfy the basic assumptions on $\R^{d_i}$ and fix
$\alpha_i>0$, $i=1,2$. For $x=(x^{(1)},x^{(2)})\in\R^{d_1}\times\R^{d_2}$, define
\begin{align}\label{sp:e:product-area}
\mathcal A_{\mathrm{prod}}f(x)^2
= &\int_{\Gamma_{\alpha_1}(x^{(1)})}
\int_{\Gamma_{\alpha_2}(x^{(2)})}
\left|f*(\varphi_{1,t_1}\otimes\varphi_{2,t_2})
       (y^{(1)},y^{(2)})\right|^2\frac{dy^{(2)}\,dt_2}{t_2^{d_2+1}}
       \frac{dy^{(1)}\,dt_1}{t_1^{d_1+1}}.
\end{align}
Write
$$
\mathfrak c_p^{\mathrm{prod}}
=\sup_{0\ne f\in L^p(\R^{d_1+d_2})}
\frac{\|f\|_{L^p(\R^{d_1+d_2})}}
     {\|\mathcal A_{\mathrm{prod}}f\|_{L^p(\R^{d_1+d_2})}}.
$$

\begin{proposition}\label{sp:p:product}
For every $2\le p<\infty$,
\begin{equation}\label{sp:e:product-constant}
\mathfrak c_p^{\mathrm{prod}}
=\mathfrak c_{p,d_1}^{\varphi_1}(\alpha_1)
 \mathfrak c_{p,d_2}^{\varphi_2}(\alpha_2).
\end{equation}
\end{proposition}

\begin{proof}
\noindent{\it Step 1.} We extend the one-parameter reverse estimate to
Hilbert-valued functions with the same constant.

For a kernel $\varphi$ on $\R^d$, {write}
$
\mathcal K=L^2\bigl(B(0,\alpha)\times(0,\infty),\,dz\,dt/t\bigr)
$ and 
$T_\varphi f(x)(z,t)=(f*\varphi_t)(x+tz).
$

With the change of variables $y=x+tz$, we have
$\left\|T_\varphi f(x)\right\|_{\mathcal K}=\mathcal A_{\varphi,\alpha}f(x)$.

Let $F=(f_1,\ldots,f_N)$ and let $\gamma_1,\ldots,\gamma_N$
be independent standard circular complex Gaussian variables, normalized by
$\mathbb E|\gamma_j|^2=1$.
Write $a_p=(\mathbb E|\gamma_1|^p)^{1/p}$. Then
$$
\mathbb E\bigg|\sum_j\gamma_jf_j(x)\bigg|^p
=a_p^p\bigg(\sum_j|f_j(x)|^2\bigg)^{p/2}.
$$
For vectors $v_j\in\mathcal K$, choose an orthonormal basis $(e_\ell)$ of $\mathcal K$.
By Minkowski's inequality with exponent $p/2$,
\begin{align*}
\bigg(\mathbb E\bigg\|\sum_j\gamma_jv_j\bigg\|_{\mathcal K}^p\bigg)^{2/p}
&\le\sum_\ell
 \bigg(\mathbb E\bigg|\sum_j\gamma_j
       \langle v_j,e_\ell\rangle\bigg|^p\bigg)^{2/p}
=a_p^2\sum_{j,\ell}|\langle v_j,e_\ell\rangle|^2
 =a_p^2\sum_j\left\|v_j\right\|_{\mathcal K}^2.
\end{align*}
Apply the scalar reverse estimate to $\sum_j\gamma_jf_j$ and integrate
in the Gaussian variables. The preceding two formulas show that
\begin{align*}
a_p\bigg\|\bigg(\sum_j|f_j|^2\bigg)^{1/2}\bigg\|_p
&\le\mathfrak c_{p,d}^{\varphi}(\alpha)
 \bigg(\mathbb E\bigg\|\sum_j\gamma_jT_\varphi f_j
                   \bigg\|_{L^p(\R^d;\mathcal K)}^p\bigg)^{1/p}
\le a_p\mathfrak c_{p,d}^{\varphi}(\alpha)
 \bigg\|\bigg(\sum_j\left\|T_\varphi f_j\right\|_{\mathcal K}^2\bigg)^{1/2}\bigg\|_p.
\end{align*}
After cancelling $a_p$, we obtain
\begin{equation}\label{sp:e:hilbert-reverse}
\left\|F\right\|_{L^p(\R^d;\mathcal H)}
\le\mathfrak c_{p,d}^{\varphi}(\alpha)
\bigg\|\bigg(\int_{\Gamma_\alpha(x)}
 \left\|(F*\varphi_t)(y)\right\|_{\mathcal H}^2
 \frac{dy\,dt}{t^{d+1}}\bigg)^{1/2}\bigg\|_{L^p_x(\R^d)}
\end{equation}
first for $\mathcal H=\mathbb C^N$.
For real-valued functions the same calculation uses real Gaussian variables.

The area function on the right is also bounded on $L^p(\R^d;\mathcal H)$.
Indeed, Jensen's inequality and the scalar forward bound in
Lemma~\ref{sp:l:continuity} yield
$$
\bigg\|\bigg(\sum_j\left\|T_\varphi f_j\right\|_{\mathcal K}^2\bigg)^{1/2}\bigg\|_p
\le\bigg(\mathbb E\bigg\|\sum_j\gamma_jT_\varphi f_j
                     \bigg\|_{L^p(\R^d;\mathcal K)}^p\bigg)^{1/p}
\le B_{p,d,\varphi,\alpha}a_p\left\|F\right\|_{L^p(\R^d;\mathcal H)}.
$$
Orthogonal projections onto finite-dimensional subspaces and monotone
convergence extend both estimates to every separable Hilbert space
$\mathcal H$. In particular, the reverse constant in
\eqref{sp:e:hilbert-reverse} is independent of its dimension.

\smallskip
\noindent{\it Step 2.} We apply the two one-parameter estimates successively.

Let $T_1,T_2$ denote the operators above acting in $x^{(1)},x^{(2)}$,
with values in $\mathcal K_1,\mathcal K_2$, respectively.
First apply the scalar estimate in $x^{(1)}$ and integrate in $x^{(2)}$.
Next, for fixed $x^{(1)}$, apply \eqref{sp:e:hilbert-reverse} in $x^{(2)}$
to the $\mathcal K_1$-valued function $T_1f$.
Thus
\begin{align*}
\left\|f\right\|_{L^p(\R^{d_1+d_2})}
&\le\mathfrak c_{p,d_1}^{\varphi_1}(\alpha_1)
     \left\|T_1f\right\|_{L^p(\R^{d_1+d_2};\mathcal K_1)}
\le\mathfrak c_{p,d_1}^{\varphi_1}(\alpha_1)
     \mathfrak c_{p,d_2}^{\varphi_2}(\alpha_2)
     \left\|\mathcal A_{\mathrm{prod}}f\right\|_{L^p(\R^{d_1+d_2})}.
\end{align*}
The two square integrals are those in \eqref{sp:e:product-area}.
The forward bounds from {\it Step 1} ensure that all these norms are finite.
We first carry out these calculations for $f\in C_c^\infty(\R^{d_1+d_2})$,
then pass to general $f\in L^p(\R^{d_1+d_2})$ by approximation.

For the opposite inequality, take
$f(x)=f_1(x^{(1)})f_2(x^{(2)})$ with nonzero $f_i\in L^p(\R^{d_i})$.
The product structure implies
$$
\mathcal A_{\mathrm{prod}}f(x)
=\mathcal A_{\varphi_1,\alpha_1}f_1(x^{(1)})
 \mathcal A_{\varphi_2,\alpha_2}f_2(x^{(2)}),
$$
and hence
$$
\frac{\|f\|_{L^p(\R^{d_1+d_2})}}
     {\|\mathcal A_{\mathrm{prod}}f\|_{L^p(\R^{d_1+d_2})}}
=\prod_{i=1}^2
 \frac{\|f_i\|_{L^p(\R^{d_i})}}
      {\|\mathcal A_{\varphi_i,\alpha_i}f_i\|_{L^p(\R^{d_i})}}.
$$
Taking the supremum over $f_1,f_2$ proves \eqref{sp:e:product-constant}.
\end{proof}

At any fixed pair of apertures, Theorem~\ref{sp:t:kernels} therefore supplies
product kernels for which the optimal orders are
$$
\sqrt p\cdot\sqrt p=p,\qquad
p\cdot\sqrt p=p^{3/2},\qquad
p\cdot p=p^2.
$$
All these kernels are independent of $p$.
For every fixed pair of kernels, the order is $p$ when both apertures are
sufficiently large. Applying Theorem~\ref{sp:t:aperture} to both factors,
we also obtain a fixed product kernel whose reverse constant has order
$p^2$ at all sufficiently small common apertures and order $p$ at all
sufficiently large common apertures.

The same proof applies to any finite number $m$ of independent parameters:
the best reverse constant is the product of the $m$ one-parameter constants.
Choosing $k$ factors of order $p$ and the remaining $m-k$ factors of order
$\sqrt p$, where $0\le k\le m$, produces the order $p^{(m+k)/2}$.

For Zygmund dilations, the scales in the last two variables are $t,st$.
The second convolution therefore depends on the scale $s$ inside the
first square integral. The following sections treat this dependence
through intermediate square functions and localized dual kernels.

\section{Dyadic models for product and Zygmund dilations}\label{sec:dyadic}

For the Zygmund area function, the one-parameter results in
Sections~\ref{sp:sec:dual}--\ref{sp:sec:dependence}
provide the first-variable estimates and lower-bound examples.
Here we recall the one- and two-parameter Haar estimates from \cite{LPW}
and obtain the three-parameter estimate through intermediate square functions.
We then use the relation between the scales $s,t,st$ to prove
an order-$p$ estimate for the Zygmund Haar model.
These estimates will be applied to Haar expansions of smooth atoms.

As in \eqref{sp:e:cover}, we cover enlarged intervals such as $3I$ and $5I$
by dyadic intervals of comparable length.
Here we use the two grids
\begin{equation}\label{eqn:D-adjacent}
\mathcal D^\varepsilon
=\left\{2^{-m}\bigl([0,1)+n+(-1)^m\varepsilon/3\bigr):
m,n\in\mathbb Z\right\},
\qquad \varepsilon\in\{0,1\}.
\end{equation}
The grid $\mathcal D^0$ is the standard dyadic grid.
These grids are adjacent: every bounded interval is contained in an interval
of comparable length from $\mathcal D^0\cup\mathcal D^1$.

\begin{lemma}\label{l:adj-int-cover}
For every bounded interval $E$ of positive length, one can choose
$$
    E\subset I\in\mathcal D^0\cup\mathcal D^1,     \qquad  |E|\leq  \ell(I)\leq   6|E|.
$$
\end{lemma}

\begin{proof}
Choose the unique $m\in\mathbb Z$ such that $3|E|<2^{-m}\le6|E|$, and let
$$
     B_m^\varepsilon=2^{-m}(\mathbb Z+(-1)^m\varepsilon/3),   \qquad  \varepsilon\in\{0,1\}.
$$
Note that the two sets of points satisfy
$$
     \operatorname{dist}(B_m^0,B_m^1)=2^{-m}/3>|E|,
$$
so at most one of $B_m^0$ and $B_m^1$ intersects $E$. Hence
$$
      \varepsilon_E:=\min\left\{\varepsilon\in\{0,1\}:E\cap B_m^\varepsilon=\emptyset  \right\}
$$
is well defined. There are unique neighboring points
$p_1,p_2=p_1+2^{-m}\in B_m^{\varepsilon_E}$ such that
$$
     E\subset(p_1,p_2)\subset[p_1,p_2)=:I\in\mathcal  D^{ \varepsilon_E},
      \qquad  |E|\le\ell(I)=2^{-m}\le6|E|.
$$
The proof is complete.
\end{proof}

\begin{remark}
Allowing $\varepsilon\in\{0,1,2\}$ in \eqref{eqn:D-adjacent} recovers the three-grid family
used in \eqref{sp:e:cover} at $d=1$.
Here we use only $\varepsilon\in\{0,1\}$, since Lemma~\ref{l:adj-int-cover} already gives the required covering
with these two grids. The choice of a grid will be made separately in each coordinate.
The $2^3=8$ possible triples are
$$
        \mathcal  D^{  \vec{ \pmb\varepsilon}  }  :=\mathcal  D^{ \varepsilon_1}\times
        \mathcal  D^{ \varepsilon_2}   \times     \mathcal  D^{ \varepsilon_3},
         \qquad    \vec{ \pmb \varepsilon} =  (\varepsilon_1,  \varepsilon_2,  \varepsilon_3 )  \in  \{0,1\}^3.
$$
Thus $I\times J\times K\in\mathcal  D^{  \vec{ \pmb \varepsilon}}$ means that
$I\in\mathcal  D^{ \varepsilon_1}$, $J\in\mathcal  D^{ \varepsilon_2}$, and $K\in\mathcal  D^{ \varepsilon_3}$.

Besides, when no distinction between $\mathcal D^0$ and $\mathcal D^1$ is needed, we denote either grid by $\mathcal D$.
\end{remark}

For a dyadic interval $I\in \mathcal D$ with children $I_{\rm l},I_{\rm r}$, recall
\begin{equation}\label{eqn:Haar}
     h_I=|I|^{-1/2}(\mathbf 1_{I_{\rm l}}-\mathbf 1_{I_{\rm r}}),  
      \qquad   h_I^1=h_I, \qquad      h_I^0=|I|^{-1/2}\mathbf 1_I.	
\end{equation}

\subsection{Product Haar estimates and intermediate square functions}\label{sub:product}

Before turning to Zygmund dilations, we state the Haar estimates in the product setting:
the one- and two-parameter bounds from  our previous work \cite{LPW} and their three-parameter extension by iteration.

First, we have the following one-parameter sharp result for the sum of Haar functions.
\begin{proposition}\label{p:haar-one}
Let $\mathcal  D^{ \varepsilon_1}$ be a dyadic grid on $\R$, and let
$(a_I)_{I\in\mathcal  D^{ \varepsilon_1}}\subset \mathbb C$ be a finitely supported family.
Define
$$
     F_1(x_1)= \sum_I a_Ih_I(x_1),    \qquad 
    S_1F_1(x_1)   = \bigg( \sum_I|a_I|^2\frac{\mathbf 1_I(x_1)}{|I|}  \bigg)^{1/2}.
$$
For $2\le p<\infty$,
\begin{equation}\label{e:haar-one}
       \left\|   F_1 \right  \|_{L^p(\R)}  \leq  C \sqrt p \,   \left\|   S_1    F_1 \right\|_{L^p(\R)},
\end{equation}
where $C$ is absolute and 
the order $\sqrt p$ is optimal as $p\to\infty$.
\end{proposition}

\begin{proof}
The estimate \eqref{e:haar-one} is Lemma~\ref{sp:l:haar} with $d=1$.
The optimality of $\sqrt p$ follows from the one-dimensional example
in Remark~\ref{rem:product-sharpness} below.
\end{proof}

Next we turn to the product-type result.

\begin{proposition}\label{p:haar-two}
Let $\mathcal  D^{ \varepsilon_1},\mathcal  D^{ \varepsilon_2}$ be dyadic grids on $\R$, and let
$(a_{I,J})_{(I,J)\in\mathcal  D^{ \varepsilon_1}\times\mathcal  D^{ \varepsilon_2}}\subset \mathbb C$
be a finitely supported family.
Define
$$
     F_2(x_1,x_2)=\sum_I\sum_J a_{I,J}h_I(x_1)h_J(x_2),
$$
$$
     S_{12} F_2(x_1,x_2)   = \bigg(   \sum_I\sum_J |a_{I,J}|^2
       \frac{   \mathbf 1_I(x_1)}{|I|}   \frac{ \mathbf 1_J(x_2) } {|J|}   \bigg)^{1/2}.
$$
For $2\le p<\infty$,
$$
 \left\|   F_2 \right\|_{L^p(\R^2)}
\le Cp \, \left\|   S_{12}F_2 \right\|_{L^p(\R^2)},
$$
where $C$ is absolute and the order $p$ is optimal as $p\to\infty$.
\end{proposition}

\begin{proof}
This is the Haar case of \cite[Proposition~1.6]{LPW}. We prove it using the intermediate square function.
 As in the proof of Lemma~\ref{sp:l:haar}, it suffices to consider real coefficients.

For fixed $x_2$,
$$
    F_2(\cdot,x_2)  =\sum_I \bigg( \sum_J a_{I,J}h_J(x_2)  \bigg)   \, h_I.
$$
Its square function in $x_1$ is
\begin{align*}
    S_1F_2(\cdot,x_2)  = & \bigg(  \sum_I \bigg|  \bigg( \sum_J a_{I,J}h_J(x_2)  \bigg) \, h_I \bigg|^2  \bigg)^{1/2}
                                    =  \bigg( \sum_I\bigg|\sum_J a_{I,J}h_J(x_2)\bigg|^2  \frac{\mathbf 1_I}{|I|}   \bigg)^{1/2}.
\end{align*}
By Proposition~\ref{p:haar-one},
$$
     \left\|  F_2(\cdot, x_2) \right\|_{L^p(\R)}  \leq C\sqrt p\, \left\|   S_1F_2(\cdot, x_2) \right\|_{L^p(\R)}.
$$
{Take the maximal intervals among the contributing $I$'s and $J$'s in each grid.
On each product of two such intervals, \cite[Proposition~4.2]{LPW} gives
the slice estimate. Summing its $p$th powers in the second grid gives, for each $x_1$,}
\begin{equation}\label{eqn:S1-2}
     \left\|   S_1 F_2(x_1, \cdot) \right\|_{L^p(\R)} \le C \sqrt p\, \left\|   S_{12}F_2(x_1, \cdot) \right\|_{L^p(\R)}.
\end{equation}
{Integrating the two slice estimates gives the stated bound.}
\end{proof}

Similarly, iteration gives the corresponding three-parameter estimate.

\begin{proposition}\label{p:haar-three}
Let $\mathcal  D^{ \varepsilon_1},\mathcal  D^{ \varepsilon_2},\mathcal  D^{ \varepsilon_3}$
be dyadic grids on $\R$, and let
$ (a_{I,J,K})_{(I,J,K)\in \mathcal  D^{ \varepsilon_1}\times\mathcal  D^{ \varepsilon_2}\times\mathcal  D^{ \varepsilon_3}}\subset \mathbb C $
be a finitely supported family.
Define
$$
    F_3(x_1,x_2,x_3)  =\sum_I\sum_J\sum_K a_{I,J,K}h_I(x_1)h_J(x_2)h_K(x_3),
$$
$$
     S_{123}F_3(x_1,x_2,x_3) = \bigg(  \sum_I \sum_J \sum_K |a_{I,J,K}|^2 \frac{ \mathbf 1_I(x_1) } {|I|}
     \frac{\mathbf 1_J(x_2)} {|J|}   \frac{ \mathbf 1_K(x_3) } {|K|}   \bigg)^{1/2}.
$$
For $2\le p<\infty$,
$$
      \left\|   F_3 \right\|_{L^p(\R^3)}  \le Cp^{3/2}\, \left\|   S_{123}F_3 \right\|_{L^p(\R^3)},
$$
where $C$ is absolute and the order $p^{3/2}$ is optimal as $p\to\infty$.
\end{proposition}

\begin{proof}
{For real coefficients this is \cite[Proposition~4.7(a)]{LPW}.
The two intermediate square functions are}
$$
     S_1F_3(x_1,x_2,x_3)   = \bigg(  \sum_I \bigg| \sum_{J}\sum_K a_{I,J,K}h_J(x_2)h_K(x_3) \bigg|^2
                                              \frac{  \mathbf 1_I(x_1) } {|I|}   \bigg)^{1/2},
$$
and
$$
         S_{12}F_3(x_1,x_2,x_3)  = \bigg(  \sum_{I}\sum_J\bigg| \sum_K a_{I, J, K}  h_K(x_3) \bigg|^2
         \frac{\mathbf 1_I(x_1)} {|I|}     \frac{\mathbf 1_J(x_2)} {|J|}  \bigg)^{1/2}.
$$
By Proposition \ref{p:haar-one},
$$
      \left\|   F_3(\cdot,x_2,x_3) \right\|_{L^p(\R)} 
      \le C\sqrt p\, \left\|   S_1 F_3(\cdot, x_2, x_3) \right\|_{L^p(\R)}.
$$
By \eqref{eqn:S1-2}, 
$$
      \left\|   S_1F_3(x_1, \cdot, x_3) \right\|_{L^p(\R)}
      \le C\sqrt p\, \left\|   S_{12} F_3(x_1, \cdot, x_3) \right\|_{L^p(\R)},
$$
{and \cite[(4.31)]{LPW}, with $(I,J)$ indexing the square sum, gives}
$$
       \left\|   S_{12}F_3(x_1, x_2, \cdot) \right\|_{L^p(\R)}
       \le C\sqrt p\, \left\|   S_{123}F_3(x_1, x_2, \cdot) \right\|_{L^p(\R)}.
$$
{Integrating in the remaining variables gives the factor $p^{3/2}$.
For complex coefficients, apply the real estimate to the real and imaginary parts.}
\end{proof}

\medskip

\begin{remark}\label{rem:product-sharpness}
For sharpness, it suffices to work with the standard dyadic grid on $[0,1)$.
Let $p\geq   2 $ be an integer, take $N=p$, and define
$$
    G_N(x)=N^{-1/2}\sum_{k=0}^{N-1}  
    \sum_{ \substack{I\in\mathcal D^0,\ I\subseteq[0,1)\\     |I|=2^{-k}}}   |I|^{1/2}h_I(x).
$$
Then $S_1G_N=\mathbf 1_{[0,1)}$, and 
$
   G_N(x)=\sqrt N$ for $  x\in  [ 0,2^{-N} )$.
Hence
$$
      \left\|   G_N \right\|_{L^p(\R)}   \geq \sqrt N\,   2^{-N/p}   =\frac 12 \sqrt p.
$$
Define
$
     F_2(x_1,x_2) = G_N(x_1) G_N(x_2)$ and 
$     F_3(x_1,x_2,x_3) = G_N(x_1)  G_N(x_2)  G_N(x_3).
$

One may verify that 
$
     S_{12}F_2=\mathbf 1_{ [0,1)^2 },   
     S_{123}F_3=\mathbf 1_{ [0,1)^3 }
$
and 
$$
      \left\|   F_2 \right\|_{L^p(\R^2)}  = \left\|   G_N \right\|_{L^p(\R)}^2
        \geq \frac{1}{4}   p, \qquad  
      \left\|   F_3 \right\|_{L^p(\R^3)}   = \left\|   G_N \right\|_{L^p(\R)}^3 \geq \frac{1}{8}  p^{3/2}.
$$
Thus the orders $\sqrt p$, $p$, and $p^{3/2}$ are optimal.
\end{remark}

\subsection{Haar estimates for Zygmund dilations}\label{sub:zyg-haar}

For $ \vec{ \pmb\varepsilon}=(\varepsilon_1,\varepsilon_2,\varepsilon_3)\in\{0,1\}^3$
and $R=I\times J\times K\in\mathcal  D^{  \vec{ \pmb\varepsilon}}$, define
$$
    h_{R,Z}^{\eta}(x)   = h_I(x_1) h_J^{\eta_2}(x_2) h_K^{\eta_3} (x_3),   \qquad 
    \eta = (\eta_2, \eta_3) \in \mathsf H:=\{( 1,0),(0,1),(1,1) \},
$$
where $h_L^0=|L|^{-1/2} \mathbf 1_L$ and $h_L^1=h_L$ are given in \eqref{eqn:Haar}.
Clearly,
\begin{itemize}
\item $\operatorname{supp}h_{R,Z}^{\eta}\subseteq\overline R$, the closure of $R$; 
\smallskip

\item $|h_{R,Z}^{\eta}|=|R|^{-1/2}\mathbf 1_R$;
\smallskip

\item cancellation in $x_1$ and jointly in $(x_2,x_3)$:
$$
\int_I h_{R,Z}^{\eta}(x)\,dx_1=0,
 \qquad 
\int_{J\times K}  h_{R,Z}^{\eta}(x)\,  dx_2\,  dx_3=0.
$$
\end{itemize}
For $j=2,3$, the $x_j$-factor is cancellative   when $\eta_j=1$.
In particular, cancellation holds separately in both $x_2$ and $x_3$   when $\eta=(1,1)$.

For $ \vec{ \pmb\varepsilon}\in\{0,1\}^3$ and $q\in\mathbb Z$, let
$$
    \mathcal D_Z^{ \vec{ \pmb\varepsilon}, q}   =\left\{ R=I \times  J \times K
    \in\mathcal  D^{  \vec{ \pmb\varepsilon} }:  \,   \ell(K)=2^q\ell(I) \ell(J)\right\}.
$$
The Zygmund grid is
$ \mathcal D_Z^{ \vec{ \pmb\varepsilon}} :=\mathcal D_Z^{ \vec{ \pmb\varepsilon},0}$. 
For $R=I\times J\times K\in\mathcal D_Z^{ \vec{ \pmb\varepsilon},q}$,
$$
q=\log_2\frac{\ell(K)}{\ell(I)\ell(J)}.
$$
Haar expansions in each variable may produce rectangles with $\ell(K)\ne\ell(I)\ell(J)$.
We therefore need a reverse estimate uniform in $q$.

\begin{proposition}\label{p:zyg-haar}
Fix $ \vec{ \pmb\varepsilon}\in\{0,1\}^3$ and $q\in\mathbb Z$.
Let $\mathcal A^{(q)}\subset\mathcal D_Z^{ \vec{ \pmb\varepsilon},q}$ be finite,
and let $a_{R,\eta}\in\mathbb C$ for $R\in\mathcal A^{(q)}$ and $\eta\in\mathsf H$.
Define
$$
X^{(q)}(x)
=\sum_{R\in\mathcal A^{(q)}}\sum_{\eta\in\mathsf H}a_{R,\eta}h_{R,Z}^{\eta}(x),
 \qquad  
S^{(q)}(x)
= \bigg( 
\sum_{R\in\mathcal A^{(q)}}\sum_{\eta\in\mathsf H}
|a_{R,\eta}|^2\frac{\mathbf 1_R(x)}{|R|}
  \bigg)^{1/2}.
$$
For $2\le p<\infty$,
\begin{equation}\label{eqn:Square-Haar-Z}
	 \left\|   X^{(q)} \right\|_{L^p(\R^3)}
\le Cp\, \left\|   S^{(q)} \right\|_{L^p(\R^3)},
\end{equation}
where  the absolute constant $C$ is   independent of $ \vec{ \pmb\varepsilon}$, $q$,
$\mathcal A^{(q)}$, and the coefficients.
\end{proposition}

\begin{proof}
Take $a_{R,\eta}=0$ for $R\notin\mathcal A^{(q)}$.
For $\eta\in\mathsf H$, define
$$
X_\eta=\sum_{R\in\mathcal A^{(q)}}a_{R,\eta}h_{R,Z}^{\eta},
 \qquad 
S_\eta(x)
= \bigg( 
\sum_{R\in\mathcal A^{(q)}}|a_{R,\eta}|^2\frac{\mathbf 1_R(x)}{|R|}
  \bigg)^{1/2}.
$$
Then $X^{(q)}(x)=\sum_{\eta\in\mathsf H}X_\eta(x)$ and
$S^{(q)}(x)^2=\sum_{\eta\in\mathsf H}S_\eta(x)^2$.
Thus it suffices to prove \eqref{eqn:Square-Haar-Z} with
$X^{(q)}$ and $S^{(q)}$ replaced by $X_\eta$ and $S_\eta$, respectively.
We will apply Proposition~\ref{p:haar-two} in each case.

{\it Case 1. }
For $\eta=(1,0)$ or $(1,1)$, fix $x_3$ and let
$$
b_{I,J}^{\eta}(x_3)
=\sum_{K:\,\ell(K)=2^q\ell(I)\ell(J)}
a_{I\times J\times K,\eta}h_K^{\eta_3}(x_3).
$$
Then
$$
X_\eta(x_1,x_2,x_3)
=\sum_I\sum_J b_{I,J}^{\eta}(x_3)h_I(x_1)h_J(x_2).
$$
For fixed $I,J$, the corresponding intervals $K$  occurring in $b_{I,J}^{\eta}$
 have the same length and are mutually disjoint. By
 $|h_K^{\eta_3}|^2=|K|^{-1}\mathbf 1_K$, we obtain
$$
|b_{I,J}^{\eta}(x_3)|^2
=\sum_{\ell(K)=2^q\ell(I)\ell(J)}
|a_{I\times J\times K,\eta}|^2\frac{\mathbf 1_K(x_3)}{|K|},
$$
and hence by definition of $S_\eta$, we have 
$$
\sum_I\sum_J |b_{I,J}^{\eta}(x_3)|^2
\frac{\mathbf 1_I(x_1)}{|I|}
\frac{\mathbf 1_J(x_2)}{|J|}
=S_\eta(x_1,x_2,x_3)^2.
$$
Applying Proposition~\ref{p:haar-two} in $(x_1,x_2)$, we obtain
$$
 \left\|   X_\eta(\cdot,\cdot,x_3) \right\|_{L^p(\R^2)}
\le Cp\, \left\|   S_\eta(\cdot,\cdot,x_3) \right\|_{L^p(\R^2)}.
$$
Consequently, 
\begin{align*}
 \left\|   X_\eta \right\|_{L^p(\R^3)}^p
= \int_\R
 \left\|   X_\eta(\cdot,\cdot,x_3) \right\|_{L^p(\R^2)}^p\,dx_3
\leq 
(Cp)^p \left\|   S_\eta \right\|_{L^p(\R^3)}^p.
\end{align*}

{\it Case 2.}
For $\eta=(0,1)$,  now fix $x_2$,  let
$$
c_{I,K}(x_2)
=\sum_{J:\, \ell(J)=2^{-q}\ell(K)/\ell(I)}
a_{I\times J\times K,(0,1)}h_J^0(x_2).
$$
Then
$$
X_{(0,1)}(x_1,x_2,x_3)
=\sum_I\sum_K c_{I,K}(x_2)h_I(x_1)h_K(x_3).
$$
Similar to   {\it Case 1}, 
$$
|c_{I,K}(x_2)|^2
=\sum_{\ell(J)=2^{-q}\ell(K)/\ell(I)}
|a_{I\times J\times K,(0,1)}|^2\frac{\mathbf 1_J(x_2)}{|J|},
$$
and
$$
\sum_I\sum_K |c_{I,K}(x_2)|^2
\frac{\mathbf 1_I(x_1)}{|I|}
\frac{\mathbf 1_K(x_3)}{|K|}
=S_{(0,1)}(x_1,x_2,x_3)^2.
$$
Proposition~\ref{p:haar-two} in $(x_1,x_3)$ implies
$$
 \left\|   X_{(0,1)}(\cdot,x_2,\cdot) \right\|_{L^p(\R^2)}
\le Cp\, \left\|   S_{(0,1)}(\cdot,x_2,\cdot) \right\|_{L^p(\R^2)}.
$$
Again,
\begin{align*}
 \left\|   X_{(0,1)} \right\|_{L^p(\R^3)}^p
 \leq (Cp)^p \left\|   S_{(0,1)} \right\|_{L^p(\R^3)}^p.
\end{align*}
 
Both bounds in {\it Case 1} and {\it Case 2} are uniform in $ \vec{ \pmb\varepsilon}$ and $q$, and we complete the proof.
\end{proof}

\begin{remark}
The restriction to a fixed $q$ in Proposition~\ref{p:zyg-haar} is essential:
 $\big (\sum_q(S^{(q)})^2\big  )^{1/2}$
may not control $\sum_qX^{(q)}$, even in $L^2$.

Indeed, take $ \vec{ \pmb\varepsilon}=(0,0,0)$ and $I=J=[0,1)$.
For $n=0,\ldots,N-1$, choose
$$
a_{I\times J\times K,(1,0)}=|K|^{1/2},
 \qquad 
K\in\mathcal D^0,\quad K\subset[0,1),\quad |K|=2^{-n},
$$
with all other coefficients zero. Then
$$
X^{(-n)}(x)=h_I(x_1)h_J(x_2)\mathbf 1_{[0,1)}(x_3),
 \qquad 
S^{(-n)}(x)=\mathbf 1_{[0,1)^3}(x).
$$
Hence
$$
\bigg  \|  \sum_{n=0}^{N-1}X^{(-n)}\bigg    \|_{L^2(\R^3)}=N,
 \qquad 
\bigg  \|   \bigg( \sum_{n=0}^{N-1}(S^{(-n)})^2  \bigg)^{1/2}
\bigg    \|_{L^2(\R^3)}=\sqrt N.
$$

Besides, if only the case $\eta=(1,1)$ is considered, the sum over $q$ satisfies the product-type estimate
in Proposition~\ref{p:haar-three}, with optimal order $p^{3/2}$.
\end{remark}

\section{Smooth atoms for Zygmund dilations}\label{sec:atoms}

We construct smooth atoms from a Calder\'on reproducing formula at the coupled scales $s,t,st$,
as in the proof of Proposition~\ref{sp:p:reverse}; see also \cite{Wilson89,WuJGA}.
The atom indexed by $R=I\times J\times K$ is supported in
$(5I)\times(5J)\times(9K)$.
We cover this set by a rectangle $Q(R)$ from one of the eight product grids, with
$\ell(Q_3(R))=\ell(Q_1(R))\ell(Q_2(R))$.
Classifying the atoms by the chosen grid, we apply the Haar estimate
from Section~\ref{sec:dyadic} separately in each class.

Under \eqref{e:supp-origin}, we choose sufficiently small dual supports
to obtain the order-$p$ bound in Theorem~\ref{t:primal}.
In Section~\ref{sec:remaining}, we use the intermediate square function
$S_{1,\alpha}^{\psi_1}f$ to treat the two support conditions separately
and determine the dual support sizes allowed by the aperture.

\subsection{The Calder\'on reproducing formula}\label{sub:calderon}

Let $\psi_1\in C_c^\infty(\R)$ be real, even and nonzero,
and let $\psi_2\in C_c^\infty(\R^2)$ be real, radial and nonzero, with
$$
\int_\R\psi_1=0 ,
 \qquad  
\int_{\R^2}\psi_2=0.
$$
Denote
\begin{equation}\label{eqn:constant-psi}
	c_{\psi_1}=\int_0^\infty|\widehat\psi_1(r)|^2\,\frac{dr}{r},
 \qquad 
c_{\psi_2}=\int_0^\infty|\widehat\psi_2(r,0)|^2\,\frac{dr}{r},
\end{equation}
where $0<c_{\psi_1},\,c_{\psi_2}<\infty$.

To construct the Calder\'on reproducing formula in the Zygmund dilation setting,
 choose real, even $\Phi_1\in C_c^\infty(\R)$ and real, radial $\Phi_2\in C_c^\infty(\R^2)$ such that
$$
\int_\R\Phi_1=0,
 \qquad  
\int_{\R^2}\Phi_2=0,
$$
and
$$
\int_0^\infty\widehat\psi_1(r)\widehat\Phi_1(r)\,\frac{dr}{r}=1,
 \qquad 
\int_0^\infty\widehat\psi_2(r,0)\widehat\Phi_2(r,0)\,\frac{dr}{r}=1.
$$
These are the one-parameter normalizations in \eqref{sp:e:dual}, in dimensions one and two.
For example, one may take $\Phi_j=\psi_j/c_{\psi_j}$, as desired.
Define $c_{\Phi_1},c_{\Phi_2}$ by the same formulas as $c_{\psi_1},c_{\psi_2}$.
For $s,t>0$, let
$$
\psi_{s,t}(x)=\frac1{s^2t^2}
\psi_1\left(\frac{x_1}{s}\right)
\psi_2\left(\frac{x_2}{t},\frac{x_3}{st}\right), \qquad  
\Phi_{s,t}(x)=\frac1{s^2t^2}
\Phi_1\left(\frac{x_1}{s}\right)
\Phi_2\left(\frac{x_2}{t},\frac{x_3}{st}\right).
$$
We use $\Phi_{s,t}$ to reconstruct $f$ from the convolutions $f*\psi_{s,t}$.

For each fixed $s>0$, radiality reduces the $t$-integral to the two-dimensional one-parameter normalization.
For $\xi_1\ne0$ and $(\xi_2,\xi_3)\ne(0,0)$, {write}
$r=t(\xi_2^2+s^2\xi_3^2)^{1/2}$. Then
\begin{align*}
\int_0^\infty\int_0^\infty
\widehat\psi_{s,t}(\xi)\widehat\Phi_{s,t}(\xi)\,\frac{dt}{t}\frac{ds}{s}
= &\int_0^\infty\widehat\psi_1(s\xi_1)\widehat\Phi_1(s\xi_1)\,
 \bigg( \int_0^\infty\widehat\psi_2(r,0)\widehat\Phi_2(r,0)\,\frac{dr}{r}  \bigg)\frac{ds}{s}\\
= &\int_0^\infty\widehat\psi_1(r)\widehat\Phi_1(r)\,\frac{dr}{r}=1.
\end{align*}
The same change of variables and Plancherel's theorem show that, for $g\in L^2(\R^3)$,
\begin{equation}\label{e:kernel-l2}
\int_0^\infty\int_0^\infty
 \left\|   g*\psi_{s,t} \right\|_{L^2(\R^3)}^2\,\frac{dt}{t}\frac{ds}{s}
=c_{\psi_1}c_{\psi_2} \left\|   g \right\|_{L^2(\R^3)}^2.
\end{equation}
The identity \eqref{e:kernel-l2} also holds with $\psi$ replaced by $\Phi$.
For $f\in L^2(\R^3)$, the corresponding reproducing formula is
\begin{equation}\label{e:calderon-z}
f(x)=\int_0^\infty\int_0^\infty\int_{\R^3}
(f*\psi_{s,t})(y)\Phi_{s,t}(x-y)\,dy\,\frac{dt}{t}\frac{ds}{s},
\end{equation}
with convergence in $L^2$.
Indeed, for $N\ge2$, let
$$
f_N(x)=\int_{N^{-1}}^N\int_{N^{-1}}^N\int_{\R^3}
(f*\psi_{s,t})(y)\Phi_{s,t}(x-y)\,dy\,\frac{dt}{t}\frac{ds}{s},
 \qquad 
E_N=(0,\infty)^2\setminus[N^{-1},N]^2.
$$
For almost every $\xi\in\R^3$,  
\begin{align*}
\widehat f_N(\xi)-\widehat f(\xi)
= &\widehat f(\xi)\,
 \bigg( 
\int_{N^{-1}}^N\int_{N^{-1}}^N
\widehat\psi_{s,t}(\xi)\widehat\Phi_{s,t}(\xi)
\,\frac{dt}{t}\frac{ds}{s}-1
  \bigg)
= -\widehat f(\xi)\iint_{E_N}
\widehat\psi_{s,t}(\xi)\widehat\Phi_{s,t}(\xi)
\,\frac{dt}{t}\frac{ds}{s},
\end{align*}
and
$$
\int_0^\infty\int_0^\infty
|\widehat\Phi_{s,t}(\xi)|^2\,\frac{dt}{t}\frac{ds}{s}
=c_{\Phi_1}c_{\Phi_2}.
$$
We estimate $f_N-f$ over the omitted scales by Plancherel and Cauchy--Schwarz:

\begin{align*}
 \left\|   f_N-f \right\|_{L^2(\R^3)}^2
= &\int_{\R^3}|\widehat f(\xi)|^2
\bigg|\iint_{E_N}
\widehat\psi_{s,t}(\xi)\widehat\Phi_{s,t}(\xi)\,
\frac{dt}{t}\frac{ds}{s}\bigg|^2\,d\xi\\
\leq \,  &c_{\Phi_1}c_{\Phi_2}
\int_{\R^3}|\widehat f(\xi)|^2
\iint_{E_N}|\widehat\psi_{s,t}(\xi)|^2
\,\frac{dt}{t}\frac{ds}{s}\,d\xi\\
= \,  &c_{\Phi_1}c_{\Phi_2}
\iint_{E_N} \left\|   f*\psi_{s,t} \right\|_{L^2(\R^3)}^2
\,\frac{dt}{t}\frac{ds}{s}
\to 0.
\end{align*}
The last integral tends to zero by \eqref{e:kernel-l2} and dominated convergence.
Thus \eqref{e:calderon-z} holds in $L^2(\R^3)$.

\smallskip

\subsection{Smooth atomic decomposition}\label{sub:atoms}

We next decompose $f$ into smooth atomic functions using \eqref{e:calderon-z}.
Choose the smallest $\delta\in2^{\mathbb Z}$ such that
$
\operatorname{supp}\Phi_1\subset[-\delta,\delta]
$ and 
$\operatorname{supp}\Phi_2\subset[-\delta,\delta]\times[-\delta^2,\delta^2].
$

Let
$
\mathscr A_0=\big \{R=I\times J\times K:
I,J,K\in\mathcal D^0,\quad \ell(K)=\ell(I)\ell(J)\big  \},
$

For each $R=I\times J\times K\in\mathscr A_0$, define
$$
\mathcal W(R)=R\times
\left[\frac{\ell(I)}{\delta},\frac{2\ell(I)}{\delta}\right)
\times\left[\frac{\ell(J)}{\delta},\frac{2\ell(J)}{\delta}\right).
$$
Due to $\delta\in 2^{\mathbb Z}$, one may verify that 
\begin{equation}\label{e:W-decomp}
	\R^3\times(0,\infty)^2= \bigcup_{R\in\mathscr A_0}\mathcal W(R) \qquad  \text{(disjoint)}.
\end{equation}
Indeed, fix $(y,s,t)\in\R^3\times(0,\infty)^2$.
There are unique integers $k,l\in\mathbb Z$ such that
$$
2^k\le\delta s<2^{k+1},
 \qquad 
2^l\le\delta t<2^{l+1}.
$$
For these $k,l$, there are unique $I,J,K\in\mathcal D^0$
such that
$$
y=(y_1,y_2,y_3)\in I\times J \times K \quad\text{with} \quad
\ell(I)=2^k,\quad
\ell(J)=2^l,\quad
\ell(K)=2^{k+l}=\ell(I)\ell(J).
$$
Hence $R=I\times J\times K\in\mathscr A_0$ and
$(y,s,t)\in\mathcal W(R)$.
The reverse inclusion is trivial and  the sets $\mathcal W(R)$ are obviously pairwise disjoint.

 Thus, for $f\in L^2(\R^3)$, we use \eqref{e:calderon-z} to write
\begin{align}\label{e:atom-sum}
f(x)
 = &\sum_{R\in\mathscr A_0}\widetilde{\mathsf{a}}_R(f)(x),
\end{align}
where
$$
\widetilde{\mathsf{a}}_R(f)(x)=\int_{\mathcal W(R)}
(f*\psi_{s,t})(y)\Phi_{s,t}(x-y)\,dy\,\frac{dt}{t}\frac{ds}{s}.
$$
We will show in Lemma \ref{l:atomic-l2} that the convergence in \eqref{e:atom-sum} is in $L^2$ norm, and consequently, a subsequence of its finite partial sums   converges to $f$ almost everywhere.

\smallskip

For an interval $I$ and $\rho\ge1$, let $\rho I$ be the closed interval with the same centre as $I$
and length $ \rho\ell(I)$.
We have the following properties.

\begin{lemma}\label{l:atom-prop}
For every $f\in L^2(\R^3)$ and $R=I\times J\times K\in\mathscr A_0$,
the function $\widetilde{\mathsf{a}}_R(f)$ is smooth and satisfies:
\begin{enumerate}
\item[(i)]
$
\operatorname{supp}\widetilde{\mathsf{a}}_R(f)\subset(5I)\times(5J)\times(9K).
$

\smallskip

\item[(ii)] For $\gamma\in\mathbb N_0^3$ with $|\gamma|\le3$,
$
|\partial^\gamma \widetilde{\mathsf{a}}_R(f)(x)|
\le C_\gamma\lambda_R(f)|R|^{-1/2}
\ell(I)^{-\gamma_1}\ell(J)^{-\gamma_2}\ell(K)^{-\gamma_3},
$
where
\begin{equation}\label{e:atom-coeff}
\lambda_R(f)= \bigg( \int_{\mathcal W(R)}
|(f*\psi_{s,t})(y)|^2\,dy\,\frac{dt}{t}\frac{ds}{s}  \bigg)^{1/2}.
\end{equation}
In particular, $ \left\|   \widetilde{\mathsf{a}}_R(f) \right\|_\infty\le C\lambda_R(f)|R|^{-1/2}$.

\smallskip

\item[(iii)] For all $x,x'\in\R^3$,
$$
\left |\widetilde{\mathsf{a}}_R(f)(x)-\widetilde{\mathsf{a}}_R(f)(x') \right|
\le C\lambda_R(f)|R|^{-1/2}
\min\left\{1,
\max\left\{
\frac{|x_1-x'_1|}{\ell(I)},
\frac{|x_2-x'_2|}{\ell(J)},
\frac{|x_3-x'_3|}{\ell(K)}
\right\}\right\}.
$$
\item[(iv)] For fixed remaining variables,
$$
\int_\R \widetilde{\mathsf{a}}_R(f) \,dx_1=0,
 \qquad 
\int_{\R^2}\widetilde{\mathsf{a}}_R(f) \,dx_2\,dx_3=0.
$$
\end{enumerate}
The constants  above are independent of $f$ and $R$.
\end{lemma}

\begin{proof}
(i) Note that
$\operatorname{supp}\Phi_{s,t}
\subset[-\delta s,\delta s]\times[-\delta t,\delta t]
\times[-\delta^2st,\delta^2st]
$,  for any $(y,s,t)\in\mathcal W(R)$, the condition $x-y\in\operatorname{supp}\Phi_{s,t}$  implies
$$
|x_1-y_1|\le\delta s\le2\ell(I),
 \qquad 
|x_2-y_2|\le\delta t\le2\ell(J),
 \qquad 
|x_3-y_3|\le\delta^2st\le4\ell(K).
$$
Since $y\in R$, we obtain (i).

\smallskip

(ii)
For $|\gamma|\le3$,  it follows from Cauchy--Schwarz and \eqref{e:atom-coeff} that
\begin{align*}
\left|\partial^\gamma \widetilde{\mathsf{a}}_R(f)(x)\right|
\leq\,   &\lambda_R(f) \,\bigg( \int_{\mathcal W(R)}
\left|\partial^\gamma\Phi_{s,t}(x-y)\right|^2\,dy\,\frac{dt}{t}\frac{ds}{s}  \bigg)^{1/2}\\
\leq \,  &(\log2)\lambda_R(f)  \,    |R|^{1/2}
\sup_{(y,s,t)\in\mathcal W(R)}  \left|\partial^\gamma\Phi_{s,t}(x-y)\right|\\
\leq \,  &C_{\Phi,\gamma}\lambda_R(f)  \,   |R|^{1/2}
\left(\frac{\ell(I)}{\delta}\right)^{-2-\gamma_1-\gamma_3}
\left(\frac{\ell(J)}{\delta}\right)^{-2-\gamma_2-\gamma_3}\\
= \,  &C_{\Phi,\gamma}\delta^{4+\gamma_1+\gamma_2+2\gamma_3}
\lambda_R(f)    \,   |R|^{-1/2}
\ell(I)^{-\gamma_1}\ell(J)^{-\gamma_2}\ell(K)^{-\gamma_3}.
\end{align*}
This proves (ii). Moreover, the same calculation applies to every derivative, so $\widetilde{\mathsf{a}}_R(f)$ is smooth.

\smallskip

(iii) As a consequence of 
  (ii),
\begin{align*}
\left|\widetilde{\mathsf{a}}_R(f)(x)-\widetilde{\mathsf{a}}_R(f)(x')  \right|
\leq &\sum_{j=1}^3   \left|x_j-x'_j\right| \left\|   \partial_j \widetilde{\mathsf{a}}_R(f) \right\|_\infty\\
\leq &C\lambda_R(f)|R|^{-1/2}
\left(\frac{|x_1-x'_1|}{\ell(I)}+
\frac{|x_2-x'_2|}{\ell(J)}+
\frac{|x_3-x'_3|}{\ell(K)}\right).
\end{align*}
We also have $|\widetilde{\mathsf{a}}_R(f)(x)-\widetilde{\mathsf{a}}_R(f)(x')|\le2 \left\|   \widetilde{\mathsf{a}}_R(f) \right\|_\infty$. Combining the two bounds proves (iii).

\smallskip
(iv)
Finally, $\int\Phi_1=\int\Phi_2=0$ implies
$$
\int_\R\Phi_{s,t}(x-y)\,dx_1=0,
 \qquad 
\int_{\R^2}\Phi_{s,t}(x-y)\,dx_2\,dx_3=0.
$$
Hence
$$
\int_\R \widetilde{\mathsf{a}}_R(f)(x)\,dx_1
=\int_{\mathcal W(R)}(f*\psi_{s,t})(y)\,
 \bigg( \int_\R\Phi_{s,t}(x-y)\,dx_1  \bigg)
\,dy\,\frac{dt}{t}\frac{ds}{s}=0,
$$
and
$$
\int_{\R^2}\widetilde{\mathsf{a}}_R(f)(x)\,dx_2\,dx_3
=\int_{\mathcal W(R)}(f*\psi_{s,t})(y) \,
 \bigg( \int_{\R^2}\Phi_{s,t}(x-y)\,dx_2\,dx_3  \bigg)
\,dy\,\frac{dt}{t}\frac{ds}{s}=0.
$$
The proof is complete.
\end{proof}

For $\lambda_R(f)>0$, let
$$
\mathsf{a}_R(f)(x)=\frac1{\lambda_R(f)}\int_{\mathcal W(R)}
(f*\psi_{s,t})(y)\Phi_{s,t}(x-y)\,dy\,\frac{dt}{t}\frac{ds}{s},
$$
and  let $\mathsf{a}_R(f)=0$ for  $\lambda_R(f)=0$.

\begin{lemma}\label{l:atomic-l2}
Let $f\in L^2(\R^3)$, and let $\lambda_R$ and $\mathsf{a}_R$ be
defined as above. Then
$$
f=\sum_{R\in\mathscr A_0}\lambda_R(f)  \mathsf{a}_R(f)
 \qquad \text{in }L^2(\R^3).
$$
More precisely, for every
finite subfamily $\mathscr A\subset\mathscr A_0$,
\begin{equation}\label{e:remainder}
	\bigg  \|  f-\sum_{R\in\mathscr A}\lambda_R  (f)   \mathsf{a}_R
  (f)   \bigg    \|_{L^2(\R^3)}
\leq
(c_{\Phi_1}c_{\Phi_2})^{1/2}\,
 \bigg( \sum_{R\in\mathscr A_0\setminus\mathscr A}
\lambda_R^2(f)  \bigg)^{1/2}.
\end{equation}
\end{lemma}

\begin{proof}
We use the duality argument from {\it Step 3} of Proposition~\ref{sp:p:reverse},
with the energy identity \eqref{e:kernel-l2}.
Write $d\nu=dy\,\frac{dt}{t}\frac{ds}{s}$.
Note that $\R^3\times(0,\infty)^2= \bigcup_{R\in\mathscr A_0}\mathcal W(R)$ up to a set of measure zero.
For a finite subfamily $\mathscr A\subset\mathscr A_0$, denote
$$
E_{\mathscr A}
=\bigcup_{R\in\mathscr A_0\setminus\mathscr A}\mathcal W(R).
$$
For any given $g\in L^2(\R^3)$ with $ \left\|   g \right\|_{L^2(\R^3)}=1$, since $\Phi_{s,t}$ is real and even, we have
$$
\int_{\R^3}\Phi_{s,t}(x-y)\overline{g(x)}\,dx
=\overline{(g*\Phi_{s,t})(y)}.
$$
Thus by \eqref{e:calderon-z},
\begin{align*}
&\bigg\langle f-\sum_{R\in\mathscr A}\lambda_R   (f)   \mathsf{a}_R   (f),g\bigg\rangle\\
= &\int_{\R^3\times(0,\infty)^2}
(f*\psi_{s,t})(y)\overline{(g*\Phi_{s,t})(y)}\,d\nu
-\sum_{R\in\mathscr A}
\int_{\mathcal W(R)}
(f*\psi_{s,t})(y)\overline{(g*\Phi_{s,t})(y)}\,d\nu\\
= &\int_{E_{\mathscr A}}
(f*\psi_{s,t})(y)\overline{(g*\Phi_{s,t})(y)}\,d\nu.
\end{align*}

Since \eqref{e:kernel-l2} also holds with $\psi$ replaced by $\Phi$, Cauchy--Schwarz gives
\begin{align*}
\bigg|
\bigg\langle f-\sum_{R\in\mathscr A}\lambda_R   (f)  \mathsf{a}_R   (f)   ,g\bigg\rangle
\bigg|
\leq \,&
 \bigg( \int_{E_{\mathscr A}}|(f*\psi_{s,t})(y)|^2\,d\nu  \bigg)^{1/2}
 \bigg( \int_{E_{\mathscr A}}|(g*\Phi_{s,t})(y)|^2\,d\nu  \bigg)^{1/2}\\
\leq\,  &
 \bigg( \sum_{R\in\mathscr A_0\setminus\mathscr A}\lambda_R^2(f)  \bigg)^{1/2}
 \bigg( \int_0^\infty\int_0^\infty
 \left\|   g*\Phi_{s,t} \right\|_{L^2(\R^3)}^2
\,\frac{dt}{t}\frac{ds}{s}  \bigg)^{1/2}\\
=\,   &(c_{\Phi_1}c_{\Phi_2})^{1/2}
 \bigg( \sum_{R\in\mathscr A_0\setminus\mathscr A}\lambda_R^2  (f)     \bigg)^{1/2}
 \left\|   g \right\|_{L^2(\R^3)}.
\end{align*}
Therefore, \eqref{e:remainder} holds.

The coefficient energy is finite by the partition \eqref{e:W-decomp} and \eqref{e:kernel-l2}:
$$
\sum_{R\in\mathscr A_0}\lambda_R^2(f)
=\int_{\R^3\times(0,\infty)^2}
|(f*\psi_{s,t})(y)|^2\,d\nu
=c_{\psi_1}c_{\psi_2} \left\|   f \right\|_{L^2(\R^3)}^2<\infty.
$$
Thus, for finite families $\mathscr A_m\subseteq\mathscr A_{m+1}$ with
$\bigcup_m\mathscr A_m=\mathscr A_0$,
$$
\lim_{m\to\infty}
\bigg\|f-\sum_{R\in\mathscr A_m}\lambda_R(f)\mathsf{a}_R(f)\bigg\|_{L^2(\R^3)}=0.
$$
The proof is complete.
\end{proof}

We next cover $(5I)\times(5J)\times(9K)$ by a dyadic Zygmund rectangle.
Applying Lemma~\ref{l:adj-int-cover} in each coordinate first gives a product of three dyadic intervals.
Its side lengths need not satisfy the Zygmund relation.
We therefore pass to ancestors in the same selected grids, choosing the first two lengths to be
$16\ell(I)$ and $16\ell(J)$ and the third to be their product $256\ell(K)$.

\begin{lemma}\label{l:zyg-cover}
For every $R=I\times J\times K\in\mathscr A_0$,
one can choose $ \vec{ \pmb\varepsilon}(R)\in\{0,1\}^3$ and
$$
Q(R)=Q_1(R)\times Q_2(R)\times Q_3(R)
\in\mathcal D_Z^{ \vec{ \pmb\varepsilon}(R)}
$$
such that
$
(5I)\times(5J)\times(9K)\subset Q(R)
$ and
$|Q(R)|=16^4|R|.
$
\end{lemma}

\begin{proof}
Applying Lemma~\ref{l:adj-int-cover} to $5I$, $5J$ and $9K$,
choose $L_j\in\mathcal  D^{ \varepsilon_j(R)}$, $j=1,2,3$, such that
$$
5I\subset L_1, \qquad  5J\subset L_2, \qquad  9K\subset L_3,
$$
and
$$
\ell(L_1)\le30\ell(I), \qquad 
\ell(L_2)\le30\ell(J), \qquad 
\ell(L_3)\le54\ell(K).
$$
Since all these lengths belong to $2^{\mathbb Z}$, we have 
$$
\ell(L_1)\le16\ell(I), \qquad 
\ell(L_2)\le16\ell(J), \qquad 
\ell(L_3)\le32\ell(K).
$$
In the respective grids $\mathcal D^{\varepsilon_j(R)}$, we take the unique dyadic ancestors
$Q_j(R)$ of $L_j$ with
$$
\ell(Q_1(R))=16\ell(I), \qquad 
\ell(Q_2(R))=16\ell(J), \qquad 
\ell(Q_3(R))=256\ell(K).
$$
Then $(5I)\times(5J)\times(9K)\subset Q(R)$.
Moreover, since $\ell(K)=\ell(I)\ell(J)$,
$$
\ell(Q_3(R))
=256\ell(K)
=(16\ell(I))(16\ell(J))
=\ell(Q_1(R))\ell(Q_2(R)).
$$
Thus $Q(R)\in\mathcal D_Z^{ \vec{ \pmb\varepsilon}(R)}$, and
$
|Q(R)|=16\cdot16\cdot256\,|R|=16^4|R|.
$
\end{proof}

Fix these choices for every $R\in\mathscr A_0$.
By Lemma~\ref{l:atom-prop} and the definition of $\mathsf{a}_R$,
$$
R\cup\operatorname{supp}\mathsf{a}_R
\subset(5I)\times(5J)\times(9K)
\subset Q(R),
 \qquad 
|Q(R)|\le C_0|R|,
 \qquad  C_0=16^4.
$$
Both $ \vec{ \pmb\varepsilon}(R)$ and $Q(R)$  are independent of $f$.

We classify the original rectangles by the selected grid of their covering rectangles:
$$
\mathscr A_{0, \vec{ \pmb\varepsilon}}
=\big \{R\in\mathscr A_0:
 \vec{ \pmb\varepsilon}(R)= \vec{ \pmb\varepsilon}\big  \},
\qquad
\mathscr A_0
=\bigcup_{ \vec{ \pmb\varepsilon}\in\{0,1\}^3}
\mathscr A_{0, \vec{ \pmb\varepsilon}}.
$$
The union is disjoint because one grid label has been fixed for each $R$.
Thus every atom occurs in exactly one class, and all $Q(R)$ in that class belong to the same dyadic Zygmund system.
Different rectangles $R$ may still have the same covering rectangle $Q(R)$; we next bound their number.

For a fixed rectangle $P=P_1\times P_2\times P_3$,
every $R=I\times J\times K$ satisfying $Q(R)=P$ has
$$
\ell(I)=\frac{\ell(P_1)}{16}, \qquad 
\ell(J)=\frac{\ell(P_2)}{16}, \qquad 
\ell(K)=\frac{\ell(P_3)}{256}.
$$
These rectangles $R$ lie in $P$ and are mutually disjoint, each with
volume $|P|/{16^4}$. Hence
\begin{align}\label{e:top-bd-preim}
\#\{R\in\mathscr A_0:Q(R)=P\}
\leq 16^4=:M_{\rm mult}.
\end{align}

We next apply Lemma~\ref{l:atom-prop} to 
$\mathsf{a}_R(f)=\widetilde{\mathsf{a}}_R(f)/\lambda_R(f)$ when $\lambda_R(f)>0$,
with $\mathsf{a}_R=0$ otherwise. Concretely,
 $\mathsf{a}_R$ is smooth and
$$
\operatorname{supp}\mathsf{a}_R(f)
\subset(5I)\times(5J)\times(9K)
\subset Q(R).
$$
For $\lambda_R(f)>0$, $\gamma\in\mathbb N_0^3$ with $|\gamma|\le3$,
and $x\in\R^3$, its derivative estimate implies
\begin{align}\label{e:prod-top-der}
|\partial^\gamma\mathsf{a}_R(f)(x)|
= \frac{|\partial^\gamma \widetilde{\mathsf{a}}_R(f)(x)|}{\lambda_R(f)}
\leq  A_\gamma|Q(R)|^{-1/2}
\prod_{j=1}^3\ell(Q_j(R))^{-\gamma_j},
\end{align}
where 
$
A_\gamma=C_\gamma  16^{2+\gamma_1+\gamma_2+2\gamma_3}
$, and $C_\gamma$ is the constant in (ii) of Lemma~\ref{l:atom-prop}.
Thus $A_\gamma$  is independent of $f$, $R$ and the selected grids.

Meanwhile,  by the  cancellation  in (iv) of  Lemma~\ref{l:atom-prop}, we have 
\begin{equation}\label{e:prod-cancel}
\int_{Q_1(R)}\mathsf{a}_R(f)(x)\,dx_1=0,
 \qquad 
\int_{Q_2(R)\times Q_3(R)}\mathsf{a}_R(f)(x)\,dx_2\,dx_3=0.
\end{equation}

When $\lambda_R(f)=0$, both the derivative estimates and the
cancellation identities hold because $\mathsf{a}_R=0$.

We call these $\mathsf{a}_R(f)(x)$ {\it normalized smooth atoms}. 

\subsection{Reverse \texorpdfstring{$L^p$}{Lp} estimates for smooth atomic sums}\label{sub:atom-bounds}

Let $2\leq p<\infty$ and $f\in L^2(\R^3)\cap L^p(\R^3)$. Let $\mathsf{a}_R$ and $Q(R)$ be the normalized atoms and enlarged 
rectangles constructed above.

For fixed $R$, by H\"older's inequality,
$$
\lambda_R(f)
\le   \|  f    \|_{L^p}
 \bigg( \int_{\mathcal W(R)}  \|  \psi_{s,t}    \|_{L^{p'}}^2\,dy\,\frac{dt}{t}\frac{ds}{s}  \bigg)^{1/2}
\le C_R  \|  f    \|_{L^p}.
$$
Thus the  definition of $\widetilde{\mathsf{a}}_R(f)$ and the proof of Lemma~\ref{l:atom-prop}
 are  valid for  $ L^p$ functions.
 
For a finite $\mathscr A\subset\mathscr A_0$, define
$$
X_{\mathscr A}(f) (x)=\sum_{R\in\mathscr A}\lambda_R(f)\mathsf{a}_R(f) (x),
 \qquad 
\mathcal S_{\mathscr A}(f)(x)^2
=\sum_{R\in\mathscr A}\lambda_R(f)^2
\frac{\mathbf 1_{Q(R)}(x)}{|Q(R)|}.
$$

The reverse estimate below controls $f$ by this discrete square function.
To compare $\mathcal S_{\mathscr A_0}(f)$ pointwise with $S_{z,\alpha}^{\psi}f$,
the enlarged supports must fit the cone at aperture $\alpha$,
as in the one-parameter comparison \eqref{sp:e:discrete-area}.

We begin with the following upper bound for $\mathcal S_{\mathscr A_0}$.

\begin{lemma}\label{l:atomic-forward}
For $2\le p<\infty$, we have
$
 \left\|   \mathcal S_{\mathscr A_0}(h) \right\|_{L^p(\R^3)}\le B_p  \|  h    \|_{L^p(\R^3)},
$ for all $ h\in L^p(\R^3),$
where $B_p$ depends only on $p$ and the fixed kernels.
\end{lemma}

\begin{proof}
For $p=2$,  combining \eqref{e:atom-coeff}, \eqref{e:W-decomp}  and \eqref{e:kernel-l2}, we have
$$
 \left\|   \mathcal S_{\mathscr A_0}(h) \right\|_{L^2}^2
=\sum_{R\in\mathscr A_0}\lambda_R(h)^2
=c_{\psi_1}c_{\psi_2}  \|  h    \|_{L^2}^2.
$$

For $2<p<\infty$, we use the weighted method of R. Fefferman and Pipher in \cite{FP}.
{Let} $r=(p/2)'\in(1,\infty)$.
For every finite $\mathscr A\subset\mathscr A_0$, by $L^{p/2}$ duality,
$$
 \left\|   \mathcal S_{\mathscr A}(h) \right\|_{L^p}^2
= \left\|   \mathcal S_{\mathscr A}(h)^2 \right\|_{L^{p/2}}
=\sup_{\substack{w\ge0\\  \|  w    \|_{L^r}=1}}
\int_{\R^3}\mathcal S_{\mathscr A}(h)^2w.
$$

Let $\mathsf M$ be the strong maximal operator on $\R^3$.
Its $L^r$ boundedness gives $N_r:=  \|  \mathsf M    \|_{L^r\to L^r}<\infty$.
To apply the weighted square-function estimate to the dual pairing above,
we majorize $w$ by an $A_1$ weight with respect to rectangles.
For $w\ge0$ with $  \|  w    \|_{L^r}=1$, define
$$
W=\sum_{k=0}^\infty\frac{\mathsf M^kw}{(2N_r)^k},
$$
then
$$
w\le W,\quad   \|  W    \|_{L^r}\le2,\quad
\mathsf MW\le2N_r(W-w)\le2N_rW.
$$
Thus $W$ is an $A_1$ weight with respect to rectangles.
For every rectangle $B$,
$$
\langle W\rangle_B:=|B|^{-1}\int_BW
\le2N_r\mathop{\rm ess\,inf}_B W,
$$
and hence
$$
\langle W\rangle_B\langle W^{-1}\rangle_B
\le
\frac{\langle W\rangle_B}{\mathop{\rm ess\,inf}_B W}
\le2N_r.
$$
The weighted estimate in \cite[Theorem~2.9]{FP}, for the kernels fixed above, says that 
$$
\int_0^\infty\int_0^\infty\int_{\R^3}
|h*\psi_{s,t}(y)|^2\langle W\rangle_{B_{s,t}(y)}\,dy\,\frac{dt}{t}\frac{ds}{s}
\le C_{p,\psi}\int_{\R^3}|h|^2W,
$$
where
$
{B_{s,t}(y)=y+\bigl((-s,s)\times(-t,t)\times(-st,st)\bigr).}
$

Let $c_\delta=1+\max\{16\delta,256\delta^2\}$.
For $(y,s,t)\in\mathcal W(R)$, the side lengths of $Q(R)$ give
$$
Q(R)\subset c_\delta B_{s,t}(y),
 \qquad 
\frac{|c_\delta B_{s,t}(y)|}{|Q(R)|}
\le 8c_\delta^3 16^{-3}\delta^{-4}=:C_\delta.
$$
Consequently,
$$
\langle W\rangle_{Q(R)}
\le C_\delta\langle W\rangle_{c_\delta B_{s,t}(y)}
\le2C_\delta N_r\mathop{\rm ess\,inf}_{c_\delta B_{s,t}(y)}W
\le2C_\delta N_r\langle W\rangle_{B_{s,t}(y)}.
$$
Therefore,
\begin{align*}
\int_{\R^3}\mathcal S_{\mathscr A}(h)^2w
\leq &\sum_{R\in\mathscr A}\int_{\mathcal W(R)}
|h*\psi_{s,t}(y)|^2\langle W\rangle_{Q(R)}\,dy\,\frac{dt}{t}\frac{ds}{s}\\
\leq &2C_\delta N_r\int_0^\infty\int_0^\infty\int_{\R^3}
|h*\psi_{s,t}(y)|^2\langle W\rangle_{B_{s,t}(y)}\,dy\,\frac{dt}{t}\frac{ds}{s}\\
= &\frac{C_\delta N_r}{4}\int_{\R^3}(S_{z,1}^{\psi}h)^2W
\le C_{p,\psi,\delta}\int_{\R^3}|h|^2W\\
\leq &C_{p,\psi,\delta}\|h\|_{L^p}^2\|W\|_{L^r}
\le2C_{p,\psi,\delta}\|h\|_{L^p}^2.
\end{align*}
Taking the supremum over $w$ with $  \|  w    \|_{L^r}=1$, we obtain the estimate for finite $\mathscr A$.
The estimate for $\mathscr A_0$ then follows by monotone convergence.
\end{proof}

\smallskip

\begin{proposition}\label{p:smooth-model}
Let $2\le p<\infty$ and $f\in L^p(\R^3)$.
For every finite $\mathscr A\subset\mathscr A_0$,
$$
 \left\|   X_{\mathscr A}(f) \right\|_{L^p(\R^3)}
\le Cp\, \left\|   \mathcal S_{\mathscr A}(f) \right\|_{L^p(\R^3)}.
$$
For every sequence $(\mathscr A_m)_{m\ge1}$ of finite families
$\mathscr A_m\subseteq\mathscr A_{m+1}\subseteq\mathscr A_0$
with $\bigcup_m\mathscr A_m=\mathscr A_0$,
$$
\lim_{m\to\infty}\left\|X_{\mathscr A_m}(f)-f\right\|_{L^p(\R^3)}=0.
$$
Moreover,
$$
  \|  f    \|_{L^p(\R^3)}
\le Cp\, \left\|   \mathcal S_{\mathscr A_0}(f) \right\|_{L^p(\R^3)}.
$$
The constant $C$ depends only on the fixed kernels.
\end{proposition}

\begin{remark}
The Haar expansion follows the proof of Lemma~\ref{sp:l:smooth}.
As in \cite[Lemma~2.2]{Wilson89} and the principal part of
\cite[Proposition~3.6]{WuJGA}, we use Haar cancellation and atom smoothness
to obtain summable decay in the relative generations.
Here cancellation in $x_1$ and jointly in $(x_2,x_3)$ gives
the three Haar types $\eta\in\mathsf H=\{(1,0),(0,1),(1,1)\}$.
For every $P\in Ch_{ \vec{ \pmb n}}(Q)$ in the local expansion, we have
$$
P\in\mathcal D_Z^{ \vec{ \pmb\varepsilon},q}
\   \    \text{with}\    \     q=n_1+n_2-n_3,
$$
so we apply Proposition~\ref{p:zyg-haar} uniformly in $q$ to obtain the factor $p$.
\end{remark}

\begin{proof}[Proof of Proposition \ref{p:smooth-model}]
Let $f\in L^p(\R^3)$ and let $\mathscr A\subset\mathscr A_0$ be finite. For $ \vec{ \pmb\varepsilon}\in\{0,1\}^3$, let
$\displaystyle
\mathscr A_{ \vec{ \pmb\varepsilon}}
=\left\{R\in\mathscr A: \vec{ \pmb\varepsilon}(R)= \vec{ \pmb\varepsilon}\right\}.
$ 
Then
$$
X_{\mathscr A}(f)
=\sum_{ \vec{ \pmb\varepsilon}\in\{0,1\}^3}X_{\mathscr A_{ \vec{ \pmb\varepsilon}}}(f),
 \qquad 
\mathcal S_{\mathscr A}(f)^2
=\sum_{ \vec{ \pmb\varepsilon}\in\{0,1\}^3}\mathcal S_{\mathscr A_{ \vec{ \pmb\varepsilon}}}(f)^2.
$$

\smallskip
\noindent{\it Step 1.}
Fix $ \vec{ \pmb\varepsilon}$.
All sums over $R$ in this step and the following {\it Step 2} are over $\mathscr A_{ \vec{ \pmb\varepsilon}}$.
For each coordinate $j$, the supporting intervals $Q_j(R)$ now belong to the single grid
$\mathcal D^{\varepsilon_j}$. We form their Haar expansions in these fixed grids,
and sum over the eight classes after obtaining the estimate for one class.
For $Q\in\mathcal D_Z^{ \vec{ \pmb\varepsilon}}$ and $ \vec{ \pmb n}\in\mathbb N_0^3$, define
$$
Ch_{ \vec{ \pmb n}}(Q)
=\left\{P\in\mathcal  D^{  \vec{ \pmb\varepsilon}}:P\subseteq Q,\
\frac{\ell(P_j)}{\ell(Q_j)}=2^{-n_j}, \    j=1,2,3\right\},
 \qquad  | \vec{ \pmb n}|=n_1+n_2+n_3.
$$
Fix $R\in\mathscr A_{ \vec{ \pmb\varepsilon}}$ and write $Q=Q(R)$.
For $n\ge0$, let
$$
\mathbb E_{j,n}^Qg(x)
=\sum_{\substack{L\in\mathcal  D^{ \varepsilon_j},\ L\subseteq Q_j\\
\ell(L)=2^{-n}\ell(Q_j)}}
\left(\frac1{|L|}\int_Lg\,dx_j\right)\mathbf 1_L(x_j),
 \qquad 
\Delta_{j,n}^Q=\mathbb E_{j,n+1}^Q-\mathbb E_{j,n}^Q.
$$
For $N\ge1$, 
$$
\mathbb E_{j,N}^Q=\mathbb E_{j,0}^Q+\sum_{n=0}^{N-1}\Delta_{j,n}^Q,
 \qquad 
\mathbb E_{j,0}^Qg=h_{Q_j}^0\int_{Q_j}g h_{Q_j}^0\,dx_j,
$$
and
$$
\Delta_{j,n}^Qg
=\sum_{\substack{L\in\mathcal  D^{ \varepsilon_j},\ L\subseteq Q_j\\
\ell(L)=2^{-n}\ell(Q_j)}}h_L\int_Lg h_L\,dx_j.
$$
Moreover, these operators in distinct coordinates commute. By \eqref{e:prod-cancel},
we have
$
\mathbb E_{1,0}^Q\mathsf{a}_R(f)=0,
$ and
$\mathbb E_{2,0}^Q\mathbb E_{3,0}^Q\mathsf{a}_R(f)=0.
$

Define $\mathsf{a}_{R,N}(f)=\mathbb E_{1,N}^Q\mathbb E_{2,N}^Q\mathbb E_{3,N}^Q\mathsf{a}_R(f)$.
Then
\begin{align*}
\mathsf{a}_{R,N}
= &\mathbb E_{2,N}^Q\mathbb E_{3,N}^Q
 \bigg( \underbrace{\mathbb E_{1,0}^Q\mathsf{a}_R (f)}_{=0}
+\sum_{n_1=0}^{N-1}\Delta_{1,n_1}^Q\mathsf{a}_R (f)  \bigg)\\
= &\sum_{n_1=0}^{N-1}\Delta_{1,n_1}^Q
 \bigg( \mathbb E_{2,0}^Q+\sum_{n_2=0}^{N-1}\Delta_{2,n_2}^Q  \bigg)\,
 \bigg( \mathbb E_{3,0}^Q+\sum_{n_3=0}^{N-1}\Delta_{3,n_3}^Q  \bigg)\, \mathsf{a}_R (f)\\
= &\sum_{n_1=0}^{N-1}\Delta_{1,n_1}^Q
\underbrace{\mathbb E_{2,0}^Q\mathbb E_{3,0}^Q\mathsf{a}_R (f)}_{=0}
+\sum_{n_1,n_2=0}^{N-1}\Delta_{1,n_1}^Q\Delta_{2,n_2}^Q\mathbb E_{3,0}^Q\mathsf{a}_R(f)\\
& \quad +\sum_{n_1,n_3=0}^{N-1}\Delta_{1,n_1}^Q\mathbb E_{2,0}^Q\Delta_{3,n_3}^Q\mathsf{a}_R(f)
+\sum_{n_1,n_2,n_3=0}^{N-1}\Delta_{1,n_1}^Q\Delta_{2,n_2}^Q\Delta_{3,n_3}^Q\mathsf{a}_R(f).
\end{align*}
The three terms have the Haar representations
\begin{align*}
\Delta_{1,n_1}^Q\Delta_{2,n_2}^Q\mathbb E_{3,0}^Q\mathsf{a}_R(f)(x)
= &\sum_{P\in Ch_{(n_1,n_2,0)}(Q)}
h_{P_1}(x_1)h_{P_2}(x_2)h_{Q_3}^0(x_3)
\int_P\mathsf{a}_R (f) (y)h_{P_1}(y_1)h_{P_2}(y_2)h_{Q_3}^0(y_3)\,dy\\
= &\sum_{P\in Ch_{(n_1,n_2,0)}(Q)}
\left\langle\mathsf{a}_R (f)  ,h_{P,Z}^{(1,0)}\right\rangle h_{P,Z}^{(1,0)}(x).
\end{align*}
Similarly,  
$$
\Delta_{1,n_1}^Q\mathbb E_{2,0}^Q\Delta_{3,n_3}^Q\mathsf{a}_R (f)(x)
=\sum_{P\in Ch_{(n_1,0,n_3)}(Q)}
h_{P_1}(x_1)h_{Q_2}^0(x_2)h_{P_3}(x_3)
\int_P\mathsf{a}_R(f)(y)h_{P_1}(y_1)h_{Q_2}^0(y_2)h_{P_3}(y_3)\,dy,
$$
and
$$
\Delta_{1,n_1}^Q\Delta_{2,n_2}^Q\Delta_{3,n_3}^Q\mathsf{a}_R(f)(x)
=\sum_{P\in Ch_{(n_1,n_2,n_3)}(Q)}
h_{P_1}(x_1)h_{P_2}(x_2)h_{P_3}(x_3)
\int_P\mathsf{a}_R(f)(y)h_{P_1}(y_1)h_{P_2}(y_2)h_{P_3}(y_3)\,dy.
$$
Therefore
$$
\mathsf{a}_{R,N}(f)=\sum_{\eta, \vec{ \pmb n}}\sum_{P\in Ch_{ \vec{ \pmb n}}(Q)}
\left \langle\mathsf{a}_R(f),h_{P,Z}^{\eta}  \right\rangle h_{P,Z}^{\eta},
$$
where
$$
\eta\in\mathsf H, \qquad  0\le n_1<N, \qquad 
0\le n_j\le\eta_j(N-1)\quad(j=2,3).
$$
In this expansion,
$
\eta_j=0$ implies $ n_j=0$ and $ P_j=Q_j$.  As a consequence, 
\begin{equation}\label{e:n-est}
n_1+\eta_2n_2+\eta_3n_3=| \vec{ \pmb n}|.
\end{equation}

\smallskip
\noindent{\it Step 2.}
Fix $\eta, \vec{ \pmb n}$ in these ranges.
For $P\in Ch_{ \vec{ \pmb n}}(Q)$, we have 
$$
|P|=2^{-| \vec{ \pmb n}|}|Q|,
 \qquad 
\frac{\ell(P_3)}{\ell(P_1)\ell(P_2)}=2^{n_1+n_2-n_3}.
$$
Thus 
\begin{equation}\label{e:Pd}
	P\in\mathcal D_Z^{ \vec{ \pmb\varepsilon},q}\     \    \text{with}\   \    q=n_1+n_2-n_3, \qquad \forall\, P\in Ch_{ \vec{
\pmb n}}(Q(R)).
\end{equation}

For $P\in Ch_{ \vec{ \pmb n}}(Q(R))$, let $c_j$ be the midpoint of $P_j$.
Haar cancellation and \eqref{e:prod-top-der} imply
\begin{align}\label{e:aux-PsiHaar}
\left|\langle\mathsf{a}_R(f),h_{P,Z}^{\eta}\rangle\right|
\leq & \left\|   \partial_1\partial_2^{\eta_2}\partial_3^{\eta_3}\mathsf{a}_R \right\|_\infty
\int_P|x_1-c_1|\prod_{j=2}^3|x_j-c_j|^{\eta_j}|h_{P,Z}^{\eta}(x)|\,dx\notag\\
\leq & \left\|   \partial_1\partial_2^{\eta_2}\partial_3^{\eta_3}\mathsf{a}_R \right\|_\infty
|P|^{1/2}\ell(P_1)\prod_{j=2}^3\ell(P_j)^{\eta_j}\notag\\
\leq &C\left(\frac{|P|}{|Q(R)|}\right)^{1/2}
\frac{\ell(P_1)}{\ell(Q_1(R))}
\prod_{j=2}^3\left(\frac{\ell(P_j)}{\ell(Q_j(R))}\right)^{\eta_j}   \notag\\
= &C2^{-| \vec{ \pmb n}|/2-n_1-\eta_2n_2-\eta_3n_3}
= C2^{-3| \vec{ \pmb n}|/2}.
\end{align}

Define
$$
X_{ \vec{ \pmb n},\eta}(f)
=\sum_R\lambda_R (f)\sum_{P\in Ch_{ \vec{ \pmb n}}(Q(R))}
\left \langle\mathsf{a}_R(f),h_{P,Z}^{\eta}  \right\rangle h_{P,Z}^{\eta},
$$
and its Haar square function is
$$
S_{ \vec{ \pmb n},\eta}(f)(x)^2
=\sum_{P\in\mathcal D_Z^{ \vec{ \pmb\varepsilon},q}}
\bigg|\sum_{R\in\mathscr A_{ \vec{ \pmb\varepsilon}}:\,P\in Ch_{ \vec{ \pmb n}}(Q(R))}
\lambda_R(f)\langle\mathsf{a}_R(f),h_{P,Z}^{\eta}\rangle\bigg|^2
\frac{\mathbf 1_P(x)}{|P|}.
$$

For fixed $P$ and $\vec{\pmb n}$, each $Q_j(R)$ must be the unique $n_j$th ancestor of $P_j$
in the fixed grid $\mathcal D^{\varepsilon_j}$.
Hence all atoms contributing to this Haar coefficient have the same covering rectangle $Q(R)$.
By \eqref{e:top-bd-preim}, their number satisfies
$$
\#\left\{R\in\mathscr A_{ \vec{ \pmb\varepsilon}}:P\in Ch_{ \vec{ \pmb n}}(Q(R))\right\}
\le M_{\rm mult}=16^4,
$$
independently of the relative generations $\vec{\pmb n}$.
Meanwhile,
$\displaystyle
\sum_{P\in Ch_{ \vec{ \pmb n}}(Q(R))}
\frac{\mathbf 1_P(x)}{|P|}
=2^{| \vec{ \pmb n}|}\frac{\mathbf 1_{Q(R)}(x)}{|Q(R)|}
$. Applying Cauchy--Schwarz and \eqref{e:aux-PsiHaar}, we obtain
\begin{align}\label{e:depth-sq}
S_{ \vec{ \pmb n},\eta}(f)(x)^2
\leq &\sum_{P\in\mathcal D_Z^{ \vec{ \pmb\varepsilon},q}}
 \bigg( \sum_{R\in\mathscr A_{ \vec{ \pmb\varepsilon}}:\,P\in Ch_{ \vec{ \pmb n}}(Q(R))}1  \bigg)\,
 \bigg( \sum_{R\in\mathscr A_{ \vec{ \pmb\varepsilon}}:\,P\in Ch_{ \vec{ \pmb n}}(Q(R))}
\lambda_R(f)^2
\left|\langle\mathsf{a}_R(f),h_{P,Z}^{\eta}\rangle\right|^2  \bigg)
\frac{\mathbf 1_P(x)}{|P|}\notag\\
\leq &C\sum_R\lambda_R(f)^2
\sum_{P\in Ch_{ \vec{ \pmb n}}(Q(R))}
\left|\langle\mathsf{a}_R(f),h_{P,Z}^{\eta}\rangle\right|^2
\frac{\mathbf 1_P(x)}{|P|}\notag\\
\leq &C2^{-3| \vec{ \pmb n}|}\sum_R\lambda_R(f)^2
\sum_{P\in Ch_{ \vec{ \pmb n}}(Q(R))}
\frac{\mathbf 1_P(x)}{|P|}\notag\\
= &C2^{-2| \vec{ \pmb n}|}\sum_R\lambda_R(f)^2
\frac{\mathbf 1_{Q(R)}(x)}{|Q(R)|}\notag\\
= &C2^{-2| \vec{ \pmb n}|}
\mathcal S_{\mathscr A_{ \vec{ \pmb\varepsilon}}}(f)(x)^2.
\end{align}

Combining \eqref{e:Pd},  Proposition~\ref{p:zyg-haar} and \eqref{e:depth-sq}, we have
\begin{align*}
 \left\|   X_{ \vec{ \pmb n},\eta}(f) \right\|_{L^p(\R^3)}
\leq &Cp\, \left\|   S_{ \vec{ \pmb n},\eta}(f) \right\|_{L^p(\R^3)}
\leq C p\,2^{-| \vec{ \pmb n}|}
 \left\|   \mathcal S_{\mathscr A_{ \vec{ \pmb\varepsilon}}}(f) \right\|_{L^p(\R^3)},
\end{align*}
with $C$ independent of $q$.

For the three Haar types, one may apply \eqref{e:n-est} to see
\begin{align*}
\sum_{\eta, \vec{ \pmb n}}2^{-| \vec{ \pmb n}|}
= &\sum_{n_1,n_2=0}^{N-1}2^{-n_1-n_2}
+\sum_{n_1,n_3=0}^{N-1}2^{-n_1-n_3}
+\sum_{n_1,n_2,n_3=0}^{N-1}2^{-n_1-n_2-n_3}
= 2\,\bigg (\sum_{k=0}^{N-1}2^{-k}  \bigg)^2
+ \bigg( \sum_{k=0}^{N-1}2^{-k}  \bigg)^3
\le16.
\end{align*}
Consequently,
\begin{align*}
\bigg  \|  \sum_R\lambda_R\mathsf{a}_{R,N}\bigg    \|_{L^p(\R^3)}
= \bigg  \|  \sum_{\eta, \vec{ \pmb n}}X_{ \vec{ \pmb n},\eta}(f)\bigg    \|_{L^p(\R^3)}
\leq \sum_{\eta, \vec{ \pmb n}}
 \left\|   X_{ \vec{ \pmb n},\eta}(f) \right\|_{L^p(\R^3)}
\leq C p\,
 \left\|   \mathcal S_{\mathscr A_{ \vec{ \pmb\varepsilon}}}(f) \right\|_{L^p(\R^3)},
\end{align*}
uniformly in $N$ and $ \vec{ \pmb\varepsilon}$.
Moreover,
$$
  \|  \mathsf{a}_{R,N}-\mathsf{a}_R    \|_\infty
\le2^{-N}\sum_{j=1}^3\ell(Q_j(R))  \|  \partial_j\mathsf{a}_R    \|_\infty
\le C2^{-N}|Q(R)|^{-1/2}.
$$
Since $\mathscr A_{ \vec{ \pmb\varepsilon}}$ is finite,
$$
\bigg  \|  \sum_R\lambda_R(\mathsf{a}_{R,N}-\mathsf{a}_R)\bigg    \|_{L^p}
\le C2^{-N}\sum_R\lambda_R|Q(R)|^{1/p-1/2}\to0.
$$
Letting $N\to\infty$ and summing over $ \vec{ \pmb\varepsilon}$, we obtain
\begin{align}\label{e:finite-atomic}
 \left\|   X_{\mathscr A}(f) \right\|_{L^p}
\leq &\sum_{ \vec{ \pmb\varepsilon}\in\{0,1\}^3}
 \left\|   X_{\mathscr A_{ \vec{ \pmb\varepsilon}}}(f) \right\|_{L^p} 
\leq Cp
\sum_{ \vec{ \pmb\varepsilon}\in\{0,1\}^3}
 \left\|   \mathcal S_{\mathscr A_{ \vec{ \pmb\varepsilon}}}(f) \right\|_{L^p} 
\leq Cp \left\|   \mathcal S_{\mathscr A}(f) \right\|_{L^p}.
\end{align}

\smallskip
\noindent{\it Step 3.}
We now use \eqref{e:finite-atomic} to prove the convergence of the atomic sums in $L^p$.
 
First let $f\in L^2(\R^3)\cap L^p(\R^3)$.
Then $\mathcal S_{\mathscr A_0}(f)\in L^p$ by Lemma~\ref{l:atomic-forward}.
Choose finite families $\mathscr A_m\subseteq\mathscr A_{m+1}$ with
$\bigcup_m\mathscr A_m=\mathscr A_0$.
Define
$$
\mathcal E_m(f)= \bigg( \sum_{R\in\mathscr A_0\setminus\mathscr A_m}
\lambda_R(f)^2\frac{\mathbf 1_{Q(R)}}{|Q(R)|}  \bigg)^{1/2}.
$$
Then
$$
0\le \mathcal E_m(f)\le\mathcal S_{\mathscr A_0}(f), \qquad 
\mathcal E_m(f)\to0\    \     \text{a.e.}, \qquad    \|  \mathcal E_m(f)    \|_{L^p}\to0.
$$
For $k>m$, \eqref{e:finite-atomic} gives
$$
  \|  X_{\mathscr A_k}(f)-X_{\mathscr A_m}(f)    \|_{L^p}
=  \|  X_{\mathscr A_k\setminus\mathscr A_m}(f)    \|_{L^p}
\le Cp\,  \|  \mathcal E_m(f)    \|_{L^p}.
$$
Thus $\lim_{m\to\infty}\|X_{\mathscr A_m}(f)-g\|_{L^p}=0$ for some $g\in L^p$.
By Lemma~\ref{l:atomic-l2}, for every $\vartheta\in C_c^\infty(\R^3)$,
$$
|\langle g-f,\vartheta\rangle|
\le   \|  g-X_{\mathscr A_m}(f)    \|_{L^p}  \|  \vartheta    \|_{L^{p'}}
+  \|  X_{\mathscr A_m}(f)-f    \|_{L^2}  \|  \vartheta    \|_{L^2}\to0.
$$
Hence $g=f$, and $\lim_{m\to\infty}\|X_{\mathscr A_m}(f)-f\|_{L^p}=0$.

Now let $f\in L^p$. By Lemma~\ref{l:atomic-forward}
and \eqref{e:finite-atomic},
$$
  \|  X_{\mathscr A}(h)    \|_{L^p}
\le Cp\,  \|  \mathcal S_{\mathscr A}(h)    \|_{L^p}
\le C_p  \|  h    \|_{L^p},
 \qquad  h\in L^p,
$$
where $C_p$ is independent of the finite family $\mathscr A$.
Choose
$$
f_k=f\mathbf 1_{[-k,k]^3}\mathbf 1_{\{|f|\le k\}}
\in L^2\cap L^p,
 \qquad    \|  f-f_k    \|_{L^p}\to0.
$$
For each fixed $k$, the above argument shows that 
$$
\lim_{m\to\infty}\left\|X_{\mathscr A_m}(f_k)-f_k\right\|_{L^p}=0.
$$
By linearity of the finite sums $X_{\mathscr A_m}$,
\begin{align*}
  \|  X_{\mathscr A_m}(f)-f    \|_{L^p}
\leq &  \|  X_{\mathscr A_m}(f-f_k)    \|_{L^p}
+  \|  X_{\mathscr A_m}(f_k)-f_k    \|_{L^p}
+  \|  f_k-f    \|_{L^p}\\
\leq &C_p  \|  f-f_k    \|_{L^p}
+  \|  X_{\mathscr A_m}(f_k)-f_k    \|_{L^p}.
\end{align*}
Hence
$$
\limsup_{m\to\infty}  \|  X_{\mathscr A_m}(f)-f    \|_{L^p}
\le C_p  \|  f-f_k    \|_{L^p}\to0
 \qquad (k\to\infty).
$$
Thus $X_{\mathscr A_m}(f)\to f$ in $L^p$.
Finally, \eqref{e:finite-atomic} gives
$$
  \|  f    \|_{L^p}
=\lim_{m\to\infty}  \|  X_{\mathscr A_m}(f)    \|_{L^p}
\le Cp\,  \|  \mathcal S_{\mathscr A_0}(f)    \|_{L^p}.
$$
The proof is complete.
\end{proof}

\begin{remark}
The finite-sum estimate \eqref{e:finite-atomic} holds with $\lambda_R(f)$ replaced by any $c_R(f)\in\mathbb C$
and $\lambda_R(f)^2$ by $|c_R(f)|^2$.
It also holds for smooth atoms supported in $Q(R)$ satisfying
\eqref{e:prod-top-der} and \eqref{e:prod-cancel}, with uniformly bounded
multiplicity of $R\mapsto Q(R)$.
The constant $C$ depends only on these bounds.
\end{remark} 

\section{Sharp order-\texorpdfstring{$p$}{p} reverse estimates for Zygmund dilations}\label{sec:remaining}

We first apply the one-parameter results to the first variable.
We then keep the square integral in that variable and estimate the remaining
variables at the coupled scales $t,st$.
Under \eqref{e:supp-origin}, these comparisons give the order-$p$ reverse estimate in Theorem~\ref{t:primal}.
The final subsection proves a general lower bound under the basic kernel assumptions.

\subsection{The first-variable estimate}\label{sub:first-estimate}

{For $x_1\in\R$, write
$\Gamma_\alpha(x_1)=\{(u,s)\in\R\times(0,\infty):|u-x_1|<\alpha s\}$.
The intermediate square function in the first variable is
\begin{equation}\label{e:intro-intermediate}
S_{1,\alpha}^{\psi_1}f(x)^2
=\int_{\Gamma_\alpha(x_1)}
|(f(\cdot,x_2,x_3)*\psi_{1,s})(u)|^2\,\frac{du\,ds}{s^2},
\end{equation}
where $\psi_{1,s}(u)=s^{-1}\psi_1(u/s)$.}
We first record the continuity needed to pass from $L^2\cap L^p$ to $L^p$.

\begin{lemma}\label{l:continuous-forward}
For every $\alpha>0$ and $2\le p<\infty$, there is a finite constant $B_{p,\psi,\alpha}$ such that
$$
 \left\|   S_{1,\alpha}^{\psi_1}h \right\|_{L^p(\R^3)}
+ \left\|   S_{z,\alpha}^{\psi}h \right\|_{L^p(\R^3)}
\le B_{p,\psi,\alpha}  \|  h    \|_{L^p(\R^3)},
 \qquad  h\in L^p(\R^3).
$$
Moreover, both square functions converge in $L^p$ whenever their arguments converge in $L^p$.
\end{lemma}

\begin{proof}
For the first estimate, apply Lemma~\ref{sp:l:continuity} with $d=1$,
$\varphi=\psi_1$ and $g=h(\cdot,\pmb x)$, and integrate the $p$th powers in
$\pmb x\in\R^2$.
We also record the corresponding case of \eqref{sp:e:energy}:
\begin{equation}\label{e:ID}
\int_\R\int_0^\infty\int_{|u-x_1|<\alpha s}
|(g*\psi_{1,s})(u)|^2\,\frac{du\,ds}{s^2}\,dx_1
=2\alpha c_{\psi_1}\|g\|_{L^2(\R)}^2,
\qquad g\in L^2(\R),
\end{equation}
where $c_{\psi_1}$ is defined in \eqref{eqn:constant-psi}.

\medskip

Next we show the second estimate
$  \|  S_{z,\alpha}^{\psi}h    \|_{L^p(\R^3)}\lesssim  \|  h    \|_{L^p(\R^3)}$.
Integrating first in the cone centre $x$ and using \eqref{e:kernel-l2}, we obtain

$$
 \left\|   S_{z,\alpha}^{\psi}h \right\|_{L^2}^2
=(2\alpha)^3c_{\psi_1}c_{\psi_2}  \|  h    \|_{L^2}^2.
$$
For $p>2$, we use the weighted argument in \cite{FP}.
Let $r=(p/2)'$ and $\beta=\max\{1,\alpha\}$.
For $w\ge0$ with $\|w\|_{L^r}=1$, take $W$ and
$N_r=\|\mathsf M\|_{L^r\to L^r}$ as in the proof of
Lemma~\ref{l:atomic-forward}, where $\mathsf M$ is the strong
maximal operator on $\R^3$.
Then
$$
w\le W, \qquad    \|  W    \|_{L^r}\le2, \qquad  \text{and} \qquad  \mathsf MW\le2N_rW.
$$
Since $\beta\ge1$,
$$
\langle W\rangle_{\beta B_{s,t}(y)}
\le2N_r\mathop{\rm ess\,inf}_{\beta B_{s,t}(y)}W
\le2N_r\langle W\rangle_{B_{s,t}(y)}.
$$
To use $L^{p/2}$ duality, first restrict the square function to $N^{-1}\le s,t\le N$, with $N\ge2$,
and denote it by $S_{z,\beta;N}^{\psi}$.
By Minkowski's and Young's inequalities and
$\|\psi_{s,t}\|_{L^1}=\|\psi_1\|_{L^1}\|\psi_2\|_{L^1}$,
\begin{align*}
 \left\|   S_{z,\beta;N}^{\psi}h \right\|_{L^p}^2
\leq &(2\beta)^3\int_{N^{-1}}^N\int_{N^{-1}}^N
  \|  h*\psi_{s,t}    \|_{L^p}^2\,\frac{dt}{t}\frac{ds}{s}
\leq (2\beta)^3(2\log N)^2
  \|  \psi_1    \|_{L^1}^2  \|  \psi_2    \|_{L^1}^2  \|  h    \|_{L^p}^2<\infty.
\end{align*}
Integrating first in $x$, we apply the weighted estimate \cite[Theorem~2.9]{FP}:
\begin{align*}
\int_{\R^3}\left(S_{z,\beta;N}^{\psi}h\right)^2w
\leq &(2\beta)^3\int_0^\infty\int_0^\infty\int_{\R^3}
|h*\psi_{s,t}(y)|^2\langle W\rangle_{\beta B_{s,t}(y)}
\,dy\,\frac{dt}{t}\frac{ds}{s}\\
\leq &2N_r(2\beta)^3\int_0^\infty\int_0^\infty\int_{\R^3}
|h*\psi_{s,t}(y)|^2\langle W\rangle_{B_{s,t}(y)}
\,dy\,\frac{dt}{t}\frac{ds}{s}\\
\leq &C_{p,\psi,\alpha}\int_{\R^3}|h|^2W 
\leq  C_{p,\psi,\alpha}  \|  h    \|_{L^p}^2  \|  W    \|_{L^r}\\
\le& 2C_{p,\psi,\alpha}  \|  h    \|_{L^p}^2.
\end{align*}
Taking the supremum over $w$ and then letting $N\to\infty$,
$$
 \left\|   S_{z,\alpha}^{\psi}h \right\|_{L^p}^2
\le \left\|   S_{z,\beta}^{\psi}h \right\|_{L^p}^2
=\lim_{N\to\infty} \left\|   S_{z,\beta;N}^{\psi}h \right\|_{L^p}^2
\le C_{p,\psi,\alpha}  \|  h    \|_{L^p}^2.
$$

Finally, let $h_n,h\in L^p(\R^3)$ with $\|h_n-h\|_{L^p(\R^3)}\to0$ as $n\to\infty$.
The triangle inequality in the defining square integrals and the two bounds above give
\begin{align*}
  \left\|   S_{1,\alpha}^{\psi_1}h_n-S_{1,\alpha}^{\psi_1}h \right\|_{L^p}
+ \left\|   S_{z,\alpha}^{\psi}h_n-S_{z,\alpha}^{\psi}h \right\|_{L^p} 
\leq & \left\|   S_{1,\alpha}^{\psi_1}(h_n-h) \right\|_{L^p}
+ \left\|   S_{z,\alpha}^{\psi}(h_n-h) \right\|_{L^p}\\
\leq &B_{p,\psi,\alpha}  \|  h_n-h    \|_{L^p}\to0,
\end{align*}
as desired.
\end{proof}

\begin{proposition}\label{prop:square-1}
Assume $0\in\operatorname{supp}(|D|\psi_1)$.
For every $\alpha>0$, $2\le p<\infty$, and $f\in L^p(\R^3)$,
\begin{equation}\label{e:reverse-one}
 \left\|   f \right\|_{L^p(\R^3)}
\le C_{\psi_1,\alpha}\sqrt p\,
 \left\|   S_{1,\alpha}^{\psi_1}f \right\|_{L^p(\R^3)},
\end{equation}
where $C_{\psi_1,\alpha}$ is independent of $p$.
\end{proposition}

\begin{proof}
For $f\in L^p(\R^3)$, Fubini's theorem gives
$f(\cdot,\pmb x)\in L^p(\R)$ for almost every $\pmb x=(x_2,x_3)$, and
$$
S_{1,\alpha}^{\psi_1}f(x_1,\pmb x)
=\mathcal A_{\psi_1,\alpha}[f(\cdot,\pmb x)](x_1).
$$
Apply the definition \eqref{sp:e:best} with $d=1$ to each such function,
take $p$th powers, and integrate in $\pmb x$. We obtain
\begin{align}\label{e:first-from-single}
\|f\|_{L^p(\R^3)}^p
= &\int_{\R^2}\|f(\cdot,\pmb x)\|_{L^p(\R)}^p\,d\pmb x
\le \bigl(\mathfrak c_{p,1}^{\psi_1}(\alpha)\bigr)^p
\int_{\R^2}
\|\mathcal A_{\psi_1,\alpha}[f(\cdot,\pmb x)]\|_{L^p(\R)}^p\,d\pmb x\notag\\
= &\bigl(\mathfrak c_{p,1}^{\psi_1}(\alpha)\bigr)^p
\|S_{1,\alpha}^{\psi_1}f\|_{L^p(\R^3)}^p.
\end{align}
The square-function norms are finite by Lemma~\ref{l:continuous-forward}.
Since $0\in\operatorname{supp}(|D|\psi_1)$, Theorem~\ref{sp:t:general}
with $d=1$ gives
$$
\mathfrak c_{p,1}^{\psi_1}(\alpha)
\le C_{\psi_1,\alpha}\sqrt p.
$$
Taking the $p$th root in \eqref{e:first-from-single} proves
\eqref{e:reverse-one}.
The same calculation also gives the estimate under
$d_1(\psi_1)<c_1\alpha$ and the general order-$p$ bound,
using the corresponding conclusions of Theorem~\ref{sp:t:general},
where $c_1$ is its constant for $d=1$.
\end{proof}

For the estimate \eqref{e:reverse-one},  the condition $0\in\operatorname{supp}(|D|\psi_1)$
in Proposition~\ref{prop:square-1} cannot be removed in general.
For every fixed aperture $\alpha>0$, we can construct a fixed kernel
for which the reverse constant is at least $cp$,
where $c>0$ is independent of $p$.  See Example \ref{r:counterexample-one}.

\subsection{The reverse estimate at the coupled scales}\label{sub:remaining-estimate}

In the previous subsection, we estimated $\|f\|_{L^p}$ in terms of
$\|S_{1,\alpha}^{\psi_1}f\|_{L^p}$.
We now compare $S_{1,\alpha}^{\psi_1}f$ with
$S_{z,\alpha}^{\psi}f$ in the remaining variables $(x_2,x_3)$.
The estimate in Proposition~\ref{prop:square-2} below,
combined with \eqref{e:reverse-one}, yields the reverse estimate
in Theorem~\ref{t:primal}, with $\sqrt p\cdot\sqrt p=p$.

For $x=(x_1,x_2,x_3)$, recall that the square function is given by
$$
S_{z,\alpha}^{\psi}f(x)^2
=\int_0^\infty\int_0^\infty
\int_{\substack{|u-x_1|<\alpha s\\   |  v_2-x_2|<\alpha t\\   |  v_3-x_3|<\alpha st}}
|f*\psi_{s,t}(u,v_2,v_3)|^2\,
\frac{du\,dv_2\,dv_3\,dt\,ds}{s^3t^3}.
$$

\begin{proposition}\label{prop:square-2}
Assume $0\in\operatorname{supp}\psi_2$.
For every $\alpha>0$, $2\le p<\infty$, and $f\in L^p(\R^3)$,
\begin{equation}\label{e:reverse-two}
     \left\|   S_{1,\alpha}^{\psi_1}f \right\|_{L^p(\R^3)}
    \le C_{\psi_2,\alpha}\sqrt p\,
     \left\|   S_{z,\alpha}^{\psi}f \right\|_{L^p(\R^3)}.
\end{equation}
The constant is independent of $p$ and $f$.
\end{proposition}

\smallskip

\begin{remark}\label{rem:square-2}
{For product square functions, Pipher \cite{Pipher86} proved exponential
integrability under a bound on the double square function. Here we adapt the
intermediate square function method of \cite[Sections~4--5]{LPW} to the coupled
scales $t,st$.} The comparison in Proposition~\ref{prop:square-2},
using the continuous square function $S_{1,\alpha}^{\psi_1}f$, plays the role of
\eqref{eqn:S1-2} in the Haar argument of Proposition~\ref{p:haar-two}.
When both support assumptions hold, we have
$$
\|f\|_{L^p(\R^3)}
\le C\sqrt p\,\|S_{1,\alpha}^{\psi_1}f\|_{L^p(\R^3)}
\le Cp\,\|S_{z,\alpha}^{\psi}f\|_{L^p(\R^3)}.
$$
For the second inequality, we decompose the first-variable convolutions
in the definition of $S_{1,\alpha}^{\psi_1}f$
into smooth atomic functions in $(x_2,x_3)$. The following observations explain
how the Zygmund scale relation is used in their Haar expansions.

\begin{itemize}
\item \textbf{The Zygmund scale relation.}
{At each first-variable scale $s$, the rectangles in $(x_2,x_3)$
have side lengths $t$ and $st$. We keep the integral over the first-variable
cone in the $L^2$ norm of each Haar coefficient.}

For $2^m\le s<2^{m+1}$, we use atoms indexed by rectangles
$R=I_2\times I_3$ with $|I_3|=2^m|I_2|$.
Their supports in $(y_2,y_3)$ lie in $(5I_2)\times(9I_3)$.
We choose covering intervals $Q_j(R)$ from the adjacent grids using Lemma~\ref{l:adj-int-cover},
and group the atoms according to the two selected grids and the two enlargement factors, where
$$
|Q_2(R)|=2^{a_2}|I_2|,
\qquad
|Q_3(R)|=2^{a_3}|I_3|.
$$
Fix the dyadic grids, $a_2,a_3$, and the relative generations
$n_2,n_3\in\mathbb N_0$.
The Haar intervals $J,K$ have lengths
$$
|J|=2^{-n_2}|Q_2(R)|,
\qquad
|K|=2^{-n_3}|Q_3(R)|.
$$
Hence
$$
\frac{|K|}{|J|}
=2^{a_3-a_2+n_2-n_3}\frac{|I_3|}{|I_2|}
=2^{m+a_3-a_2+n_2-n_3},
$$
as in \eqref{e:second-ratio}.
For fixed $m,J$, the intervals $K$ have the same length and are disjoint.
The coefficients for each $m$ vanish unless $2^m\le s<2^{m+1}$,
so coefficients for different $m$ have disjoint supports in $s$.
Thus there are no cross terms when we square the sum over $m,K$.

\smallskip
  \item \textbf{One Haar estimate at fixed generations.}
Here we only have cancellation in $(x_2,x_3)$ jointly,
so we include all three Haar types:
$$
h_Jh_K^0,\qquad h_J^0h_K,\qquad h_Jh_K.
$$
These are product Haar functions, but their scales are not independent.
For each type, we first fix the generations $n_2,n_3$.
By the Zygmund scale relation and the disjointness above,
we apply the $L^2$-valued one-dimensional Haar estimate
only once, with a factor $\sqrt p$:
for $h_Jh_K^0$ and $h_Jh_K$, we apply it in $y_2$;
for $h_J^0h_K$, we apply it in $y_3$.

\smallskip
\item \textbf{Summing over the generations.}
In the Haar expansion, the factors with superscript $0$
are $h_{Q_2(R)}^0$ and $h_{Q_3(R)}^0$.
Thus, for $\eta\in\{(1,0),(0,1),(1,1)\}$ and $j=2,3$,
we have $n_j=0$ whenever $\eta_j=0$, and hence
$$
\eta_2n_2+\eta_3n_3=n_2+n_3.
$$
By \eqref{e:second-depth} and the Haar estimate,
the total factor after summing over the three types and their allowed generations is
$$
C\sum_{\substack{\eta\in\{(1,0),(0,1),(1,1)\}\\
n_2,n_3\ge0,\ (1-\eta_j)n_j=0\ (j=2,3)}}
\underbrace{2^{-n_2-n_3}}_{\text{smoothness}}
\,\underbrace{\sqrt p}_{\text{one Haar estimate}}
=C\sqrt p\,(2+2+4)
=8C\sqrt p.
$$
The constant $C$ is independent of $p$.
\end{itemize}
{For atomic decompositions in the flag setting, see \cite{HLLW}.
In the flag concentration estimates of \cite{HLWW}, the two scales are $t$
and $\max\{t,s\}$.}
\end{remark}

\smallskip

\begin{remark}
For the estimate \eqref{e:reverse-two}, the condition
$0\in\operatorname{supp}\psi_2$
in Proposition~\ref{prop:square-2} cannot be removed in general,
even when $0\in\operatorname{supp}(|D|\psi_1)$.
For every fixed aperture $\alpha>0$,
Example~\ref{r:counterexample-two} constructs a fixed kernel pair with
$$
0\in\operatorname{supp}(|D|\psi_1),
\qquad
0\notin\operatorname{supp}\psi_2,
$$
and functions $f_p\in\mathcal S(\R^3)$ such that
$$
\frac{\|S_{1,\alpha}^{\psi_1}f_p\|_{L^p(\R^3)}}
{\|S_{z,\alpha}^{\psi}f_p\|_{L^p(\R^3)}}
\ge c_{\psi,\alpha}p,
\qquad 2\le p<\infty.
$$
\end{remark}

\begin{remark}\label{rem:support-why}
 The different forms of the support assumptions in Propositions~\ref{prop:square-1} and \ref{prop:square-2}
come from the radial measures $dr$ and $r\,dr$ in dimensions one and two, respectively.
For real $\chi_1\in C_c^\infty(\R)$ and $\chi_2\in C_c^\infty(\R^2)$,
with $\chi_1$ even and $\chi_2$ radial, we have 
$$
\int_0^\infty\widehat\psi_1(r)\widehat{(-\chi_1'')}(r)\,\frac{dr}{r}
=  4\pi^2\int_0^\infty r\widehat\psi_1(r)\widehat\chi_1(r)\,dr 
=  2\pi^2\int_\R(|D|\psi_1)\chi_1,
$$
while
\begin{equation}\label{e:pair-two}
\int_0^\infty\widehat\psi_2(r,0)\widehat{(-\Delta\chi_2)}(r,0)\,\frac{dr}{r}
=  4\pi^2\int_0^\infty r\widehat\psi_2(r,0)\widehat\chi_2(r,0)\,dr 
=  2\pi\int_{\R^2}\psi_2\chi_2.
\end{equation}
Taking $\chi_1=\zeta_1|D|\psi_1$ and $\chi_2=\zeta_2\psi_2$,
the respective support assumptions allow
$$
\int_\R\zeta_1(|D|\psi_1)^2>0, \qquad 
\int_{\R^2}\zeta_2\psi_2^2>0
$$
with nonnegative cutoffs supported in arbitrarily small neighborhoods of $0$.
Thus, for any fixed aperture $\alpha>0$, we can choose the normalized dual kernels
with support radius $\rho$ satisfying $16\rho<\alpha$
in the one-parameter argument \eqref{sp:e:discrete-area} at $d=1$,
and $54\rho<\alpha$ in Proposition~\ref{prop:square-2}.
\end{remark}

\begin{proof}[Proof of Proposition~\ref{prop:square-2}]

Take $\rho\in2^{\mathbb Z}$ with $54\rho<\alpha$.

\noindent{\it Step 1.} 
{First let $f\in L^2(\R^3)\cap L^p(\R^3)$. Since $0\in\operatorname{supp}\psi_2$,
choose a nonnegative radial $\zeta\in C_c^\infty(B(0,\rho))$ with
$\int_{\R^2}\zeta\psi_2^2>0$. As in Remark~\ref{rem:support-why}, let
$$
\Phi^{(2)}=-a^{-1}\Delta(\zeta\psi_2),
\qquad a=2\pi\int_{\R^2}\zeta\psi_2^2.
$$}
Then $\Phi^{(2)}$ is real and radial, with
$$
\operatorname{supp}\Phi^{(2)}\subset B(0,\rho) \qquad \text{and} \qquad  
\int_{\R^2}\Phi^{(2)}=0.
$$
By \eqref{e:pair-two},
\begin{equation}\label{e:dual-second}
\int_0^\infty\widehat\psi_2(r,0)\widehat{\Phi^{(2)}}(r,0)\,\frac{dr}{r}
=\frac{2\pi}{a}\int_{\R^2}\zeta\psi_2^2=1.
\end{equation}
For $(x_2,x_3)\in\R^2$, define
\begin{equation}\label{eqn:psi2}
\psi_{s,t}^{(2)}(x_2,x_3)
=\frac1{st^2}\psi_2\left(\frac{x_2}{t},\frac{x_3}{st}\right),
 \qquad 
\Phi_{s,t}^{(2)}(x_2,x_3)
=\frac1{st^2}\Phi^{(2)}\left(\frac{x_2}{t},\frac{x_3}{st}\right).
\end{equation}
By radiality and \eqref{e:dual-second}, we have for $s>0$ and $\pmb\xi\ne0$ that
$$
\int_0^\infty\widehat\psi_2(t\xi_2,st\xi_3)
\widehat{\Phi^{(2)}}(t\xi_2,st\xi_3)\,\frac{dt}{t}
=\int_0^\infty\widehat\psi_2(r,0)
\widehat{\Phi^{(2)}}(r,0)\,\frac{dr}{r}=1.
$$

\smallskip

\noindent{\it Step 2.}  
We decompose $f(\cdot,x_2,x_3)*\psi_{1,s}$ into smooth atoms in $(x_2,x_3)$
and compare their discrete square function with $S_{z,\alpha}^{\psi}f$.
{For $x_1\in\R$ and $m\in\mathbb Z$, write
$$
\Gamma_{\alpha,m}(x_1)
=\{(u,s)\in\Gamma_\alpha(x_1):2^m\le s<2^{m+1}\}.
$$}
We have
$$
(f*\psi_{s,t})(u,v_2,v_3)
=\int_{\R^2}(f(\cdot,y_2,y_3)*\psi_{1,s})(u)
\psi_{s,t}^{(2)}(v_2-y_2,v_3-y_3)\,dy_2\,dy_3,
$$
and
$$
\int_{\Gamma_\alpha(x_1)}
\left|(f(\cdot,y_2,y_3)*\psi_{1,s})(u)\right|^2
\,\frac{du\,ds}{s^2}
=S_{1,\alpha}^{\psi_1}f(x_1,y_2,y_3)^2.
$$
 
Write $\pmb y=(y_2,y_3)$ and $\pmb v=(v_2,v_3)$.
By the first-variable $L^2$ identity in the proof of Lemma~\ref{l:continuous-forward}, we have
$$
\int_\R\int_{\R^2}\int_{\Gamma_\alpha(x_1)}
\left|(f(\cdot,\pmb y)*\psi_{1,s})(u)\right|^2
\,\frac{du\,ds}{s^2}\,d\pmb y\,dx_1
=2\alpha c_{\psi_1}\left\|f\right\|_{L^2(\R^3)}^2<\infty.
$$

Thus,
\begin{equation}\label{eqn:aux-finite}
\int_{\R^2}\int_{\Gamma_\alpha(x_1)}
\left|(f(\cdot,\pmb y)*\psi_{1,s})(u)\right|^2
\,\frac{du\,ds}{s^2}\,d\pmb y<\infty
 \qquad \text{for a.e. }x_1.
\end{equation}
We fix such an $x_1$ now.   For $s>0$ and $\pmb\xi=(\xi_2,\xi_3)\ne0$,  it  follows from the definition of $\psi_{s,t}^{(2)}$ in \eqref{eqn:psi2} that 
$
\widehat{\psi_{s,t}^{(2)}}(\pmb\xi)
=\widehat\psi_2(t\xi_2,st\xi_3).
$
Moreover,
by radiality and the substitution $r=t(\xi_2^2+s^2\xi_3^2)^{1/2}$,
$$
\int_0^\infty|\widehat\psi_2(t\xi_2,st\xi_3)|^2\,\frac{dt}{t}
=\int_0^\infty|\widehat\psi_2(r,0)|^2\,\frac{dr}{r}
=c_{\psi_2},
$$
where $c_{\psi_2}$ is defined in \eqref{eqn:constant-psi}.
For this fixed $x_1$, Plancherel's identity in the last two variables identifies
the square integral of $f*\psi_{s,t}$ with the finite first-variable square integral in \eqref{eqn:aux-finite}:
\begin{align}\label{e:second-energy}
&\int_{\Gamma_\alpha(x_1)}\int_0^\infty\int_{\R^2}
|(f*\psi_{s,t})(u,\pmb v)|^2
\,d\pmb v\,\frac{dt}{t}\,\frac{du\,ds}{s^2}\notag\\
= &c_{\psi_2}\int_{\R^2}\int_{\Gamma_\alpha(x_1)}
\left|(f(\cdot,\pmb y)*\psi_{1,s})(u)\right|^2
\,\frac{du\,ds}{s^2}\,d\pmb y<\infty.
\end{align}

Recall that $\rho\in2^{\mathbb Z}$.  In the rest of this proof, $R$ denotes a dyadic rectangle in $\R^2$.
For $m,k\in\mathbb Z$, define
$$
\mathscr R_{m,k}
=\left\{R=I_2\times I_3:I_2,I_3\in\mathcal D^0,\
|I_2|=\rho2^k,\ |I_3|=\rho2^{m+k}\right\},
 \qquad 
\mathscr R=\bigcup_{m,k\in\mathbb Z}\mathscr R_{m,k}.
$$
Define
$$
\mathcal W^{(2)}(R)
=\left\{(u,s,\pmb v,t)\in\R^5:
(u,s)\in\Gamma_{\alpha,m}(x_1),\
\pmb v\in R,\ 2^k\le t<2^{k+1}\right\}\quad\text{for}\quad R\in\mathscr R_{m,k}.
$$
Furthermore, one may verify that 
\begin{equation}\label{W2R}
\Gamma_\alpha(x_1)\times\R^2\times(0,\infty)
=\bigcup_{R\in\mathscr R}\mathcal W^{(2)}(R)
\quad\text{(disjoint)}.
 \end{equation}

For $R\in\mathscr R_{m,k}$, define
$$
\lambda_R^{(2)}(f)
= \bigg( \int_{\mathcal W^{(2)}(R)}
|(f*\psi_{s,t})(u,\pmb v)|^2
\,\frac{du\,ds\,d\pmb v\,dt}{s^2t}  \bigg)^{1/2},
$$
and
$$
\widetilde{\mathsf{a}}_R^{(2)}(f)(u,s;\pmb y)
=\mathbf 1_{\Gamma_{\alpha,m}(x_1)}(u,s)
\int_{2^k}^{2^{k+1}}\int_R
(f*\psi_{s,t})(u,\pmb v)\Phi_{s,t}^{(2)}(\pmb y-\pmb v)
\,d\pmb v\,\frac{dt}{t}.
 $$
The corresponding Calder\'on reproducing formula is
\begin{align}\label{CRF2}
(f(\cdot,\pmb y)*\psi_{1,s})(u)
= &\int_0^\infty\int_{\R^2}
(f*\psi_{s,t})(u,\pmb v)\Phi_{s,t}^{(2)}(\pmb y-\pmb v)
\,d\pmb v\,\frac{dt}{t}\notag\\[2pt]
= &\sum_{R\in\mathscr R}\widetilde{\mathsf{a}}_R^{(2)}(f)(u,s;\pmb y),
 \qquad  (u,s)\in\Gamma_\alpha(x_1),
\end{align}
converging in
$ 
L^2\left(\Gamma_\alpha(x_1)\times\R^2,
\frac{du\,ds}{s^2}\,d\pmb y\right)
$, which will be  
 verified in {\it Step 4}.

Let
$$
\mathsf{a}_R^{(2)}(f)=\frac{\widetilde{\mathsf{a}}_R^{(2)}(f)}{\lambda_R^{(2)}(f)}
$$
with $\mathsf{a}_R^{(2)}(f)=0$ when $\lambda_R^{(2)}(f)=0$.
By \eqref{e:second-energy} and \eqref{W2R},
$$
\sum_{R\in\mathscr R}\lambda_R^{(2)}(f)^2<\infty.
$$
If $\Phi_{s,t}^{(2)}(\pmb y-\pmb v)\ne0$, with $2^m\le s<2^{m+1}$ and $2^k\le t<2^{k+1}$, then
$$
|y_2-v_2|\le\rho t<2|I_2|,
\qquad
|y_3-v_3|\le\rho st<4|I_3|.
$$
Since $\pmb v\in I_2\times I_3$, we obtain
$$
\operatorname{supp}\mathsf{a}_R^{(2)}(f)
\subseteq\overline{\Gamma_{\alpha,m}(x_1)}
\times(5I_2)\times(9I_3),
\qquad
\int_{\R^2}\mathsf{a}_R^{(2)}(f)(u,s;\pmb y)\,d\pmb y=0.
$$
The enlarged intervals need not lie in $\mathcal D^0$.
Apply Lemma~\ref{l:adj-int-cover} separately to $5I_2$ and $9I_3$, and fix the resulting covering intervals
and their grid labels for each $R$. Then $Q(R)=Q_2(R)\times Q_3(R)$ satisfies
$$
(5I_2)\times(9I_3)\subset Q(R),
 \qquad 
Q_j(R)\in\mathcal  D^{ \varepsilon_j(R)},
\quad \varepsilon_j(R)\in\{0,1\},\quad j=2,3,
$$
and
$$
5|I_2|\le|Q_2(R)|\le30|I_2|,
 \qquad 
9|I_3|\le|Q_3(R)|\le54|I_3|.
$$

Since the interval lengths are powers of two, the possible enlargements are
$$
|Q_2(R)|=2^{a_2(R)}\rho2^k, \qquad
|Q_3(R)|=2^{a_3(R)}\rho2^{m+k}, \qquad
a_2(R)\in\{3,4\},\quad a_3(R)\in\{4,5\}.
$$
Both the grid labels and the enlargement factors may depend on $R$.
We fix them by passing to the following classes. For
$\vec{\pmb\varepsilon}=(\varepsilon_2,\varepsilon_3)\in\{0,1\}^2$
and $\vec{\pmb a}=(a_2,a_3)\in\{3,4\}\times\{4,5\}$, define
$$
\mathscr R_{ \vec{ \pmb\varepsilon}, \vec{ \pmb a}}
=\left\{R\in\mathscr R:
\varepsilon_j(R)=\varepsilon_j,\
|Q_j(R)|=2^{a_j}|I_j|,\ j=2,3\right\}.
$$
These classes give a disjoint decomposition of $\mathscr R$, with at most $4\cdot4=16$ nonempty classes.
Within each class, the covering intervals belong to one fixed pair of grids,
and the ratios $|Q_j(R)|/|I_j|=2^{a_j}$ are fixed.
For $R\in\mathscr R_{ \vec{ \pmb\varepsilon}, \vec{ \pmb a}}\cap\mathscr R_{m,k}$,
$$
\frac{|Q_3(R)|}{|Q_2(R)|}=2^{m+a_3-a_2},
$$
and, writing $Q(R)=Q_2\times Q_3$, we have
$$
2^k=\frac{|Q_2|}{2^{a_2}\rho},
 \qquad 
2^m=2^{a_2-a_3}\frac{|Q_3|}{|Q_2|}.
$$
Hence, within each fixed class $\mathscr R_{ \vec{ \pmb\varepsilon}, \vec{ \pmb a}}$, all $R$ assigned to the same $Q$
belong to the same $\mathscr R_{m,k}$.
Since these $R$ are disjoint and   $R\subset Q$,
\begin{equation}\label{multi-23}
\#\left\{R\in\mathscr R_{ \vec{ \pmb\varepsilon}, \vec{ \pmb a}}:Q(R)=Q\right\}
\le\frac{|Q|}{\rho^22^{m+2k}}
=2^{a_2+a_3}\le C.
 \end{equation}

We next estimate  derivatives of $\mathsf{a}_R^{(2)}(f)$,
which will be used to bound its Haar coefficients in {\it Step 3}.
For $\gamma_2,\gamma_3\in\mathbb N_0$ with $\gamma_2+\gamma_3\le2$, scaling gives
\begin{align*}
\left\|\partial_2^{\gamma_2}\partial_3^{\gamma_3}\Phi_{s,t}^{(2)}\right\|_{L^2(\R^2)}^2
= &s^{-1-2\gamma_3}t^{-2-2\gamma_2-2\gamma_3}
\left\|\partial_2^{\gamma_2}\partial_3^{\gamma_3}\Phi^{(2)}\right\|_{L^2(\R^2)}^2
\leq C_{\psi_2,\alpha}s^{-1-2\gamma_3}t^{-2-2\gamma_2-2\gamma_3}.
\end{align*}

For $\lambda_R^{(2)}(f)\neq 0$, by Cauchy--Schwarz again, 
\begin{align*}
 \int_{\Gamma_\alpha(x_1)}
\left|\partial_2^{\gamma_2}\partial_3^{\gamma_3}
\mathsf{a}_R^{(2)}(f)(u,s;\pmb y)\right|^2\,\frac{du\,ds}{s^2} 
\leq &\frac1{\lambda_R^{(2)}(f)^2}
\int_{\Gamma_{\alpha,m}(x_1)}
 \bigg( \int_{2^k}^{2^{k+1}}\int_R
|(f*\psi_{s,t})(u,\pmb v)|^2\,d\pmb v\,\frac{dt}{t}  \bigg)\\
&\quad \times  \bigg( \int_{2^k}^{2^{k+1}}\int_R
\left|\partial_2^{\gamma_2}\partial_3^{\gamma_3}
\Phi_{s,t}^{(2)}(\pmb y-\pmb v)\right|^2
\,d\pmb v\,\frac{dt}{t}  \bigg)\,\frac{du\,ds}{s^2}\\
\leq &\sup_{2^m\le s<2^{m+1}}
\int_{2^k}^{2^{k+1}}
 \left\|   \partial_2^{\gamma_2}\partial_3^{\gamma_3}
\Phi_{s,t}^{(2)} \right\|_{L^2(\R^2)}^2\,\frac{dt}{t}\\
\leq &C2^{-m(1+2\gamma_3)}
\int_{2^k}^{2^{k+1}}t^{-3-2\gamma_2-2\gamma_3}\,dt\\
\leq &C|Q(R)|^{-1}|Q_2(R)|^{-2\gamma_2}|Q_3(R)|^{-2\gamma_3},
\end{align*}
where $C$ depends only on $\psi_2$ and $\alpha$.
If $\lambda_R^{(2)}(f)=0$, then $\mathsf{a}_R^{(2)}(f)=0$, so the same bound holds.

For $\pmb y\in Q(R)$ and $(u,s,\pmb v,t)\in\mathcal W^{(2)}(R)$,
$$
|v_2-y_2|\le30\rho2^k<\alpha t,
 \qquad 
|v_3-y_3|\le54\rho2^{m+k}<\alpha st,
$$
and
$$
\frac1{|Q(R)|}
\le\frac1{\rho^22^{m+2k}}
\le\frac8{\rho^2st^2}.
$$
For a finite  $\mathscr F\subset\mathscr R$, define
$$
\mathcal S_{\mathscr F}^{(2)}(f)(\pmb y)^2
=\sum_{R\in\mathscr F}\lambda_R^{(2)}(f)^2
\frac{\mathbf 1_{Q(R)}(\pmb y)}{|Q(R)|}.
$$
By the disjointness of $\mathcal W^{(2)}(R)$,
\begin{align}\label{e:second-square}
\mathcal S_{\mathscr F}^{(2)}(f)(\pmb y)^2
\leq &\frac8{\rho^2}
\sum_{R\in\mathscr F}\mathbf 1_{Q(R)}(\pmb y)
\int_{\mathcal W^{(2)}(R)}
|(f*\psi_{s,t})(u,\pmb v)|^2
\,\frac{du\,ds\,d\pmb v\,dt}{s^3t^3}\notag\\
\leq &\frac8{\rho^2}\int_{\Gamma_\alpha(x_1)}\int_0^\infty
\int_{\substack{|v_2-y_2|<\alpha t\\  |  v_3-y_3|<\alpha st}}
|(f*\psi_{s,t})(u,\pmb v)|^2
\,\frac{d\pmb v\,dt}{st^3}\,\frac{du\,ds}{s^2}
= \frac8{\rho^2}S_{z,\alpha}^{\psi}f(x_1,\pmb y)^2.
\end{align}

\noindent{\it Step 3.}  
We now expand the atoms $\mathsf{a}_R^{(2)}(f)$ in terms of Haar functions.
 Using the Zygmund scale relation, we will apply the $L^2$-valued form of
Proposition~\ref{p:haar-one} to obtain the factor $\sqrt p$.

To obtain this $L^2$-valued form, we rewrite \eqref{eqn:S1-2}
in terms of coefficients depending on $(u,s)$.
Let $v_J\in L^2(\Gamma_\alpha(x_1),du\,ds/s^2)$ be finitely many functions,
with $J$ in one dyadic system.
Their real linear combinations form a finite-dimensional space.
By Gram--Schmidt, write
$$
v_J=\sum_{\ell=1}^D v_{J,\ell}e_\ell,
\qquad v_{J,\ell}\in\R,
$$
where
$$
\operatorname{Re}\int_{\Gamma_\alpha(x_1)}
e_\ell(u,s)\overline{e_k(u,s)}\,\frac{du\,ds}{s^2}
=\delta_{\ell k}.
$$
Choose dyadic intervals $I_D\subsetneq\cdots\subsetneq I_1$
and an interior point $x_0$ of $I_D$.
Take
$$
a_{I_\ell,J}=|I_\ell|^{1/2}v_{J,\ell},
 \qquad  \ell=1,\ldots,D,
$$
with all other coefficients zero, and apply \eqref{eqn:S1-2} at $x_0$.
Since $\mathbf 1_{I_\ell}(x_0)=1$ for every $\ell$, by
orthonormality and \eqref{eqn:S1-2}, we obtain 
\begin{align}\label{e:haar-hilbert}
\bigg  \|   \bigg( \int_{\Gamma_\alpha(x_1)}
\bigg|\sum_Jv_J(u,s)h_J(x)\bigg|^2
\,\frac{du\,ds}{s^2}  \bigg)^{1/2}\bigg    \|_{L^p(\R)} 
= \,  &\bigg  \|   \bigg( \sum_{\ell=1}^D
\bigg|\sum_Jv_{J,\ell}h_J(x)\bigg|^2  \bigg)^{1/2}\bigg    \|_{L^p(\R)}\notag\\
= \,  &\bigg  \|   \bigg( \sum_{\ell=1}^D
\bigg|\sum_Ja_{I_\ell,J}h_J(x)\bigg|^2
\frac{\mathbf 1_{I_\ell}(x_0)}{|I_\ell|}
  \bigg)^{1/2}\bigg    \|_{L^p(\R)}\notag\\
\leq  \,  &C\sqrt p\,
\bigg  \|   \bigg( \sum_{\ell=1}^D\sum_J
|a_{I_\ell,J}|^2
\frac{\mathbf 1_{I_\ell}(x_0)}{|I_\ell|}
\frac{\mathbf 1_J(x)}{|J|}
  \bigg)^{1/2}\bigg    \|_{L^p(\R)}\notag\\
=\,   &C\sqrt p\,
\bigg  \|   \bigg( \sum_J
 \bigg( \int_{\Gamma_\alpha(x_1)}|v_J(u,s)|^2\,\frac{du\,ds}{s^2}  \bigg)
\frac{\mathbf 1_J(x)}{|J|}
  \bigg)^{1/2}\bigg    \|_{L^p(\R)}.
\end{align}
The constant is independent of $D,x_0$ and $x_1$.
The same proof applies to any $L^2$ space, with the same constant.

Fix a class $\mathscr R_{ \vec{ \pmb\varepsilon}, \vec{ \pmb a}}$
and a finite family $\mathscr F$ of original rectangles in this class.
All Haar intervals in $y_2$ will belong to $\mathcal D^{\varepsilon_2}$,
and all Haar intervals in $y_3$ to $\mathcal D^{\varepsilon_3}$.
Write $ \vec{ \pmb n}=(n_2,n_3)$ and $| \vec{ \pmb n}|=n_2+n_3$.
For $Q=Q_2\times Q_3$, let
$$
Ch_{ \vec{ \pmb n}}(Q)
=\Big \{J\times K:J\in\mathcal  D^{ \varepsilon_2},\
K\in\mathcal  D^{ \varepsilon_3},\
J\subseteq Q_2,\ K\subseteq Q_3,\
|J|=2^{-n_2}|Q_2|,\ |K|=2^{-n_3}|Q_3|\,\Big \}.
$$
Use the local averages $\mathbb E_{j,n}^Q$ and differences $\Delta_{j,n}^Q$
from Proposition~\ref{p:smooth-model}, acting in $y_j$.
For $Q=Q(R)$,
\begin{align}\label{H-expan2}
\mathsf{a}_R^{(2)}(f)
= &\underbrace{\mathbb E_{2,0}^Q\mathbb E_{3,0}^Q\mathsf{a}_R^{(2)}(f)}_{=0}
+\sum_{n_2=0}^{\infty}\Delta_{2,n_2}^Q\mathbb E_{3,0}^Q\mathsf{a}_R^{(2)}(f)\notag\\
&\quad  +\sum_{n_3=0}^{\infty}\mathbb E_{2,0}^Q\Delta_{3,n_3}^Q\mathsf{a}_R^{(2)}(f)
+\sum_{n_2,n_3=0}^{\infty}\Delta_{2,n_2}^Q\Delta_{3,n_3}^Q\mathsf{a}_R^{(2)}(f)\notag\\
= &\sum_{\eta, \vec{ \pmb n}}
\sum_{J\times K\in Ch_{ \vec{ \pmb n}}(Q(R))}
c_{R,J\times K}^{\eta}(f)
h_J^{\eta_2}  h_K^{\eta_3},
\end{align}
and 
$$
\lim_{N\to\infty}
\left\|\mathbb E_{2,N}^Q\mathbb E_{3,N}^Q\mathsf{a}_R^{(2)}(f)
-\mathsf{a}_R^{(2)}(f)\right\|_{
L^2  (\Gamma_\alpha(x_1)\times\R^2,\frac{du\,ds}{s^2}\,d\pmb y)}=0,
$$
where
$$
c_{R,J\times K}^{\eta}(f)(u,s)
=\int_{Q(R)}\mathsf{a}_R^{(2)}(f)(u,s;\pmb y)
h_J^{\eta_2}(y_2)h_K^{\eta_3}(y_3)\,d\pmb y,
$$
with
$$
\eta\in\mathsf H=\{(1,0),(0,1),(1,1)\}, \qquad 
n_j\in\mathbb N_0,\quad (1-\eta_j)n_j=0\quad(j=2,3).
$$

In particular, $\eta_j=0$ implies $n_j=0$. Consequently,
$$
\eta_2n_2+\eta_3n_3=| \vec{ \pmb n}|.
$$
This, combined with the cancellation of Haar functions and {\it Step 2}, yields   
\begin{align*}
& \bigg( \int_{\Gamma_\alpha(x_1)}
|c_{R,J\times K}^{\eta}(f)(u,s)|^2\,\frac{du\,ds}{s^2}  \bigg)^{1/2}\\
\leq \,  &|J\times K|^{1/2}|J|^{\eta_2}|K|^{\eta_3}
\sup_{\pmb y\in Q(R)}
 \bigg( \int_{\Gamma_\alpha(x_1)}
\left|\partial_2^{\eta_2}\partial_3^{\eta_3}
\mathsf{a}_R^{(2)}(f)(u,s;\pmb y)\right|^2
\,\frac{du\,ds}{s^2}  \bigg)^{1/2}\\
\leq  \,   &C\left(\frac{|J\times K|}{|Q(R)|}\right)^{1/2}
\left(\frac{|J|}{|Q_2(R)|}\right)^{\eta_2}
\left(\frac{|K|}{|Q_3(R)|}\right)^{\eta_3}\\
= \,  &C2^{-| \vec{ \pmb n}|/2-\eta_2n_2-\eta_3n_3}
=C2^{-3| \vec{ \pmb n}|/2}.
\end{align*}
 Then
\begin{align*}
& \bigg( \int_{\Gamma_\alpha(x_1)}
\bigg|\sum_{J\times K\in Ch_{ \vec{ \pmb n}}(Q(R))}
c_{R,J\times K}^{\eta}(f)(u,s)
 h_J^{\eta_2}(y_2)h_K^{\eta_3}(y_3)\bigg|^2
\,\frac{du\,ds}{s^2}  \bigg)^{1/2}\\
\leq \,&C2^{-3| \vec{ \pmb n}|/2}
\sum_{J\times K\in Ch_{ \vec{ \pmb n}}(Q(R))}
\frac{\mathbf 1_{J\times K}(\pmb y)}{|J\times K|^{1/2}}
=C2^{-| \vec{ \pmb n}|}\frac{\mathbf 1_{Q(R)}(\pmb y)}{|Q(R)|^{1/2}}.
\end{align*}
Since $\sum_{\eta, \vec{ \pmb n}}2^{-| \vec{ \pmb n}|}<\infty$,
this bound and the $L^2$ convergence in \eqref{H-expan2} give
$$
\lim_{N\to\infty}
\left\|\mathbb E_{2,N}^Q\mathbb E_{3,N}^Q\mathsf{a}_R^{(2)}(f)
-\mathsf{a}_R^{(2)}(f)\right\|_{
L^p  (\R^2;L^2  (\Gamma_\alpha(x_1),\frac{du\,ds}{s^2})    )}=0.
$$

Fix $\eta, \vec{ \pmb n}$ and define
$$
b_{m,J,K}(f)(u,s)
=\sum_{k\in\mathbb Z}
\sum_{\substack{R\in\mathscr F\cap\mathscr R_{m,k}\\
J\times K\in Ch_{ \vec{ \pmb n}}(Q(R))}}
\lambda_R^{(2)}(f)c_{R,J\times K}^{\eta}(f)(u,s).
$$
Then $b_{m,J,K}(f)=0$ outside $\Gamma_{\alpha,m}(x_1)$.
Whenever $b_{m,J,K}(f)\not\equiv0$, we also have
\begin{equation}\label{e:second-ratio}
\frac{|K|}{|J|}=2^{m+a_3-a_2+n_2-n_3}.
\end{equation}
Let
$$
X_{ \vec{ \pmb n},\eta}(f)(u,s;\pmb y)
=\sum_{m,J,K}b_{m,J,K}(f)(u,s)
h_J^{\eta_2}(y_2)h_K^{\eta_3}(y_3),
$$
and
$$
S_{ \vec{ \pmb n},\eta}(f)(\pmb y)^2
=\sum_{m,J,K}
 \bigg( \int_{\Gamma_{\alpha,m}(x_1)}
|b_{m,J,K}(f)(u,s)|^2\,\frac{du\,ds}{s^2}  \bigg)
\frac{\mathbf 1_{J\times K}(\pmb y)}{|J\times K|}.
$$
For fixed $J,K$ and $\vec{\pmb n}$, the covering intervals $Q_2(R)$ and $Q_3(R)$ are the unique
$n_2$th and $n_3$th ancestors of $J$ and $K$ in the fixed grids.
Thus all rectangles contributing to the same Haar coefficient have the same $Q(R)$, and \eqref{multi-23} gives
$$
\#\left\{R\in\mathscr F:J\times K\in Ch_{\vec{\pmb n}}(Q(R))\right\}
\le 2^{a_2+a_3}\le 2^9.
$$
Applying Cauchy--Schwarz to this finite sum, we obtain
\begin{align}\label{e:second-depth}
S_{ \vec{ \pmb n},\eta}(f)(\pmb y)^2
\leq &C\sum_{R\in\mathscr F}\lambda_R^{(2)}(f)^2
\sum_{J\times K\in Ch_{ \vec{ \pmb n}}(Q(R))}
 \bigg( \int_{\Gamma_\alpha(x_1)}
|c_{R,J\times K}^{\eta}(f)(u,s)|^2\,\frac{du\,ds}{s^2}  \bigg)
\frac{\mathbf 1_{J\times K}(\pmb y)}{|J\times K|}\notag\\
\leq &C2^{-3| \vec{ \pmb n}|}
\sum_{R\in\mathscr F}\lambda_R^{(2)}(f)^2
\sum_{J\times K\in Ch_{ \vec{ \pmb n}}(Q(R))}
\frac{\mathbf 1_{J\times K}(\pmb y)}{|J\times K|}
= C2^{-2| \vec{ \pmb n}|}\mathcal S_{\mathscr F}^{(2)}(f)(\pmb y)^2.
\end{align}

We now use \eqref{e:second-ratio} to apply \eqref{e:haar-hilbert} in just one coordinate.
For $\eta\in\left\{(1,0),(1,1)\right\}$, namely, when $\eta_2=1$,  write
$$
X_{ \vec{ \pmb n},\eta}(f)(u,s;\pmb y)
=\sum_Jh_J(y_2)
\sum_{m,K}b_{m,J,K}(f)(u,s)h_K^{\eta_3}(y_3).
$$
The sets $\Gamma_{\alpha,m}(x_1)$ are disjoint.
For fixed $m,J$, all relevant $K$ have the same length by
\eqref{e:second-ratio}, so distinct $K$ are disjoint.
Thus
\begin{align*}
&\sum_J\frac{\mathbf 1_J(y_2)}{|J|}
\int_{\Gamma_\alpha(x_1)}
\bigg|\sum_{m,K}b_{m,J,K}(f)(u,s)h_K^{\eta_3}(y_3)\bigg|^2
\,\frac{du\,ds}{s^2}\\
= &\sum_{J,m}\frac{\mathbf 1_J(y_2)}{|J|}
\int_{\Gamma_{\alpha,m}(x_1)}
\bigg|\sum_{|K|=2^{m+a_3-a_2+n_2-n_3}|J|}
b_{m,J,K}(f)(u,s)h_K^{\eta_3}(y_3)\bigg|^2
\,\frac{du\,ds}{s^2}\\
= &\sum_{m,J,K}
 \bigg( \int_{\Gamma_{\alpha,m}(x_1)}
|b_{m,J,K}(f)(u,s)|^2\,\frac{du\,ds}{s^2}  \bigg)
\frac{\mathbf 1_{J\times K}(\pmb y)}{|J\times K|}
= S_{ \vec{ \pmb n},\eta}(f)(\pmb y)^2.
\end{align*}
Applying \eqref{e:haar-hilbert} in $y_2$ and integrating the $p$th powers in $y_3$,
$$
\bigg  \|   \bigg( \int_{\Gamma_\alpha(x_1)}
|X_{ \vec{ \pmb n},\eta}(f)(u,s;\cdot)|^2
\,\frac{du\,ds}{s^2}  \bigg)^{1/2}\bigg    \|_{L^p(\R^2)}
\le C\sqrt p\, \left\|   S_{ \vec{ \pmb n},\eta}(f) \right\|_{L^p(\R^2)}.
$$

For $\eta=(0,1)$, write
$$
X_{ \vec{ \pmb n},(0,1)}(f)(u,s;\pmb y)
=\sum_Kh_K(y_3)
\sum_{m,J}b_{m,J,K}(f)(u,s)h_J^0(y_2).
$$
By \eqref{e:second-ratio}, for fixed $m,K$, the relevant $J$ have length
$$
|J|=2^{-m+a_2-a_3-n_2+n_3}|K|
$$
and are disjoint. The same calculation gives
\begin{align*}
\sum_K\frac{\mathbf 1_K(y_3)}{|K|}
\int_{\Gamma_\alpha(x_1)}
\bigg|\sum_{m,J}b_{m,J,K}(f)(u,s)h_J^0(y_2)\bigg|^2
\,\frac{du\,ds}{s^2}
= S_{ \vec{ \pmb n},(0,1)}(f)(\pmb y)^2.
\end{align*}
Applying \eqref{e:haar-hilbert} in $y_3$ and integrating the $p$th powers in $y_2$,
$$
\bigg  \|   \bigg( \int_{\Gamma_\alpha(x_1)}
|X_{ \vec{ \pmb n},(0,1)}(f)(u,s;\cdot)|^2
\,\frac{du\,ds}{s^2}  \bigg)^{1/2}\bigg    \|_{L^p(\R^2)}
\le C\sqrt p\, \left\|   S_{ \vec{ \pmb n},(0,1)}(f) \right\|_{L^p(\R^2)}.
$$
Hence this estimate holds for all $\eta\in \mathsf H$.
Combined with \eqref{e:second-depth},
$$
\bigg  \|   \bigg( \int_{\Gamma_\alpha(x_1)}
|X_{ \vec{ \pmb n},\eta}(f)(u,s;\cdot)|^2
\,\frac{du\,ds}{s^2}  \bigg)^{1/2}\bigg    \|_{L^p(\R^2)}
\le C\sqrt p\,2^{-| \vec{ \pmb n}|}
 \left\|   \mathcal S_{\mathscr F}^{(2)}(f) \right\|_{L^p(\R^2)}.
$$
For this finite $\mathscr F$, the Haar expansions imply
$$
\sum_{R\in\mathscr F}\widetilde{\mathsf{a}}_R^{(2)}(f)
=\sum_{R\in\mathscr F}\lambda_R^{(2)}(f)\mathsf{a}_R^{(2)}(f)
=\sum_{\eta, \vec{ \pmb n}}X_{ \vec{ \pmb n},\eta}(f).
$$
Since
$
\sum_{\eta, \vec{ \pmb n}}2^{-| \vec{ \pmb n}|}
=2\sum_{n=0}^{\infty}2^{-n}
+\big (\sum_{n=0}^{\infty}2^{-n}\big)^2\le C,
$
by Minkowski's inequality we have, for $p\geq 2$,
\begin{align*}
&\bigg  \|   \bigg( \int_{\Gamma_\alpha(x_1)}
\bigg|\sum_{R\in\mathscr F}\widetilde{\mathsf{a}}_R^{(2)}(f)(u,s;\cdot)\bigg|^2
\,\frac{du\,ds}{s^2}  \bigg)^{1/2}\bigg    \|_{L^p(\R^2)}\\
\leq &\sum_{\eta, \vec{ \pmb n}}
\bigg  \|   \bigg( \int_{\Gamma_\alpha(x_1)}
|X_{ \vec{ \pmb n},\eta}(f)(u,s;\cdot)|^2
\,\frac{du\,ds}{s^2}  \bigg)^{1/2}\bigg    \|_{L^p(\R^2)}
\leq  C\sqrt p\, \left\|   \mathcal S_{\mathscr F}^{(2)}(f) \right\|_{L^p(\R^2)}.
\end{align*}

For a general finite $\widetilde{\mathscr F}\subset\mathscr R$, split its atomic sum according to the classes
$\widetilde{\mathscr F}\cap\mathscr R_{\vec{\pmb\varepsilon},\vec{\pmb a}}$.
By the triangle inequality, the estimate just proved in each class, and \eqref{e:second-square}, we have
\begin{align}\label{e:second-finite}
&\bigg  \|   \bigg( \int_{\Gamma_\alpha(x_1)}
\bigg|\sum_{R\in\widetilde{\mathscr F}}\widetilde{\mathsf{a}}_R^{(2)}(f)(u,s;\cdot)\bigg|^2
\,\frac{du\,ds}{s^2}  \bigg)^{1/2}\bigg    \|_{L^p(\R^2)}\notag\\
\leq &C\sqrt p\sum_{ \vec{ \pmb\varepsilon}, \vec{ \pmb a}}
 \left\|   \mathcal S_{\widetilde{\mathscr F}\cap
\mathscr R_{ \vec{ \pmb\varepsilon}, \vec{ \pmb a}}}^{(2)}(f)
 \right\|_{L^p(\R^2)}
\leq C\sqrt p\,
 \left\|   S_{z,\alpha}^{\psi}f(x_1,\cdot) \right\|_{L^p(\R^2)}.
\end{align}

\noindent{\it Step 4.}  
 We now verify the $L^2$ reproducing formula \eqref{CRF2}, stated in {\it Step 2}, and then
 use this convergence and \eqref{e:second-finite} to estimate
$S_{1,\alpha}^{\psi_1}f(x_1,\cdot)$.

{Write $d\mu=du\,ds/s^2$ on $\Gamma=\Gamma_\alpha(x_1)$.
For $G\in L^2(\Gamma\times\R^2,d\mu\,d\pmb y)$, Plancherel in $\pmb y$ gives
$$
\int_\Gamma\int_0^\infty
\big\|G(u,s;\cdot)*\Phi_{s,t}^{(2)}\big\|_{L^2(\R^2)}^2
\,\frac{dt}{t}\,d\mu
=c_{\Phi^{(2)}}\|G\|_{L^2(\Gamma\times\R^2)}^2,
$$
where $c_{\Phi^{(2)}}=\int_0^\infty|\widehat{\Phi^{(2)}}(r,0)|^2\,dr/r$
by the radial substitution in \eqref{e:second-energy}.}
Choose finite families $\widetilde{\mathscr F}_\ell\subseteq\widetilde{\mathscr F}_{\ell+1}$
with $\bigcup_\ell\widetilde{\mathscr F}_\ell=\mathscr R$.
By \eqref{e:dual-second}, \eqref{e:second-energy} and the partition \eqref{W2R},
the duality argument in Lemma~\ref{l:atomic-l2} gives
\begin{align*}
&\bigg(\int_{\R^2}\int_{\Gamma_\alpha(x_1)}
\bigg|(f(\cdot,\pmb y)*\psi_{1,s})(u)
-\sum_{R\in\widetilde{\mathscr F}_\ell}\widetilde{\mathsf{a}}_R^{(2)}(f)(u,s;\pmb y)\bigg|^2
\,\frac{du\,ds}{s^2}\,d\pmb y\bigg)^{1/2}
\leq c_{\Phi^{(2)}}^{1/2}
\bigg(\sum_{R\notin\widetilde{\mathscr F}_\ell}\lambda_R^{(2)}(f)^2\bigg)^{1/2}\to0.
\end{align*}
This proves \eqref{CRF2}. We may pass to a subsequence whose $L^2(\Gamma_\alpha(x_1),du\,ds/s^2)$ norms
converge almost everywhere in $\pmb y$.
Fatou's lemma and \eqref{e:second-finite} therefore give
$$
\|S_{1,\alpha}^{\psi_1}f(x_1,\cdot)\|_{L^p(\R^2)}
\le C\sqrt p\,\|S_{z,\alpha}^{\psi}f(x_1,\cdot)\|_{L^p(\R^2)}.
$$
Integrating the $p$th powers in $x_1$ proves \eqref{e:reverse-two}
 for $f\in L^2(\R^3)\cap L^p(\R^3)$.

By density of $L^2\cap L^p$ in $L^p$ and the continuity of both square functions
in Lemma~\ref{l:continuous-forward}, \eqref{e:reverse-two} holds for every
$f\in L^p(\R^3)$, with the same constant $C_{\psi_2,\alpha}\sqrt p$.
\end{proof}

 \smallskip

\begin{proof}[Proof of the reverse estimate in Theorem~\ref{t:primal}]
Under \eqref{e:supp-origin}, Propositions~\ref{prop:square-1} and \ref{prop:square-2} show that
$$
\|f\|_{L^p(\R^3)}
\le C_{\psi_1,\alpha}\sqrt p\,\|S_{1,\alpha}^{\psi_1}f\|_{L^p(\R^3)}
\le C_{\psi,\alpha}p\,\|S_{z,\alpha}^{\psi}f\|_{L^p(\R^3)}.
$$
The sharpness of $p$ will be proved in Subsection~\ref{sec:sharp-p}.
\end{proof}

\subsection{A general lower bound for the reverse constant}\label{sec:sharp-p}

For $2\le p<\infty$ and $\alpha>0$, write
$$
\mathfrak C_p^\psi(\alpha)
=\sup_{0\ne f\in L^p(\R^3)}
\frac{\|f\|_{L^p(\R^3)}}
{\|S_{z,\alpha}^{\psi}f\|_{L^p(\R^3)}}.
$$
In this subsection, we assume only the basic conditions on $\psi_1,\psi_2$,
without requiring either support condition in \eqref{e:supp-origin}.
For every fixed kernel pair and every $\alpha>0$, we prove
$$
c_{\psi,\alpha}p\le\mathfrak C_p^\psi(\alpha),
\qquad 2\le p<\infty,
$$
where $c_{\psi,\alpha}>0$ is independent of $p$.
In particular, this lower bound also holds under \eqref{e:supp-origin}.
This, combined with the reverse estimate proved in the preceding subsection,
shows that the factor $p$ in Theorem~\ref{t:primal} is optimal.

We combine the one-dimensional case of Lemma~\ref{sp:l:lower}
with a limit obtained by dilating the third variable.
For a real-valued  and even kernel $ \kappa \in C_c^\infty(\R)$ with $\int_\R \kappa =0$, define
$$
\mathcal A_{ \kappa ,\alpha}g(x_1)^2
=\int_0^\infty\int_{|u-x_1|<\alpha s}
|(g* \kappa _s)(u)|^2\,\frac{du\,ds}{s^2},
 \qquad \text{where} \qquad 
 \kappa _s(x_1)=\frac1s \kappa \left(\frac{x_1}{s}\right).
$$

For the transfer argument below, differentiate $j$ times in $r$
and integrate by parts $M$ times. For $r\ne0$, this gives
$$
\widehat \kappa ^{(j)}(r)
=(-2\pi i)^j\int_\R z^j \kappa (z)e^{-2\pi irz}\,dz
=\frac{(-2\pi i)^j}{(2\pi ir)^M}
\int_\R\frac{d^M}{dz^M}\big (z^j \kappa (z)\big  )e^{-2\pi irz}\,dz.
$$
Using the first expression for $|r|\le1$ and the second for $|r|>1$,
\begin{equation}\label{e:kappa-fourier-decay}
|\widehat \kappa ^{(j)}(r)|
\le C_{ \kappa ,j,M}(1+|r|)^{-M},
 \qquad  r\in\R,\quad j\in\mathbb N_0,\quad M\ge1.
\end{equation}

\medskip

The next lemma gives a product limit for the Zygmund area function.

\begin{lemma}\label{l:lower-transfer}
Fix $\alpha>0$ and $2\le p<\infty$, and let
\begin{equation}\label{kappa}
 \kappa (x_2)=\int_\R\psi_2(x_2,x_3)\,dx_3.
 \end{equation}
 Let $g_1,g_2,\chi\in\mathcal S(\R)$ be nonzero, with
$
\operatorname{supp}\widehat g_1
\cup\operatorname{supp}\widehat g_2
\cup\operatorname{supp}\widehat\chi
\Subset\R\setminus\{0\}.
$

For $L\ge1$, define
$
f_L(x)=g_1(x_1)g_2(x_2)\chi(x_3/L).
$
Then $f_L\in\mathcal S(\R^3)$, and
\begin{equation}\label{e:lower-transfer}
L^{-1/p} \left\|   S_{z,\alpha}^{\psi}f_L \right\|_{L^p(\R^3)}
\to
(2\alpha)^{1/2}  \|  \chi    \|_{L^p(\R)}
 \left\|   \mathcal A_{\psi_1,\alpha}g_1 \right\|_{L^p(\R)}
 \left\|   \mathcal A_{ \kappa ,\alpha}g_2 \right\|_{L^p(\R)}
 \qquad (L\to\infty).
\end{equation}
\end{lemma}

\begin{proof} 
Our idea is to take $L$ large so that $\chi(x_3/L)$ changes slowly in $x_3$.
After rescaling, the area function tends to $\sqrt{2\alpha}|\chi|$
times the product of the two one-dimensional area functions.

\smallskip
\noindent{\it Step 1.}
We first show that, for fixed $s,t>0$,
$$
(f_L*\psi_{s,t})(x_1,x_2,Lx_3)
\to 
(g_1*\psi_{1,s})(x_1)(g_2* \kappa _t)(x_2)\chi(x_3)
\quad\text{in }L^p(\R^3)
 \qquad (L\to\infty).
$$

Let $\displaystyle q=\frac{st}{L}$ and $\pmb z=(z_2,z_3)$. Then
\begin{align*}
\mathcal T_{s,t,L}(x):=&(f_L*\psi_{s,t})(x_1,x_2,Lx_3)
= (g_1*\psi_{1,s})(x_1)
\int_{\R^2}g_2(x_2-tz_2)\chi(x_3-qz_3)
\psi_2(\pmb z)\,d\pmb z.
\end{align*}
Since $q\to0$ as $L\to\infty$, the formal limit of
$\mathcal T_{s,t,L}$ is
\begin{align*}
\mathcal T_{s,t,\infty}(x)&=(g_1*\psi_{1,s})(x_1)
\int_{\R^2}g_2(x_2-tz_2)
\psi_2(\pmb z)\,d\pmb z \,  \chi(x_3 )\\
&=  (g_1*\psi_{1,s})(x_1)
\int_\R g_2(x_2-tz_2) \kappa (z_2)\,dz_2\,     \chi(x_3 )\\
&=(g_1*\psi_{1,s})(x_1)\,   (g_2* \kappa _t)(x_2)\,  \chi(x_3),
 \end{align*}
 where $ \kappa $ is given in \eqref{kappa} and $ \kappa _t(\cdot)=t^{-1} \kappa (\cdot/t)$.
 
To justify this limit, apply Minkowski's inequality:
\begin{align*}
&  \|  \mathcal T_{s,t,L}-\mathcal T_{s,t,\infty}    \|_{L^p(\R^3)}
\leq      \|  g_1*\psi_{1,s}    \|_{L^p(\R)}  \|  g_2    \|_{L^p(\R)}
\int_{\R^2}|\psi_2(\pmb z)|
  \|  \chi(\cdot-qz_3)-\chi    \|_{L^p(\R)}\,d\pmb z.
\end{align*}
For every fixed $z_3$, as $q\to 0$ (i.e., $L\to \infty$),
$$
  \|  \chi(\cdot-qz_3)-\chi    \|_{L^p(\R)}\to0,
 \qquad 
  \|  \chi(\cdot-qz_3)-\chi    \|_{L^p(\R)}
\le2  \|  \chi    \|_{L^p(\R)}.
$$
Hence, it follows from the dominated convergence theorem that
\begin{equation}\label{e:trans-conv}
  \|  \mathcal T_{s,t,L}-\mathcal T_{s,t,\infty}    \|_{L^p(\R^3)}
\to0
 \qquad \text{for every fixed }s,t>0.
\end{equation}

\smallskip

\noindent{\it Step 2.}
The Fourier supports of the fixed functions $g_1,g_2,\chi$
stay a positive distance from zero.
We use this condition and the rapid decay of
$\widehat\psi_1,\widehat\psi_2$ to control large scales.
Together with $\widehat\psi_1(0)=\widehat\psi_2(0,0)=0$,
we will prove
$$
\|\mathcal T_{s,t,L}\|_{L^p(\R^3)}
\le C_p\min\{s,s^{-1}\}\min\{t,t^{-1}\},
\qquad s,t>0,\quad L\ge1.
$$
Note that
$
\widehat{\mathcal T_{s,t,L}}(\xi)
=\widehat g_1(\xi_1)\widehat g_2(\xi_2)\widehat\chi(\xi_3)
\widehat\psi_1(s\xi_1)\widehat\psi_2(t\xi_2,q\xi_3),
$
and we can choose $0<c<C$, independent of $s,t,L$, such that
$$
\operatorname{supp}\widehat{\mathcal T_{s,t,L}}
\subset
\operatorname{supp}\widehat g_1
\times\operatorname{supp}\widehat g_2
\times\operatorname{supp}\widehat\chi   \subset
\left\{\xi\in\R^3:c\le|\xi_j|\le C,\ j=1,2,3\right\}.
$$
Consequently, 
\begin{align}\label{e:trans-ibp}
(1+4\pi^2|x|^2)^2|\mathcal T_{s,t,L}(x)|
= &\bigg|\int_{\R^3}(1-\Delta_\xi)^2
\widehat{\mathcal T_{s,t,L}}(\xi)e^{2\pi ix\cdot\xi}\,d\xi\bigg|
\leq C\max_{|\beta|\le4}
 \left\|   \partial^\beta\widehat{\mathcal T_{s,t,L}} \right\|_{L^\infty}\notag\\
\leq &C\max_{\substack{i,j,k\ge0\\i+j+k\le4}}
\sup_{\substack{c\le|\xi_\ell|\le C\\\ell=1,2,3}}
s^i\left|\widehat\psi_1^{(i)}(s\xi_1)\right|
t^jq^k
\left|(\partial_2^j\partial_3^k
\widehat\psi_2)(t\xi_2,q\xi_3)\right|,
\end{align}
where $ 
 \displaystyle      
\widehat\psi_1^{(i)}(r)
:=\frac{d^i}{dr^i}\widehat\psi_1(r)$.

By $\widehat\psi_1(0)=0$ and  \eqref{e:kappa-fourier-decay}    with $M=8$,
\begin{align}\label{e:trans-ibp-I}
s^i|\widehat\psi_1^{(i)}(s\xi_1)|
\leq\, &C
\begin{cases}
s,&0<s\le1,\ i=0,\\[3pt]
s^i,&0<s\le1,\ 1\le i\le4,\\[3pt]
s^{i-8},&s>1,\ 0\le i\le4,
\end{cases}  \notag\\
\leq \, &C\min\{s,s^{-3}\}.
\end{align}

Similarly, since $\widehat\psi_2(0,0)=0$, for $j,k\in\mathbb N_0$
with $j+k\le4$,
\begin{align*}
 \left|\partial_{\xi_2}^j\partial_{\xi_3}^k
\widehat\psi_2(t\xi_2,q\xi_3)\right| 
= &t^jq^k
\left|(\partial_2^j\partial_3^k\widehat\psi_2)(t\xi_2,q\xi_3)\right|\\
\leq &C
\begin{cases}
t+q,&t+q\le1,\ j+k=0,\\[2pt]
t^jq^k,&t+q\le1,\ 1\le j+k\le4,\\[2pt]
t^jq^k(1+\sqrt{t^2\xi_2^2+q^2\xi_3^2})^{-8},
&t+q>1.
\end{cases}
\end{align*}
Noting
$$
t^jq^k\le(t+q)^{j+k},
 \qquad 
\sqrt{t^2\xi_2^2+q^2\xi_3^2}\ge\frac c{\sqrt2}(t+q),
$$
we obtain
\begin{align}\label{e:trans-ibp-II}
t^jq^k
\left|(\partial_2^j\partial_3^k\widehat\psi_2)(t\xi_2,q\xi_3)\right|
\leq &C
\begin{cases}
t+q,&t+q\le1,\\[2pt]
(t+q)^{-4},&t+q>1,
\end{cases}  \notag\\
\leq &C(1+t+q)^{-3}\min\{1,t+q\}   \notag\\
\leq & C(1+t)^{-3}  (1+s)   \min\{1,t\},
\end{align}
where the last inequality follows from $\min\{1,t+q\}\le(1+s)\min\{1,t\}$, since $q=st/L$ and $L\ge1$.
 
Combining \eqref{e:trans-ibp-I} and \eqref{e:trans-ibp-II}, we have 
\begin{align*}
s^i|\widehat\psi_1^{(i)}(s\xi_1)|
t^jq^k
\left|(\partial_2^j\partial_3^k\widehat\psi_2)(t\xi_2,q\xi_3)\right|
&\leq C(1+s)\min\{s,s^{-3}\}(1+t)^{-3}\min\{1,t\}\\
&\leq C\min\{s,s^{-1}\}\min\{t,t^{-1}\}.
\end{align*}
Substituting this bound into \eqref{e:trans-ibp}, we obtain
$$
|\mathcal T_{s,t,L}(x)|
\le C(1+|x|)^{-4}\min\{s,s^{-1}\}\min\{t,t^{-1}\},
$$
consequently,
\begin{align*}
  \|  \mathcal T_{s,t,L}    \|_{L^p(\R^3)}
\leq &C\min\{s,s^{-1}\}\min\{t,t^{-1}\}
\left(\int_{\R^3}(1+|x|)^{-4p}\,dx\right)^{1/p}
\leq C_p\min\{s,s^{-1}\}\min\{t,t^{-1}\},
\end{align*}
which proves the claimed estimate.
Passing to the limit in \eqref{e:trans-conv}, we also have
$$
  \|  \mathcal T_{s,t,\infty}    \|_{L^p(\R^3)}
=\lim_{L\to\infty}  \|  \mathcal T_{s,t,L}    \|_{L^p(\R^3)}
\le C_p\min\{s,s^{-1}\}\min\{t,t^{-1}\}.
$$
Therefore,
\begin{equation}\label{e:trans-bd}
  \|  \mathcal T_{s,t,L}    \|_{L^p(\R^3)}
+  \|  \mathcal T_{s,t,\infty}    \|_{L^p(\R^3)}
\le C_p\min\{s,s^{-1}\}\min\{t,t^{-1}\}.
\end{equation}

\smallskip

\noindent{\it Step 3.}
We now pass from the fixed-scale convergence to the area functions by expressing the cones
in the fixed parameter region $(-\alpha,\alpha)^3$.
By definition,
\begin{align*}
 S_{z,\alpha}^{\psi}f_L(x_1,x_2,Lx_3)^2 
= &\int_0^\infty\int_0^\infty
\int_{\substack{|u-x_1|<\alpha s\\
|v_2-x_2|<\alpha t\\
|v_3-Lx_3|<\alpha st}}
|(f_L*\psi_{s,t})(u,\pmb v)|^2
\,\frac{du\,d\pmb v\,ds\,dt}{s^3t^3}.
\end{align*}
For $(a_1,a_2,a_3)\in(-\alpha,\alpha)^3$,  let
$
u=x_1+sa_1, \ 
v_2=x_2+ta_2, \ 
v_3=Lx_3+sta_3.
$
Then
$$
du\,d\pmb v=s^2t^2\,da_1\,da_2\,da_3,
 \qquad 
\frac{v_3}{L}=x_3+qa_3,
$$
and hence
\begin{equation}\label{e:trans-cone}
S_{z,\alpha}^{\psi}f_L(x_1,x_2,Lx_3)^2
=\int_0^\infty\int_0^\infty\int_{(-\alpha,\alpha)^3}
\left|\mathcal T_{s,t,L}
\big (x+(sa_1,ta_2,qa_3)\big  )\right|^2
\,da_1\,da_2\,da_3\,\frac{ds}{s}\frac{dt}{t}.
\end{equation}

Note that the limit  $\mathcal T_{s,t,\infty}$ in  {\it Step 1} satisfies
$$
\mathcal T_{s,t,\infty}\big (x+(sa_1,ta_2,0)\big  )
=(g_1*\psi_{1,s})(x_1+sa_1)(g_2* \kappa _t)(x_2+ta_2)\chi(x_3),
$$
which is independent of $a_3$. Integrating in $a_3$ first, we obtain
\begin{align*}
&\int_0^\infty\int_0^\infty\int_{(-\alpha,\alpha)^3}
\left|\mathcal T_{s,t,\infty}
\big (x+(sa_1,ta_2,0)\big  )\right|^2
\,da_1\,da_2\,da_3\,\frac{ds}{s}\frac{dt}{t}\\
= &2\alpha|\chi(x_3)|^2
\mathcal A_{\psi_1,\alpha}g_1(x_1)^2
\mathcal A_{ \kappa ,\alpha}g_2(x_2)^2.
\end{align*}

{For $a=(a_1,a_2,a_3)\in(-\alpha,\alpha)^3$, write
$$
 E_L(s,t,a)=
 \big\|\mathcal T_{s,t,L}\big(\cdot+(sa_1,ta_2,qa_3)\big)
 -\mathcal T_{s,t,\infty}\big(\cdot+(sa_1,ta_2,0)\big)\big\|_p.
$$
For fixed $s,t,a$, we have $q=st/L\to0$. By \eqref{e:trans-conv}
and continuity of translations in $L^p$,
\begin{equation}\label{e:trans-shift}
 E_L(s,t,a)\le
 \|\mathcal T_{s,t,L}-\mathcal T_{s,t,\infty}\|_p
 +\|\mathcal T_{s,t,\infty}(\cdot+(0,0,qa_3))-\mathcal T_{s,t,\infty}\|_p
 \to 0.
\end{equation}
By \eqref{e:trans-bd}, $E_L(s,t,a)\le C_p\min\{s,s^{-1}\}\min\{t,t^{-1}\}$
uniformly in $L$ and $a$. Since $p\ge2$, the triangle inequality in the cone
integral and Minkowski's inequality give
\begin{align*}
 &\big\|S_{z,\alpha}^{\psi}f_L(x_1,x_2,Lx_3)
 -\sqrt{2\alpha}\,|\chi(x_3)|\,
 \mathcal A_{\psi_1,\alpha}g_1(x_1)\mathcal A_{\kappa,\alpha}g_2(x_2)
 \big\|_{L^p(dx)}\\
 &\qquad\le
 \bigg(\int_0^\infty\int_0^\infty\int_{(-\alpha,\alpha)^3}
 E_L(s,t,a)^2\,da\,\frac{ds}{s}\frac{dt}{t}\bigg)^{1/2}
 \to0,
\end{align*}
where \eqref{e:trans-shift} and dominated convergence give the limit, since
$\int_0^\infty\min\{r,r^{-1}\}^2\,dr/r=1$.
Changing variables in $x_3$ on the left and taking the norm of the product
gives \eqref{e:lower-transfer}.}
\end{proof}

Combining the case $d=1$ of Lemma~\ref{sp:l:lower} with Lemma~\ref{l:lower-transfer},
we obtain the following estimate
without requiring the support conditions  $0\in\operatorname{supp}(|D|\psi_1)$ and $
0\in\operatorname{supp}\psi_2$   in \eqref{e:supp-origin}.

\begin{proposition}\label{p:lower-bound}
Under the basic assumptions on $\psi_1$ and $\psi_2$,
for every $\alpha>0$ and $2\le p<\infty$,
$$
c_{\psi,\alpha}p\le\mathfrak C_p^\psi(\alpha),
$$
where $c_{\psi,\alpha}>0$ is independent of $p$.
\end{proposition}

\begin{proof}
Note that
the kernel $ \kappa $ in \eqref{kappa} is real and even, with
$$
 \kappa \in C_c^\infty(\R), \qquad 
\int_\R \kappa =\int_{\R^2}\psi_2=0, \qquad 
\widehat \kappa (r)=\widehat\psi_2(r,0).
$$
If $ \kappa =0$, radiality would force $\widehat\psi_2=0$
on $\R^2$, contrary to $\psi_2\ne0$.
Thus $ \kappa $ satisfies the assumptions of Lemma~\ref{sp:l:lower} with $d=1$.

Choose $g_1,g_2$ by Lemma~\ref{sp:l:lower} with $d=1$
for the kernels $\psi_1, \kappa $, respectively,
and take $\chi,f_L$ as in Lemma~\ref{l:lower-transfer}.
Since
$$
  \|  f_L    \|_{L^p(\R^3)}
=L^{1/p}  \|  g_1    \|_{L^p(\R)}
  \|  g_2    \|_{L^p(\R)}  \|  \chi    \|_{L^p(\R)},
$$
it follows from \eqref{e:lower-transfer} that
\begin{align*}
\mathfrak C_p^\psi(\alpha)
\geq &\lim_{L\to\infty}
\frac{  \|  f_L    \|_{L^p(\R^3)}}
{  \|  S_{z,\alpha}^{\psi}  f_L    \|_{  L^p(\R^3) }   }  
=  \frac {1} {\sqrt{2\alpha}}
\frac{  \|  g_1    \|_{L^p(\R)}}
{  \|  \mathcal A_{\psi_1,\alpha}g_1    \|_{L^p(\R)}}
\frac{  \|  g_2    \|_{L^p(\R)}}
{  \|  \mathcal A_{ \kappa ,\alpha}g_2    \|_{L^p(\R)}} 
\geq  c_{\psi,\alpha}p.
\end{align*}
The proof is complete.
\end{proof}

In particular, the lower bound for $\mathfrak C_p^\psi(\alpha)$
also holds under \eqref{e:supp-origin}.
Hence one may combine Theorem~\ref{t:primal}
and Proposition~\ref{p:lower-bound} to obtain
$$
c_{\psi,\alpha}p
\le\mathfrak C_p^\psi(\alpha)
\le C_{\psi,\alpha}p.
$$
That is, the factor $p$ in Theorem~\ref{t:primal} is optimal.

\section{Support distances and removal of the first condition}\label{sec:support-distance}

\subsection{Support distances and large apertures}\label{sub:large-ap}

Under the conditions in \eqref{e:supp-origin}, the reverse estimate in
Theorem~\ref{t:primal} holds for every $\alpha>0$.
In this section, we prove Theorem~\ref{t:large-ap}: for each fixed kernel pair,
the same estimate holds without requiring \eqref{e:supp-origin}
when $\alpha$ is sufficiently large.

 To this end, we will prove the two $\sqrt p$ estimates under
$d_1(\psi_1)<c_*\alpha$ and $d_2(\psi_2)<c_*\alpha$, respectively;
see Proposition~\ref{p:support-distance}.
These conditions compare $\alpha$ with the support distances
$d_1,d_2$ defined in the Introduction. Recall that
$$
d_1(\psi_1)=\operatorname{dist}\big (0,\operatorname{supp}(|D|\psi_1)\big  ),
\qquad
d_2(\psi_2)=\operatorname{dist}\big (0,\operatorname{supp}\psi_2\big  ).
$$
Thus
$$
0\in\operatorname{supp}(|D|\psi_1)
\quad\text{if and only if}\quad d_1(\psi_1)=0,
$$
and
$$
0\in\operatorname{supp}\psi_2
\quad\text{if and only if}\quad d_2(\psi_2)=0.
$$

Recall the one-dimensional kernel
$\kappa(u)=\int_\R\psi_2(u,v)\,dv$ in \eqref{kappa}.
{Proposition~\ref{sp:p:marginal} relates $|D|\kappa$ to $\psi_2$ in dimension two.}

\begin{lemma}\label{l:marginal-supp}
For every $\varepsilon>0$,
$$
\psi_2=0\text{ on }B(0,\varepsilon)
\qquad\text{if and only if}   \qquad 
|D|\kappa=0\text{ on }(-\varepsilon,\varepsilon).
$$
In particular,
\begin{equation}\label{e:marginal-distance}
d_2(\psi_2)=d_1(\kappa).
\end{equation}
\end{lemma}

\begin{proof}
The kernel $\psi_2$ satisfies the basic assumptions in dimension two,
and $\kappa$ is its integral in the last coordinate.
Equation~\eqref{sp:e:spherical-fourier} becomes
$$
\psi_2(\pmb x)
=\frac12\int_0^{2\pi}
(|D|\kappa)(x_2\cos\theta+x_3\sin\theta)\,d\theta.
$$
Since $|D|^{2-2}\psi_2=\psi_2$, Proposition~\ref{sp:p:marginal}
gives the stated equivalence and \eqref{e:marginal-distance}.
\end{proof}

Remark~\ref{rem:support-why} explains the construction of normalized dual
kernels with arbitrarily small support under \eqref{e:supp-origin}.
For a fixed support radius, the necessary and sufficient conditions can be expressed in terms of $d_1,d_2$.
The following lemma specializes Proposition~\ref{sp:p:local-dual} to $d=1$ and $d=2$,
with the normalization written explicitly for each kernel.

\begin{lemma}\label{l:loc-pairing}
Fix $\varepsilon>0$.
\begin{itemize}
\item[(i)] There is a real even $\Phi^{(1)}\in C_c^\infty((-\varepsilon,\varepsilon))$ such that
$$
\int_\R\Phi^{(1)}=0,
\qquad
\int_0^\infty\widehat\psi_1(r)\widehat{\Phi^{(1)}}(r)\,\frac{dr}{r}=1
$$
if and only if $d_1(\psi_1)<\varepsilon$.
\smallskip
\item[(ii)] There is a real radial $\Phi^{(2)}\in C_c^\infty(B(0,\varepsilon))$ such that
$$
\int_{\R^2}\Phi^{(2)}=0,
\qquad
\int_0^\infty\widehat\psi_2(r,0)\widehat{\Phi^{(2)}}(r,0)\,\frac{dr}{r}=1
$$
if and only if $d_2(\psi_2)<\varepsilon$.
\end{itemize}
\end{lemma}

\smallskip

\begin{remark}
For a fixed radius $\rho$, Lemma~\ref{l:loc-pairing}
requires only $d_j(\psi_j)<\rho$ for the corresponding normalized dual kernel.
The pointwise comparisons \eqref{sp:e:discrete-area}, with $d=1$,
and \eqref{e:second-square} require $16\rho<\alpha$ and $54\rho<\alpha$,
respectively.
Thus \eqref{e:supp-origin} allows every $\alpha>0$.
More generally, compared with Propositions~\ref{prop:square-1}
and \ref{prop:square-2}, Proposition~\ref{p:support-distance}
below establishes the two $\sqrt p$ estimates under the weaker conditions
$$
d_1(\psi_1)<c_*\alpha,
\qquad
d_2(\psi_2)<c_*\alpha,
$$
respectively, for an absolute constant $c_*>0$.
For fixed $\psi_1,\psi_2$, these conditions hold for sufficiently large $\alpha$.

Examples~\ref{r:counterexample-one} and
\ref{r:counterexample-two} show that, if the respective support conditions
in \eqref{e:supp-origin} are removed, the $\sqrt p$ estimates may fail
at a fixed $\alpha$; they do not contradict the estimates
for sufficiently large $\alpha$ with $\psi_1,\psi_2$ fixed.
\end{remark}

\smallskip

 \begin{proof}[Proof of Lemma~\ref{l:loc-pairing}]
Apply Proposition~\ref{sp:p:local-dual} with
$(d,\varphi)=(1,\psi_1)$ and $(d,\varphi)=(2,\psi_2)$, respectively.
The two constants in the pairing formula \eqref{sp:e:pair} are
$4\pi^2/\sigma_0=2\pi^2$ and $4\pi^2/\sigma_1=2\pi$.
Thus, for real even $\chi_1$ and real radial $\chi_2$ supported in the
respective balls of radius $\varepsilon$,
$$
\int_0^\infty\widehat\psi_1(r)\widehat{-\chi_1''}(r)\,\frac{dr}{r}
=2\pi^2\int_\R(|D|\psi_1)\chi_1,
$$
and
$$
\int_0^\infty\widehat\psi_2(r,0)\widehat{-\Delta\chi_2}(r,0)\,\frac{dr}{r}
=2\pi\int_{\R^2}\psi_2\chi_2.
$$
{Under the respective distance conditions, choose nonnegative smooth cutoffs
$\zeta_1,\zeta_2$ supported in these balls, even on $\R$ and radial on $\R^2$,
so that the denominators below are positive.}
The normalized dual kernels are
$$
\Phi^{(1)}
=-\frac{(\zeta_1|D|\psi_1)''}
{2\pi^2\int_\R\zeta_1(|D|\psi_1)^2},
\qquad
\Phi^{(2)}
=-\frac{\Delta(\zeta_2\psi_2)}
{2\pi\int_{\R^2}\zeta_2\psi_2^2}.
$$
The necessity of the distance conditions is the converse assertion in
Proposition~\ref{sp:p:local-dual}.
\end{proof}

{We now use these support distances to bound the two reverse constants at a fixed aperture.}

\begin{proposition}\label{p:support-distance}
There is an absolute constant $c_*>0$ such that, for every $\alpha>0$,
$2\le p<\infty$ and $f\in L^p(\R^3)$,
\begin{itemize}
\item[(i)] if $d_1(\psi_1)<c_*\alpha$, then
$$
\|f\|_{L^p(\R^3)}
\le C_{\psi_1,\alpha}\sqrt p\,
\|S_{1,\alpha}^{\psi_1}f\|_{L^p(\R^3)};
$$
\item[(ii)] if $d_2(\psi_2)<c_*\alpha$, then
$$
\|S_{1,\alpha}^{\psi_1}f\|_{L^p(\R^3)}
\le C_{\psi_2,\alpha}\sqrt p\,
\|S_{z,\alpha}^{\psi}f\|_{L^p(\R^3)}.
$$
\end{itemize}
The constants are independent of $p$ and $f$.
\end{proposition}

\begin{proof}
Fix $0<c_*<\min\{c_1,1/108\}$, where $c_1$ is the constant in Theorem~\ref{sp:t:general} for $d=1$,
and choose $\rho\in2^{\mathbb Z}$ such that
$$
c_*\alpha\le\rho<2c_*\alpha.
$$
First let $f\in L^2(\R^3)\cap L^p(\R^3)$.

\smallskip
\noindent(i).
Under $d_1(\psi_1)<c_*\alpha<c_1\alpha$, Theorem~\ref{sp:t:general}
with $d=1$ gives
$$
\mathfrak c_{p,1}^{\psi_1}(\alpha)\le C_{\psi_1,\alpha}\sqrt p.
$$
Substituting this bound into \eqref{e:first-from-single} and taking
the $p$th root proves (i).

\smallskip
\noindent(ii).
Similarly, under $d_2(\psi_2)<c_*\alpha\le\rho$,
(ii) of Lemma~\ref{l:loc-pairing} provides a real radial
$\Phi^{(2)}\in C_c^\infty(B(0,\rho))$
with zero integral and the required normalization.
We construct the atoms as in Proposition~\ref{prop:square-2}.
Since $54\rho<\alpha$, the calculation in \eqref{e:second-square}
bounds their discrete square function by
$C_{\psi_2,\alpha}S_{z,\alpha}^{\psi}f$.
The Haar estimate and reconstruction in that proof establish (ii).

Finally, by the approximation argument using Lemma~\ref{l:continuous-forward}, 
both estimates hold for every $f\in L^p(\R^3)$.
\end{proof}

\begin{proof}[Proof of Theorem \ref{t:large-ap}]
Since $|D|\psi_1\ne0$ and $\psi_2\ne0$, both support distances are finite.
Choose
$$
\alpha_\psi>c_*^{-1}\max\{d_1(\psi_1),d_2(\psi_2)\}.
$$
For $\alpha\ge\alpha_\psi$, Proposition~\ref{p:support-distance} yields
$$
\|f\|_{L^p(\R^3)}
\le C_{\psi_1,\alpha}\sqrt p\,
\|S_{1,\alpha}^{\psi_1}f\|_{L^p(\R^3)}
\le C_{\psi,\alpha}p\,
\|S_{z,\alpha}^{\psi}f\|_{L^p(\R^3)}.
$$
Combining this with Proposition~\ref{p:lower-bound} proves the two-sided estimate.
\end{proof}

\subsection{Removing the first distance condition}\label{sub:first-free}

We now fix any $\alpha>0$.
In Theorem~\ref{t:first-free}, we remove the first distance condition $d_1(\psi_1)<c_*\alpha$
from Proposition~\ref{p:support-distance},
while $d_2(\psi_2)<c_*\alpha$ is still assumed.
The factor $p^{3/2}$ in Theorem~\ref{t:first-free} is optimal over the corresponding kernel class;
we will prove this in Section~\ref{sec:sharp-aperture}.

The first-variable estimate follows directly from Theorem~\ref{sp:t:general}
with $d=1$, by the slice calculation \eqref{e:first-from-single}.
We also state the corresponding aperture comparison for $S_{1,\alpha}^{\psi_1}$;
its sharpness is proved in Remark~\ref{r:aperture-sharp}.

\begin{lemma}\label{l:aperture-one}
For $0<\alpha\le\beta$, $2\le p<\infty$ and $f\in L^p(\R^3)$,
\begin{equation}\label{e:aperture-one}
\left\|S_{1,\beta}^{\psi_1}f\right\|_{L^p(\R^3)}
\le C\left(\frac\beta\alpha\right)^{1/2}\sqrt p\,
\left\|S_{1,\alpha}^{\psi_1}f\right\|_{L^p(\R^3)},
\end{equation}
where the constant $C$ is absolute.
\end{lemma}

\begin{proof}
For almost every $\pmb x=(x_2,x_3)$, apply Lemma~\ref{sp:l:aperture}
with $d=1$ to $f(\cdot,\pmb x)$.
Taking $p$th powers and integrating in $\pmb x$, we obtain
\begin{align*}
\|S_{1,\beta}^{\psi_1}f\|_{L^p(\R^3)}^p
&=\int_{\R^2}
\|\mathcal A_{\psi_1,\beta}[f(\cdot,\pmb x)]\|_{L^p(\R)}^p\,d\pmb x\\
&\le \left[C(\beta/\alpha)^{1/2}\sqrt p\right]^p
\int_{\R^2}
\|\mathcal A_{\psi_1,\alpha}[f(\cdot,\pmb x)]\|_{L^p(\R)}^p\,d\pmb x\\
&=\left[C(\beta/\alpha)^{1/2}\sqrt p\right]^p
\|S_{1,\alpha}^{\psi_1}f\|_{L^p(\R^3)}^p.
\end{align*}
Taking the $p$th root proves \eqref{e:aperture-one}.
\end{proof}

We now apply the general one-parameter bound at the fixed aperture $\alpha$.

\begin{proposition}\label{p:first-general}
For every $\alpha>0$, $2\le p<\infty$ and $f\in L^p(\R^3)$,
\begin{equation}\label{e:first-general}
\left\|f\right\|_{L^p(\R^3)}
\le C_{\psi_1,\alpha}p\,
\left\|S_{1,\alpha}^{\psi_1}f\right\|_{L^p(\R^3)}.
\end{equation}
The constant is independent of $p$ and $f$.
No support condition on $|D|\psi_1$ is required.
\end{proposition}

 The factor $p$ in \eqref{e:first-general} is optimal over the full kernel class,
as shown by Example~\ref{r:counterexample-one}.

\begin{proof}[Proof of Proposition \ref{p:first-general}]
Theorem~\ref{sp:t:general} with $d=1$ gives
$\mathfrak c_{p,1}^{\psi_1}(\alpha)\le C_{\psi_1,\alpha}p$.
Taking the $p$th root in \eqref{e:first-from-single}, we obtain
$$
\|f\|_{L^p(\R^3)}
\le\mathfrak c_{p,1}^{\psi_1}(\alpha)\|S_{1,\alpha}^{\psi_1}f\|_{L^p(\R^3)}
\le C_{\psi_1,\alpha}p\,\|S_{1,\alpha}^{\psi_1}f\|_{L^p(\R^3)}.
$$
The proof is complete.
\end{proof}

\begin{proof}[Proof of Theorem~\ref{t:first-free}]
By Proposition~\ref{p:first-general} and (ii) of Proposition~\ref{p:support-distance},
\begin{align*}
\left\|f\right\|_{L^p(\R^3)}
\leq &C_{\psi_1,\alpha}p\,
\left\|S_{1,\alpha}^{\psi_1}f\right\|_{L^p(\R^3)}
\leq  C_{\psi,\alpha}p^{3/2}\,
\left\|S_{z,\alpha}^{\psi}f\right\|_{L^p(\R^3)}.
\end{align*}
The proof is complete.
\end{proof}

\section{Removing the second distance condition}\label{sec:second-free}

We remove the second distance condition $d_2(\psi_2)<c_*\alpha$
from Proposition~\ref{p:support-distance}.
We obtain a $p^{3/2}$ estimate when $d_1(\psi_1)<c_*\alpha$ is still assumed,
and a $p^2$ estimate when this condition is removed as well.
Both orders are optimal over the corresponding kernel classes;
see Section~\ref{sec:sharp-aperture}.

To this end, we first prove the following estimate without either distance condition.
Combining it with (i) of Proposition~\ref{p:support-distance} proves
(i) of Theorem~\ref{t:second-free}.
Combining it with Proposition~\ref{p:first-general} proves
(ii) of Theorem~\ref{t:second-free}.

\begin{proposition}\label{p:second-general}
For every $\alpha>0$, $2\le p<\infty$ and $f\in L^p(\R^3)$,
\begin{equation}\label{e:second-general}
\left\|S_{1,\alpha}^{\psi_1}f\right\|_{L^p(\R^3)}
\le C_{\psi_2,\alpha}p\,
\left\|S_{z,\alpha}^{\psi}f\right\|_{L^p(\R^3)}.
\end{equation}
The constant is independent of $p$ and $f$.
No additional support condition on either kernel is required.
\end{proposition}

To prove \eqref{e:second-general}, we first reconstruct
$(f(\cdot,\pmb y)*\psi_{1,s})(u)$ from
$(f*\psi_{s,t})(u,\pmb v)$.

The proof of Proposition~\ref{prop:square-2}
needs a normalized radial dual kernel supported in $B(0,\rho)$
with $54\rho<\alpha$.
By (ii) of Lemma~\ref{l:loc-pairing}, this requires
$d_2(\psi_2)<\rho<\alpha/54$, which may not hold here.
We instead use the two strip kernels in Lemma~\ref{l:strip-duals}.

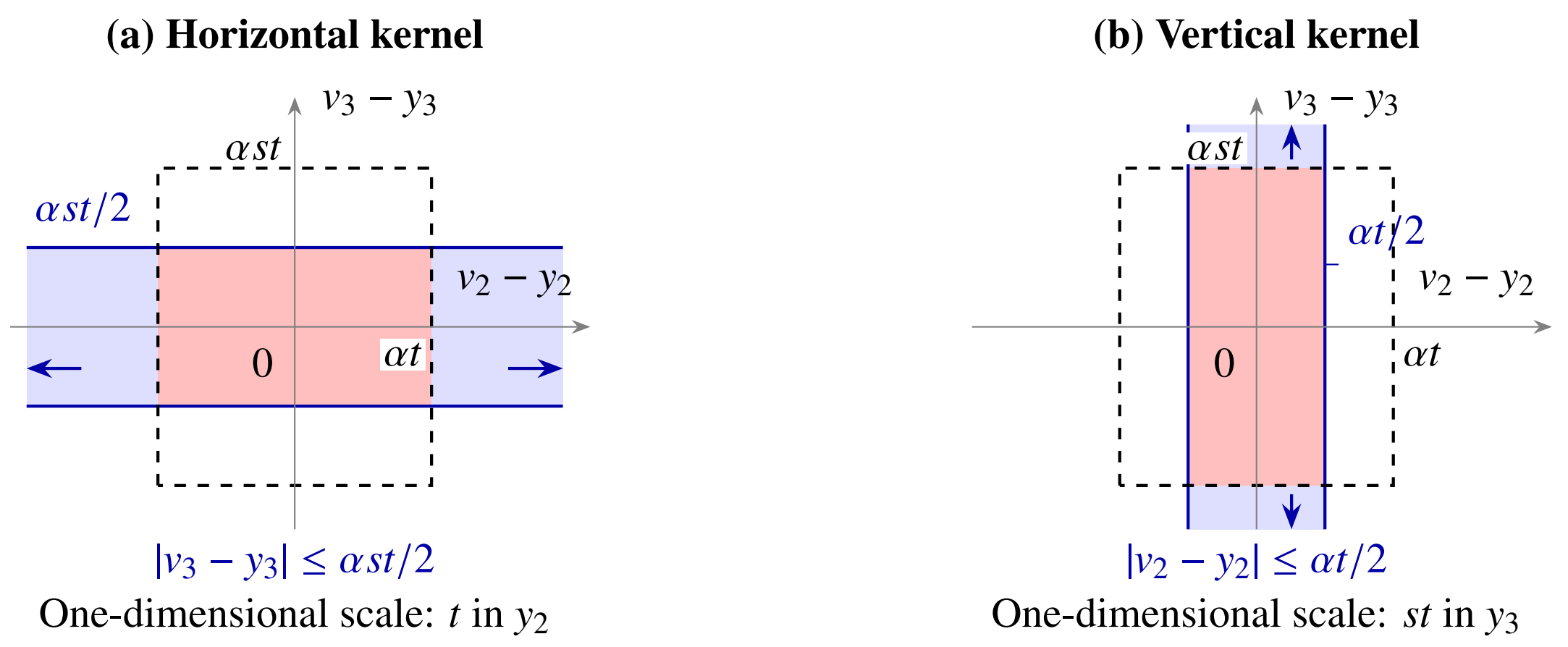
\begin{figure}[!htbp]
\centering
\begin{minipage}[t]{0.48\textwidth}
\centering
\begin{tikzpicture}[x=1cm,y=1cm,>=Stealth,
every node/.style={font=\small,text=black}]
\fill[blue!13] (-2.45,-.725) rectangle (2.45,.725);
\fill[red!25] (-1.25,-.725) rectangle (1.25,.725);
\draw[blue!65!black,thick] (-2.45,-.725)--(2.45,-.725);
\draw[blue!65!black,thick] (-2.45,.725)--(2.45,.725);
\draw[->,blue!65!black,thick] (-1.95,-.38)--(-2.45,-.38);
\draw[->,blue!65!black,thick] (1.95,-.38)--(2.45,-.38);
\draw[dashed,black,thick] (-1.25,-1.45) rectangle (1.25,1.45);
\draw[->,gray] (-2.6,0)--(2.7,0);
\draw[->,gray] (0,-1.85)--(0,2.1);
\node[anchor=south east] at (2.7,.12) {$v_2-y_2$};
\node[anchor=west] at (.10,2.05) {$v_3-y_3$};
\node[anchor=north east] at (-.05,-.05) {$0$};
\node[anchor=north east,fill=white,inner sep=1pt] at (1.20,-.10) {$\alpha t$};
\node[anchor=south east,fill=white,inner sep=1pt] at (-.08,1.48) {$\alpha st$};
\node[anchor=south east,text=blue!65!black] at (-1.35,.76) {$\alpha st/2$};
\node[font=\small\bfseries] at (0,2.65) {(a) Horizontal kernel};
\node[text=blue!65!black] at (0,-2.15) {$|v_3-y_3|\le\alpha st/2$};
\node at (0,-2.65) {One-dimensional scale: $t$ in $y_2$};
\end{tikzpicture}
\end{minipage}\hfill
\begin{minipage}[t]{0.48\textwidth}
\centering
\begin{tikzpicture}[x=1cm,y=1cm,>=Stealth,
every node/.style={font=\small,text=black}]
\fill[blue!13] (-.625,-1.85) rectangle (.625,1.85);
\fill[red!25] (-.625,-1.45) rectangle (.625,1.45);
\draw[blue!65!black,thick] (-.625,-1.85)--(-.625,1.85);
\draw[blue!65!black,thick] (.625,-1.85)--(.625,1.85);
\draw[->,blue!65!black,thick] (.32,1.53)--(.32,1.85);
\draw[->,blue!65!black,thick] (.32,-1.53)--(.32,-1.85);
\draw[dashed,black,thick] (-1.25,-1.45) rectangle (1.25,1.45);
\draw[->,gray] (-2.6,0)--(2.7,0);
\draw[->,gray] (0,-1.85)--(0,2.1);
\node[anchor=south east] at (2.7,.12) {$v_2-y_2$};
\node[anchor=west] at (.10,2.05) {$v_3-y_3$};
\node[anchor=north east] at (-.05,-.05) {$0$};
\node[anchor=north west,fill=white,inner sep=1pt] at (1.30,-.10) {$\alpha t$};
\node[anchor=south east,fill=white,inner sep=1pt] at (-.08,1.48) {$\alpha st$};
\node[anchor=south west,text=blue!65!black] at (.69,.55) {$\alpha t/2$};
\draw[blue!65!black] (.625,.57)--(.75,.57);
\node[font=\small\bfseries] at (0,2.65) {(b) Vertical kernel};
\node[text=blue!65!black] at (0,-2.15) {$|v_2-y_2|\le\alpha t/2$};
\node at (0,-2.65) {One-dimensional scale: $st$ in $y_3$};
\end{tikzpicture}
\end{minipage}
\caption[The two strip kernels at scales $t,st$]{{The dashed rectangles
show the cone restrictions in the $(v_2-y_2,v_3-y_3)$ plane. The blue strips contain
the supports of the dilated kernels. Their intersections with the rectangles are red.}}
\label{fig:strips}
\end{figure}

\begin{itemize}
\item \textbf{Two strip kernels for the reproducing formula.}
Fix $\alpha>0$.
We construct $\Phi_{\mathrm{hor}}^{(2)},\Phi_{\mathrm{ver}}^{(2)}$ with
$$
\operatorname{supp}\Phi_{\mathrm{hor}}^{(2)}
\subset\R\times[-\alpha/2,\alpha/2],
\qquad
\operatorname{supp}\Phi_{\mathrm{ver}}^{(2)}
\subset[-\alpha/2,\alpha/2]\times\R.
$$
For every $\pmb\theta=(\theta_2,\theta_3)\in\mathbb S^1$ (i.e., $\theta_2^2+\theta_3^2=1$), 
these kernels satisfy
\begin{equation}\label{Phi-h-v}
	\int_0^\infty\widehat\psi_2(r,0)
\widehat{\Phi_{\mathrm{hor}}^{(2)}}(r\pmb\theta)\,\frac{dr}{r}
 =\theta_2^2, \qquad
\int_0^\infty\widehat\psi_2(r,0)
\widehat{\Phi_{\mathrm{ver}}^{(2)}}(r\pmb\theta)\,\frac{dr}{r}
 =\theta_3^2.
\end{equation}
Thus the two kernels together take the place of $\Phi^{(2)}$
in the Calder\'on reproducing formula; see \eqref{e:strip-repro}.

\smallskip

\item \textbf{Matching the cone restrictions.}
The dilation $(x_2,x_3)\mapsto(tx_2,stx_3)$ gives
$$
\Phi_{j;s,t}^{(2)}(\pmb y-\pmb v)
=\frac1{st^2}\Phi_j^{(2)}
\left(\frac{y_2-v_2}{t},\frac{y_3-v_3}{st}\right),
\qquad j\in\{\mathrm{hor},\mathrm{ver}\}.
$$
By the strip supports above, if
$\Phi_{\mathrm{hor};s,t}^{(2)}(\pmb y-\pmb v)\ne0$, then
$$
\left|\frac{y_3-v_3}{st}\right|\le\alpha/2,
\qquad
|v_3-y_3|\le\alpha st/2<\alpha st.
$$
Similarly, if $\Phi_{\mathrm{ver};s,t}^{(2)}(\pmb y-\pmb v)\ne0$, then
$$
\left|\frac{y_2-v_2}{t}\right|\le\alpha/2,
\qquad
|v_2-y_2|\le\alpha t/2<\alpha t.
$$
Each support ensures one cone restriction in $S_{z,\alpha}^{\psi}f$;
see Figure~\ref{fig:strips}.

\smallskip
\item \textbf{The extra $\sqrt p$ in the long direction.}
The kernels have integral zero in their long directions:
$$
\int_\R\Phi_{\mathrm{hor}}^{(2)}(x_2,x_3)\,dx_2=0,
\qquad
\int_\R\Phi_{\mathrm{ver}}^{(2)}(x_2,x_3)\,dx_3=0.
$$
In each long direction, we first decompose the kernel into
compactly supported kernels with integral zero.
For each kernel, the Haar estimate uses an area function whose aperture
depends on its support radius.
If this aperture exceeds $\alpha$, we compare it with $\alpha$:
$$
\underbrace{\sqrt p}_{L^2\text{-valued Haar estimate}}
\;\times\;
\underbrace{\sqrt p}_{\text{aperture comparison}}=p.
$$
The kernel decay allows us to sum these bounds;
see \eqref{e:molecular-one}.
At aperture $\alpha$, both restrictions
$|v_2-y_2|<\alpha t$ and $|v_3-y_3|<\alpha st$ hold,
so the resulting area function is bounded by $S_{z,\alpha}^{\psi}f$.
We add the horizontal and vertical estimates.
The long-direction scale is $t$ for the horizontal kernel and $st$ for the vertical kernel.
\end{itemize}

\subsection{Dual kernels supported in strips}\label{sub:strip-duals}

We first construct the horizontal kernel, then obtain the vertical kernel
by interchanging $x_2$ and $x_3$.

\begin{lemma}\label{l:strip-duals}
For every $\alpha>0$, there is a real-valued smooth function
$\Phi_{\mathrm{hor}}^{(2)}\in L^1(\R^2)$ such that
$$
\operatorname{supp}\Phi_{\mathrm{hor}}^{(2)}
\subset\R\times[-\alpha/2,\alpha/2],
\qquad
\int_\R\Phi_{\mathrm{hor}}^{(2)}(x_2,x_3)\,dx_2=0.
$$
For some $m\in\mathbb N_0$, every $a,b\in\mathbb N_0$
and $(x_2,x_3)\in\R^2$,
\begin{equation}\label{e:strip-decay}
\left|\partial_2^a\partial_3^b\Phi_{\mathrm{hor}}^{(2)}(x_2,x_3)\right|
\le C_{a,b}\mathbf1_{\{|x_3|<\alpha/2\}}
(1+|x_2|)^{-3-a}\log^m(e+|x_2|).
\end{equation}
Define
$$
\Phi_{\mathrm{ver}}^{(2)}(x_2,x_3)
=\Phi_{\mathrm{hor}}^{(2)}(x_3,x_2).
$$
The same support, cancellation and derivative bounds hold for
$\Phi_{\mathrm{ver}}^{(2)}$ with $x_2,x_3$ interchanged.

Moreover, for every $\pmb\theta=(\theta_2,\theta_3)\in\mathbb S^1$,
\begin{equation}\label{e:strip-dual}
\int_0^\infty\widehat\psi_2(r,0)
\left(\widehat{\Phi_{\mathrm{hor}}^{(2)}}(r\pmb\theta)
+\widehat{\Phi_{\mathrm{ver}}^{(2)}}(r\pmb\theta)\right)
\,\frac{dr}{r}=1,
\end{equation}
and
\begin{equation}\label{e:strip-ray}
\sup_{\pmb\theta\in\mathbb S^1}\int_0^\infty
\left(\,  \Big |\widehat{\Phi_{\mathrm{hor}}^{(2)}}(r\pmb\theta)\Big |^2
+\Big |\widehat{\Phi_{\mathrm{ver}}^{(2)}}(r\pmb\theta)\Big |^2\right)
\,\frac{dr}{r}<\infty.
\end{equation}
The kernels and $m$ depend only on $\psi_2$ and $\alpha$;
$C_{a,b}$ depends only on $\psi_2,\alpha,a,b$.
\end{lemma}

The estimate \eqref{e:strip-ray} will be used to prove the
$L^2$ bound \eqref{e:strip-recon-l2} and the $L^2$ convergence
of the two reconstruction integrals below.

\begin{proof}[Proof of Lemma~\ref{l:strip-duals}]
\smallskip
\noindent{\it Step 1.}
We begin by constructing the horizontal kernel.

{Choose a fixed real radial bump of integral one in $B(0,1)$
and dilate it to obtain $\zeta\in C_c^\infty(B(0,\alpha/2))$ of integral one.}
To keep the kernel zero for $|x_3|\ge\alpha/2$,
we multiply $\widehat\zeta(\pmb\xi)$ by a function of $\xi_2$ only.

For $\pmb\xi=r\pmb\theta$, we have $\xi_2^2=r^2\theta_2^2$.
Hence we first consider the multiplier $\xi_2^2$
to obtain $\theta_2^2$ in \eqref{Phi-h-v}.
By radiality of $\zeta$,
$$
\int_0^\infty\widehat\psi_2(r,0)
(r\theta_2)^2\widehat\zeta(r\pmb\theta)\,\frac{dr}{r}
=\theta_2^2\int_0^\infty
\widehat\psi_2(r,0)\widehat\zeta(r,0)r^2\,\frac{dr}{r}.
$$
If the integral on the right is nonzero, we can divide by it to obtain the required normalization.
To also cover the case when it vanishes, define

$$
\mathcal M(z)=\int_0^\infty
\widehat\psi_2(r,0)\widehat\zeta(r,0)r^z\,\frac{dr}{r},
\qquad \operatorname{Re}z>-2.
$$
Since $\psi_2$ is radial and $\int_{\R^2}\psi_2=0$, we have
$$
\widehat\psi_2(r,0)
=\int_{\R^2}\psi_2(\pmb x)
\big (\cos(2\pi r x_2)-1\big  )\,d\pmb x
=O(r^2)
\qquad(r\to0).
$$
The Fourier transforms $\widehat\psi_2$ and $\widehat\zeta$
decay rapidly at infinity.
For a compact set $K\subset\{z:\operatorname{Re}z>-2\}$,
choose $\delta>0$ such that
$\operatorname{Re}z\ge-2+\delta$ for all $z\in K$.
Then, for every $j\in\mathbb N_0$,
$$
\sup_{z\in K}
\left|
\widehat\psi_2(r,0)\widehat\zeta(r,0)
r^{z-1}(\log r)^j
\right|
\le C_{K,j}
\begin{cases}
r^{-1+\delta}|\log r|^j,&0<r<1,\\[2pt]
r^{-2}(\log r)^j,&r\ge1.
\end{cases}
$$
The right-hand side is integrable in $r$, so differentiation under the integral shows that
$\mathcal M$ is analytic on $\operatorname{Re}z>-2$, with
$$
\mathcal M^{(j)}(z)
=\int_0^\infty
\widehat\psi_2(r,0)\widehat\zeta(r,0)
r^z(\log r)^j\,\frac{dr}{r}.
$$

{To see that $\mathcal M\not\equiv0$, let
$F(v)=\widehat\psi_2(e^v,0)\widehat\zeta(e^v,0)$. Then $F\in L^1(\R)$ and
$$
\widehat F(\tau)=\mathcal M(-2\pi i\tau).
$$
Since $\psi_2$ is nonzero, radial and compactly supported,
$r\mapsto\widehat\psi_2(r,0)$ is real analytic and cannot vanish on an interval.
Also $\widehat\zeta(0)=1$, so $F\ne0$. Fourier uniqueness gives $\mathcal M\not\equiv0$.}

{By analyticity,}
\begin{equation}\label{mD}
m:=\min\left\{j\in\mathbb N_0:\mathcal M^{(j)}(2)\ne0\right\}
\qquad \text{and}\qquad
\mathsf D:=\mathcal M^{(m)}(2)\in\R\setminus\{0\}
\end{equation}
are well defined. To normalize the horizontal kernel, let
\begin{equation}\label{e:strip-symbol}
\widehat{\Phi_{\mathrm{hor}}^{(2)}}(\pmb\xi)
=
\begin{cases}
\mathsf D^{-1}\xi_2^2(\log|\xi_2|)^m\widehat\zeta(\pmb\xi),
&\xi_2\ne0,\\[2pt]
0,&\xi_2=0.
\end{cases}
\end{equation}

\smallskip

\noindent{\it Step 2.}
By partial Fourier inversion, \eqref{e:strip-symbol} becomes
$$
\Phi_{\mathrm{hor}}^{(2)}(x_2,x_3)
=\mathsf D^{-1}\int_\R\xi^2(\log|\xi|)^m
\bigg(\int_\R\zeta(v_2,x_3)e^{-2\pi i v_2\xi}\,dv_2\bigg)
e^{2\pi i x_2\xi}\,d\xi.
$$
Thus $\Phi_{\mathrm{hor}}^{(2)}$ is real, smooth and even in each variable,
with
$
\operatorname{supp}\Phi_{\mathrm{hor}}^{(2)}
\subset\R\times[-\alpha/2,\alpha/2].
$

For $a,b\in\mathbb N_0$,     
$$
\partial_2^a\partial_3^b\Phi_{\mathrm{hor}}^{(2)}(\pmb x)
=\int_\R b_{a,b}(\xi,x_3)e^{2\pi i x_2\xi}\,d\xi,
$$
where
$$
b_{a,b}(\xi,x_3)
=\mathsf D^{-1}(2\pi i\xi)^a\xi^2(\log|\xi|)^m
\int_\R\partial_3^b\zeta(v_2,x_3)e^{-2\pi i v_2\xi}\,dv_2.
$$
By the product rule, we obtain
$$
|\partial_\xi^j b_{a,b}(\xi,x_3)|
\le C_{a,b} \, |\xi|^{a+2-j} \left(1+\left|\log|\xi|\right|\right)^m,
\qquad 0<|\xi|\le1,\quad 0\le j\le a+4.
$$
For $|\xi|\ge1$, these derivatives decay rapidly. All these bounds are uniform in $x_3$.

For $|x_2|\ge2$, split the integral at $|\xi|=|x_2|^{-1}$ and integrate by parts $a+4$ times
on each of the two unbounded intervals. This gives
\begin{align*}
\left|\partial_2^a\partial_3^b\Phi_{\mathrm{hor}}^{(2)}(\pmb x)\right|
\leq &\int_{|\xi|\le |x_2|^{-1}}|b_{a,b}(\xi,x_3)|\,d\xi
 +\sum_{j=0}^{a+3}(2\pi|x_2|)^{-j-1}
\left(
|\partial_\xi^j b_{a,b}(|x_2|^{-1},x_3)|
+|\partial_\xi^j b_{a,b}(-|x_2|^{-1},x_3)|
\right)\\
&   \qquad +(2\pi|x_2|)^{-a-4}\int_{|\xi|\ge |x_2|^{-1}}
\left |\partial_\xi^{a+4}b_{a,b}(\xi,x_3)  \right|   d\xi\\
\leq &C_{a,b}\bigg[
\int_0^{1/|x_2|}\xi^{a+2}(1+|\log\xi|)^m\,d\xi   +(1+\log|x_2|)^m
\sum_{j=0}^{a+3}|x_2|^{-j-1}|x_2|^{j-a-2}\\
&\qquad +|x_2|^{-a-4}\,
\bigg(1+\int_{1/|x_2|}^1\xi^{-2}(1+|\log\xi|)^m\,d\xi\bigg)
\bigg]\\
\leq &C_{a,b}|x_2|^{-a-3}(1+\log|x_2|)^m.
\end{align*}
For $|x_2|<2$,
$$
\left|\partial_2^a\partial_3^b\Phi_{\mathrm{hor}}^{(2)}(\pmb x)\right|
\le\int_\R|b_{a,b}(\xi,x_3)|\,d\xi
\le C_{a,b}.
$$
Together with the vanishing for $|x_3|\ge\alpha/2$,
we obtain \eqref{e:strip-decay}.

Taking $a=b=0$, we have $\Phi_{\mathrm{hor}}^{(2)}\in L^1(\R^2)$
and $\Phi_{\mathrm{hor}}^{(2)}(\cdot,x_3)\in L^1(\R)$ for every $x_3$.
By partial Fourier transformation,
\begin{align*}
\int_\R\Phi_{\mathrm{hor}}^{(2)}(x_2,x_3)\,dx_2
= & \mathsf  D^{-1}\lim_{\xi\to0}\xi^2(\log|\xi|)^m
\int_\R\zeta(v_2,x_3)e^{-2\pi i v_2\xi}\,dv_2
=0.
\end{align*}

\smallskip
\noindent{\it Step 3.}
We now verify the normalization and the square integral bound for
$\Phi_{\mathrm{hor}}^{(2)}$.
For $\theta_2\ne0$, by radiality and $\mathcal M^{(j)}(2)=0$ for $j<m$, we have
\begin{align*}
\int_0^\infty\widehat\psi_2(r,0)
\widehat{\Phi_{\mathrm{hor}}^{(2)}}(r\pmb\theta)\,\frac{dr}{r}
= &\frac{\theta_2^2}{\mathsf D}\int_0^\infty
\widehat\psi_2(r,0)\widehat\zeta(r,0)r^2
(\log r+\log|\theta_2|)^m\,\frac{dr}{r}\\
= &\frac{\theta_2^2}{\mathsf D}\sum_{j=0}^m\binom mj
\, (\log|\theta_2|)^{m-j}\mathcal M^{(j)}(2)
=\theta_2^2.
\end{align*}
For $\theta_2=0$, both sides are zero.

By \eqref{e:strip-symbol} and noting
$\sup_{0<a\le1}a^2|\log a|^m<\infty$, we have
\begin{align*}
\sup_{\pmb\theta\in\mathbb S^1}
\left|\widehat{\Phi_{\mathrm{hor}}^{(2)}}(r\pmb\theta)\right|
\leq &Cr^2|\widehat\zeta(r,0)|
\sup_{0<a\le1}a^2|\log r+\log a|^m
\leq Cr^2(1+|\log r|)^m|\widehat\zeta(r,0)|.
\end{align*}
Hence
\begin{align*}
\sup_{\pmb\theta\in\mathbb S^1}\int_0^\infty
\Big |\widehat{\Phi_{\mathrm{hor}}^{(2)}}(r\pmb\theta)\Big |^2
\,\frac{dr}{r}
\leq  C\int_0^\infty
r^3(1+|\log r|)^{2m}|\widehat\zeta(r,0)|^2\,dr<\infty.
\end{align*}

Define $\Phi_{\mathrm{ver}}^{(2)}$ by interchanging $x_2,x_3$. Then
$
\widehat{\Phi_{\mathrm{ver}}^{(2)}}(\xi_2,\xi_3)
=\widehat{\Phi_{\mathrm{hor}}^{(2)}}(\xi_3,\xi_2).
$
The support, cancellation and derivative bounds for
$\Phi_{\mathrm{ver}}^{(2)}$ follow from interchanging $x_2,x_3$.
Similarly,
$$
\int_0^\infty\widehat\psi_2(r,0) \,
\widehat{\Phi_{\mathrm{ver}}^{(2)}}(r\pmb\theta)\,\frac{dr}{r}
=\theta_3^2.
$$
This proves the two identities in \eqref{Phi-h-v}.
Since $\pmb\theta=(\theta_2,\theta_3)\in\mathbb S^1$, their sum is $1$, which proves \eqref{e:strip-dual}.

The same square integral bound holds for $\Phi_{\mathrm{ver}}^{(2)}$
by interchanging coordinates. Hence \eqref{e:strip-ray} holds. 
\end{proof}

\begin{remark}
Since $\psi_2$ and $\zeta$ are real and radial, we have
$$
\mathcal M(2)
=\frac1{2\pi}\int_{\R^2}
\widehat\psi_2(\pmb\xi)\widehat\zeta(\pmb\xi)\,d\pmb\xi
=\frac1{2\pi}\int_{\R^2}\psi_2(\pmb x)\zeta(\pmb x)\,d\pmb x.
$$
In particular, $\mathcal M(2)=0$ if $d_2(\psi_2)\ge\alpha/2$,
since $\operatorname{supp}\zeta\subset B(0,\alpha/2)$.
The logarithmic factor in \eqref{e:strip-symbol} comes from
$$
\left.\frac{d^m}{dz^m}|\xi_2|^z\right|_{z=2}
=\xi_2^2\left(\log|\xi_2| \right)^m,
\qquad \xi_2\ne0.
$$
Using the first nonzero derivative of $\mathcal M$ at $2$
allows normalization even when $\mathcal M(2)=0$.
The lower derivatives vanish, so the terms containing
$\log|\theta_2|$ disappear from the normalization integral.

\end{remark}

\subsection{Reconstruction with the two strip kernels}\label{sub:strip-repro}

We now prove Proposition~\ref{p:second-general}.
Fix $\alpha>0$ and $2\le p<\infty$.
First let $f\in L^2(\R^3)\cap L^p(\R^3)$.
Fix $x_1$ for which \eqref{eqn:aux-finite}, \eqref{e:second-energy} hold
and $S_{z,\alpha}^{\psi}f(x_1,\cdot)\in L^p(\R^2)$.
These conditions hold for almost every $x_1$ by those identities
and Lemma~\ref{l:continuous-forward}.
Choose the kernels from Lemma~\ref{l:strip-duals} for this $\alpha$.
All constants below depend only on $\psi_2$ and $\alpha$.

We first reconstruct $(f(\cdot,\pmb y)*\psi_{1,s})(u)$
using the two strip kernels.
Let
$$
\mathcal R_j(f)(u,s;\pmb y)
=\int_0^\infty\int_{\R^2}(f*\psi_{s,t})(u,\pmb v)
\Phi_{j;s,t}^{(2)}(\pmb y-\pmb v)\,d\pmb v\,\frac{dt}{t},\qquad j\in\{\mathrm{hor},\mathrm{ver}\},
$$
where 
$$
\Phi_{j;s,t}^{(2)}(\pmb y)
=\frac1{st^2}\Phi_j^{(2)}\left(\frac{y_2}{t},\frac{y_3}{st}\right).
$$
{We interpret the $t$-integral as an $L^2$ limit.
For each $s>0$ and $\pmb\xi\ne0$, the substitution
$r=t|(\xi_2,s\xi_3)|$ in \eqref{e:strip-ray} gives
$$
\int_0^\infty
|\widehat{\Phi_j^{(2)}}(t\xi_2,st\xi_3)|^2\,\frac{dt}{t}\le C.
$$
For $0<a<b<\infty$, Cauchy--Schwarz in $t$ and Plancherel in $\pmb y$
give the estimate below with both sides restricted to $a<t<b$. By \eqref{e:second-energy}
the truncations are Cauchy in $L^2$. In the limit,}
\begin{align}\label{e:strip-recon-l2}
&\int_{\R^2}\int_{\Gamma_\alpha(x_1)}
|\mathcal R_j(f)(u,s;\pmb y)|^2\,\frac{du\,ds}{s^2}\,d\pmb y\notag\\
\leq &C\int_{\Gamma_\alpha(x_1)}\int_0^\infty\int_{\R^2}
|(f*\psi_{s,t})(u,\pmb v)|^2\,d\pmb v\,\frac{dt}{t}\,\frac{du\,ds}{s^2}<\infty.
\end{align}
{In particular,}

\begin{equation}\label{e:strip-recon-limit}
\lim_{N\to\infty}
\bigg\|
\int_{N^{-1}}^N\int_{\R^2}(f*\psi_{s,t})(u,\pmb v)
\Phi_{j;s,t}^{(2)}(\pmb y-\pmb v)\,d\pmb v\,\frac{dt}{t}
-\mathcal R_j(f)(u,s;\pmb y)
\bigg\|_{L^2   (\Gamma_\alpha(x_1)\times\R^2,\frac{du\,ds}{s^2}\,d\pmb y   )}=0.
\end{equation}

For $s>0$ and $\pmb\xi\ne0$,  note that  \eqref{e:strip-dual} with 
$r=t(\xi_2^2+s^2\xi_3^2)^{1/2}$ becomes
$$
\int_0^\infty\widehat\psi_2(t\xi_2,st\xi_3)
\left(\widehat{\Phi_{\mathrm{hor}}^{(2)}}(t\xi_2,st\xi_3)
+\widehat{\Phi_{\mathrm{ver}}^{(2)}}(t\xi_2,st\xi_3)\right)\,\frac{dt}{t}=1.
$$
The integral is absolutely convergent by Cauchy--Schwarz,
\eqref{e:strip-ray}, and $c_{\psi_2}<\infty$. Hence,
 by Fourier transformation in $\pmb y$, we have 
\begin{equation}\label{e:strip-repro}
(f(\cdot,\pmb y)*\psi_{1,s})(u)
=\mathcal R_{\mathrm{hor}}(f)(u,s;\pmb y)
+\mathcal R_{\mathrm{ver}}(f)(u,s;\pmb y)\qquad \text{in}  \   \   L^2(\Gamma_\alpha(x_1)\times\R^2,\frac{du\,ds}{s^2}\,d\pmb y).
\end{equation}

We will estimate the horizontal and vertical terms in the following two subsections, respectively.

\subsection{Estimate for the horizontal term}\label{sub:horizontal}

{In the truncated horizontal integral, write $r=(y_3-v_3)/(st)$. Then}
\begin{align*}
\Phi_{\mathrm{hor};s,t}^{(2)}(\pmb y-\pmb v)
\left|\frac{\partial v_3}{\partial r}\right|\,dr
= &\frac1{st^2}\Phi_{\mathrm{hor}}^{(2)}
\left(\frac{y_2-v_2}{t},r\right)st\,dr
= \frac1t\Phi_{\mathrm{hor}}^{(2)}
\left(\frac{y_2-v_2}{t},r\right)\,dr,
\end{align*}
{Passing to the $L^2$ limit in \eqref{e:strip-recon-limit} gives}
\begin{align}\label{e:horizontal-recon}
\mathcal R_{\mathrm{hor}}(f)(u,s;\pmb y)
= &\int_{-\alpha/2}^{\alpha/2}\int_0^\infty\int_\R
(f*\psi_{s,t})(u,v_2,y_3-str)
 \frac1t\Phi_{\mathrm{hor}}^{(2)}
\left(\frac{y_2-v_2}{t},r\right)\,dv_2\,\frac{dt}{t}\,dr.
\end{align}

For each fixed $y_3$,  let 
\begin{equation}\label{e:coeff-h}
F(v_2,t;u,s,r):=(f*\psi_{s,t})(u,v_2,y_3-str).
\end{equation}
The one-dimensional reconstruction in $y_2$ has scale $t$:
$$
\mathcal R_{\Phi_{\mathrm{hor}}^{(2)}}F(y_2;u,s,r)
=\int_0^\infty\int_\R F(v_2,t;u,s,r)
\frac1t\Phi_{\mathrm{hor}}^{(2)}\left(\frac{y_2-v_2}{t},r\right)\,dv_2\,\frac{dt}{t}.
$$
We keep $(u,s,r)$ in the $L^2$ norm and define, for $\beta>0$,
\begin{equation}\label{e:compact-area}
\mathcal A_\beta F(y_2)^2
=\int_0^\infty\int_{|v_2-y_2|<\beta t}
\int_{\Gamma_\alpha(x_1)}\int_{-\alpha/2}^{\alpha/2}
|F(v_2,t;u,s,r)|^2\,dr\,\frac{du\,ds}{s^2}\,\frac{dv_2\,dt}{t^2}.   
\end{equation}

For fixed $y_3$, changing variables from $r$ to $v_3=y_3-str$, we have 
\begin{align}\label{e:horizontal-area}
\mathcal A_\alpha F(y_2)^2
= &\int_{\Gamma_\alpha(x_1)}\int_{-\alpha/2}^{\alpha/2}
\int_0^\infty\int_{|v_2-y_2|<\alpha t}
|(f*\psi_{s,t})(u,v_2,y_3-str)|^2
\,\frac{dv_2\,dt}{t^2}\,dr\,\frac{du\,ds}{s^2}\notag\\
= &\int_{\Gamma_\alpha(x_1)}\int_0^\infty
\int_{\substack{|v_2-y_2|<\alpha t\\|v_3-y_3|<\alpha st/2}}
|(f*\psi_{s,t})(u,\pmb v)|^2
\,\frac{d\pmb v\,dt}{st^3}\,\frac{du\,ds}{s^2}\notag\\
\leq &S_{z,\alpha}^{\psi}f(x_1,\pmb y)^2.
\end{align}
 To verify the $L^2$ condition on $F$, integrate first in $y_3$, keeping $s,t,r$ fixed,
{and write} $v_3=y_3-str$, so that $dy_3=dv_3$:
\begin{align}\label{e:strip-coeff}
&\int_\R
\left\|(f*\psi_{s,t})(u,v_2,y_3-str)\right\|_{L^2(dv_2\,dt/t\,du\,ds/s^2\,dr)}^2
\,dy_3\notag\\
= &\int_{\Gamma_\alpha(x_1)}\int_0^\infty
\int_{-\alpha/2}^{\alpha/2}\int_{\R^2}
|(f*\psi_{s,t})(u,\pmb v)|^2
\,d\pmb v\,dr\,\frac{dt}{t}\,\frac{du\,ds}{s^2}\notag\\
= &\alpha\int_{\Gamma_\alpha(x_1)}\int_0^\infty\int_{\R^2}
|(f*\psi_{s,t})(u,\pmb v)|^2
\,d\pmb v\,\frac{dt}{t}\,\frac{du\,ds}{s^2}<\infty,
\end{align}
where the finiteness follows from \eqref{e:second-energy}.  Thus, for almost every $y_3$, $F$ satisfies the two assumptions
in Lemma~\ref{l:molecular-one} below:
$$
\mathcal A_\alpha F\in L^p(\R)
\qquad  \text{and}\qquad 
 F\in L^2\Big (dv_2\,\frac{dt}{t}\,\frac{du\,ds}{s^2}\,dr\Big ).
$$

We next estimate the horizontal term
$\mathcal R_{\mathrm{hor}}(f)$ in \eqref{e:horizontal-recon}.
To construct smooth atoms and expand them in Haar functions
as in the proof of Lemma~\ref{sp:l:smooth},
we first decompose $\Phi_{\mathrm{hor}}^{(2)}(\cdot,r)$
into compactly supported kernels, each with integral zero.

\smallskip
\noindent{\bf {Decomposition of the strip kernels.}}   \   \   
To preserve the zero integral when localizing in the long direction,
we cut off a primitive and then differentiate. Define
$$
P(x_2,r)=\int_{-\infty}^{x_2}\Phi_{\mathrm{hor}}^{(2)}(v,r)\,dv,\qquad x_2\in\R,\     \       |r|<\alpha/2.
$$
This, combined with $\int_\R\Phi_{\mathrm{hor}}^{(2)}(v,r)\,dv=0$, yields
$$
P(x_2,r)=-\int_{x_2}^\infty\Phi_{\mathrm{hor}}^{(2)}(v,r)\,dv.
$$
Using the first formula for $x_2\le0$ and the second for $x_2>0$, together with \eqref{e:strip-decay}, we obtain
\begin{equation}\label{e:P-bound}
|P(x_2,r)|
\leq C_{\psi_2,\alpha}
\int_{|x_2|}^\infty(1+v)^{-3}\log^m(e+v)\,dv
\leq   C_{\psi_2,\alpha}
(1+|x_2|)^{-2}\log^m(e+|x_2|),
\end{equation}
where $m$ is defined in \eqref{mD}.

Choose $\vartheta\in C_c^\infty((-2,2))$ satisfying $\vartheta\equiv 1$  on $[-1,1]$. Let
$$
\eta_0(x_2)=\vartheta(x_2),
\qquad
\eta_k(x_2)=\vartheta(2^{-k}x_2)-\vartheta(2^{1-k}x_2)
\quad(k\ge1),
$$
and
$$
\widetilde{\mathsf{b}}_{k,r}(x_2)
=\frac{d}{dx_2}\big (\eta_k(x_2)P(x_2,r)\big  ),
\qquad k\ge0,   \     \       |r|<\alpha/2.
$$
Then
$
\operatorname{supp}\widetilde{\mathsf{b}}_{k,r}\subset[-2^{k+1},2^{k+1}]$ and 
$\int_\R\widetilde{\mathsf{b}}_{k,r}=0$.

{By the product rule,
$$
\widetilde{\mathsf b}_{k,r}
=\eta_k\Phi_{\mathrm{hor}}^{(2)}(\cdot,r)+\eta_k'P(\cdot,r).
$$}
Note that for $k\ge1$, $\operatorname{supp}\eta_k
\subset\{x_2:2^{k-1}\le|x_2|\le2^{k+1}\}$, and 
$$
\bigg\|\frac{d^j\eta_k}{dx_2^j}\bigg\|_{L^\infty}
\le C2^{-jk},
\quad j=0,1,2.
$$
It follows from \eqref{e:strip-decay} and \eqref{e:P-bound} that
$$
\bigg\|\frac{d^\ell\, \widetilde{\mathsf b}_{k,r}}{dx_2^\ell}\bigg\|_{L^\infty}
\le C_{\psi_2,\alpha}(1+k)^m2^{-(3+\ell)k},
\qquad k\ge0,\quad \ell=0,1,
$$
uniformly for $|r|<\alpha/2$. For $k=0$, this follows from the boundedness of $\Phi_{\mathrm{hor}}^{(2)}$, its derivative and $P$.

Since $\sum_{k=0}^N\eta_k=\vartheta(2^{-N}\cdot)$,
$$
\Phi_{\mathrm{hor}}^{(2)}(\cdot,r)-\sum_{k=0}^N\widetilde{\mathsf{b}}_{k,r}
=(1-\vartheta(2^{-N}\cdot))\Phi_{\mathrm{hor}}^{(2)}(\cdot,r)
-2^{-N}\vartheta'(2^{-N}\cdot)P(\cdot,r).
$$
By \eqref{e:strip-decay} and \eqref{e:P-bound} again,
\begin{align}\label{e:Phi-h-converg}
&\sup_{|r|<\alpha/2}
\bigg\|\Phi_{\mathrm{hor}}^{(2)}(\cdot,r)
-\sum_{k=0}^N\widetilde{\mathsf{b}}_{k,r}\bigg\|_{L^1} \notag \\
\leq &C_{\psi_2,\alpha}\bigg(
\int_{2^N}^\infty(1+x_2)^{-3}\log^m(e+x_2)\,dx_2
+2^{-N}\int_{2^N}^{2^{N+1}}
(1+x_2)^{-2}\log^m(e+x_2)\,dx_2
\bigg) \notag \\
\leq &C_{\psi_2,\alpha}(1+N)^m2^{-2N}\to0
\qquad(N\to\infty).
\end{align}

To obtain the size and derivative bounds below, we normalize
$\widetilde{\mathsf{b}}_{k,r}$ as follows:
$$
\mathsf{b}_{k,r}=c_k^{-1}\,\widetilde{\mathsf{b}}_{k,r},
$$
where $c_k=C_{\psi_2,\alpha}(1+k)^m2^{-2k}$, with
$C_{\psi_2,\alpha}$ large enough that
$$
2^{k+1}\|\mathsf{b}_{k,r}\|_{L^\infty}
+2^{2k+2}\Big \|\frac{d\mathsf{b}_{k,r}}{dx_2}\Big \|_{L^\infty}\le1,
\qquad \forall\, k\ge0,\quad |r|<\alpha/2.
$$
Also,
\begin{equation}\label{e:ck}
\sum_{k=0}^\infty c_k
=C_{\psi_2,\alpha}
\sum_{k=0}^\infty(1+k)^m2^{-2k}<\infty.
\end{equation}
It follows from the $L^1$ convergence in \eqref{e:Phi-h-converg}
and $c_k\mathsf{b}_{k,r}=\widetilde{\mathsf{b}}_{k,r}$ that
\begin{equation}\label{e:strip-compact}
\Phi_{\mathrm{hor}}^{(2)}(\cdot,r)
=\sum_{k=0}^\infty c_k\mathsf{b}_{k,r}
\quad\text{in }L^1(\R),  \   \    \text{ uniformly for }|r|<\alpha/2.
\end{equation}

The $k$th term of \eqref{e:strip-compact} gives the following term
in \eqref{e:horizontal-recon}:
$$
c_k\int_{-\alpha/2}^{\alpha/2}
\bigg(
\int_0^\infty\int_\R
F(v_2,t;u,s,r)\,
\frac1t\mathsf{b}_{k,r}\left(\frac{y_2-v_2}{t}\right)
\,dv_2\,\frac{dt}{t}
\bigg)\,dr.
$$
The next lemma estimates the integral in $(v_2,t)$,
taking the $L^2$ norm in $(u,s,r)$ and then the $L^p$ norm in $y_2$.
 We will apply it with
$$
\mathsf b_r:=\mathsf b_{k,r},
\qquad \mathsf L:=2^{k+1}.
$$
The required support, zero integral, and size and derivative bounds
for $\mathsf b_{k,r}$ were proved above.
 We will then sum with the coefficients $c_k$ in \eqref{e:strip-compact}
and integrate in $r$ to estimate the horizontal term
$\mathcal R_{\mathrm{hor}}(f)$.

\begin{lemma}\label{l:molecular-one}
Fix $\alpha>0$, $x_1\in\R$ and $\mathsf L\in2^{\mathbb N_0}$.
Let $\mathsf{b}_r(y_2)$ be real and smooth in
$(y_2,r)\in\R\times(-\alpha/2,\alpha/2)$.
For every $|r|<\alpha/2$, assume
$$
\operatorname{supp}\mathsf{b}_r(\cdot)\subset[-\mathsf L,\mathsf L],
\qquad
\int_\R\mathsf{b}_r(y_2)\,dy_2=0,
\qquad
\mathsf L\|\mathsf{b}_r\|_{L^\infty}
+\mathsf L^2\|\partial_{y_2}\mathsf{b}_r\|_{L^\infty}\le1.
$$
For a measurable $F(v_2,t;u,s,r)$ satisfying
$$
\int_{\Gamma_\alpha(x_1)}\int_{-\alpha/2}^{\alpha/2}
\int_0^\infty\int_\R |F(v_2,t;u,s,r)|^2
\,dv_2\,\frac{dt}{t}\,dr\,\frac{du\,ds}{s^2}<\infty,
$$
define
$$
\mathcal R_{\mathsf{b}} F(y_2;u,s,r)
=\int_0^\infty\int_\R F(v_2,t;u,s,r)
\frac1t\mathsf{b}_r\left(\frac{y_2-v_2}{t}\right)\,dv_2\,\frac{dt}{t}
$$
in $L^2\left(dy_2\,\frac{du\,ds}{s^2}\,dr\right)$.
For  $2\le p<\infty$ and $\mathcal A_\alpha F\in L^p(\R)$,
where $\mathcal A_\alpha F$ is defined in \eqref{e:compact-area}, we have
\begin{equation}\label{e:compact-recon}
\bigg\|\bigg(\int_{\Gamma_\alpha(x_1)}\int_{-\alpha/2}^{\alpha/2}
|\mathcal R_{\mathsf{b}} F(\cdot;u,s,r)|^2\,dr\,\frac{du\,ds}{s^2}
\bigg)^{1/2}\bigg\|_{L^p(\R)}
\le Cp\max\left\{\mathsf L^{-1/2},\alpha^{-1/2}\right\}\|\mathcal A_\alpha F\|_{L^p(\R)}.
\end{equation}
The constant $C$ is absolute.
\end{lemma}

\begin{remark}
For fixed $\alpha$, the bound in \eqref{e:compact-recon} is uniform in $\mathsf L\ge1$.
 We will use this bound to sum the estimates for the compactly supported kernels
in \eqref{e:strip-compact}.

When $\alpha\ge30\mathsf L$, aperture comparison is unnecessary.
The proof below then gives the stronger estimate
$$
\|\mathcal R_{\mathsf b}F\|_{L^p  \big (dy_2;L^2  \big(\frac{du\,ds}{s^2}\,dr \big) \, \big)}
\le C\sqrt p\,\mathsf L^{-1/2}\|\mathcal A_\alpha F\|_{L^p(\R)}.
$$
\end{remark}

\smallskip

\begin{proof}[Proof of Lemma~\ref{l:molecular-one}]
We use the smooth-atom and Haar arguments in
Lemma~\ref{sp:l:smooth},
with the $L^2$ norm taken in $(u,s,r)$.

\noindent{\it Step 1.}
We begin by decomposing $\mathcal R_{\mathsf b}F$ into smooth atoms in $y_2$.

Concretely,
for $I\in\mathcal D^0$, let $\displaystyle\mathcal W_{\mathsf L}(I)=I\times \Big [\frac{|I|}{\mathsf L},\frac{2|I|}{\mathsf L}\Big )$,
$$
\lambda_I(F)^2
=\int_{\mathcal W_{\mathsf L}(I)}
\int_{\Gamma_\alpha(x_1)}\int_{-\alpha/2}^{\alpha/2}
|F(v_2,t;u,s,r)|^2\,dr\,\frac{du\,ds}{s^2}\,dv_2\,\frac{dt}{t},
$$
and
$$
\mathsf{a}_I(F)(y_2;u,s,r)
=\frac1{\lambda_I(F)}
\int_{\mathcal W_{\mathsf L}(I)}F(v_2,t;u,s,r)
\frac1t\mathsf{b}_r\left(\frac{y_2-v_2}{t}\right)\,dv_2\,\frac{dt}{t}
$$
when $\lambda_I(F)>0$, and $\mathsf{a}_I(F)=0$ otherwise.
These atoms have integral zero in $y_2$.
For $(v_2,t)\in\mathcal W_{\mathsf L}(I)$, the kernel vanishes unless
$|y_2-v_2|\le\mathsf L t<2|I|$.
Thus their supports in $y_2$ lie in the concentric interval $5I$.
To use Haar functions, we cover this interval by an interval in one of the two adjacent grids.
Apply the choice made in the proof of Lemma~\ref{l:adj-int-cover} to $5I$.
Since $|I|$ is a power of two, the selected scale satisfies
$$
15|I|<|Q(I)|\le30|I|,
\qquad\text{hence}\qquad |Q(I)|=16|I|.
$$
Fix the selected index $\varepsilon(I)\in\{0,1\}$ and interval $Q(I)$, so that
$$
5I\subset Q(I)\in\mathcal D^{\varepsilon(I)},
\qquad |Q(I)|=16|I|.
$$
These choices depend on $I$ only.
All intervals $I$ assigned to a fixed $Q$ have length $|Q|/16$, lie in $Q$,
and are mutually disjoint. Therefore
$$
\#\{I\in\mathcal D^0:Q(I)=Q\}\le16.
$$
For $\ell=0,1$, Cauchy--Schwarz and the normalization of $\mathsf{b}_r$ show that
\begin{align*}
\int_{\Gamma_\alpha(x_1)}\int_{-\alpha/2}^{\alpha/2}
|\partial_{y_2}^\ell\mathsf{a}_I(F)(y_2;u,s,r)|^2
\,dr\,\frac{du\,ds}{s^2}
\leq &\sup_{|r|<\alpha/2}\int_{\mathcal W_{\mathsf L}(I)}
t^{-2-2\ell}
\left|\left(\partial_{y_2}^\ell\mathsf{b}_r\right)
\left(\frac{y_2-v_2}{t}\right)\right|^2\,dv_2\,\frac{dt}{t}\\
\leq &\mathsf L^{-2-2\ell}|I|
\int_{|I|/\mathsf L}^{2|I|/\mathsf L}t^{-3-2\ell}\,dt\\
\le & C|I|^{-1-2\ell}
\le C|Q(I)|^{-1-2\ell}.
\end{align*}

\smallskip
\noindent{\it Step 2.}
For each $I$, we expand the atom in the Haar functions of the selected grid $\mathcal D^{\varepsilon(I)}$.
To estimate a finite sum, fix $\mathscr I\subset\mathcal D^0$ and write
$$
\mathscr I=\mathscr I_0\cup\mathscr I_1\quad\text{(disjoint)},
\qquad
\mathscr I_\varepsilon=\{I\in\mathscr I:\varepsilon(I)=\varepsilon\}.
$$
Within $\mathscr I_\varepsilon$, the supporting intervals $Q(I)$ belong to one grid;
we estimate the Haar sums at each fixed generation relative to these intervals.

For $J\in Ch_n(Q(I))$,  by Haar cancellation and the derivative bound  again, we obtain
$$
\bigg(\int_{\Gamma_\alpha(x_1)}\int_{-\alpha/2}^{\alpha/2}
|\langle\mathsf{a}_I(F)(\cdot;u,s,r),h_J\rangle|^2
\,dr\,\frac{du\,ds}{s^2}\bigg)^{1/2}
\le C\frac{|J|^{3/2}}{|Q(I)|^{3/2}}
=C2^{-3n/2}.
$$
As in the proof of Lemma~\ref{sp:l:smooth}, the Haar expansion of $\mathsf a_I(F)$ converges to it in
$L^2\big(dy_2\,\frac{du\,ds}{s^2}\,dr\big)$.
Since the intervals $J\in Ch_n(Q(I))$ for fixed $n$ are disjoint,
\begin{align*}
&\bigg(\int_{\Gamma_\alpha(x_1)}\int_{-\alpha/2}^{\alpha/2}
\bigg|\sum_{J\in Ch_n(Q(I))}
\langle\mathsf{a}_I(F)(\cdot;u,s,r),h_J\rangle h_J(y_2)\bigg|^2
\,dr\,\frac{du\,ds}{s^2}\bigg)^{1/2}\\
\leq &C2^{-3n/2}
\sum_{J\in Ch_n(Q(I))}
\frac{\mathbf1_J(y_2)}{|J|^{1/2}}
=C2^{-n}\frac{\mathbf1_{Q(I)}(y_2)}{|Q(I)|^{1/2}}.
\end{align*}
Thus
$$
\lim_{N\to\infty}
\bigg\|
\sum_{n=0}^{N}\sum_{J\in Ch_n(Q(I))}
\left\langle\mathsf a_I(F)(\cdot;u,s,r),h_J\right\rangle h_J(y_2)
-\mathsf a_I(F)(y_2;u,s,r)
\bigg\|_{L^p   (dy_2;L^2   (\frac{du\,ds}{s^2}\,dr  )  )}=0.
$$

The covering intervals also determine the aperture in the discrete square-function estimate.
For $y_2\in Q(I)$ and $(v_2,t)\in\mathcal W_{\mathsf L}(I)$, we have
$$
|v_2-y_2|\le|Q(I)|=16|I|<30\mathsf L t,
\qquad
\frac1{|Q(I)|}\le\frac1{|I|}\le\frac2{\mathsf L t}.
$$
Since the regions $\mathcal W_{\mathsf L}(I)$ are disjoint, it follows that
$$
\sum_{I\in\mathscr I}\lambda_I(F)^2
\frac{\mathbf1_{Q(I)}(y_2)}{|Q(I)|}
\le\frac2{\mathsf L}\mathcal A_{30\mathsf L}F(y_2)^2.
$$
To compare the aperture $30\mathsf L$ with $\alpha$, first note that when $30\mathsf L\le\alpha
$, then  $\mathcal A_{30\mathsf L}F\le\mathcal A_\alpha F$, and so it suffices to consider  the case $30\mathsf L>\alpha$.
For $p=2$, integrating first in $y_2$ gives
$$
\|\mathcal A_{30\mathsf L}F\|_{L^2(\R)}^2
=\frac{30\mathsf L}{\alpha}
\|\mathcal A_\alpha F\|_{L^2(\R)}^2.
$$

For $p>2$, we apply the weighted argument in the proof of
Lemma~\ref{sp:l:aperture} to the nonnegative function
$$
\int_{\Gamma_\alpha(x_1)}\int_{-\alpha/2}^{\alpha/2}
|F(v_2,t;u,s,r)|^2\,dr\,\frac{du\,ds}{s^2}.
$$
Only the aperture in $y_2$ changes; the integration domains in $(u,s,r)$ are fixed.
For $w\ge0$ with $\|w\|_{L^{(p/2)'}(\R)}=1$, the same calculation gives
\begin{align*}
\int_\R\mathcal A_{30\mathsf L}F(y_2)^2w(y_2)\,dy_2
\leq &\frac{30\mathsf L}{\alpha}
\int_\R\mathcal A_\alpha F(y_2)^2\mathsf Mw(y_2)\,dy_2\\
\leq &\frac{30\mathsf L}{\alpha}\|\mathsf Mw\|_{L^{(p/2)'}(\R)}
\|\mathcal A_\alpha F\|_{L^p(\R)}^2
\le Cp\frac{\mathsf L}{\alpha}\|\mathcal A_\alpha F\|_{L^p(\R)}^2,
\end{align*}
where $\mathsf M$ is the uncentred Hardy--Littlewood maximal operator on $\R$.
Taking the supremum over $w$ and then the square root, and including
the cases $p=2$ and $30\mathsf L\le\alpha$, we obtain
$$
\|\mathcal A_{30\mathsf L}F\|_{L^p(\R)}
\le C\sqrt p\,
\max\bigg \{1,\bigg (\frac{\mathsf L}{\alpha}\bigg )^{1/2}\bigg \}
\|\mathcal A_\alpha F\|_{L^p(\R)}.
$$

Fix $\varepsilon\in\{0,1\}$ and a relative generation $n$.
For a Haar interval $J\in\mathcal D^\varepsilon$, all atoms from $\mathscr I_\varepsilon$
contributing to its coefficient have $Q(I)$ equal to the unique $n$th ancestor of $J$ in this grid.
The multiplicity bound in {\it Step 1} therefore gives
$$
\#\{I\in\mathscr I_\varepsilon:J\in Ch_n(Q(I))\}\le16.
$$
We can now collect their coefficients and apply Cauchy--Schwarz as in the proof of Lemma~\ref{sp:l:smooth}.
Apply \eqref{e:haar-hilbert} to this sum in $\mathcal D^\varepsilon$, with the $L^2$ norm in $(u,s,r)$.
The triangle inequality over the two classes and then the sum over the relative generations give
\begin{align}\label{e:compact-recon-one}
&\bigg\|\bigg(\int_{\Gamma_\alpha(x_1)}\int_{-\alpha/2}^{\alpha/2}
\bigg|\sum_{I\in\mathscr I}\lambda_I(F)\mathsf{a}_I(F)(\cdot;u,s,r)\bigg|^2
\,dr\,\frac{du\,ds}{s^2}\bigg)^{1/2}\bigg\|_{L^p(\R)}\notag\\
\leq &C\sqrt p\sum_{n=0}^\infty2^{-n}
\bigg\|\bigg(\sum_{I\in\mathscr I}\lambda_I(F)^2
\frac{\mathbf1_{Q(I)}}{|Q(I)|}\bigg)^{1/2}\bigg\|_{L^p(\R)}\notag\\
\leq &C\sqrt p\,\mathsf L^{-1/2}\|\mathcal A_{30\mathsf L}F\|_{L^p(\R)}\notag\\
\leq &Cp\,\mathsf L^{-1/2}
\max\left\{1,\left(\frac{\mathsf L}{\alpha}\right)^{1/2}\right\}
\|\mathcal A_\alpha F\|_{L^p(\R)}\notag\\
= &Cp\max\{\mathsf L^{-1/2},\alpha^{-1/2}\}
\|\mathcal A_\alpha F\|_{L^p(\R)}.
\end{align}

\smallskip
\noindent{\it Step 3.}
We prove that  
$\sum_{I\in\mathscr I}\lambda_I(F)\mathsf a_I(F)$
in \eqref{e:compact-recon-one}    converges to $\mathcal R_{\mathsf b}F$
in $L^2\big (dy_2\,\frac{du\,ds}{s^2}\,dr\big )$
as the finite families $\mathscr I$ increase to $\mathcal D^0$.  For $\xi_2\ne0$,
\begin{align*}
|\widehat{\mathsf{b}}_r(\xi_2)|
\leq &\min\left\{
2\pi|\xi_2|\int_\R|y_2\mathsf{b}_r(y_2)|\,dy_2,\,
\frac{\|\partial_{y_2}\mathsf{b}_r\|_{L^1}}{2\pi|\xi_2|}
\right\}
\leq C\min\{\mathsf L|\xi_2|,(\mathsf L|\xi_2|)^{-1}\}.
\end{align*}
Hence
\begin{equation}\label{e:FT-b}
\sup_{\xi_2\ne0,\,   |r|<\alpha/2}
\int_0^\infty|\widehat{\mathsf{b}}_r(t\xi_2)|^2\,\frac{dt}{t}
\le C\int_0^\infty\min\{t,t^{-1}\}^2\,\frac{dt}{t}
=C.
\end{equation}
By \eqref{e:FT-b} and the duality argument in Lemma~\ref{l:atomic-l2},
the integrals defining $\mathcal R_{\mathsf b}F$, truncated to $N^{-1}\le t\le N$,
converge in $L^2(dy_2\,du\,ds/s^2\,dr)$, and
\begin{equation}\label{e:compact-l2}
\|\mathcal R_{\mathsf b}F\|_{L^2  (dy_2\,\frac{du\,ds}{s^2}\,dr)}
\le C\,\|F\|_{L^2  (dv_2\,\frac{dt}{t}\,\frac{du\,ds}{s^2}\,dr)},
\end{equation}
where the constant $C$ is independent of $\mathsf L$.

Choose finite families $\mathscr I_\ell\subseteq\mathscr I_{\ell+1}$
with $\bigcup_\ell\mathscr I_\ell=\mathcal D^0$.
Since the rectangles $\mathcal W_{\mathsf L}(I)$ partition $\R\times(0,\infty)$,
one may apply \eqref{e:compact-l2} to yield
$$
\bigg\|\mathcal R_{\mathsf b}F
-\sum_{I\in\mathscr I_\ell}\lambda_I(F)\mathsf a_I(F)
\bigg\|_{L^2   (dy_2\,\frac{du\,ds}{s^2}\,dr  )}
\le C\, \bigg(\sum_{I\notin\mathscr I_\ell}\lambda_I(F)^2\bigg)^{1/2}
\to0.
$$
Again, combining Fatou's lemma along a subsequence and \eqref{e:compact-recon-one}, we obtain
\eqref{e:compact-recon}, as in {\it Step 4} of
Proposition~\ref{prop:square-2}.
\end{proof}

\noindent{\bf{Summing the compact kernels.}}   \   \   
For $x_1,p,F$ as in Lemma~\ref{l:molecular-one}, take
$\mathsf b_r=\mathsf b_{k,r}$ and $\mathsf L=2^{k+1}$.
The bound in \eqref{e:FT-b} is independent of $k$, so \eqref{e:strip-compact} implies
\begin{align*}
&\sup_{|r|<\alpha/2,\,\xi_2\ne0}
\bigg(\int_0^\infty\bigg|
\widehat{\Phi_{\mathrm{hor}}^{(2)}(\cdot,r)}(t\xi_2)
-\sum_{k=0}^Nc_k\widehat{\mathsf b}_{k,r}(t\xi_2)
\bigg|^2\,\frac{dt}{t}\bigg)^{1/2}
\leq C\sum_{k>N}c_k\to0
\qquad(N\to\infty).
\end{align*}
The duality argument for \eqref{e:compact-l2}, applied to the full kernel and to these tails, gives
\begin{equation}\label{e:molecular-l2}
\|\mathcal R_{\Phi_{\mathrm{hor}}^{(2)}}F\|_{L^2  (dy_2\,\frac{du\,ds}{s^2}\,dr  )}
\le C_{\psi_2,\alpha}
\|F\|_{L^2   (dv_2\,\frac{dt}{t}\,\frac{du\,ds}{s^2}\,dr  )},
\end{equation}
and
\begin{equation}\label{e:converg-r}
\lim_{N\to\infty}
\bigg\|
\sum_{k=0}^N c_k\mathcal R_{\mathsf{b}_k}F
-\mathcal R_{\Phi_{\mathrm{hor}}^{(2)}}F
\bigg\|_{L^2  (dy_2\,\frac{du\,ds}{s^2}\,dr  )}=0.
\end{equation}
By \eqref{e:compact-recon} with $\mathsf L=2^{k+1}$,  
and Fatou's lemma along a subsequence in \eqref{e:converg-r}, we obtain
\begin{align}\label{e:molecular-one}
\bigg\|\bigg(\int_{\Gamma_\alpha(x_1)}\int_{-\alpha/2}^{\alpha/2}
|\mathcal R_{\Phi_{\mathrm{hor}}^{(2)}}F(\cdot;u,s,r)|^2
\,dr\,\frac{du\,ds}{s^2}\bigg)^{1/2}\bigg\|_{L^p(\R)}
\leq &Cp\sum_{k=0}^\infty c_k\max\{2^{-(k+1)/2},\alpha^{-1/2}\}
\|\mathcal A_\alpha F\|_{L^p(\R)}\notag\\
\leq &Cp\max\{1,\alpha^{-1/2}\}
\left(\sum_{k=0}^\infty c_k\right)
\|\mathcal A_\alpha F\|_{L^p(\R)}\notag\\
\leq &C_{\psi_2,\alpha}\max\{1,\alpha^{-1/2}\}\,p\,
\|\mathcal A_\alpha F\|_{L^p(\R)},
\end{align}
where the last inequality follows from \eqref{e:ck}.
The constant is independent of $p,x_1,F$.

\medskip

\noindent{\bf{The horizontal estimate.}}  \   \   
Now we return to $F$ defined in \eqref{e:coeff-h}. By \eqref{e:horizontal-recon},
$$
\mathcal R_{\mathrm{hor}}(f)(u,s;\pmb y)
=\int_{-\alpha/2}^{\alpha/2}
\big[\mathcal R_{\Phi_{\mathrm{hor}}^{(2)}}F\big](y_2;u,s,r)\,dr.
$$
The equality holds in $L^2\big(\Gamma_\alpha(x_1)\times\R^2,\frac{du\,ds}{s^2}\,d\pmb y\big)$
by \eqref{e:molecular-l2}, \eqref{e:strip-coeff} and the limit in \eqref{e:strip-recon-limit}.
It follows from Cauchy--Schwarz on $(-\alpha/2,\alpha/2)$ that
$$
\left|\mathcal R_{\mathrm{hor}}(f)(u,s;\pmb y)\right|^2
\le\alpha\int_{-\alpha/2}^{\alpha/2}
\left|[\mathcal R_{\Phi_{\mathrm{hor}}^{(2)}}F](y_2;u,s,r)\right|^2\,dr.
$$
Then by \eqref{e:molecular-one} and
\eqref{e:horizontal-area}, we obtain
\begin{align}\label{e:horizontal-bound}
&\bigg\|\bigg(\int_{\Gamma_\alpha(x_1)}
\left|\mathcal R_{\mathrm{hor}}(f)(u,s;\cdot)\right|^2
\,\frac{du\,ds}{s^2}\bigg)^{1/2}\bigg\|_{L^p(\R^2)}^p\notag\\
\leq &\alpha^{p/2}\int_\R
\bigg\|\bigg(\int_{\Gamma_\alpha(x_1)}\int_{-\alpha/2}^{\alpha/2}
\left|\mathcal R_{\Phi_{\mathrm{hor}}^{(2)}}F(\cdot;u,s,r) \right|^2
\,dr\,\frac{du\,ds}{s^2}\bigg)^{1/2}\bigg\|_{L^p(dy_2)}^p\,dy_3\notag\\
\leq &\left(C_{\psi_2,\alpha}p\max\left\{1,\alpha^{1/2}\right\}\right)^p
\int_\R\|\mathcal A_\alpha F\|_{L^p(dy_2)}^p\,dy_3\notag\\
\leq &(C_{\psi_2,\alpha}p)^p
\|S_{z,\alpha}^{\psi}f(x_1,\cdot)\|_{L^p(\R^2)}^p.
\end{align}

\subsection{The vertical term and the reverse estimate}\label{sub:vertical}

For the vertical strip, the scale in the long direction is $st$.
For fixed $s,y_2$, change variables from $(v_2,t)$ to $(r,\tau)$ by
$$
r=\frac{y_2-v_2}{t},
\qquad \tau=st,
\qquad
v_2=y_2-\frac{\tau r}{s},
\qquad t=\frac{\tau}{s}.
$$
Then
$$
\left|\frac{\partial(v_2,t)}{\partial(r,\tau)}\right|
=\frac{\tau}{s^2},
\qquad
\frac{dv_2\,dt}{st^3}
=\frac{dr\,d\tau}{\tau^2}.
$$
For each fixed $y_2$, define
$$
F(v_3,\tau;u,s,r)
:=(f*\psi_{s,\tau/s})(u,y_2-\tau r/s,v_3).
$$
We use the one-dimensional reconstruction in $y_3$ with scale $\tau$,
keeping $(u,s,r)$ in the $L^2$ norm.

Since $d\tau/\tau=dt/t$, it follows from translation in $y_2$ and \eqref{e:second-energy} that
$$
\int_\R\|F\|_{L^2  (dv_3\,\frac{d\tau}{\tau}\,\frac{du\,ds}{s^2}\,dr  )}^2\,dy_2
=\alpha\int_{\Gamma_\alpha(x_1)}\int_0^\infty\int_{\R^2}
|(f*\psi_{s,t})(u,\pmb v)|^2\,d\pmb v\,\frac{dt}{t}\,\frac{du\,ds}{s^2}<\infty.
$$

Now keep $y_2$ fixed.
Replacing $(y_2,v_2,t)$ by $(y_3,v_3,\tau)$ in
\eqref{e:compact-area}, we obtain
\begin{align*}
\mathcal A_\alpha F(y_3)^2
= &\int_{\Gamma_\alpha(x_1)}\int_{-\alpha/2}^{\alpha/2}
\int_0^\infty\int_{|v_3-y_3|<\alpha\tau}
|(f*\psi_{s,\tau/s})(u,y_2-\tau r/s,v_3)|^2
\,\frac{dv_3\,d\tau}{\tau^2}\,dr\,\frac{du\,ds}{s^2}\\
= &\int_{\Gamma_\alpha(x_1)}\int_0^\infty
\int_{\substack{|v_2-y_2|<\alpha t/2\\
|v_3-y_3|<\alpha st}}
|(f*\psi_{s,t})(u,\pmb v)|^2
\,\frac{d\pmb v\,dt}{st^3}\,\frac{du\,ds}{s^2}
\leq S_{z,\alpha}^{\psi}f(x_1,\pmb y)^2.
\end{align*}
Hence,  for almost every $y_2$,
$$
F\in L^2\Big (dv_3\,\frac{d\tau}{\tau}\,
\frac{du\,ds}{s^2}\,dr\Big ) 
\qquad\text{and}\qquad 
\mathcal A_\alpha F\in L^p(dy_3).
$$
Thus \eqref{e:molecular-one} also holds with
$(y_2,v_2,t)$ replaced by $(y_3,v_3,\tau)$.

{Under $\tau=st$, the truncation $N^{-1}\le t\le N$ gives
$F_N=F\mathbf1_{\{s/N\le\tau\le sN\}}$. For each $s,\tau>0$ the indicator tends to one.
The energy bound above gives $F_N\to F$ in $L^2$, also after integration in $y_2$.
By \eqref{e:molecular-l2} and Cauchy--Schwarz in $r$, the reconstructions
converge in $L^2$. Since $\Phi_{\mathrm{ver}}^{(2)}(r,\cdot)
=\Phi_{\mathrm{hor}}^{(2)}(\cdot,r)$, \eqref{e:strip-recon-limit} gives}
\begin{align*}
\mathcal R_{\mathrm{ver}}(f)(u,s;\pmb y)
= &\int_{-\alpha/2}^{\alpha/2}\int_0^\infty\int_\R
F(v_3,\tau;u,s,r)
\frac1{\tau}\Phi_{\mathrm{hor}}^{(2)}
\left(\frac{y_3-v_3}{\tau},r\right)
\,dv_3\,\frac{d\tau}{\tau}\,dr\\
= &\int_{-\alpha/2}^{\alpha/2}
[\mathcal R_{\Phi_{\mathrm{hor}}^{(2)}}F](y_3;u,s,r)\,dr.
\end{align*}
As in \eqref{e:horizontal-bound}, by Cauchy--Schwarz and \eqref{e:molecular-one}, we have
\begin{align}\label{e:vertical-bound}
&\bigg\|\bigg(\int_{\Gamma_\alpha(x_1)}
\left|\mathcal R_{\mathrm{ver}}(f)(u,s;\cdot)  \right|^2
\,\frac{du\,ds}{s^2}\bigg)^{1/2}\bigg\|_{L^p(\R^2)} 
\leq   C_{\psi_2,\alpha}p\,
\|S_{z,\alpha}^{\psi}f(x_1,\cdot)\|_{L^p(\R^2)}.
\end{align}

\medskip

\begin{proof}[Completion of the proof of Proposition~\ref{p:second-general}]
For $f\in L^2(\R^3)\cap L^p(\R^3)$ and almost every $x_1$,
by \eqref{e:strip-repro},
\eqref{e:horizontal-bound} and \eqref{e:vertical-bound}, we obtain 
\begin{align*}
\|S_{1,\alpha}^{\psi_1}f(x_1,\cdot)\|_{L^p(\R^2)}
=\, &\bigg\|\bigg(\int_{\Gamma_\alpha(x_1)}
\big|\mathcal R_{\mathrm{hor}}(f)(u,s;\cdot)
+\mathcal R_{\mathrm{ver}}(f)(u,s;\cdot) \big|^2
\,\frac{du\,ds}{s^2}\bigg)^{1/2}\bigg\|_{L^p(\R^2)}\\
\leq\, &\sum_{j\in\{\mathrm{hor},\mathrm{ver}\}}
\bigg\|\bigg(\int_{\Gamma_\alpha(x_1)}
|\mathcal R_j(f)(u,s;\cdot)|^2
\,\frac{du\,ds}{s^2}\bigg)^{1/2}\bigg\|_{L^p(\R^2)}\\
\leq \, &C_{\psi_2,\alpha}p\,
\|S_{z,\alpha}^{\psi}f(x_1,\cdot)\|_{L^p(\R^2)}.
\end{align*}
Integrating the $p$th powers in $x_1$, we obtain  
\eqref{e:second-general} for $f\in L^2(\R^3)\cap L^p(\R^3)$.

Using Lemma~\ref{l:continuous-forward}, we obtain
\eqref{e:second-general} for every $f\in L^p(\R^3)$
by approximation.
\end{proof}

\medskip

\begin{proof}[Proof of Theorem~\ref{t:second-free}]
Under $d_1(\psi_1)<c_*\alpha$, (i) of Proposition~\ref{p:support-distance}
and Proposition~\ref{p:second-general} give
$$
\|f\|_{L^p}
\le C_{\psi_1,\alpha}\sqrt p\,\|S_{1,\alpha}^{\psi_1}f\|_{L^p}
\le C_{\psi,\alpha}p^{3/2}\,\|S_{z,\alpha}^{\psi}f\|_{L^p}.
$$
Furthermore, without this condition, by Proposition~\ref{p:first-general},
$$
\|f\|_{L^p}
\le C_{\psi_1,\alpha}p\,\|S_{1,\alpha}^{\psi_1}f\|_{L^p}
\le C_{\psi,\alpha}p^2\,\|S_{z,\alpha}^{\psi}f\|_{L^p}. 
$$
The proof is complete.
\end{proof}

\begin{remark}
The general $p^2$ bound in (ii) of Theorem~\ref{t:second-free} also follows from
Theorem~\ref{t:large-ap} and a comparison of the two apertures.
Fix $0<\alpha\le\beta$, $2\le p<\infty$ and $f\in L^p(\R^3)$. For $y=(y_1,y_2,y_3)\in\R^3$, let
$$
B_\gamma(y;s,t)
=y+(-\gamma s,\gamma s)\times(-\gamma t,\gamma t)\times(-\gamma st,\gamma st),
\qquad \gamma\in\{\alpha,\beta\}.
$$
Let $M_Z$ be the uncentred maximal operator over axis-parallel boxes with side lengths $s,t,st$.
The box centred at $y$ with side lengths $c_\beta s,c_\beta t,c_\beta^2st$,
where $c_\beta=\max\{2\beta,\sqrt{2\beta}\}$, contains $B_\beta(y;s,t)$.
Thus, for $w\ge0$, averaging over $B_\alpha(y;s,t)$ yields
$$
\int_{B_\beta(y;s,t)}w(x)\,dx
\le\frac{c_\beta^4}{8\alpha^3}
\int_{B_\alpha(y;s,t)}M_Zw(x)\,dx.
$$
We use the following consequence of C\'ordoba's theorem
\cite[Theorem~(b), p.~256]{Cor}:
$$
\|M_Zw\|_{L^q(\R^3)}
\le C\max\{1,(q-1)^{-2}\}\|w\|_{L^q(\R^3)},
\qquad 1<q<\infty,
$$
{as shown in \eqref{e:mult-maximal}.}

For $p>2$, let $q=(p/2)'$ and take $w\ge0$ with
$\|w\|_{L^q}=1$.
Apply the preceding estimate inside the area integral, then use H\"older's inequality:
\begin{align*}
\int_{\R^3}(S_{z,\beta}^{\psi}f)^2w
\leq &C_{\alpha,\beta}
\int_{\R^3}(S_{z,\alpha}^{\psi}f)^2M_Zw\\
\leq &C_{\alpha,\beta}\|M_Zw\|_{L^q}
\|S_{z,\alpha}^{\psi}f\|_{L^p}^2\\
\leq &C_{\alpha,\beta}\max\{1,(q-1)^{-2}\}
\|S_{z,\alpha}^{\psi}f\|_{L^p}^2\\
\leq &C_{\alpha,\beta}p^2
\|S_{z,\alpha}^{\psi}f\|_{L^p}^2,
\end{align*}
since $q-1=2/(p-2)$.
For $p=2$, integration in the cone centre $x$ gives
$\|S_{z,\beta}^{\psi}f\|_2^2=(\beta/\alpha)^3\|S_{z,\alpha}^{\psi}f\|_2^2$.
Taking the supremum over $w$ and including $p=2$, we obtain
\begin{equation}\label{e:aperture-zyg}
\|S_{z,\beta}^{\psi}f\|_{L^p}
\le C_{\alpha,\beta}p\,\|S_{z,\alpha}^{\psi}f\|_{L^p},
\qquad 2\le p<\infty.
\end{equation}
For fixed $\psi_1,\psi_2$ and $\alpha$, choose $\beta\ge\alpha$ sufficiently large for
Theorem~\ref{t:large-ap}. Combining it with \eqref{e:aperture-zyg}, we obtain
$$
\|f\|_{L^p}
\le C_{\psi,\beta}p\,\|S_{z,\beta}^{\psi}f\|_{L^p}
\le C_{\psi,\alpha}p^2\,\|S_{z,\alpha}^{\psi}f\|_{L^p}.
$$
Here $\beta$ is independent of $p$.
This argument proves (ii) of Theorem~\ref{t:second-free}, but does not give
\eqref{e:second-general}, which is used for the $p^{3/2}$ bound in (i).
\end{remark}

\section{Sharpness and comparison of reverse constants}\label{sec:sharp-aperture}

\subsection{Kernel dependence at a fixed aperture}\label{sub:sharp-kernels}

We now prove the optimality of $p^{3/2}$ in Theorem~\ref{t:first-free}
and (i) of Theorem~\ref{t:second-free}, and of $p^2$ in
(ii) of Theorem~\ref{t:second-free}.
To this end, we prove Theorem~\ref{t:sharp-cases}
by constructing the three fixed kernel pairs separately.

For $g_1,g_2$ as in Lemma~\ref{l:lower-transfer}, we have
$$
\mathfrak C_p^\psi(\alpha)
\ge\frac1{\sqrt{2\alpha}}
\frac{\|g_1\|_{L^p}}{\|\mathcal A_{\psi_1,\alpha}g_1\|_{L^p}}
\frac{\|g_2\|_{L^p}}{\|\mathcal A_{\kappa,\alpha}g_2\|_{L^p}},
$$
where $\kappa(x_2)=\int_\R\psi_2(x_2,x_3)\,dx_3$ is given in \eqref{kappa}.
By Lemma~\ref{sp:l:lower} with $d=1$, we can choose $g_1,g_2$
so that each ratio is at least $c\sqrt p$.
Their product gives the lower bound of order $p$.
To obtain the lower bounds of order $p^{3/2}$ and $p^2$,
we will make one or both ratios at least $cp$, respectively.

For this purpose, we use
$$
h(u)=\vartheta(u)\log|u|,
$$
where $\vartheta\in C_c^\infty(\R)$ and $\vartheta\equiv1$ on $[-1,1]$.
Then $\|h\|_{L^p}\ge cp$.
 We will choose kernels whose convolution with $\log|\cdot|$ vanishes near zero,
and use this cancellation to bound the area-function norm independently of $p$.
The approximation in Lemma~\ref{sp:l:log-test} then gives the test functions
required by Lemma~\ref{l:lower-transfer}, with norm ratio at least $cp$.

\begin{lemma}\label{l:log-test}
Fix $0<a<A$ and a nonzero radial function
$0\le\zeta\in C_c^\infty(\{\pmb z\in\R^2:a<|\pmb z|<A\})$.
Let $\omega=\Delta\zeta$ and
\begin{equation}\label{kappa0}
\kappa_0(u)=\int_\R\omega(u,v)\,dv.
\end{equation}
Then $\omega\in C_c^\infty(\R^2)$ is nonzero, real-valued and radial,
and $\kappa_0\in C_c^\infty(\R)$ is nonzero, real-valued and even.
Moreover,
$$
\int_{\R^2}\omega=0,
\qquad
\int_\R\kappa_0=0,
$$
and
$$
d_1(\kappa_0)=d_2(\omega)\ge a.
$$
For every $0<\alpha<a$ and $2\le p<\infty$, there is
$g_p\in\mathcal S(\R)$ with
$\operatorname{supp}\widehat g_p\Subset\R\setminus\{0\}$
such that
\begin{equation}\label{e:log-test}
\|g_p\|_{L^p(\R)}\ge cp,
\qquad
0<\|\mathcal A_{\kappa_0,\alpha}g_p\|_{L^p(\R)}\le C.
\end{equation}
Here $c,C>0$ are independent of $p$.
\end{lemma}

The Fourier support condition on $g_p$ allows its use in Lemma~\ref{l:lower-transfer}.
We will apply it with $\psi_1=\kappa_0$, with $\psi_2=\omega$, or with both choices
in the proof of Theorem~\ref{t:sharp-cases} below.

\begin{proof}[Proof of Lemma~\ref{l:log-test}]
Apply Lemma~\ref{sp:l:log-kernel} with $m=1$.
For $d=2$, its kernel $\varphi$ is $\omega=\Delta\zeta$;
for $d=1$, its kernel is $\kappa_0$.
Thus both kernels satisfy the stated smoothness, symmetry and cancellation
conditions and are nonzero.
Since $\omega$ is supported in $\{a<|\pmb z|<A\}$,
Lemma~\ref{l:marginal-supp} gives
$$
d_1(\kappa_0)=d_2(\omega)\ge a.
$$
Lemma~\ref{sp:l:log-test} with $d=1$ and $\varphi=\kappa_0$
now provides $g_p\in\mathcal S_0(\R)$ satisfying \eqref{e:log-test}.
\end{proof}

As a consequence, using the test functions in Lemma~\ref{l:log-test},
we show that the support conditions in
Propositions~\ref{prop:square-1} and \ref{prop:square-2}
cannot be removed in general.

\begin{example}\label{r:counterexample-one}
Fix $0<a<A$, and take $\omega,\kappa_0$ from
Lemma~\ref{l:log-test}.
Let $\psi_1=\kappa_0$. Then $d_1(\psi_1)\ge a>0$, so
$0\notin\operatorname{supp}(|D|\psi_1)$.

For $0<\alpha<a$ and $2\le p<\infty$, choose $g_p$ as in
\eqref{e:log-test}, and let
$$
f_p(x)=g_p(x_1)b(x_2,x_3),
$$
where $0\ne b\in\mathcal S(\R^2)$ is fixed.
By the definition of $S_{1,\alpha}^{\psi_1}$,
$$
S_{1,\alpha}^{\psi_1}f_p(x)
=\mathcal A_{\kappa_0,\alpha}g_p(x_1)|b(x_2,x_3)|.
$$
Thus \eqref{e:log-test} gives
$$
\frac{\|f_p\|_{L^p(\R^3)}}
{\|S_{1,\alpha}^{\psi_1}f_p\|_{L^p(\R^3)}}
=\frac{\|g_p\|_{L^p(\R)}}
{\|\mathcal A_{\kappa_0,\alpha}g_p\|_{L^p(\R)}}
\ge c_{\psi_1,\alpha}p.
$$
Here $\psi_1$ is independent of $\alpha$ and $p$,
and $c_{\psi_1,\alpha}>0$ is independent of $p$.
\end{example}

\begin{example}\label{r:counterexample-two}
Keep $\omega,\kappa_0$ from Example~\ref{r:counterexample-one}.
Choose a nonzero real $\eta_0\in C_c^\infty(\R)$, and let
$$
\widetilde\eta_0(x_1)=\eta_0(-x_1),
\qquad
\psi_1=-(\eta_0*\widetilde\eta_0)'',
\qquad
\psi_2=\omega.
$$
Then $\psi_1\in C_c^\infty(\R)$ is real and even, with
$\int_\R\psi_1=0$, and
$$
(|D|\psi_1)(0)
=4\pi^2\int_\R|\xi_1|^3|\widehat\eta_0(\xi_1)|^2\,d\xi_1>0,
\qquad
d_2(\psi_2)\ge a>0.
$$
Thus $0\in\operatorname{supp}(|D|\psi_1)$, whereas
$0\notin\operatorname{supp}\psi_2$.

For $0<\alpha<a$ and $2\le p<\infty$, choose $g_p$ as in
\eqref{e:log-test}.
Fix nonzero $g,\chi\in\mathcal S(\R)$ with Fourier supports
compactly contained in $\R\setminus\{0\}$, and let
$$
f_{p,L}(x)=g(x_1)g_p(x_2)\chi(x_3/L),
\qquad L\ge1.
$$
By the definition of $S_{1,\alpha}^{\psi_1}$,
$$
L^{-1/p}\|S_{1,\alpha}^{\psi_1}f_{p,L}\|_{L^p(\R^3)}
=\|\mathcal A_{\psi_1,\alpha}g\|_{L^p(\R)}
\|g_p\|_{L^p(\R)}\|\chi\|_{L^p(\R)}.
$$
Since $\psi_2=\omega$, \eqref{kappa} and \eqref{kappa0} give
$\kappa=\kappa_0$.
Combined with Lemma~\ref{l:lower-transfer} and \eqref{e:log-test},
we obtain
$$
\lim_{L\to\infty}
\frac{\|S_{1,\alpha}^{\psi_1}f_{p,L}\|_{L^p(\R^3)}}
{\|S_{z,\alpha}^{\psi}f_{p,L}\|_{L^p(\R^3)}}
=\frac{\|g_p\|_{L^p(\R)}}
{\sqrt{2\alpha}\|\mathcal A_{\kappa_0,\alpha}g_p\|_{L^p(\R)}}
\ge 2c_{\psi,\alpha}p.
$$
For each $p$, choose $L_p\ge1$ large enough that
$f_p=f_{p,L_p}$ satisfies
$$
\|S_{1,\alpha}^{\psi_1}f_p\|_{L^p(\R^3)}
\ge c_{\psi,\alpha}p\,
\|S_{z,\alpha}^{\psi}f_p\|_{L^p(\R^3)}.
$$
Here $\psi_1,\psi_2$ are independent of $\alpha$ and $p$,
and $c_{\psi,\alpha}>0$ is independent of $p$.
\end{example}

\smallskip

\begin{proof}[Proof of Theorem~\ref{t:sharp-cases}]
Fix $A>a$, and take $\omega,\kappa_0$ from Lemma~\ref{l:log-test}.
To obtain kernels satisfying the support conditions, choose a nonzero
real radial $\eta\in C_c^\infty(\R^2)$, and let
$$
\gamma_2=-\Delta(\eta*\eta),
\qquad
\gamma_1(u)=\int_\R\gamma_2(u,v)\,dv.
$$
Both kernels are real, smooth, compactly supported and have integral zero.
Since
$$
\gamma_2(0)=\int_{\R^2}|\nabla\eta|^2>0,
$$
radiality and Lemma~\ref{l:marginal-supp} show that $\gamma_1\ne0$ and
$$
d_1(\gamma_1)=d_2(\gamma_2)=0.
$$

Fix $0<\alpha<a$ and $2\le p<\infty$.
For $g_1,g_2$ as in Lemma~\ref{l:lower-transfer}, that lemma gives
\begin{equation}\label{e:product-test}
\mathfrak C_p^\psi(\alpha)
\ge\frac1{\sqrt{2\alpha}}
\frac{\|g_1\|_{L^p}}{\|\mathcal A_{\psi_1,\alpha}g_1\|_{L^p}}
\frac{\|g_2\|_{L^p}}{\|\mathcal A_{\kappa,\alpha}g_2\|_{L^p}},
\qquad
\kappa(x_2)=\int_\R\psi_2(x_2,x_3)\,dx_3.
\end{equation}

\smallskip
\noindent(i).
Take $(\psi_1,\psi_2)=(\kappa_0,\gamma_2)$. Then
$$
d_1(\psi_1)\ge a,\qquad d_2(\psi_2)=0,\qquad \kappa=\gamma_1.
$$
Choose $g_1$ by Lemma~\ref{l:log-test} and $g_2$ by
Lemma~\ref{sp:l:lower} with $d=1$ and $\varphi=\gamma_1$.
By \eqref{e:product-test} and Theorem~\ref{t:first-free},
$$
c_{\psi,\alpha}p\sqrt p
\le\mathfrak C_p^\psi(\alpha)
\le C_{\psi,\alpha}p^{3/2}.
$$

\smallskip
\noindent(ii).
Take $(\psi_1,\psi_2)=(\gamma_1,\omega)$. Then
$$
d_1(\psi_1)=0,\qquad d_2(\psi_2)\ge a,\qquad \kappa=\kappa_0.
$$
Choose $g_1$ by Lemma~\ref{sp:l:lower} with $d=1$ and $\varphi=\gamma_1$,
and $g_2$ by Lemma~\ref{l:log-test}.
By \eqref{e:product-test} and (i) of Theorem~\ref{t:second-free},
$$
c_{\psi,\alpha}\sqrt p\,p
\le\mathfrak C_p^\psi(\alpha)
\le C_{\psi,\alpha}p^{3/2}.
$$

\smallskip
\noindent(iii).
Take $(\psi_1,\psi_2)=(\kappa_0,\omega)$. Then
$$
d_1(\psi_1)=d_2(\psi_2)\ge a,\qquad \kappa=\kappa_0.
$$
Choose $g_1=g_2=g_p$ from Lemma~\ref{l:log-test}.
By \eqref{e:product-test} and (ii) of Theorem~\ref{t:second-free},
$$
c_{\psi,\alpha}p^2
\le\mathfrak C_p^\psi(\alpha)
\le C_{\psi,\alpha}p^2.
$$
All three kernel pairs are independent of $\alpha$ and $p$,
and the positive constants are independent of $p$.
\end{proof}

For any fixed aperture $\alpha>0$, choose $a>\alpha$ in the preceding construction.
Then $(\gamma_1,\gamma_2)$ has optimal order $p$ by Theorem~\ref{t:primal},
while $(\kappa_0,\gamma_2)$ and $(\gamma_1,\omega)$ have optimal order $p^{3/2}$,
and $(\kappa_0,\omega)$ has optimal order $p^2$.
Thus changing either kernel can change the optimal order while the aperture remains fixed.

\subsection{Aperture dependence for a fixed kernel pair}

Combining Theorem~\ref{t:large-ap} with
(iii) of Theorem~\ref{t:sharp-cases}, we obtain a fixed kernel pair
for which the optimal order is $p^2$ at small apertures
and $p$ at sufficiently large apertures.
Thus the optimal order depends on the aperture even when the kernels are fixed.

\begin{corollary}\label{c:aperture-change}
There is a fixed kernel pair $\psi_1,\psi_2$ and numbers
$0<\alpha_-<\alpha_+$ such that, for every $2\le p<\infty$,
$$
c_{\psi,\alpha}p^2\le\mathfrak C_p^\psi(\alpha)
\le C_{\psi,\alpha}p^2
\qquad( 0<\alpha\le\alpha_-),
$$
and
$$
c_{\psi,\alpha}p\le\mathfrak C_p^\psi(\alpha)
\le C_{\psi,\alpha}p
\qquad(\alpha\ge\alpha_+).
$$
The positive constants are independent of $p$.
\end{corollary}

\begin{proof}
Choose the pair in (iii) of Theorem~\ref{t:sharp-cases}, and let
$$
\alpha_-=a/2,
\qquad \alpha_+=\max\{2a,\alpha_\psi\},
$$
where $\alpha_\psi$ is from Theorem~\ref{t:large-ap}.
Since $\alpha_-<a$, (iii) of Theorem~\ref{t:sharp-cases} gives
the first estimate for $0<\alpha\le\alpha_-$.
The second estimate follows from Theorem~\ref{t:large-ap},
since $\alpha_+\ge\alpha_\psi$.
\end{proof}

\begin{remark}\label{r:aperture-sharp}
The factor $\sqrt p$ in Lemma~\ref{l:aperture-one} is necessary in general,
even for a fixed kernel and two fixed apertures.
Take $\kappa_0$ from Lemma~\ref{l:log-test}, fix $0<\alpha<a$, and choose
$$
\beta>\max\{\alpha,d_1(\kappa_0)/c_*\}.
$$
Since $c_*<c_1$, Corollary~\ref{sp:c:endpoint} with $d=1$ applies.
For its fixed test function $h$, it gives
$$
\frac{\|\mathcal A_{\kappa_0,\beta}h\|_{L^p(\R)}}
{\|\mathcal A_{\kappa_0,\alpha}h\|_{L^p(\R)}}\ge c\sqrt p,
\qquad
\mathcal A_{\kappa_0,\alpha}h\in L^\infty(\R),
\qquad
\mathcal A_{\kappa_0,\beta}h\notin L^\infty(\R).
$$
For any nonzero $\chi\in C_c^\infty(\R^2)$, let
$f(x_1,\pmb x)=h(x_1)\chi(\pmb x)$.
Since
$$
S_{1,\gamma}^{\kappa_0}f(x_1,\pmb x)
=|\chi(\pmb x)|\,\mathcal A_{\kappa_0,\gamma}h(x_1),
\qquad \gamma\in\{\alpha,\beta\},
$$
the same norm ratio proves the sharpness of \eqref{e:aperture-one}.
\end{remark}
\subsection{Comparison with one-parameter constants}

Let $\psi_1\in C_c^\infty(\R)$ and $\psi_2\in C_c^\infty(\R^2)$
satisfy the basic assumptions.
The intermediate square function $S_{1,\alpha}^{\psi_1}f$ allows us
to estimate $\|f\|_{L^p}$ by $\|S_{z,\alpha}^{\psi}f\|_{L^p}$ in two steps.
By Lemma~\ref{l:lower-transfer}, we obtained lower bounds from
one-dimensional test functions for the kernels $\psi_1$ and $\kappa$,
where
$$
\kappa(u)=\int_\R\psi_2(u,v)\,dv.
$$

We use the one-dimensional best constants $\mathfrak c_{p,1}^{\varphi_1}(\alpha)$
from \eqref{sp:e:best}.
In particular, \eqref{e:first-from-single} applies with $\varphi_1=\psi_1$.

Theorems~\ref{sp:t:general}--\ref{sp:t:aperture} describe the one-parameter bounds
and their dependence on the kernel and the aperture.
We now compare
$\mathfrak C_p^\psi(\alpha)$, defined in \eqref{e:reverse-constant},
with $\mathfrak c_{p,1}^{\psi_1}(\alpha)\mathfrak c_{p,1}^\kappa(\alpha)$.

\begin{proposition}\label{p:constants}
Let $\psi_1,\psi_2$ satisfy the basic assumptions.
For every $\alpha>0$ and $2\le p<\infty$,
\begin{equation}\label{e:constants-general}
\frac{\mathfrak c_{p,1}^{\psi_1}(\alpha)\mathfrak c_{p,1}^\kappa(\alpha)}
{\sqrt{2\alpha}}
\le\mathfrak C_p^\psi(\alpha)
\le C_{\psi_2,\alpha}p\,\mathfrak c_{p,1}^{\psi_1}(\alpha).
\end{equation}
Moreover,
\begin{itemize}
\item[(i)] if $d_2(\psi_2)<c_*\alpha$, then
\begin{equation}\label{e:constants-local}
\mathfrak C_p^\psi(\alpha)
\simeq_{\psi,\alpha}
\mathfrak c_{p,1}^{\psi_1}(\alpha)\mathfrak c_{p,1}^\kappa(\alpha)
\simeq_{\psi,\alpha}\sqrt p\,\mathfrak c_{p,1}^{\psi_1}(\alpha);
\end{equation}

\item[(ii)] if $0<\alpha\le d_2(\psi_2)$ and
$$
\int_{\R^2}\psi_2(\pmb y)\log|\pmb y|\,d\pmb y=0,
$$
then
\begin{equation}\label{e:constants-gap}
\mathfrak C_p^\psi(\alpha)
\simeq_{\psi,\alpha}
\mathfrak c_{p,1}^{\psi_1}(\alpha)\mathfrak c_{p,1}^\kappa(\alpha)
\simeq_{\psi,\alpha}p\,\mathfrak c_{p,1}^{\psi_1}(\alpha).
\end{equation}
\end{itemize}
All comparison constants are independent of $p$.
\end{proposition}

\begin{proof}
Fix $\alpha>0$ and $2\le p<\infty$.
By the approximation following Lemma~\ref{sp:l:continuity}, the supremum
in $\mathfrak c_{p,1}^{\varphi_1}(\alpha)$ can be taken over
$0\ne g\in\mathcal S_0(\R)$.

For such $g_1,g_2$, Lemma~\ref{l:lower-transfer} implies
$$
\mathfrak C_p^\psi(\alpha)
\ge\frac1{\sqrt{2\alpha}}
\frac{\|g_1\|_{L^p(\R)}}
{\|\mathcal A_{\psi_1,\alpha}g_1\|_{L^p(\R)}}
\frac{\|g_2\|_{L^p(\R)}}
{\|\mathcal A_{\kappa,\alpha}g_2\|_{L^p(\R)}}.
$$
Taking the supremum over $g_1$ and $g_2$ separately proves
the lower bound in \eqref{e:constants-general}.

For the upper bound, apply the definition of
$\mathfrak c_{p,1}^{\psi_1}(\alpha)$ to $f(\cdot,x_2,x_3)$
and then use Proposition~\ref{p:second-general}:
\begin{align*}
\|f\|_{L^p(\R^3)}
\leq &\mathfrak c_{p,1}^{\psi_1}(\alpha)
\|S_{1,\alpha}^{\psi_1}f\|_{L^p(\R^3)} 
\leq  C_{\psi_2,\alpha}p\,\mathfrak c_{p,1}^{\psi_1}(\alpha)
\|S_{z,\alpha}^{\psi}f\|_{L^p(\R^3)}.
\end{align*}
This  proves  \eqref{e:constants-general}.

\smallskip

\noindent(i). Under the assumption $d_2(\psi_2)<c_*\alpha$, Lemma~\ref{l:marginal-supp} gives
$$
d_1(\kappa)=d_2(\psi_2)<c_*\alpha.
$$
Since $c_*<c_1$, Theorem~\ref{sp:t:general} with $d=1$ gives
$$
\mathfrak c_{p,1}^\kappa(\alpha)
\simeq_{\psi_2,\alpha}\sqrt p.
$$
By (ii) of Proposition~\ref{p:support-distance},
$$
\|f\|_{L^p(\R^3)}
\le\mathfrak c_{p,1}^{\psi_1}(\alpha)
\|S_{1,\alpha}^{\psi_1}f\|_{L^p(\R^3)}
\le C_{\psi_2,\alpha}\sqrt p\,
\mathfrak c_{p,1}^{\psi_1}(\alpha)
\|S_{z,\alpha}^{\psi}f\|_{L^p(\R^3)}.
$$
Combining this estimate with \eqref{e:constants-general},
we obtain \eqref{e:constants-local}.

\smallskip
\noindent(ii).
To obtain the lower bound of order $p$ for $\mathfrak c_{p,1}^\kappa(\alpha)$,
we again use the logarithmic test. Let
$$
U(u)=\int_\R\kappa(v)\log|u-v|\,dv.
$$
Since $d_2(\psi_2)\ge\alpha$, it follows from radiality and \eqref{e:log}
that, for $|u|<\alpha$,
\begin{align*}
U(u)
= &2\pi\int_\alpha^\infty
\psi_2(r,0)r\log(r/2)\,dr
= \int_{\R^2}\psi_2(\pmb y)\log|\pmb y|\,d\pmb y
-(\log2)\int_{\R^2}\psi_2(\pmb y)\,d\pmb y
=0.
\end{align*}

Choose $\vartheta$ as in {\it Step 1} of Lemma~\ref{sp:l:log-test}, and let
$$
h(u)=\vartheta(u)\log|u|.
$$
Since $\kappa\in C_c^\infty(\R)$ is even, has integral zero and satisfies
$U=0$ on $(-\alpha,\alpha)$, the estimate \eqref{sp:e:one-log-area}
from that step applies to the present kernel and gives
$$
\|\mathcal A_{\kappa,\alpha}h\|_{L^\infty}\le C_{\psi_2,\alpha}.
$$
By the $L^2$ identity \eqref{e:ID}, with $\psi_1=\kappa$,
and interpolation,
$$
0<\|\mathcal A_{\kappa,\alpha}h\|_{L^p}
\le
\|\mathcal A_{\kappa,\alpha}h\|_{L^\infty}^{1-2/p}
\bigl(2\alpha c_\kappa\|h\|_{L^2}^2\bigr)^{1/p}
\le C_{\psi_2,\alpha},
\qquad 2\le p<\infty.
$$
Since $\|h\|_{L^p}\ge cp$ by \eqref{sp:e:log-size},
the general one-dimensional upper bound implies
$$
c_{\psi_2,\alpha}p
\le
\frac{\|h\|_{L^p}}
{\|\mathcal A_{\kappa,\alpha}h\|_{L^p}}
\le\mathfrak c_{p,1}^\kappa(\alpha)
\le C_{\psi_2,\alpha}p.
$$
{This includes $\alpha=d_2(\psi_2)$, since \eqref{sp:e:one-log-area}
requires $U=0$ only on the open interval $(-\alpha,\alpha)$.}
Combining these bounds with \eqref{e:constants-general}
proves \eqref{e:constants-gap}.
\end{proof}

\section{Poisson and heat area functions for Zygmund dilations}\label{sec:semigroups}

Here we prove Theorem~\ref{t:semigroup} for the Poisson and heat
area functions under the coupled scales $s,t,st$.
We use the kernels $\varphi_\nu^{(d)}$ and the localized dual kernels
from Subsection~\ref{sp:sub:semigroups}, with $d=1,2$.
Proposition~\ref{sp:p:semigroups} supplies the first-variable reverse estimate
and the one-dimensional test functions.
The remaining arguments concern the coupled scales in the last two variables.
We first consider the Poisson square function.

In dimensions one and two, differentiating the Poisson kernel gives
$$
\varphi_{\mathrm P}^{(1)}(x_1)
=\frac{2(1-4\pi^2x_1^2)}{(1+4\pi^2x_1^2)^2},
$$
and
$$
\varphi_{\mathrm P}^{(2)}(\pmb y)
=\frac{2\pi[2-4\pi^2(y_2^2+y_3^2)]}{[1+4\pi^2(y_2^2+y_3^2)]^{5/2}},
\qquad \pmb y=(y_2,y_3).
$$

As before, for $s,t>0$, let
$$
\varphi_{\mathrm P,s}^{(1)}(x_1)
=\frac1s\varphi_{\mathrm P}^{(1)}\left(\frac{x_1}{s}\right),
\qquad
\varphi_{\mathrm P,s,t}^{(2)}(y_2,y_3)
=\frac1{st^2}\varphi_{\mathrm P}^{(2)}
\left(\frac{y_2}{t},\frac{y_3}{st}\right).
$$
For $f\in\mathcal S(\R^3)$, the Poisson operators in
Theorem~\ref{t:semigroup} satisfy
\begin{equation}\label{e:Mb}
\widehat{Q_{s,t}^{\mathrm P}f}(\xi):=[M_{\mathrm P}(\xi)](s,t)\widehat f(\xi)
=b_{\mathrm P}(s|\xi_1|)b_{\mathrm P}\big (t\sqrt{\xi_2^2+s^2\xi_3^2}\big  )\widehat f(\xi).
\end{equation}
Fourier inversion then gives
$$
Q_{s,t}^{\mathrm P}f=f*K_{s,t}^{\mathrm P},
\qquad\text{where}\  \   
K_{s,t}^{\mathrm P}(x)=\frac1{s^2t^2}
\varphi_{\mathrm P}^{(1)}\left(\frac{x_1}{s}\right)
\varphi_{\mathrm P}^{(2)}\left(\frac{x_2}{t},\frac{x_3}{st}\right).
$$
Thus the three spatial scales are $s,t,st$.
A change of variables yields
$$
\|K_{s,t}^{\mathrm P}\|_{L^1(\R^3)}
=\|\varphi_{\mathrm P}^{(1)}\|_{L^1(\R)}
\|\varphi_{\mathrm P}^{(2)}\|_{L^1(\R^2)}<\infty.
$$
The convolution formula therefore defines $Q_{s,t}^{\mathrm P}f$
for every $f\in L^p(\R^3)$, $1\le p\le\infty$.
Recall that the Poisson area function in \eqref{e:sg-area} is
$$
S_{z,\alpha}^{\mathrm P}f(x)^2
=\int_0^\infty\int_0^\infty
\int_{\substack{|u-x_1|<\alpha s\\
|v_2-x_2|<\alpha t\\|v_3-x_3|<\alpha st}}
|Q_{s,t}^{\mathrm P}f(u,v_2,v_3)|^2
\,\frac{du\,dv_2\,dv_3\,dt\,ds}{s^3t^3}.
$$

 The upper $L^p$ estimates below will allow us to extend the results from
$L^2\cap L^p$ to $L^p$. 
 We will use compactly supported dual kernels and the intermediate
square function to prove the two $\sqrt p$ estimates. 
 The Poisson semigroup derivative kernels are not compactly supported.
To transfer the one-parameter sharpness from Proposition~\ref{sp:p:semigroups},
we will verify the bounds involving the coupled scales in Subsection~\ref{sub:sg-sharp}. 
 For the heat case, we will give the definitions and the parts of the proof that differ
in Subsection~\ref{sub:sg-heat}. 

\subsection{Upper \texorpdfstring{$L^p$}{Lp} estimates for the Poisson square function}\label{sub:sg-forward}

For fixed $(x_2,x_3)$, the intermediate square function is
$$
S_{1,\alpha}^{\varphi_{\mathrm P}^{(1)}}f(x_1,x_2,x_3)^2
=\int_0^\infty\int_{|u-x_1|<\alpha s}
\left|(f(\cdot,x_2,x_3)*\varphi_{\mathrm P,s}^{(1)})(u)\right|^2
\,\frac{du\,ds}{s^2}.
$$

\begin{lemma}\label{l:sg-forward}
For every $\alpha>0$, $2\le p<\infty$, and $f\in L^p(\R^3)$,
$$
\big \|S_{1,\alpha}^{\varphi_{\mathrm P}^{(1)}}f\big \|_{L^p(\R^3)}
+\big \|S_{z,\alpha}^{\mathrm P}f  \big \|_{L^p(\R^3)}
\le B_{\mathrm P,\alpha,p}\|f\|_{L^p(\R^3)}.
$$
Both square functions are continuous under $L^p$ convergence.
\end{lemma}

\begin{proof}

Fix $2\le p<\infty$.
In {\it Steps 1 and 2}, assume
$f\in L^2(\R^3)\cap L^p(\R^3)$.
 We will prove the bounds for every $f\in L^p(\R^3)$
and prove continuity under $L^p$ convergence  in {\it Step 3}. 
Recall $c_{\mathrm P}=1/4$ from \eqref{sp:e:sg-energy-constants}.

\smallskip
\noindent{\it Step 1.} We begin by showing the upper $L^p$  bounds for $S_{z,\alpha}^{\mathrm P}$, $2\leq p<\infty$.

\textbf{Case 1. $p=2$.}\    \   
For $\xi_1\ne0$ and $(\xi_2,\xi_3)\ne(0,0)$, substituting  
$v=t\sqrt{\xi_2^2+s^2\xi_3^2}$ and then $r=s|\xi_1|$, we have 
$$
\int_0^\infty\int_0^\infty|[M_{\mathrm P}(\xi)](s,t)|^2\,\frac{dt}{t}\frac{ds}{s}
=c_{\mathrm P}\int_0^\infty b_{\mathrm P}(s|\xi_1|)^2\,\frac{ds}{s}=c_{\mathrm P}^2,
$$
where $[M_{\mathrm P}(\xi)](s,t)\widehat{f}(\xi)=\widehat{Q_{s,t}^{\mathrm P}f}(\xi)$, as defined in \eqref{e:Mb}.
By Plancherel and Tonelli's theorem, for every $f\in L^2(\R^3)$,
\begin{align}\label{e:sg-energy}
\int_0^\infty\int_0^\infty\|Q_{s,t}^{\mathrm P}f\|_{L^2(\R^3)}^2
\,\frac{dt}{t}\frac{ds}{s}
=c_{\mathrm P}^2\|f\|_{L^2(\R^3)}^2\qquad\text{and}\qquad
\big \|S_{z,\alpha}^{\mathrm P}f  \big \|_{L^2(\R^3)}^2
=(2\alpha)^3c_{\mathrm P}^2\|f\|_{L^2(\R^3)}^2.
\end{align}
 
\smallskip

\textbf{Case 2. $2<p<\infty$.}\   \  
We first prove
\begin{equation}\label{e:sg-vertical}
\bigg\|\bigg(
\int_0^\infty\int_0^\infty
|Q_{s,t}^{\mathrm P}f(\cdot)|^2\,\frac{ds}{s}\frac{dt}{t}
\bigg)^{1/2}\bigg\|_{L^p(\R^3)}
\le C_{p,\mathrm P}\|f\|_{L^p(\R^3)}.
\end{equation}
Note that on the left-hand side, we take the $L^2$ norm in $(s,t)$
 and then the $L^p$ norm in $x$.
We will use the Hilbert-space case of 
\cite[Theorem~4.5 and Remark~4.6 ~ii)]{SW}
with the multiplier
$$
M_{\mathrm P}:\xi\longmapsto M_{\mathrm P}(\xi)\in
L^2\bigg((0,\infty)^2,\frac{ds}{s}\frac{dt}{t}\bigg).
$$
We start by showing the required derivative bounds for $M_{\mathrm P}$.
Let
$$
r=s|\xi_1|,
\qquad v=t\sqrt{\xi_2^2+s^2\xi_3^2},
\qquad D_j=\xi_j\partial_{\xi_j}\quad(j=1,2,3).
$$
For $\xi_1\xi_2\xi_3\ne0$, direct computation gives
\begin{itemize}
	\item $\displaystyle [(D_1M_{\mathrm P})(\xi)](s,t)  =rb_{\mathrm P}'(r)b_{\mathrm P}(v)$,
	
	\medskip
	
	\item $\displaystyle [(D_2M_{\mathrm P})(\xi)](s,t)  =b_{\mathrm P}(r)\frac{\xi_2^2}{\xi_2^2+s^2\xi_3^2}\,vb_{\mathrm P}'(v) $,
	
	\medskip
	
	\item $\displaystyle  [(D_3M_{\mathrm P})(\xi)](s,t)  =b_{\mathrm P}(r)\frac{s^2\xi_3^2}{\xi_2^2+s^2\xi_3^2}\,vb_{\mathrm P}'(v)$,
	
	\medskip
	
	\item $\displaystyle [(D_2D_3M_{\mathrm P})(\xi)](s,t)  =b_{\mathrm P}(r)\frac{\xi_2^2s^2\xi_3^2}{(\xi_2^2+s^2\xi_3^2)^2}
\left (v^2b_{\mathrm P}''(v)-vb_{\mathrm P}'(v)  \right)$.
\end{itemize}
In the last three formulas, only $b_{\mathrm P}(r)$ depends on $\xi_1$.
Since $D_1b_{\mathrm P}(r)=rb_{\mathrm P}'(r)$, we obtain the formulas for
$D_1D_2M_{\mathrm P}$, $D_1D_3M_{\mathrm P}$ and $D_1D_2D_3M_{\mathrm P}$
by replacing $b_{\mathrm P}(r)$ with $rb_{\mathrm P}'(r)$.
All three fractions lie in $[0,1]$.

By $dt/t=dv/v$ first and then $ds/s=dr/r$, we obtain
\begin{align}\label{e:sg-mixed}
&\max_{\varepsilon_j\in\{0,1\}}\sup_{\xi_1\xi_2\xi_3\ne0}
\int_0^\infty\int_0^\infty
\left|\left[(D_1^{\varepsilon_1}D_2^{\varepsilon_2}D_3^{\varepsilon_3}M_{\mathrm P})(\xi)\right](s,t)\right|^2
\,\frac{dt}{t}\frac{ds}{s}\notag\\
\leq&\ \bigg(\int_0^\infty\big (b_{\mathrm P}(r)^2+|rb_{\mathrm P}'(r)|^2\big  )\,\frac{dr}{r}\bigg)\,
\bigg(\int_0^\infty
\big (b_{\mathrm P}(v)^2+|vb_{\mathrm P}'(v)|^2+|v^2b_{\mathrm P}''(v)-vb_{\mathrm P}'(v)|^2\big  )\,\frac{dv}{v}\bigg)<\infty.
\end{align}
The last integrals are finite, since $b_{\mathrm P}(r)=re^{-r}$ and
$$
rb_{\mathrm P}'(r)=r(1-r)e^{-r},
\qquad r^2b_{\mathrm P}''(r)-rb_{\mathrm P}'(r)=r(r^2-r-1)e^{-r}.
$$

To obtain a smooth multiplier on $\R^3$, we multiply $M_{\mathrm P}$
by a smooth function that is zero near $\xi_j=0$, $j=1,2,3$.
Choose $0\le\chi\le1$ in $C_c^\infty(\R)$ with $\chi=1$ on $[-1,1]$, and let
$$
\theta_N(\xi)=\prod_{j=1}^3\left (1-\chi(N\xi_j)\right ),
\qquad N\ge2.
$$
Then $\theta_NM_{\mathrm P}$ is smooth on $\R^3$.
For $\varepsilon\in\{0,1\}^3$,
$$
\big|\xi^\varepsilon\partial^\varepsilon\theta_N(\xi)\big|
=\prod_{\varepsilon_j=0}\left (1-\chi(N\xi_j)\right )
\prod_{\varepsilon_j=1}  \left |N\xi_j\,\chi'(N\xi_j)\right |
\le C_\chi,
$$
uniformly in $N,\xi$.

For $\xi_1\xi_2\xi_3\ne0$, by \eqref{e:sg-mixed}, we have 
\begin{align}\label{e:n0}
\big\|\xi^\varepsilon\partial^\varepsilon(\theta_NM_{\mathrm P})(\xi)
\big\|_{L^2\left(\frac{ds}{s}\frac{dt}{t}\right)}
= &\bigg\|\sum_{\delta\le\varepsilon}
\left (\xi^\delta\partial^\delta\theta_N(\xi)\right)
\left (\xi^{\varepsilon-\delta}
\partial^{\varepsilon-\delta}M_{\mathrm P}(\xi)\right)
\bigg\|_{L^2\left(\frac{ds}{s}\frac{dt}{t}\right)} \notag\\
\leq &\sum_{\delta\le\varepsilon}
\left |\xi^\delta\partial^\delta\theta_N(\xi)  \right|
\left \|(D_1^{\varepsilon_1-\delta_1}
D_2^{\varepsilon_2-\delta_2}
D_3^{\varepsilon_3-\delta_3}M_{\mathrm P})(\xi)
\right \|_{L^2\left(\frac{ds}{s}\frac{dt}{t}\right)}
\le C_{\chi,\mathrm P}.
\end{align}
Since $\chi=1$ on $[-1,1]$,
$$
(\theta_NM_{\mathrm P})(\xi)=0,
\qquad \min_{1\le j\le3}|\xi_j|<N^{-1}.
$$
Thus its derivatives also vanish near the coordinate planes.
Together with \eqref{e:n0}, this gives
$$
\sup_{N\ge2}\sup_{\xi\in\R^3}
\big\|\xi^\varepsilon\partial^\varepsilon(\theta_NM_{\mathrm P})(\xi)
\big\|_{L^2\left(\frac{ds}{s}\frac{dt}{t}\right)}
\le C_{\chi,\mathrm P},
\qquad \varepsilon\in\{0,1\}^3.
$$
By the Hilbert-space case of
\cite[Theorem~4.5 and Remark~4.6~ii)]{SW}, we therefore have
$$
\bigg\|\bigg(\int_0^\infty\int_0^\infty
\left|\mathcal F^{-1}[\theta_N\widehat{Q_{s,t}^{\mathrm P}f}]\right|^2
\,\frac{dt}{t}\frac{ds}{s}\bigg)^{1/2}\bigg\|_{L^p(\R^3)}
\le C_{\mathrm P,p}\|f\|_{L^p(\R^3)},
$$
with a constant independent of $N$.
By Plancherel,
\begin{align*}
&\left\|\mathcal F^{-1}[\theta_N\widehat{Q_{s,t}^{\mathrm P}f}]
-Q_{s,t}^{\mathrm P}f\right\|_{L^2\left(dx\,\frac{ds}{s}\frac{dt}{t}\right)}^2
= c_{\mathrm P}^2\int_{\R^3}|\theta_N(\xi)-1|^2|\widehat f(\xi)|^2\,d\xi
\to 0
\qquad(N\to\infty),
\end{align*}
where the limit follows from the dominated convergence theorem.
Along a subsequence, this convergence holds in
$L^2(\frac{ds}{s}\frac{dt}{t})$ for almost every $x$.
Fatou's lemma then yields \eqref{e:sg-vertical}.

\smallskip

We now use \eqref{e:sg-vertical} to bound $\|S_{z,\alpha}^{\mathrm P}f\|_{L^p}$.
For $w\ge0$,  
\begin{align*}
\int_{\R^3}|S_{z,\alpha}^{\mathrm P}f(x)|^2w(x)\,dx
=&\ \int_{\R^3}\int_0^\infty\int_0^\infty|Q_{s,t}^{\mathrm P}f(y)|^2
\bigg(\int_{\substack{|x_1-y_1|<\alpha s\\
|x_2-y_2|<\alpha t\\|x_3-y_3|<\alpha st}}w(x)\,dx\bigg)
\frac{dt\,ds}{s^3t^3}\,dy\\
\leq&\ (2\alpha)^3\int_{\R^3}
\bigg(\int_0^\infty\int_0^\infty|Q_{s,t}^{\mathrm P}f(y)|^2\,\frac{dt}{t}\frac{ds}{s}\bigg)\,
\mathsf Mw(y)\,dy,
\end{align*}
where $\mathsf M$ is the strong maximal operator in $\R^3$.

For $2<p<\infty$ and $\|w\|_{L^{(p/2)'}}=1$, H\"older's inequality,
\eqref{e:sg-vertical} and the strong maximal bound imply
$$
\int_{\R^3}|S_{z,\alpha}^{\mathrm P}f|^2w
\le (2\alpha)^3C_{\mathrm P,p}^2\|f\|_{L^p}^2\|\mathsf Mw\|_{L^{(p/2)'}}
\le C_{\mathrm P,\alpha,p}\|f\|_{L^p}^2.
$$
Taking the supremum over these $w$ proves the bound for $S_{z,\alpha}^{\mathrm P}f$ for $2<p<\infty$.

\smallskip
\noindent{\it Step 2.}
For the intermediate square function, \eqref{sp:e:sg-forward} with $d=1$ gives
$$
\big \|S_{1,\alpha}^{\varphi_{\mathrm P}^{(1)}}f(\cdot,x_2,x_3)\big \|_{L^p(\R)}
\le B_{p,1,\alpha}\|f(\cdot,x_2,x_3)\|_{L^p(\R)}
$$
for almost every $(x_2,x_3)$.
Integrating the $p$th powers in $(x_2,x_3)$, we obtain
$$
\big \|S_{1,\alpha}^{\varphi_{\mathrm P}^{(1)}}f\big \|_{L^p(\R^3)}
\le B_{p,1,\alpha}\|f\|_{L^p(\R^3)}.
$$
Its continuity under $L^p$ convergence follows from the pointwise triangle inequality
and this bound.

\smallskip
\noindent{\it Step 3.}
We extend the bounds to every $f\in L^p(\R^3)$ by the approximation argument
in {\it Step 1} of Proposition~\ref{sp:p:semigroups}.
Indeed, the kernels of $Q_{s,t}^{\mathrm P}$ have uniformly bounded $L^1$ norms.
If $f_j\in L^2\cap L^p$ and $\|f_j-f\|_p\to0$, Young's inequality gives
$$
\int_{N^{-1}}^N\int_{N^{-1}}^N
\|Q_{s,t}^{\mathrm P}(f_j-f)\|_{L^p(\R^3)}^p\,\frac{dt}{t}\frac{ds}{s}
\le C_{\mathrm P}^p(2\log N)^2\|f_j-f\|_p^p\to0,
\qquad N\ge2.
$$
The same argument applies to the first-variable convolutions.
A common subsequence converges almost everywhere in the spatial and scale variables;
Fatou's lemma then gives both area-function bounds and \eqref{e:sg-vertical} for $f$.
The pointwise triangle inequality, together with these bounds, proves
continuity of both area functions under $L^p$ convergence.
\end{proof}

\subsection{Dual kernels and the Poisson reverse estimate}\label{sub:sg-reverse}

Fix $\rho\in2^{\mathbb Z}$ with $54\rho<\alpha$.
Take the dual kernels $\Phi_{\mathrm P}^{(1)}$ and $\Phi_{\mathrm P}^{(2)}$
from \eqref{sp:e:sg-dual}, using this $\rho$ in dimensions one and two.
Their support, cancellation and nonnegative Fourier transforms were verified there.
By \eqref{sp:e:sg-norm},
\begin{equation}\label{e:sg-duals}
\int_0^\infty b_{\mathrm P}(r)\widehat{\Phi_{\mathrm P}^{(1)}}(r)\,\frac{dr}{r}
=\int_0^\infty b_{\mathrm P}(r)\widehat{\Phi_{\mathrm P}^{(2)}}(r,0)\,\frac{dr}{r}=1.
\end{equation}
We use the dilates
$$
\Phi_{\mathrm P,s}^{(1)}(x_1)=\frac1s\Phi_{\mathrm P}^{(1)}\left(\frac{x_1}{s}\right),
\qquad
\Phi_{\mathrm P,s,t}^{(2)}(y_2,y_3)=\frac1{st^2}
\Phi_{\mathrm P}^{(2)}\left(\frac{y_2}{t},\frac{y_3}{st}\right).
$$
Both dual kernels depend only on $\alpha$ and are fixed independently of $p$.

\begin{proof}[Proof of the Poisson reverse estimate in Theorem~\ref{t:semigroup}]
We apply Proposition~\ref{sp:p:semigroups} in the first variable
and the atom argument of Proposition~\ref{prop:square-2} in the remaining variables.
The dual kernels determine the supports and derivative bounds of the smooth atoms.
For the second step, we verify the reproducing formula and energy identity at the coupled scales.

\smallskip
\noindent{\it Step 1.}
By \eqref{sp:e:sg-reverse} with $d=1$,
$$
\|h\|_{L^p(\R)}
\le C_{\mathrm P,\alpha}\sqrt p\,
\|\mathcal A_{\varphi_{\mathrm P}^{(1)},\alpha}h\|_{L^p(\R)},
\qquad h\in L^p(\R).
$$
Write $\pmb y=(y_2,y_3)$ and $\pmb v=(v_2,v_3)$.
Taking $h=f(\cdot,\pmb y)$ and integrating the $p$th powers in $\pmb y$, we obtain
\begin{equation}\label{e:sg-first}
\|f\|_{L^p(\R^3)}\le C_{\mathrm P,\alpha}\sqrt p\,
\|S_{1,\alpha}^{\varphi_{\mathrm P}^{(1)}}f\|_{L^p(\R^3)},
\qquad f\in L^2(\R^3)\cap L^p(\R^3).
\end{equation}
The same slice argument proves \eqref{e:sg-first} for every $f\in L^p(\R^3)$.
For the $L^2$ reconstruction below, we continue with
$f\in L^2(\R^3)\cap L^p(\R^3)$, and will extend 
 the final estimate to all $f\in L^p(\R^3)$
at the end of the proof.

\smallskip
\noindent{\it Step 2.}
We now reconstruct $(f(\cdot,\pmb y)*\varphi_{\mathrm P,s}^{(1)})(u)$ from
$$
Q_{s,t}^{\mathrm P}f(u,\pmb y)
=\int_{\R^2}(f(\cdot,\pmb v)*\varphi_{\mathrm P,s}^{(1)})(u)
\varphi_{\mathrm P,s,t}^{(2)}(\pmb y-\pmb v)\,d\pmb v.
$$
By radiality and \eqref{e:sg-duals}, the substitution
$r=t\sqrt{\xi_2^2+s^2\xi_3^2}$ gives
$$
\int_0^\infty b_{\mathrm P}\big (t\sqrt{\xi_2^2+s^2\xi_3^2}\big  )
\widehat{\Phi_{\mathrm P}^{(2)}}(t\xi_2,st\xi_3)\,\frac{dt}{t}=1,
\qquad s>0,\quad (\xi_2,\xi_3)\ne(0,0).
$$
For $f\in L^2(\R^3)\cap L^p(\R^3)$, the same Plancherel calculation
as in \eqref{e:second-energy} gives
\begin{align}\label{e:sg-second-energy}
&\int_{\Gamma_\alpha(x_1)}\int_0^\infty\int_{\R^2}
|Q_{s,t}^{\mathrm P}f(u,\pmb v)|^2
\,d\pmb v\,\frac{dt}{t}\frac{du\,ds}{s^2}\notag\\
= &c_{\mathrm P}\int_{\R^2}\int_{\Gamma_\alpha(x_1)}
|(f(\cdot,\pmb y)*\varphi_{\mathrm P,s}^{(1)})(u)|^2
\,\frac{du\,ds}{s^2}\,d\pmb y<\infty
\end{align}
for almost every $x_1$, since integration in $x_1$ gives
$2\alpha c_{\mathrm P}^2\|f\|_{L^2(\R^3)}^2<\infty$.

Fix such an $x_1$.
By \eqref{e:sg-duals} and the nonnegativity of
$b_{\mathrm P}$ and $\widehat{\Phi_{\mathrm P}^{(2)}}$,
the truncated reproducing multipliers lie in $[0,1]$
and converge to $1$ for $(\xi_2,\xi_3)\ne(0,0)$.
Combining Plancherel, dominated convergence and
\eqref{e:sg-second-energy}, we obtain
$$
\lim_{N\to\infty}
\bigg\|\int_{N^{-1}}^N\int_{\R^2}
Q_{s,t}^{\mathrm P}f(u,\pmb v)
\Phi_{\mathrm P,s,t}^{(2)}(\pmb y-\pmb v)
\,d\pmb v\,\frac{dt}{t}
-(f(\cdot,\pmb y)*\varphi_{\mathrm P,s}^{(1)})(u)
\bigg\|_{L^2\left(\Gamma_\alpha(x_1)\times\R^2,
\frac{du\,ds}{s^2}\,d\pmb y\right)}=0.
$$
Besides,
$$
\int_0^\infty
|\widehat{\Phi_{\mathrm P}^{(2)}}(r,0)|^2\,\frac{dr}{r}<\infty,
$$
 which we will use for the atomic reconstruction below,
as in {\it Step 4} of Proposition~\ref{prop:square-2}. 

Take $\mathscr R$, $\mathcal W^{(2)}(R)$ and
$Q(R)=Q_2(R)\times Q_3(R)$ as in {\it Step 2} of
Proposition~\ref{prop:square-2}.
For $R=I_2\times I_3\in\mathscr R_{m,k}$, define
$$
\lambda_R^{(2)}(f)^2
=\int_{\mathcal W^{(2)}(R)}
|Q_{s,t}^{\mathrm P}f(u,\pmb v)|^2
\,\frac{du\,ds\,d\pmb v\,dt}{s^2t},
$$
and
$$
\widetilde{\mathsf a}_R^{(2)}(f)(u,s;\pmb y)
=\mathbf1_{\Gamma_{\alpha,m}(x_1)}(u,s)
\int_{2^k}^{2^{k+1}}\int_R
Q_{s,t}^{\mathrm P}f(u,\pmb v)
\Phi_{\mathrm P,s,t}^{(2)}(\pmb y-\pmb v)
\,d\pmb v\,\frac{dt}{t}.
$$
Let
$$
\mathsf a_R^{(2)}(f)
=\frac{\widetilde{\mathsf a}_R^{(2)}(f)}{\lambda_R^{(2)}(f)},
$$
with $\mathsf a_R^{(2)}(f)=0$ when $\lambda_R^{(2)}(f)=0$.

Since $\operatorname{supp}\Phi_{\mathrm P}^{(2)}\subset B(0,\rho)$
and $\int_{\R^2}\Phi_{\mathrm P}^{(2)}=0$,
$$
\operatorname{supp}\mathsf a_R^{(2)}(f)(u,s;\cdot)
\subset(5I_2)\times(9I_3)\subset Q(R),
\qquad
\int_{\R^2}\mathsf a_R^{(2)}(f)(u,s;\pmb y)\,d\pmb y=0.
$$
For $\gamma_2,\gamma_3\in\mathbb N_0$ with $\gamma_2+\gamma_3\le2$,
the substitution $\pmb y=(tz_2,stz_3)$ gives $d\pmb y=st^2\,d\pmb z$, and hence
$$
\left\|\partial_{y_2}^{\gamma_2}\partial_{y_3}^{\gamma_3}
\Phi_{\mathrm P,s,t}^{(2)} \right\|_{L^2(\R^2)}^2
=s^{-1-2\gamma_3}t^{-2-2\gamma_2-2\gamma_3}
\left \|\partial_{y_2}^{\gamma_2}\partial_{y_3}^{\gamma_3}
\Phi_{\mathrm P}^{(2)}  \right\|_{L^2(\R^2)}^2.
$$
Thus, by Cauchy--Schwarz as in {\it Step 2} of
Proposition~\ref{prop:square-2},
\begin{align*}
 \int_{\Gamma_\alpha(x_1)}
\left |\partial_{y_2}^{\gamma_2}\partial_{y_3}^{\gamma_3}
\mathsf a_R^{(2)}(f)(u,s;\pmb y)  \right|^2\,\frac{du\,ds}{s^2} 
\leq &\sup_{2^m\le s<2^{m+1}}
\int_{2^k}^{2^{k+1}}
\left \|\partial_{y_2}^{\gamma_2}\partial_{y_3}^{\gamma_3}
\Phi_{\mathrm P,s,t}^{(2)}   \right \|_{L^2(\R^2)}^2\,\frac{dt}{t}\\
\leq &C_{\mathrm P,\alpha}
|Q(R)|^{-1}|Q_2(R)|^{-2\gamma_2}|Q_3(R)|^{-2\gamma_3}.
\end{align*}

For a finite family $\mathscr F\subset\mathscr R$, use the partition from {\it Step 2} of
Proposition~\ref{prop:square-2}:
$$
\mathscr F
=\bigcup_{\vec{\pmb\varepsilon},\vec{\pmb a}}
(\mathscr F\cap\mathscr R_{\vec{\pmb\varepsilon},\vec{\pmb a}})
\quad\text{(disjoint)},
$$
where $\vec{\pmb\varepsilon}\in\{0,1\}^2$ and
$\vec{\pmb a}\in\{3,4\}\times\{4,5\}$.
In each class, the two selected grids and the ratios $|Q_j(R)|/|I_j|=2^{a_j}$ are fixed.
Thus the multiplicity bound \eqref{multi-23} applies to the present atoms as well.
The factor $\mathbf1_{\Gamma_{\alpha,m}(x_1)}$ keeps the restriction $2^m\le s<2^{m+1}$.
Within each class, the Haar coefficients therefore satisfy the same scale relation
\eqref{e:second-ratio} and estimate \eqref{e:second-depth}.

Since $54\rho<\alpha$, the Haar estimate \eqref{e:second-finite}
and the square-function comparison \eqref{e:second-square} ensure
for every finite $\mathscr F\subset\mathscr R$,
\begin{align*}
 \bigg\|\bigg(\int_{\Gamma_\alpha(x_1)}
\bigg|\sum_{R\in\mathscr F}
\widetilde{\mathsf a}_R^{(2)}(f)(u,s;\cdot)\bigg|^2
\,\frac{du\,ds}{s^2}\bigg)^{1/2}\bigg\|_{L^p(\R^2)} 
\leq &C_{\mathrm P,\alpha}\sqrt p\,
\bigg\|\bigg(\sum_{R\in\mathscr F}
\lambda_R^{(2)}(f)^2
\frac{\mathbf1_{Q(R)}}{|Q(R)|}\bigg)^{1/2}\bigg\|_{L^p(\R^2)}\\
\leq &C_{\mathrm P,\alpha}\sqrt p\,
\left\|S_{z,\alpha}^{\mathrm P}f(x_1,\cdot) \right\|_{L^p(\R^2)}.
\end{align*}

By \eqref{W2R} and \eqref{e:sg-second-energy},
$$
\sum_{R\in\mathscr R}\lambda_R^{(2)}(f)^2
=\int_{\Gamma_\alpha(x_1)}\int_0^\infty\int_{\R^2}
 \left |Q_{s,t}^{\mathrm P}f(u,\pmb v)  \right|^2
\,d\pmb v\,\frac{dt}{t}\frac{du\,ds}{s^2}<\infty.
$$
As in {\it Step 4} of Proposition~\ref{prop:square-2},
the $L^2$ reconstruction and Fatou's lemma give
$$
\big \|S_{1,\alpha}^{\varphi_{\mathrm P}^{(1)}}f(x_1,\cdot)  \big \|_{L^p(\R^2)}
\le C_{\mathrm P,\alpha}\sqrt p\,
\big \|S_{z,\alpha}^{\mathrm P}f(x_1,\cdot) \big \|_{L^p(\R^2)}.
$$
This, combined with \eqref{e:sg-first},  yields
$$
\|f\|_{L^p(\R^3)}
\le C_{\mathrm P,\alpha}\sqrt p\,
\big \|S_{1,\alpha}^{\varphi_{\mathrm P}^{(1)}}f\big \|_{L^p(\R^3)}
\le C_{\mathrm P,\alpha}p\,
\big \|S_{z,\alpha}^{\mathrm P}f \big \|_{L^p(\R^3)},
\qquad f\in L^2(\R^3)\cap L^p(\R^3).
$$
For $f\in L^p(\R^3)$, choose
$f_j\in L^2(\R^3)\cap L^p(\R^3)$ with
$\|f_j-f\|_{L^p}\to0$.
By the $L^p$ continuity in Lemma~\ref{l:sg-forward},
letting $j\to\infty$ gives the same estimate for $f$.
\end{proof}

\subsection{Sharpness for the Poisson square function}\label{sub:sg-sharp}

We prove that the factor $p$ in the Poisson reverse estimate is optimal.
For this, we use the one-dimensional test functions and the dilation
in $x_3$ from Subsection~\ref{sec:sharp-p}.
In the construction of Lemma~\ref{l:lower-transfer}, we now replace
$\psi_1,\psi_2$ by
$\varphi_{\mathrm P}^{(1)},\varphi_{\mathrm P}^{(2)}$, respectively.
The kernel $\kappa$ in \eqref{kappa} then becomes
$$
\kappa(y_2)
=\int_\R\varphi_{\mathrm P}^{(2)}(y_2,y_3)\,dy_3
=\varphi_{\mathrm P}^{(1)}(y_2),
$$
 since taking the Fourier transform in $y_2$ gives
$b_{\mathrm P}(|\xi_2|)$ for both the integral and $\varphi_{\mathrm P}^{(1)}$.

The one-dimensional construction was verified in Proposition~\ref{sp:p:semigroups}.
We now check the derivative bounds and convergence needed to transfer it to the coupled scales.

\begin{proof}[Proof of sharpness for the Poisson square function]
Recall  $b_{\mathrm P}(r)=re^{-r}$ and note  $\kappa=\varphi_{\mathrm P}^{(1)}$ now.

\noindent{\it Step 1.}
By \eqref{sp:e:sg-lower} with $d=1$ and $\nu=\mathrm P$, there is a nonzero
$g_p\in\mathcal S_0(\R)$ such that
$$
0<\|\mathcal A_{\kappa,\alpha}g_p\|_{L^p(\R)}<\infty,
\qquad
\|g_p\|_{L^p(\R)}\ge c_{\mathrm P,\alpha}\sqrt p\,
\|\mathcal A_{\kappa,\alpha}g_p\|_{L^p(\R)}.
$$
\smallskip
\noindent{\it Step 2.}
Choose a nonzero
$\chi\in\mathcal S(\R)$ with $\operatorname{supp}\widehat\chi\Subset\R\setminus\{0\}$, and let
$$
F(x)=g_p(x_1)g_p(x_2)\chi(x_3),
\qquad
f_L(x)=g_p(x_1)g_p(x_2)\chi(x_3/L),\qquad L\ge1.
$$
Denote
$$
q=st/L,\qquad
\mathcal T_{s,t,L}(x)=Q_{s,t}^{\mathrm P}f_L(x_1,x_2,Lx_3).
$$
Then
$$
\widehat{\mathcal T_{s,t,L}}(\xi)=m_{s,t,L}(\xi)\widehat F(\xi),
\qquad
m_{s,t,L}(\xi)=b_{\mathrm P}(s|\xi_1|)\widehat{\varphi_{\mathrm P}^{(2)}}(t\xi_2,q\xi_3).
$$
For $\pmb\zeta=(\zeta_2,\zeta_3)\ne0$, use \eqref{sp:e:sg-symbol} with $d=2$ and $\nu=\mathrm P$.
On $\operatorname{supp}\widehat F\Subset\{\xi_1\xi_2\xi_3\ne0\}$,
$|(t\xi_2,q\xi_3)|\simeq_F t+q$.
It follows from \eqref{sp:e:sg-symbol} with $M=8$ that, for $j+k\le4$,
\begin{align*}
\left |\partial_{\xi_2}^j\partial_{\xi_3}^k
[\widehat{\varphi_{\mathrm P}^{(2)}}(t\xi_2,q\xi_3)]  \right|
=&\ t^jq^k  \left |\partial_{\zeta_2}^j\partial_{\zeta_3}^k
\widehat{\varphi_{\mathrm P}^{(2)}}(t\xi_2,q\xi_3)  \right|\\
\leq&\ C_{\mathrm P,F}
\begin{cases}
t^jq^k(t+q)^{1-j-k},&t+q\le1,\\[2pt]
t^jq^k(t+q)^{-8},&t+q>1
\end{cases}\\
\leq&\ C_{\mathrm P,F}(1+t+q)^{-4}\min\{1,t+q\}\\[2pt]
\leq&\ C_{\mathrm P,F}(1+s)(1+t)^{-4}\min\{1,t\},
\end{align*}
where $t^jq^k\le(t+q)^{j+k}$ and $t\le t+q\le t(1+s)$.
 Using \eqref{sp:e:sg-annular} for derivatives of order at most $4$ with $M=6$, we therefore obtain
\begin{align*}
\max_{|\beta|\le4}\|\partial^\beta(m_{s,t,L}\widehat F)\|_{L^\infty}
\leq&\ C_{\mathrm P,F}(1+s)\min\{s,s^{-6}\}(1+t)^{-4}\min\{1,t\}
\leq C_{\mathrm P,F}\min\{s,s^{-1}\}\min\{t,t^{-1}\}.
\end{align*}
Then, using  integration by parts in \eqref{e:trans-ibp} , we obtain
$$
|\mathcal T_{s,t,L}(x)|
\le C_{\mathrm P,F}(1+|x|)^{-4}\min\{s,s^{-1}\}\min\{t,t^{-1}\},
$$
and hence
\begin{equation}\label{e:sg-transfer-bound}
\|\mathcal T_{s,t,L}\|_{L^p(\R^3)}
\le C_{\mathrm P,F,p}\min\{s,s^{-1}\}\min\{t,t^{-1}\},
\qquad s,t>0,\quad L\ge1,
\end{equation}
where the constant is independent of $s,t,L$.

For fixed $s,t>0$, we have $q\to0$ as $L\to\infty$.
Since $|\xi_2|$ is bounded below by a positive constant on
$\operatorname{supp}\widehat F$,
$$
\lim_{L\to\infty}
\sup_{\xi\in\operatorname{supp}\widehat F}
\left|
\partial_\xi^\gamma
\left[
m_{s,t,L}(\xi)
-b_{\mathrm P}(s|\xi_1|)b_{\mathrm P}(t|\xi_2|)
\right]
\right|=0,
\qquad |\gamma|\le4.
$$
By \eqref{e:trans-ibp}, we obtain
$$
\lim_{L\to\infty}
\left\|\mathcal T_{s,t,L}
-(g_p*\kappa_s)(x_1)(g_p*\kappa_t)(x_2)\chi(x_3)
\right\|_{L^p(\R^3)}=0.
$$
\smallskip
\noindent{\it Step 3.}
The change of variables in \eqref{e:trans-cone} gives
$$
S_{z,\alpha}^{\mathrm P}f_L(x_1,x_2,Lx_3)^2
=\int_{(-\alpha,\alpha)^3}\int_0^\infty\int_0^\infty
\left|\mathcal T_{s,t,L}(x_1+sa_1,x_2+ta_2,x_3+qa_3)\right|^2
\,\frac{dt}{t}\frac{ds}{s}\,da_1\,da_2\,da_3.
$$
For fixed $s,t>0$ and $a_1,a_2,a_3\in(-\alpha,\alpha)$,
we have $q\to0$ as $L\to\infty$.
The $L^p$ limit above and continuity of translations give
\begin{align*}
\lim_{L\to\infty}\bigg\|
&\mathcal T_{s,t,L}(x_1+sa_1,x_2+ta_2,x_3+qa_3)
-(g_p*\kappa_s)(x_1+sa_1)
(g_p*\kappa_t)(x_2+ta_2)\chi(x_3)
\bigg\|_{L^p(\R^3)}=0.
\end{align*}
By \eqref{e:sg-transfer-bound} and this limit,
the $L^p$ norm of the difference is bounded by
$$
C_{\mathrm P,F,p}\min\{s,s^{-1}\}\min\{t,t^{-1}\},
$$
uniformly in $L,a_1,a_2,a_3$, and
$$
\int_0^\infty\int_0^\infty
\min\{s,s^{-1}\}^2\min\{t,t^{-1}\}^2
\,\frac{dt}{t}\frac{ds}{s}=1.
$$
Combining Minkowski's inequality for $p\ge2$ and dominated convergence,
as in {\it Step 3} of Lemma~\ref{l:lower-transfer}, we obtain
$$
\lim_{L\to\infty}
\bigg\|S_{z,\alpha}^{\mathrm P}f_L(x_1,x_2,Lx_3)
-\sqrt{2\alpha}\,|\chi(x_3)|
\mathcal A_{\kappa,\alpha}g_p(x_1)
\mathcal A_{\kappa,\alpha}g_p(x_2)
\bigg\|_{L^p(\R^3)}=0.
$$
{The same bound \eqref{e:sg-transfer-bound} gives
$\|S_{z,\alpha}^{\mathrm P}f_L\|_{L^p}<\infty$ for each $L$.}

Since
$
L^{-1/p}\|f_L\|_{L^p(\R^3)}
=\|g_p\|_{L^p(\R)}^2\|\chi\|_{L^p(\R)},
$
we obtain
$$
\lim_{L\to\infty}
\frac{\|f_L\|_{L^p(\R^3)}}{\|S_{z,\alpha}^{\mathrm P}f_L\|_{L^p(\R^3)}}
=\frac1{\sqrt{2\alpha}}
\bigg(
\frac{\|g_p\|_{L^p(\R)}}
{\|\mathcal A_{\kappa,\alpha}g_p\|_{L^p(\R)}}
\bigg)^2
\ge c_{\mathrm P,\alpha}p.
$$
For each $p$, choose $L_p\ge1$ large enough that
$f_p:=f_{L_p}$ satisfies 
$$
\|f_p\|_{L^p(\R^3)}
\ge c_{\mathrm P,\alpha}p\,
\|S_{z,\alpha}^{\mathrm P}f_p\|_{L^p(\R^3)},
\qquad
0<\|S_{z,\alpha}^{\mathrm P}f_p\|_{L^p(\R^3)}<\infty,
$$
where $\|S_{z,\alpha}^{\mathrm P}f_p\|_{L^p(\R^3)}>0$
follows from \eqref{e:sg-energy}.
Additionally, 
$$
\widehat f_p(\xi)
=L_p\widehat g_p(\xi_1)\widehat g_p(\xi_2)\widehat\chi(L_p\xi_3),
\qquad
\operatorname{supp}\widehat f_p\Subset\{\xi_1\xi_2\xi_3\ne0\}.
$$
This proves the sharpness.
\end{proof}

\subsection{The heat semigroup case}\label{sub:sg-heat}

For the heat semigroup, take $b_{\mathrm H}(r)=r^2e^{-r^2}$ and the kernels
$\varphi_{\mathrm H}^{(1)},\varphi_{\mathrm H}^{(2)}$ from \eqref{sp:e:sg-kernels}.
Their cancellation and decay were verified in Subsection~\ref{sp:sub:semigroups}.
The heat semigroup derivative kernels are
$$
\varphi_{\mathrm H}^{(1)}(x_1)
=\sqrt\pi\,\Big(\frac12-\pi^2x_1^2\Big)e^{-\pi^2x_1^2},
$$
and
$$
\varphi_{\mathrm H}^{(2)}(\pmb y)
=\pi[1-\pi^2(y_2^2+y_3^2)]e^{-\pi^2(y_2^2+y_3^2)},
\qquad \pmb y=(y_2,y_3).
$$

Again, for $s,t>0$, denote
$$
\varphi_{\mathrm H,s}^{(1)}(x_1)
=\frac1s\varphi_{\mathrm H}^{(1)}\left(\frac{x_1}{s}\right),
\qquad
\varphi_{\mathrm H,s,t}^{(2)}(y_2,y_3)
=\frac1{st^2}\varphi_{\mathrm H}^{(2)}
\left(\frac{y_2}{t},\frac{y_3}{st}\right).
$$
For $f\in L^p(\R^3)$, $1\le p\le\infty$, define
$$
Q_{s,t}^{\mathrm H}f=f*K_{s,t}^{\mathrm H},
\qquad\text{where}\qquad
K_{s,t}^{\mathrm H}(x)
=\frac1{s^2t^2}\varphi_{\mathrm H}^{(1)}\left(\frac{x_1}{s}\right)
\varphi_{\mathrm H}^{(2)}\left(\frac{x_2}{t},\frac{x_3}{st}\right).
$$
Then 
$$
\widehat{Q_{s,t}^{\mathrm H}f}(\xi)
=b_{\mathrm H}(s|\xi_1|)
b_{\mathrm H}\big (t\sqrt{\xi_2^2+s^2\xi_3^2}\big  )\widehat f(\xi),\qquad f\in\mathcal S(\R^3).
$$

Here the corresponding area functions are
$$
S_{z,\alpha}^{\mathrm H}f(x)^2
=\int_0^\infty\int_0^\infty
\int_{\substack{|u-x_1|<\alpha s\\
|v_2-x_2|<\alpha t\\|v_3-x_3|<\alpha st}}
|Q_{s,t}^{\mathrm H}f(u,v_2,v_3)|^2
\,\frac{du\,dv_2\,dv_3\,dt\,ds}{s^3t^3},
$$
and, for fixed $(x_2,x_3)$, 
$$
S_{1,\alpha}^{\varphi_{\mathrm H}^{(1)}}f(x_1,x_2,x_3)^2
=\int_0^\infty\int_{|u-x_1|<\alpha s}
\left|(f(\cdot,x_2,x_3)*\varphi_{\mathrm H,s}^{(1)})(u)\right|^2
\,\frac{du\,ds}{s^2}.
$$

The argument for the heat reverse estimate in Theorem~\ref{t:semigroup}
follows the Poisson case; we verify the changes below.

For the upper $L^p$ estimates, recall $c_{\mathrm H}=1/8$ from \eqref{sp:e:sg-energy-constants}.
Also,
$$
rb_{\mathrm H}'(r)=2r^2(1-r^2)e^{-r^2},
\qquad r^2b_{\mathrm H}''(r)-rb_{\mathrm H}'(r)=4r^4(r^2-2)e^{-r^2}.
$$
Thus 
$$b_{\mathrm H}(r),\  \    rb_{\mathrm H}'(r),\  \    r^2b_{\mathrm H}''(r)-rb_{\mathrm H}'(r)\in L^2(\frac{dr}{r}).
$$
Applying  \eqref{e:sg-energy} with $b_{\mathrm P}$ replaced by $b_{\mathrm H}$, we obtain 
\begin{equation}\label{e:heat-energy}
\|S_{z,\alpha}^{\mathrm H}f\|_{L^2(\R^3)}^2
=(2\alpha)^3c_{\mathrm H}^2\|f\|_{L^2(\R^3)}^2,
\qquad f\in L^2(\R^3).
\end{equation}
As in the proof of Lemma~\ref{l:sg-forward}, we also have for $f\in L^p(\R^3)$ with $2\le p<\infty$ that 
\begin{equation}\label{e:heat-vertical}
\bigg\|\bigg(\int_0^\infty\int_0^\infty
|Q_{s,t}^{\mathrm H}f|^2\,\frac{dt}{t}\frac{ds}{s}\bigg)^{1/2}\bigg\|_{L^p(\R^3)}
\le C_{\mathrm H,p}\|f\|_{L^p(\R^3)}, 
\end{equation}
and 
$$
\big\|S_{1,\alpha}^{\varphi_{\mathrm H}^{(1)}}f \big \|_{L^p(\R^3)}
+\big\|S_{z,\alpha}^{\mathrm H}f \big\|_{L^p(\R^3)}
\le B_{\mathrm H,\alpha,p}\|f\|_{L^p(\R^3)}.
$$

For the reverse $L^p$ estimate, use the same $\rho$ and the dual kernels
$\Phi_{\mathrm H}^{(1)},\Phi_{\mathrm H}^{(2)}$ from \eqref{sp:e:sg-dual}.
Taking the same $\eta_1,\eta_2$ as in the Poisson case, these satisfy
$$
\Phi_{\mathrm H}^{(1)}=\frac{c_{\mathrm P,1}}{c_{\mathrm H,1}}\Phi_{\mathrm P}^{(1)},
\qquad
\Phi_{\mathrm H}^{(2)}=\frac{c_{\mathrm P,2}}{c_{\mathrm H,2}}\Phi_{\mathrm P}^{(2)}.
$$
Their normalizations are \eqref{sp:e:sg-norm} with $\nu=\mathrm H$ and $d=1,2$.
Consequently, the $L^2$ reproducing formulas there hold with
$\varphi_{\mathrm H}^{(d)}$ and $\Phi_{\mathrm H}^{(d)}$, $d=1,2$.
The first-variable estimate follows from \eqref{sp:e:sg-reverse} with $d=1$.
For the remaining variables, the proof in Subsection~\ref{sub:sg-reverse} applies with these dual kernels:
the support radius is still $\rho$, the truncated reproducing multipliers lie in $[0,1]$,
and the energy constant is $c_{\mathrm H}$.
The derivative norms of $\Phi_{\mathrm H}^{(2)}$ differ from those of $\Phi_{\mathrm P}^{(2)}$
by the fixed factor $c_{\mathrm P,2}/c_{\mathrm H,2}$.
Thus we obtain
$$
\|f\|_{L^p(\R^3)}
\le C_{\mathrm H,\alpha}\sqrt p\,
\big \|S_{1,\alpha}^{\varphi_{\mathrm H}^{(1)}}f \big \|_{L^p(\R^3)}
\le C_{\mathrm H,\alpha}p\,  \big \|S_{z,\alpha}^{\mathrm H}f \big \|_{L^p(\R^3)},
\qquad f\in L^2(\R^3)\cap L^p(\R^3),
$$
which extends to each $f\in L^p(\R^3)$ because
$S_{1,\alpha}^{\varphi_{\mathrm H}^{(1)}}f$ and $S_{z,\alpha}^{\mathrm H}f$
are continuous under $L^p$ convergence.

It remains to prove that the order $p$ in the reverse estimate for $S_{z,\alpha}^{\mathrm H}$ is sharp.
Let $\kappa=\varphi_{\mathrm H}^{(1)}$.
As in the Poisson case,
$$
\int_\R\varphi_{\mathrm H}^{(2)}(y_2,y_3)\,dy_3=\kappa(y_2).
$$
Take $g_p$ from \eqref{sp:e:sg-lower} with $\nu=\mathrm H$ and $d=1$,
and use the test functions $f_L$ from {\it Step 2} of Subsection~\ref{sub:sg-sharp}.
For $\nu=\mathrm H$, \eqref{sp:e:sg-annular} gives the one-dimensional derivative bounds,
and \eqref{sp:e:sg-symbol} gives the two-dimensional bounds.
The  above derivative estimates  yield  the heat version of \eqref{e:sg-transfer-bound}, uniformly in $L\ge1$.
The argument in {\it Steps 2 and 3} of Subsection~\ref{sub:sg-sharp} thus ensures that 
for every $2\le p<\infty$, there exists a function $f_p\in\mathcal S(\R^3)$ such that
$$
\operatorname{supp}\widehat f_p\Subset  \left\{\xi_1\xi_2\xi_3\ne0 \right\},
\qquad
0<\|S_{z,\alpha}^{\mathrm H}f_p\|_{L^p}<\infty,
$$
and
$$
\|f_p\|_{L^p}\ge c_{\mathrm H,\alpha}p\,\|S_{z,\alpha}^{\mathrm H}f_p\|_{L^p}.
$$
Here $\|S_{z,\alpha}^{\mathrm H}f_p\|_{L^p}>0$ follows  from \eqref{e:heat-energy}.
We complete  the proof of the heat case.

\section{Zygmund multiplier estimates}\label{sec:mult}

In this section, we prove the $p^2$ estimate for Zygmund multipliers in
Theorem~\ref{t:multipliers} and show that this order is optimal.
R. Fefferman and Pipher anticipated the order $p^2$ in \cite[p.~360]{FP}.

Recall that the Fourier multiplier $T_m$ is defined by
$$
\widehat{T_mf}(\xi)=m(\xi)\widehat f(\xi).
$$
Following \cite[p.~349]{FP}, in our coordinate notation, let
$$
A^{x_1}=\left\{\xi\in\R^3:
1/2<|\xi_1|\le1,\quad 1/2<|(\xi_2,\xi_3)|\le1\right\},
$$
and, by interchanging $\xi_1$ and $\xi_2$, let
$$
A^{x_2}=\left\{\xi\in\R^3:
1/2<|\xi_2|\le1,\quad 1/2<|(\xi_1,\xi_3)|\le1\right\}.
$$
Note that interchanging $x_1$ and $x_2$, together with $s$ and $t$,
keeps the form of the Zygmund dilations unchanged.
There is no such symmetry between $x_3$ and either of the first two coordinates,
since the third scale is $st$.

For a positive integer $N$, the seminorm in \eqref{e:mult-snorm} is
$$
\|m\|_{\mathcal M_{Z,x_1}^{N}}
=\max_{0\le a,b,c\le N}\sup_{r,q>0}\sup_{\xi\in A^{x_1}}
\left|\partial_1^a\partial_2^b\partial_3^c
[m(r\xi_1,q\xi_2,rq\xi_3)]\right|.
$$
The derivatives are taken with respect to $\xi=(\xi_1,\xi_2,\xi_3)$,
with $r,q$ fixed.
For the second coordinate, define
$$
\|m\|_{\mathcal M_{Z,x_2}^{N}}
=\|m(\xi_2,\xi_1,\xi_3)\|_{\mathcal M_{Z,x_1}^{N}}.
$$

We will prove that there is a positive integer $N$ such that, if
\begin{equation}\label{e:mult-assum}
m\in C^{3N}\left(\left\{\xi\in\R^3:
\xi_1\ne0,\ (\xi_2,\xi_3)\ne(0,0)\right\}\right),
\qquad
\|m\|_{\mathcal M_{Z,x_1}^{N}}<\infty,
\end{equation}
then
\begin{equation}\label{e:mult-upper}
\|T_mf\|_{L^p(\R^3)}
\le Cp^2\|m\|_{\mathcal M_{Z,x_1}^{N}}\|f\|_{L^p(\R^3)},
\qquad 2\le p<\infty,\quad f\in\mathcal S(\R^3).
\end{equation}
The integer $N$ and the constant $C$ are independent of $p$, $m$ and $f$.
We will then extend $T_m$ to $L^p(\R^3)$ by density.

To prove \eqref{e:mult-upper}, we combine our order-$p$ reverse estimate with the weighted
square-function argument in \cite[Theorem~4.4 and Lemmas~4.6--4.7]{FP}.
We give the details to show how the bound depends on
$\|m\|_{\mathcal M_{Z,x_1}^{N}}$.

For each positive integer $N$, choose a nonzero real-valued even kernel
$\sigma_1\in C_c^\infty(\R)$ supported in $(-1,1)$.
Also choose a nonzero real radial kernel
$\sigma_2\in C_c^\infty(\R^2)$ supported in $B(0,1)\subset\R^2$.
Assume both kernels satisfy
\begin{equation}\label{e:mult-moments}
\begin{cases}\displaystyle 
\int_\R x_1^j\sigma_1(x_1)\,dx_1
= 0,\qquad j=0,\ldots,N+1, \\[10pt]  \displaystyle 
\int_{\R^2}x_2^{\beta_2}x_3^{\beta_3}\sigma_2(x_2,x_3)\,dx_2\,dx_3
= 0,\qquad \beta_2,\beta_3\in\mathbb N_0,\quad \beta_2+\beta_3\le N+1.
\end{cases}
\end{equation}
Such kernels exist. Choose a nonzero real even function
$\eta_1\in C_c^\infty((-1,1))$ and a nonzero real radial function
$\eta_2\in C_c^\infty(B(0,1))$. Then
\begin{equation}\label{e:mult-kernel-choice}
\sigma_1(x_1)
=\frac{d^{2N+2}\eta_1}{dx_1^{2N+2}}(x_1),
\qquad
\sigma_2(x_2,x_3)=\Delta^{N+1}\eta_2(x_2,x_3)
\end{equation}
are nonzero and have the required supports and symmetry.
The moment conditions in \eqref{e:mult-moments} follow by integration by parts, due to  $2N+2>N+1$.

Let
$$
\Psi_1(x_1)=(\sigma_1*\sigma_1)(x_1),
\qquad x_1\in\R,
$$
and
$$
\Psi_2(x_2,x_3)=(\sigma_2*\sigma_2)(x_2,x_3),
\qquad (x_2,x_3)\in\R^2.
$$
For $x=(x_1,x_2,x_3)\in\R^3$, define
$$
\sigma(x)=\sigma_1(x_1)\sigma_2(x_2,x_3),
\qquad
\Psi(x)=\Psi_1(x_1)\Psi_2(x_2,x_3).
$$
For $s,t>0$, their Zygmund dilations are
$$
\sigma_{s,t}(x)=\frac1{s^2t^2}
\sigma_1\left(\frac{x_1}{s}\right)
\sigma_2\left(\frac{x_2}{t},\frac{x_3}{st}\right),\qquad
\Psi_{s,t}(x)=\frac1{s^2t^2}
\Psi_1\left(\frac{x_1}{s}\right)
\Psi_2\left(\frac{x_2}{t},\frac{x_3}{st}\right).
$$
Since $\sigma*\sigma=\Psi$, dilation gives
$$
(\sigma_{s,t}*\sigma_{s,t})(x)
=\frac1{s^2t^2}(\sigma*\sigma)
\left(\frac{x_1}{s},\frac{x_2}{t},\frac{x_3}{st}\right)
=\Psi_{s,t}(x).
$$

{Since $\Psi_{s,t}=\sigma_{s,t}*\sigma_{s,t}$,
$$
(T_mf)*\Psi_{s,t}=(f*\sigma_{s,t})*(T_m\sigma_{s,t}).
$$
To estimate the last convolution, we use the weight}
$$
W_{s,t}(x)=\frac1{s^2t^2}
\left(1+\frac{|x_1|}{s}\right)^{-4}
\left(1+\frac{x_2^2}{t^2}+\frac{x_3^2}{s^2t^2}\right)^{-4},
\qquad s,t>0,\quad x\in\R^3,
$$
which decays at the scales $s,t,st$.
With this weight, define
$$
g_z^{*,\sigma}f(x)^2
=\int_0^\infty\int_0^\infty\int_{\R^3}
|f*\sigma_{s,t}(x+u)|^2W_{s,t}(u)\,du\,\frac{dt}{t}\frac{ds}{s}.
$$
We compare $S_{z,1}^{\Psi}(T_mf)$ with
$g_z^{*,\sigma}f$ and prove the $L^p$ bound for
$g_z^{*,\sigma}f$ needed in \eqref{e:mult-upper} as follows.

\begin{lemma}\label{l:mult-gstar}
There is a positive integer $N$ such that, for every $m$ satisfying \eqref{e:mult-assum},
every pair of kernels $\sigma_1,\sigma_2$ as above satisfying \eqref{e:mult-moments} for this $N$,
and every $f\in\mathcal S(\R^3)$, we have
\begin{equation}\label{e:mult-ptw}
S_{z,1}^{\Psi}(T_mf)(x)
\le C_\sigma\|m\|_{\mathcal M_{Z,x_1}^{N}}g_z^{*,\sigma}f(x).
\end{equation}
Moreover, for every $2\le p<\infty$ and $f\in L^p(\R^3)$,
\begin{equation}\label{e:mult-gstar}
\|g_z^{*,\sigma}f\|_{L^p(\R^3)}\le C_\sigma p\|f\|_{L^p(\R^3)}.
\end{equation}
The constants are independent of $p$, $m$ and $f$.
\end{lemma}

The estimate \eqref{e:mult-ptw} also holds for every $f\in L^p(\R^3)$ with $2\le p<\infty$; see Remark~\ref{rem:aeconverg}.

\begin{proof}[Proof of Lemma \ref{l:mult-gstar}]
\noindent{\it Step 1.}
We first estimate $T_m\sigma$ for a positive integer $N$,
with $\sigma_1,\sigma_2$ chosen as in \eqref{e:mult-moments}.

For any $r,q>0$ and $\zeta\in A^{x_1}$,  it follows from the chain rule and \eqref{e:mult-snorm} that
for $\gamma\in\mathbb N_0^3$ and $|\gamma|\le N$,
\begin{align}\label{e:aux-der1}
r^{\gamma_1+\gamma_3}q^{\gamma_2+\gamma_3}
|\partial^\gamma m(r\zeta_1,q\zeta_2,rq\zeta_3)|
= &\left|\partial_\zeta^\gamma[m(r\zeta_1,q\zeta_2,rq\zeta_3)]\right|
\leq \|m\|_{\mathcal M_{Z,x_1}^{N}}.
\end{align}
For each $\xi$ with $\xi_1\ne0$ and $(\xi_2,\xi_3)\ne(0,0)$, choose
$$
r_\xi=\frac43|\xi_1|,
\qquad
q_\xi=\frac43\, \bigg(\xi_2^2+\frac{\xi_3^2}{r_\xi^2}\bigg)^{1/2},
\qquad
\zeta_\xi=\bigg (\frac{\xi_1}{r_\xi},\frac{\xi_2}{q_\xi},\frac{\xi_3}{r_\xi q_\xi}\bigg).
$$
Then
$$
\frac{|\xi_1|}{r_\xi}
=\bigg(\frac{\xi_2^2}{q_\xi^2}
+\frac{\xi_3^2}{r_\xi^2q_\xi^2}\bigg)^{1/2}
=\frac34,
$$
and so $\zeta_\xi\in A^{x_1}$.
By \eqref{e:aux-der1} with $r=r_\xi$, $q=q_\xi$ and $\zeta=\zeta_\xi$, we have
$$
\left|\partial^\gamma m(\xi)\right|
\le r_\xi^{-\gamma_1-\gamma_3}q_\xi^{-\gamma_2-\gamma_3}
\|m\|_{\mathcal M_{Z,x_1}^{N}},\qquad |\gamma|\le N.
$$
In particular, taking $\gamma=0$, we obtain
$\|m\|_{L^\infty}\le\|m\|_{\mathcal M_{Z,x_1}^{N}}$,
hence  $T_m$ is bounded on $L^2(\R^3)$ from  Plancherel.

Note that  $r_\xi=4|\xi_1|/3$ and
$$
q_\xi\ge\frac{4|(\xi_2,\xi_3)|}{3\max\{1,r_\xi\}}
\ge\frac{c|(\xi_2,\xi_3)|}{1+|\xi_1|},
$$
thus
\begin{equation}\label{e:mult-symbol-der}
 \left|\partial^\gamma m(\xi)\right|
\le C_N\|m\|_{\mathcal M_{Z,x_1}^{N}}
|\xi_1|^{-\gamma_1-\gamma_3}
(1+|\xi_1|)^{\gamma_2+\gamma_3}
|(\xi_2,\xi_3)|^{-\gamma_2-\gamma_3},
\qquad |\gamma|\le N.
\end{equation}

By \eqref{e:mult-moments}, we have
$
\partial_{\xi_1}^j\widehat\sigma_1(0)
=0$ for
$0\le j\le N+1$, and 
$\partial_{\xi_2}^{\beta_2}\partial_{\xi_3}^{\beta_3}
\widehat\sigma_2(0,0)   =0$ for $\beta_2+\beta_3\le N+1$. 
By Taylor's formula near zero and the rapid decay of
$\widehat\sigma_1$ and $\widehat\sigma_2$ at infinity,
we have for every $M>0$,
$$
\left|\partial_{\xi_1}^j\widehat\sigma_1(\xi_1)\right|
\le C_{\sigma,N,M}
\min\{|\xi_1|,1\}^{N+2-j}(1+|\xi_1|)^{-M},
\qquad 0\le j\le N,
$$
and
$$
\left|\partial_{\xi_2}^{j_2}\partial_{\xi_3}^{j_3}
\widehat\sigma_2(\xi_2,\xi_3)\right|
\le C_{\sigma,N,M}
\min\{|(\xi_2,\xi_3)|,1\}^{N+2-j_2-j_3}
(1+|(\xi_2,\xi_3)|)^{-M},
\qquad j_2,j_3\in\mathbb N_0,\quad j_2+j_3\le N.
$$

Taking $M=N+5$ in these two estimates and using
\eqref{e:mult-symbol-der}, we obtain, for
$\beta,\gamma\in\mathbb N_0^3$ with $\gamma\le\beta$ and $|\beta|\le N$,
\begin{align*}
&\left|\partial^\gamma m(\xi)
\partial_1^{\beta_1-\gamma_1}\widehat\sigma_1(\xi_1)
\partial_2^{\beta_2-\gamma_2}\partial_3^{\beta_3-\gamma_3}
\widehat\sigma_2(\xi_2,\xi_3)\right|\\
\leq &C_{\sigma,N}\|m\|_{\mathcal M_{Z,x_1}^{N}}
\min\{|\xi_1|,1\}^{N+2-\beta_1-\gamma_3}
\min\{|(\xi_2,\xi_3)|,1\}^{N+2-\beta_2-\beta_3}
(1+|\xi_1|)^{-5}(1+|(\xi_2,\xi_3)|)^{-5}\\
\leq &C_{\sigma,N}\|m\|_{\mathcal M_{Z,x_1}^{N}}
\min\{|\xi_1|,1\}^{N+2-|\beta|}
\min\{|(\xi_2,\xi_3)|,1\}^{N+2-|\beta|}
(1+|\xi|)^{-5},
\end{align*}
whenever $\xi_1\ne0$ and $(\xi_2,\xi_3)\ne(0,0)$.
Here we used $\gamma_3\le\beta_3$ and $\gamma_2+\gamma_3\le N$.

Let $F(\xi)=m(\xi)\widehat\sigma(\xi)$ for
$\xi_1\ne0$ and $(\xi_2,\xi_3)\ne(0,0)$.
By the product rule and the estimate above,
\begin{align}\label{e:mult-der}
|\partial^\beta F(\xi)|
\leq &\sum_{\gamma\le\beta}\binom{\beta}{\gamma}
\left|\partial^\gamma m(\xi)
\partial_1^{\beta_1-\gamma_1}\widehat\sigma_1(\xi_1)
\partial_2^{\beta_2-\gamma_2}\partial_3^{\beta_3-\gamma_3}
\widehat\sigma_2(\xi_2,\xi_3)\right|\notag\\
\leq &C_{\sigma,N}\|m\|_{\mathcal M_{Z,x_1}^{N}}
\min\{|\xi_1|,1\}^{N+2-|\beta|}
\min\{|(\xi_2,\xi_3)|,1\}^{N+2-|\beta|}
(1+|\xi|)^{-5},
\qquad |\beta|\le N.
\end{align}

{Extend $F$ and its derivatives of order at most $N$ by zero
to $\{\xi_1=0\}\cup\{(\xi_2,\xi_3)=(0,0)\}$. By \eqref{e:mult-der},
these extensions are continuous, and}
$$
\frac{|\partial^\beta F(\xi)|}{|\xi-\xi^0|}
\le C_{\sigma,N}\|m\|_{\mathcal M_{Z,x_1}^{N}}
|\xi-\xi^0|^{N+1-|\beta|}
\longrightarrow0
\qquad(\xi\to\xi^0),\quad |\beta|<N,
$$
whenever $\xi_1^0=0$ or $(\xi_2^0,\xi_3^0)=(0,0)$.
Hence $F\in C^N(\R^3)$, and \eqref{e:mult-der} also implies
$$
\sum_{|\beta|\le N}\|\partial^\beta F\|_{L^1(\R^3)}
\le C_{\sigma,N}\|m\|_{\mathcal M_{Z,x_1}^{N}}.
$$
Since $(1+|u|)^N\le C_N(1+\sum_{j=1}^3|u_j|^N)$,
by Fourier inversion and integration by parts, we have
\begin{align}\label{e:takingN}
(1+|u|)^N|T_m\sigma(u)|
= &(1+|u|)^N|\mathcal F^{-1}F(u)|\notag  \\
\leq &C_N\bigg(
\|F\|_{L^1(\R^3)}
+\sum_{j=1}^3\|\partial_{\xi_j}^NF\|_{L^1(\R^3)}
\bigg)    
\leq C_{\sigma,N}\|m\|_{\mathcal M_{Z,x_1}^{N}}.
\end{align}

By the definition of $W_{1,1}$, we have
$$
W_{1,1}(u)^{-1}
=(1+|u_1|)^4(1+u_2^2+u_3^2)^4
\le(1+|u|)^{4+8}
=(1+|u|)^{12}.
$$
Taking $N=12$ in \eqref{e:takingN}, we have
\begin{equation}\label{e:mult-decay2}
|T_m\sigma(u)|
\le C_\sigma\|m\|_{\mathcal M_{Z,x_1}^{N}}
(1+|u|)^{-12}
\le C_\sigma\|m\|_{\mathcal M_{Z,x_1}^{N}}W_{1,1}(u).
\end{equation}

For fixed $s,t>0$, let
$$
m_{s,t}(\xi)
=m\left(\frac{\xi_1}{s},\frac{\xi_2}{t},
\frac{\xi_3}{st}\right).
$$
For every $\tilde r,\tilde q>0$ and $\zeta\in A^{x_1}$, we have
$$
m_{s,t}(\tilde r\zeta_1,\tilde q\zeta_2,\tilde r\tilde q\zeta_3)
=m\left(\frac{\tilde r}{s}\zeta_1,\frac{\tilde q}{t}\zeta_2,
\frac{\tilde r}{s}\frac{\tilde q}{t}\zeta_3\right).
$$
Since $(\tilde r/s,\tilde q/t)$ ranges over $(0,\infty)^2$
as $(\tilde r,\tilde q)$ ranges over $(0,\infty)^2$,
by \eqref{e:mult-snorm} we have
$$
\|m_{s,t}\|_{\mathcal M_{Z,x_1}^{N}}
=\|m\|_{\mathcal M_{Z,x_1}^{N}}.
$$
Applying \eqref{e:mult-decay2} to the rescaled multiplier, we obtain
\begin{align}\label{e:mult-kernel}
|T_m\sigma_{s,t}(u)|
= &\frac1{s^2t^2}
\left|(T_{m_{s,t}}\sigma)
\left(\frac{u_1}{s},\frac{u_2}{t},\frac{u_3}{st}\right)\right|
\leq \frac{C_\sigma}{s^2t^2}
\|m_{s,t}\|_{\mathcal M_{Z,x_1}^{N}}
W_{1,1}\left(\frac{u_1}{s},\frac{u_2}{t},\frac{u_3}{st}\right)\notag\\
= &C_\sigma\|m\|_{\mathcal M_{Z,x_1}^{N}}W_{s,t}(u),
\end{align}
where  $C_\sigma$ is independent of $s,t$ and $m$.

\smallskip

\noindent{\it Step 2.}
Let $f\in\mathcal S(\R^3)$.
Since
$$
\int_{\R^3}W_{s,t}(u)\,du
=\int_{\R^3}W_{1,1}(u)\,du
=\frac{2\pi}{9},
$$
Cauchy--Schwarz and \eqref{e:mult-kernel} give
\begin{align*}
|(T_mf)*\Psi_{s,t}(y)|^2
= &|(f*\sigma_{s,t})*(T_m\sigma_{s,t})(y)|^2\\
\leq &C_\sigma\|m\|_{\mathcal M_{Z,x_1}^{N}}^2
\bigg(\int_{\R^3}|f*\sigma_{s,t}(v)|
W_{s,t}(y-v)\,dv\bigg)^2\\
\leq &C_\sigma\|m\|_{\mathcal M_{Z,x_1}^{N}}^2
\bigg(\int_{\R^3}|f*\sigma_{s,t}(v)|^2
W_{s,t}(y-v)\,dv\bigg)
\bigg(\int_{\R^3}W_{s,t}(u)\,du\bigg)\\
\leq &C_\sigma\|m\|_{\mathcal M_{Z,x_1}^{N}}^2
\int_{\R^3}|f*\sigma_{s,t}(v)|^2W_{s,t}(y-v)\,dv.
\end{align*}
For $|y_1-x_1|<s$, $|y_2-x_2|<t$ and $|y_3-x_3|<st$,
by  the definition of $W_{s,t}$, we have
$$
W_{s,t}(y-v)\le C W_{s,t}(x-v),
\qquad v\in\R^3.
$$
Integrating over the cone section, whose volume is
$$
\int_{\substack{|y_1-x_1|<s\\
|y_2-x_2|<t\\|y_3-x_3|<st}}dy
=(2s)(2t)(2st)=8s^2t^2,
$$
we obtain
\begin{align*}
S_{z,1}^{\Psi}(T_mf)(x)^2
= &\int_0^\infty\int_0^\infty
\int_{\substack{|y_1-x_1|<s\\
|y_2-x_2|<t\\|y_3-x_3|<st}}
|(T_mf)*\Psi_{s,t}(y)|^2
\,\frac{dy\,dt\,ds}{s^3t^3}\\
\leq &C_\sigma\|m\|_{\mathcal M_{Z,x_1}^{N}}^2
\int_0^\infty\int_0^\infty\int_{\R^3}
|f*\sigma_{s,t}(v)|^2W_{s,t}(x-v)
\,dv\,\frac{dt}{t}\frac{ds}{s}\\
= &C_\sigma\|m\|_{\mathcal M_{Z,x_1}^{N}}^2
g_z^{*,\sigma}f(x)^2,
\end{align*}
where the last equality follows from the change of variables
$u=v-x$ and the identity $W_{s,t}(-u)=W_{s,t}(u)$.
This proves \eqref{e:mult-ptw}.

\smallskip

\noindent{\it Step 3.}
We next prove a weighted $L^2$ estimate for $g_z^{*,\sigma}$.
For $w\in L^1_{\mathrm{loc}}(\R^3)$, define
$$
M_Zw(x)=\sup_{\substack{R\ni x,\ s,t>0\\
\ell(R_1)=s,\ \ell(R_2)=t,\ \ell(R_3)=st}}
\frac1{|R|}\int_R|w(y)|\,dy,
$$
where $R=R_1\times R_2\times R_3$ ranges over rectangles with sides parallel to the coordinate axes.
For $2<p<\infty$, $f\in L^p(\R^3)$ and
nonzero $0\le w\in C_c^\infty(\R^3)$, we claim that 
\begin{equation}\label{e:mult-weighted}
\int_{\R^3}(g_z^{*,\sigma}f)^2w
\le C_\sigma\int_{\R^3}|f|^2M_Zw.
\end{equation}
Indeed, 
by Tonelli's theorem and the evenness of $W_{s,t}$,
$$
\int_{\R^3}(g_z^{*,\sigma}f)^2w
=\int_0^\infty\int_0^\infty\int_{\R^3}
|f*\sigma_{s,t}(u)|^2(w*W_{s,t})(u)\,du\,\frac{dt}{t}\frac{ds}{s}.
$$
To estimate $w*W_{s,t}$, note 
$$
W_{s,t}(u)
\le\frac{C}{s^2t^2}\sum_{j,k\ge0}2^{-4j-8k}
\mathbf1_{\{|u_1|<2^js,\ |u_2|<2^kt,\ |u_3|<2^kst\}},
$$
thus
\begin{align}\label{e:mult-weight}
(w*W_{s,t})(x)
\leq &\frac{C}{s^2t^2}
\sum_{j,k\ge0}2^{-4j-8k}
\int_{\substack{|u_1|<2^js\\
|u_2|<2^kt\\|u_3|<2^kst}}
w(x-u)\,du\notag\\
\leq &\frac{C}{s^2t^2}
\sum_{j,k\ge0}2^{-4j-8k}
\int_{\substack{|u_1|<2^{j+1}s\\
|u_2|<2^kt\\
|u_3|<(2^{j+1}s)(2^{k+1}t)}}
w(x-u)\,du\notag\\
\leq &\frac{C}{s^2t^2}
\sum_{j,k\ge0}2^{-4j-8k}
(2^{j+2}s)(2^{k+1}t)(2^{j+k+3}st)
M_Zw(x)\notag\\
\leq &C\sum_{j,k\ge0}2^{-2j-6k}M_Zw(x)
\le C M_Zw(x).
\end{align}

Recall that A. C\'ordoba \cite[Theorem~(b), p.~256]{Cor} proved a weak $L\log L$
estimate for rectangles with side lengths $s,t,\Phi(s,t)$,
where $\Phi$ is positive and nondecreasing in each variable.
Taking $\Phi(s,t)=st$, we have
\begin{equation}\label{e:mult-weak}
|\{M_Zw>\lambda\}|
\le C\int_{\R^3}\frac{|w(x)|}{\lambda}
\left(1+\log^+\frac{|w(x)|}{\lambda}\right)\,dx,
\qquad \lambda>0.
\end{equation}

For $k\in\mathbb Z$, let $O_k=\{M_Zw>2^k\}$.
{Since $w$ is bounded and compactly supported, \eqref{e:mult-weak} gives $|O_k|<\infty$. Hence}
$$
\|f\mathbf1_{O_k}\|_{L^2(\R^3)}
\le |O_k|^{1/2-1/p}\|f\|_{L^p(\R^3)}<\infty.
$$

To localize the convolution to the sets $O_k$, for $u\in\R^3$ and $s,t>0$, let
$$
R_{s,t}(u)
=u+[-2s,2s]\times[-2t,2t]\times[-8st,8st].
$$
By the definition of $W_{s,t}$,
$$
W_{s,t}(u-y)\le C W_{s,t}(v-y),
\qquad v\in R_{s,t}(u),\quad y\in\R^3.
$$
Combined with \eqref{e:mult-weight}, we have
\begin{align}\label{e:mult-local}
(w*W_{s,t})(u)
\leq &C\inf_{v\in R_{s,t}(u)}(w*W_{s,t})(v)
\le C\inf_{v\in R_{s,t}(u)}M_Zw(v)\notag\\
\leq &C\sum_{\{k:\,2^k<\inf_{v\in R_{s,t}(u)}M_Zw(v)\}}2^k 
\leq C\sum_{\{k:\,R_{s,t}(u)\subset O_k\}}2^k.
\end{align}

If $R_{s,t}(u)\subset O_k$, then
$$
u-\operatorname{supp}\sigma_{s,t}
\subset u+[-s,s]\times[-t,t]\times[-st,st]
\subset R_{s,t}(u)\subset O_k,
$$
and so
$$
f*\sigma_{s,t}(u)
=(f\mathbf1_{O_k})*\sigma_{s,t}(u), \qquad R_{s,t}(u)\subset O_k.
$$
We apply \eqref{e:mult-local} and
\eqref{e:kernel-l2} with $\psi_i=\sigma_i$ to obtain
\begin{align*}
\int_{\R^3}(g_z^{*,\sigma}f)^2w
= &\int_0^\infty\int_0^\infty\int_{\R^3}
|f*\sigma_{s,t}(u)|^2(w*W_{s,t})(u)
\,du\,\frac{dt}{t}\frac{ds}{s}\\
\leq &C\sum_{k\in\mathbb Z}2^k
\int_{\{(u,s,t):\,R_{s,t}(u)\subset O_k\}}
|(f\mathbf1_{O_k})*\sigma_{s,t}(u)|^2
\,du\,\frac{dt}{t}\frac{ds}{s}\\
\leq &C\sum_{k\in\mathbb Z}2^k
\int_0^\infty\int_0^\infty
\|(f\mathbf1_{O_k})*\sigma_{s,t}\|_{L^2(\R^3)}^2
\,\frac{dt}{t}\frac{ds}{s}\\
= &C_\sigma\sum_{k\in\mathbb Z}2^k\int_{O_k}|f(x)|^2\,dx\\
= &C_\sigma\int_{\R^3}|f(x)|^2
\sum_{\{k:\,2^k<M_Zw(x)\}}2^k\,dx\\
\leq &2C_\sigma\int_{\R^3}|f(x)|^2M_Zw(x)\,dx.
\end{align*}
This is the weighted estimate \eqref{e:mult-weighted}.

\smallskip

\noindent{\it Step 4.}
We estimate $\|M_Z\|_{L^q\to L^q}$ to use \eqref{e:mult-weighted} with $q=(p/2)'$.
First let $1<q\le2$ and $0\le w\in L^q(\R^3)$.
By the definition of $M_Z$, for every $\lambda>0$ and $x\in\R^3$,
we have
$$
M_Z(w\mathbf1_{\{w\le\lambda/2\}})(x)
\le\|w\mathbf1_{\{w\le\lambda/2\}}\|_{L^\infty(\R^3)}
\le\lambda/2.
$$
Consequently, 
\begin{align*}
M_Zw(x)
\leq &M_Z(w\mathbf1_{\{w\le\lambda/2\}})(x)
+M_Z(w\mathbf1_{\{w>\lambda/2\}})(x)
\leq \lambda/2
+M_Z(w\mathbf1_{\{w>\lambda/2\}})(x),
\end{align*}
and so
$$
\{x:M_Zw(x)>\lambda\}
\subset
\{x:M_Z(w\mathbf1_{\{w>\lambda/2\}})(x)>\lambda/2\}.
$$
Applying \eqref{e:mult-weak} to this level-set inclusion and integrating in $\lambda$, we obtain
\begin{align*}
\|M_Zw\|_{L^q}^q
= &q\int_0^\infty
\lambda^{q-1}|\{M_Zw>\lambda\}|\,d\lambda\\
\leq &2Cq\int_{\R^3}w(x)
\int_0^{2w(x)}\lambda^{q-2}
\left(1+\log\frac{2w(x)}{\lambda}\right)
\,d\lambda\,dx\\
= &Cq2^q\|w\|_{L^q}^q
\int_0^1\tau^{q-2}(1-\log\tau)\,d\tau
\leq \frac{C}{(q-1)^2}\|w\|_{L^q}^q.
\end{align*}
Thus
$$
\|M_Zw\|_{L^q}
\le C(q-1)^{-2/q}\|w\|_{L^q}
\le C(q-1)^{-2}\|w\|_{L^q},
\qquad 1<q\le2.
$$
For $q>2$,  we have by  H\"older's inequality that
$$
M_Zw(x)\le\bigl[M_Z(|w|^{q/2})(x)\bigr]^{2/q}.
$$
Therefore, for $q>2$,
\begin{align*}
\|M_Zw\|_{L^q}
= &\|(M_Zw)^{q/2}\|_{L^2}^{2/q}
\leq \|M_Z(|w|^{q/2})\|_{L^2}^{2/q}
\leq C^{2/q}\||w|^{q/2}\|_{L^2}^{2/q}
\le C\|w\|_{L^q}.
\end{align*}

Combining the two cases, we obtain
\begin{equation}\label{e:mult-maximal}
\|M_Z\|_{L^q\to L^q}
\le C\max\{1,(q-1)^{-2}\},
\qquad 1<q<\infty,
\end{equation}
with $C$ independent of $q$.
Consequently, 
for $2<p<\infty$ and $f\in L^p(\R^3)$, let $q=(p/2)'$.
By \eqref{e:mult-weighted} and
\eqref{e:mult-maximal},  we have for $0\le w\in C_c^\infty(\R^3)$ that 
\begin{align*}
\int_{\R^3}(g_z^{*,\sigma}f)^2w
\leq &C_\sigma\int_{\R^3}|f|^2M_Zw
\leq C_\sigma\|f\|_{L^p}^2\|M_Zw\|_{L^q}
\leq C_\sigma p^2\|f\|_{L^p}^2\|w\|_{L^q},
\end{align*}
since $q-1=2/(p-2)$. Therefore,
$$
\|g_z^{*,\sigma}f\|_{L^p}^2
\le C_\sigma p^2\|f\|_{L^p}^2.
$$

For $p=2$, by   \eqref{e:kernel-l2} again, 
$$
\|g_z^{*,\sigma}f\|_{L^2}^2
=\bigg(\int_{\R^3}W_{1,1}(u)\,du\bigg )\,
c_{\sigma_1}c_{\sigma_2}\|f\|_{L^2}^2.
$$
Thus \eqref{e:mult-gstar} holds for every $2\le p<\infty$.
\end{proof}

\begin{proof}[Proof of Theorem~\ref{t:multipliers}]
\noindent{\it Step 1.}
{Take $N=12$ and fix the kernels $\sigma_1,\sigma_2$ in \eqref{e:mult-kernel-choice}.
As in \cite[Lemma~4.6]{FP}, let $\Psi_1=\sigma_1*\sigma_1$ and $\Psi_2=\sigma_2*\sigma_2$.
Since $\sigma_1$ is real and even,
$$
(|D|\Psi_1)(0)=\int_{\R}|\xi_1|\,|\widehat\sigma_1(\xi_1)|^2\,d\xi_1>0.
$$
Likewise, $\Psi_2(0)=\|\sigma_2\|_2^2>0$. Thus $\Psi_1,\Psi_2$ satisfy \eqref{e:supp-origin}.
Under these conditions, the reverse estimate has order $p$.\footnote{At aperture one,
this also implies the order-$p^{3/2}$ reverse bound in \cite[Theorem~4.2]{FP}
for these convolution-square kernels.
For the full kernel class in \cite[Definition~2.8]{FP},
(iii) of Theorem~\ref{t:sharp-cases} allows order $p^2$ even at aperture one.}

Let $f\in\mathcal S(\R^3)$. Since $m$ is bounded, $m\widehat f\in L^1\cap L^2$.
Fourier inversion and Plancherel give $T_mf\in L^\infty\cap L^2$, hence $T_mf\in L^p$ for $2\le p<\infty$.
Theorem~\ref{t:primal} at aperture one and Lemma~\ref{l:mult-gstar} now give
\begin{align*}
\|T_m f\|_{L^p}
&\le C_\Psi p\,\|S_{z,1}^{\Psi}(T_m f)\|_{L^p}
\le C_\sigma p\,\|m\|_{\mathcal M_{Z,x_1}^{12}}\,\|g_z^{*,\sigma}f\|_{L^p}
\le C p^2\,\|m\|_{\mathcal M_{Z,x_1}^{12}}\,\|f\|_{L^p}.
\end{align*}
The kernels are fixed independently of $p$ and $m$. Density gives \eqref{e:mult-upper} on $L^p(\R^3)$.}

\smallskip
\noindent{\it Step 2.}
To prove optimality, recall the fixed symbol
$$
m_\star(\xi)=
\begin{cases}
\displaystyle\frac{\xi_1\xi_2}{((\xi_1\xi_2)^2+\xi_3^2)^{1/2}},
 &(\xi_1\xi_2,\xi_3)\ne(0,0),\\[4pt]
0, &(\xi_1\xi_2,\xi_3)=(0,0).
\end{cases}
$$
We will show that $m_\star$ satisfies the multiplier assumptions and
$\|T_{m_\star}\|_{L^p\to L^p}\ge cp^2$ for $p\ge2$.
The symbol is smooth when $\xi_1\ne0$ and
$(\xi_2,\xi_3)\ne(0,0)$.
Moreover,
\begin{equation}\label{e:m-dilation}
m_\star(r\xi_1,q\xi_2,rq\xi_3)=m_\star(\xi), \quad \forall\,  r,q>0.
\end{equation}
On $\overline{A^{x_1}}$, we have
$$
(\xi_1\xi_2)^2+\xi_3^2
\ge\frac14(\xi_2^2+\xi_3^2)\ge\frac1{16}.
$$
Thus $m_\star$ is smooth on a neighborhood of this compact set.
By \eqref{e:mult-snorm} and \eqref{e:m-dilation}, for every  $N\geq 1$,  we have
$$
\|m_\star\|_{\mathcal M_{Z,x_1}^{N}}
=\max_{0\le a,b,c\le N}\sup_{\xi\in A^{x_1}}
\left|\partial_1^a\partial_2^b\partial_3^c m_\star(\xi)\right|<\infty.
$$

Let $H$ denote the Hilbert transform on $\R$.
Fix $p\ge2$, let $\varepsilon_p=e^{-2p}/4$, and choose $h_p\in C_c^\infty((0,1))$
with $0\le h_p\le1$ and $h_p=1$ on $[\varepsilon_p,1/2]$.
For $-\varepsilon_p<x_1<0$, we have
\begin{align*}
|Hh_p(x_1)|
\geq &\frac1\pi\int_{\varepsilon_p}^{1/2}\frac{dy_1}{y_1+\varepsilon_p}
=\frac1\pi\log\frac{1/2+\varepsilon_p}{2\varepsilon_p}
\geq 
\frac{2p}{\pi}.
\end{align*}
Since $0<\|h_p\|_{L^p}\le1$, we obtain
$$
\frac{\|Hh_p\|_{L^p}}{\|h_p\|_{L^p}}
\ge\frac{2p}{\pi}\varepsilon_p^{1/p}
\ge\frac{e^{-2}}{\pi}p.
$$

To use the Hilbert-transform lower bound in both $x_1$ and $x_2$,
choose a fixed nonzero $\Theta\in\mathcal S(\R)$ and let
$$
f_{p,L}(x)=h_p(x_1)h_p(x_2)\Theta(x_3/L),\qquad L\ge1.
$$
Then $f_{p,L}\in\mathcal S(\R^3)$, and
$
\widehat f_{p,L}(\xi)
=L\widehat h_p(\xi_1)\widehat h_p(\xi_2)\widehat\Theta(L\xi_3).
$
By Fourier inversion and a change of variable in the third frequency coordinate,
we have
$$
(T_{m_\star}f_{p,L})(x_1,x_2,Lx_3)
=\int_{\R^3}m_\star(\xi_1,\xi_2,\xi_3/L)
\widehat h_p(\xi_1)\widehat h_p(\xi_2)\widehat\Theta(\xi_3)
e^{2\pi ix\cdot\xi}\,d\xi.
$$
For $\xi_1\xi_2\ne0$,
$$
\lim_{L\to\infty}m_\star(\xi_1,\xi_2,\xi_3/L)
=\frac{\xi_1\xi_2}{|\xi_1\xi_2|}
=\operatorname{sgn}(\xi_1)\operatorname{sgn}(\xi_2).
$$
Since $|m_\star|\le1$ and
$\widehat h_p,\widehat\Theta\in L^1(\R)$,
we obtain by dominated convergence
\begin{align*}
\lim_{L\to\infty}(T_{m_\star}f_{p,L})(x_1,x_2,Lx_3)
= &\int_{\R^3}\operatorname{sgn}(\xi_1)\operatorname{sgn}(\xi_2)
\widehat h_p(\xi_1)\widehat h_p(\xi_2)\widehat\Theta(\xi_3)
e^{2\pi ix\cdot\xi}\,d\xi\\
= &-Hh_p(x_1)Hh_p(x_2)\Theta(x_3),
\qquad x\in\R^3.
\end{align*}
By Fatou's lemma and the change of variable $y_3=Lx_3$,
we have
\begin{align*}
\|Hh_p\|_{L^p(\R)}^2\|\Theta\|_{L^p(\R)}
\leq &\liminf_{L\to\infty}
\|(T_{m_\star}f_{p,L})(x_1,x_2,Lx_3)\|_{L^p(dx)}\\
= &\liminf_{L\to\infty}
L^{-1/p}\|T_{m_\star}f_{p,L}\|_{L^p(\R^3)}\\
\leq &\|T_{m_\star}\|_{L^p\to L^p}
\liminf_{L\to\infty}L^{-1/p}\|f_{p,L}\|_{L^p(\R^3)}\\
= &\|T_{m_\star}\|_{L^p\to L^p}
\|h_p\|_{L^p(\R)}^2\|\Theta\|_{L^p(\R)}.
\end{align*}
Since $\|h_p\|_{L^p(\R)}^2\|\Theta\|_{L^p(\R)}>0$, division by this quantity gives
$$
\|T_{m_\star}\|_{L^p\to L^p}
\ge
\left(\frac{\|Hh_p\|_{L^p(\R)}}{\|h_p\|_{L^p(\R)}}\right)^2
\ge\frac{e^{-4}}{\pi^2}p^2,
\qquad p\ge2.
$$
This proves the optimality of $p^2$.
\end{proof}

\begin{remark}\label{rem:aeconverg}
The  estimate \eqref{e:mult-ptw} also holds almost everywhere for
$f\in L^p(\R^3)$, $2\le p<\infty$.
Indeed, choose $f_j\in\mathcal S(\R^3)$ with $\|f_j-f\|_{L^p}\to0$.
Combining Theorem~\ref{t:multipliers}, Lemma~\ref{l:continuous-forward} and \eqref{e:mult-gstar},
$$
\lim_{j\to\infty}\left(
\|S_{z,1}^{\Psi}(T_mf_j)-S_{z,1}^{\Psi}(T_mf)\|_{L^p}
+\|g_z^{*,\sigma}f_j-g_z^{*,\sigma}f\|_{L^p}\right)=0.
$$
{Along a common subsequence both square functions converge almost everywhere.
Passing to the limit in \eqref{e:mult-ptw} gives the estimate for $f$.}
\end{remark}

\vskip .8cm

\noindent{\bf Acknowledgements:} \par
\noindent
We would like to thank Jill Pipher for all the discussions on the paper during her visits
to Macquarie University in 2024, 2025 and 2026, after we completed our joint work with Jill Pipher on multiparameter Hoeffding inequalities in the tensor product setting \cite{LPW} in 2024.
We recognized that the pointwise intermediate square functions used there
provided a natural starting point for the present project on Zygmund dilations.

Li is supported by the Australian Research Council (DP220100285 and DP260100485).
Wu is supported by the National Natural Science Foundation of China (Grant No. 12201002).
This work was initiated during Wu's postdoctoral fellowship at Macquarie University in 2023--2024.
He thanks Macquarie for its hospitality and academic support.

Li also thanks Anhui University for its hospitality during his visit to Wu in Dec 2025.
During this visit, Li and Wu discussed square functions associated with Zygmund dilations,
as well as related flag and twisted multiparameter problems.


\bigskip

\begin{thebibliography} {99}

\bibitem{CWW}
S.-Y. A. Chang, J.M. Wilson and T.H. Wolff,
\emph{Some weighted norm inequalities concerning the Schr\"odinger operators},
Comment. Math. Helv. 60  (1985), no.~2, 217--246.

\bibitem{CDWY}
P. Chen, X. T. Duong, L.C. Wu and L.X. Yan,
\emph{Exponential-square integrability, weighted inequalities for the square functions associated to operators,
and applications},
Int. Math. Res. Not. IMRN 2021  (2021), no.~23, 18057--18117.

\bibitem{CLWY}
P. Chen, X. Lin, L.C. Wu and L.X. Yan,
\emph{On maximal functions generated by H\"ormander-type spectral multipliers},
J. Funct. Anal. 290 (2026), no.~2, Paper No.~111238.

\bibitem{Cor}
A. C\'ordoba,
\emph{Maximal functions: a proof of a conjecture of A. Zygmund},
Bull. Amer. Math. Soc.  (N.S.)  1  (1979), no.~1, 255--257.

\bibitem{FP}
R. Fefferman and J. Pipher,
\emph{Multiparameter operators and sharp weighted inequalities},
Amer. J. Math. 119  (1997), 337--369.

\bibitem{GHS}
L. Grafakos, P. Honz\'ik and A. Seeger,
\emph{On maximal functions for Mikhlin--H\"ormander multipliers},
Adv. Math. 204 (2006), no.~2, 363--378.

\bibitem{GM}
L. Grafakos and S. Montgomery-Smith,
Best constants for uncentred maximal functions,
Bull. London Math. Soc. \textbf{29} (1997), no.~1, 60--64.

\bibitem{GRSY}
S. Guo, J. Roos, A. Seeger and P.-L. Yung,
\emph{A maximal function for families of Hilbert transforms along homogeneous curves},
Math. Ann. 377 (2020), no.~1--2, 69--114.

\bibitem{GuthDecoupling}
L. Guth,
\emph{18.118 Decoupling, Lecture 3: Decoupling for the Paraboloid I},
MIT lecture notes, 
\url{https://math.mit.edu/~lguth/Math118/DecLect3.pdf}.

\bibitem{Han94}
Y. Han,
\emph{Calder\'on-type reproducing formula and the $Tb$ theorem},
Rev. Mat. Iberoamericana 10 (1994), no.~1, 51--91.

\bibitem{HLLW}
Y. Han, M.-Y. Lee, J. Li and B. D. Wick,
\emph{Maximal functions, Littlewood--Paley theory, Riesz transforms and atomic decomposition in the
multi-parameter flag setting},
Mem. Amer. Math. Soc. 279 (2022), no.~1373, 102 pp.

\bibitem{HLLT}
Y. Han, J. Li, C.-C. Lin and C. Tan,
\emph{Singular integrals associated with Zygmund dilations},
J. Geom. Anal. 29 (2019), no.~3, 2410--2455.

\bibitem{HLWW}
Y. Han, J. Li, B. D. Wick and L.C. Wu,
\emph{A theory of stratified concentration for Euclidean flag square functions}, preprint, 2026.

\bibitem{Hoeffding}
W. Hoeffding,
\emph{Probability inequalities for sums of bounded random variables},
J. Amer. Statist. Assoc. 58 (1963), no.~301, 13--30.

\bibitem{HLMV}
T. Hyt\"onen, K. Li, H. Martikainen and E. Vuorinen,
\emph{Multiresolution analysis and Zygmund dilations},
Amer. J. Math. 148  (2026), no.~2, 313--387.

\bibitem{JMZ}
B. Jessen, J. Marcinkiewicz and A. Zygmund,
\emph{Note on the differentiability of multiple integrals},
Fund. Math. 25  (1935), 217--234.

\bibitem{JN}
F. John and L. Nirenberg,
\emph{On functions of bounded mean oscillation},
Comm. Pure Appl. Math. 14 (1961), no.~3, 415--426.

\bibitem{LPW}
J. Li, J. Pipher and L.C. Wu,
\emph{Chernoff--Hoeffding inequalities, almost orthogonality and bounded
multi-parameter square functions}, preprint, 2024.

\bibitem{LPWCWW}
J. Li, J. Pipher and L.C. Wu,
\emph{Chang--Wilson--Wolff estimates for Schr\"odinger operators:
zero-energy obstruction and full-energy estimates},
preprint, 2026.

\bibitem{Pipher86}
J. Pipher,
\emph{Bounded double square functions},
Ann. Inst. Fourier  (Grenoble)  36  (1986), no.~2, 69--82.

\bibitem{Saff}
E. B. Saff,
\emph{Logarithmic potential theory with applications to approximation theory},
Surv. Approx. Theory 5 (2010), 165--200.

\bibitem{SVBMO}
L. Slavin and V. Vasyunin,
\emph{Sharp $L^p$ estimates on BMO},
Indiana Univ. Math. J. 61 (2012), no.~3, 1051--1110.

\bibitem{Sor}
F. Soria,
\emph{Examples and counterexamples to a conjecture in the theory of differentiation of integrals},
Ann. of Math.  (2)  123  (1986), no.~1, 1--9.

\bibitem{Stein}
E.M. Stein,
\emph{Harmonic analysis: real-variable methods, orthogonality, and oscillatory integrals},
Princeton Mathematical Series, vol.~43, Princeton University Press, Princeton, NJ, 1993.

\bibitem{SW}
\v Z. \v Strkalj and L. Weis,
\emph{On operator-valued Fourier multiplier theorems},
Trans. Amer. Math. Soc. 359 (2007), no.~8, 3529--3547.

\bibitem{Vershynin}
R. Vershynin,
\emph{High-dimensional probability: An introduction with applications in data science},
Cambridge Series in Statistical and Probabilistic Mathematics, vol.~47,
Cambridge University Press, Cambridge, 2018.

\bibitem{Wilson89}
J.M. Wilson,
\emph{Weighted norm inequalities for the continuous square function},
Trans. Amer. Math. Soc. 314 (1989), no.~2, 661--692.

\bibitem{Wilson08}
M. Wilson,
\emph{Weighted Littlewood-Paley Theory and Exponential-Square Integrability},
Lecture Notes in Mathematics, vol.~1924, Springer, Berlin, 2008.

\bibitem{WuJGA}
L.C. Wu,
\emph{From Hoeffding's inequality to Chang--Wilson--Wolff theorem---a unified approach},
J. Geom. Anal. 35  (2025), Paper No.~261, 18 pp.

\end{thebibliography}
\end{document}